\documentclass[a4paper]{amsart}

\usepackage{amsmath, stmaryrd}
\usepackage{amsfonts}
\usepackage{amssymb}
\usepackage{graphicx}
\usepackage{color}
\usepackage{xfrac}
\usepackage[shortlabels]{enumitem}

\usepackage{calligra}
\DeclareMathAlphabet{\mathcalligra}{T1}{calligra}{m}{n}

\newtheorem{lemma}{Lemma}[section]
\newtheorem{proposition}[lemma]{Proposition}
\newtheorem{corollary}[lemma]{Corollary}
\newtheorem{example}[lemma]{Example}

\theoremstyle{remark}
\newtheorem{remark}[lemma]{Remark}

\usepackage{tikz}
\usepackage{tikz-cd}
\tikzset{cong/.style={draw=none,edge node={node [sloped, allow upside down, auto=false]{$\cong$}}},
         Isom/.style={draw=none,every to/.append style={edge node={node [sloped, allow upside down, auto=false]{$\cong$}}}}}
\usetikzlibrary{arrows}

\renewcommand{\imath}{\mathrm{i}}

\newcommand{\CC}{\mathcal{C}}
\newcommand{\CL}{\mathcal{L}}
\newcommand{\CQ}{\mathcal{Q}}
\newcommand{\CF}{\mathcal{F}}
\newcommand{\CS}{\mathcal{S}}
\newcommand{\CM}{\mathcal{M}}
\newcommand{\CR}{\mathcal{R}}

\newcommand{\XM}{\hbox{${\buildrel \times\over {\mathcal M}}$}}
\newcommand{\xtriangleleft}{\hbox{${\buildrel \times\over \triangleleft}$}}

\newcommand{\xPsi}{\hbox{${\buildrel \times\over \Psi}$}}

\newcommand{\C}{\mathbb{C}}
\newcommand{\R}{\mathbb{R}}

\newcommand{\Z}{\mathbb{Z}}

\newcommand{\sign}{\mathrm{sign}}

\newcommand{\del}{\partial}

\newcommand{\extd}{\mathrm{d}}

\newcommand{\isom}{{\cong}}
\newcommand{\eps}{{\epsilon}}
\newcommand{\tens}{\mathop{{\otimes}}}
\newcommand{\la}{{\triangleright}}
\newcommand{\ra}{{\triangleleft}}

\newcommand{\Ad}{\mathrm{ Ad}}

\newcommand{\id}{\mathrm{id}}

\newcommand{\<}{\langle}
\renewcommand{\>}{\rangle}
\newcommand{\End}{\mathrm{ End}}
\newcommand{\from}{{\longleftarrow}}

\newcommand{\lcross}{{>\!\!\!\triangleleft}}

\newcommand{\lbiprod}{{>\!\!\!\triangleleft\kern-.33em\cdot}}
\newcommand{\rbiprod}{{\cdot\kern-.33em\triangleright\!\!\!<}}
\newcommand{\rcocross}{{\blacktriangleright\!\!<}}
\newcommand{\lcocross}{{>\!\!\blacktriangleleft}}

 \newcommand{\cg}{{\mathfrak{g}}}
 \newcommand{\metric}{{\mathfrak{g}}}

\allowdisplaybreaks

\makeatletter

\renewcommand\subsubsection{\@secnumfont}{\bfseries}%
\renewcommand\subsubsection{\@startsection{subsubsection}{3}
  \z@{.5\linespacing\@plus.7\linespacing}{-.5em}%
  {\bfseries\itshape}}
  
  \makeatother

\begin{document}

\author{ Shahn Majid and  Leo Sean McCormack}
\address{School of Mathematical Sciences\\ Queen Mary University of London \\ Mile End Rd, London E1 4NS }
\email{ s.majid@qmul.ac.uk, l.s.mccormack@qmul.ac.uk}
\thanks{Ver. 1.06}

\title{Quantum geometric Wigner construction for $D(G)$ and braided racks}
	
\begin{abstract} The quantum double $D(G)=\C(G)\lcross \C G$  of a finite group plays an important role in the Kitaev model for quantum computing, as well as in associated TQFT's, as a kind of Poincar\'e group. We interpret the known construction of its irreps, which are quasiparticles for the model, in a geometric manner strictly analogous to the Wigner construction for the usual Poincar\'e group of $\R^{1,3}$. Irreps are labelled by pairs $(\CC, \pi)$, where $\CC$ is a conjugacy class in the role of a mass-shell, and $\pi$ is a representation of the isotropy group $C_G$ in the role of spin. The geometric picture entails $D^\vee(G)\to \C(C_G)\rcocross \C G$ as a quantum homogeneous bundle where the base is $G/C_G$, and $D^\vee(G)\to  \C(G)$ as another homogeneous bundle where the base is the group algebra $\C G$ as noncommutative spacetime. Analysis of the latter leads to a duality whereby the differential calculus and solutions of the wave equation on $\C G$  are governed by irreps and conjugacy classes of $G$ respectively, while the same picture on $\C(G)$ is governed by the reversed data.  Quasiparticles as irreps of $D(G)$ also turn out to classify irreducible bicovariant differential structures $\Omega^1_{\CC, \pi}$ on $D^\vee(G)$ and these in turn correspond to braided-Lie algebras $\CL_{\CC, \pi}$ in the braided category of $G$-crossed modules, which we call `braided racks' and study. We show under mild assumptions that $U(\mathcal{L}_{\CC,\pi})$ quotients to a braided Hopf algebra $B_{\CC,\pi}$ related by transmutation to a coquasitriangular Hopf algebra $H_{\CC,\pi}$.  \end{abstract}
\maketitle

\section{Introduction}\label{secintro}

For a finite (or discrete) group $G$, Whitehead \cite{Whi} essentially invented what is now called the braided-category of crossed $G$-modules. An object here is a vector space $V$ on which $G$ acts and which is also $G$-graded, i.e. $V=\oplus_{g\in G} V_g$ such that  $h\la: V_g\to  V_{h g h^{-1}}$ for all $h\in G$. The latter says that if $v$ has grade $|v|$ then $|h\la v|=h|v|h^{-1}$. Two such crossed modules have a braiding
\begin{equation}\label{GcrossPsi} \Psi_{V,W}(v\tens w)=|v|\la w\tens v, \end{equation}
if $v$ has grade $|v|\in G$. Moreover, this braided category is a simple example of how quantum groups generate braided categories, namely there is a quantum group $D(G)=\C (G)\lcross \C  G$, where $G$ acts by the dual adjoint action and the coalgebra is the tensor product of the coalgebras of $\C (G)$, the functions on $G$, and $\C  G$ the group algebra of $G$.  The latter is the vector space with basis labelled by $G$ and the product of $G$ extended linearly. The quantum group $D(G)$ is an easy example of a general Hopf algebra construction \cite{Dri} known as the Drinfeld double $D(H)$ (for a finite-dimensional Hopf algebra $H$), and more generally its braided category of representations is the category of $H$-crossed modules, a.k.a. Yetter-Drinfeld modules \cite{Yet}  and an example of the Drinfeld-Majid centre of a monoidal category \cite{Ma:rep}. This quantum group and braided category are important as  key structures behind the Kitaev model \cite{Kit} for fault-tolerant quantum computing, see \cite{Kit,CowMa} for recent work. One can also view the associated TQFT as a baby version of 2+1 quantum gravity, where the latter is based on $D(U(su_2))= C[SU_2]\lcross U(su_2)$ for the Lie group $SU_2$ rather than a finite group, or $D(U_q(su_2))= B_q[SU_2]\lbiprod U_q(su_2)$ if we include a cosmological constant, see \cite{MaSch} for more background here. 

In both points of view, $D(G)$ plays the role of `Poincar\'e group' in that its representations are the quasiparticles of the Kitaev model among the lattice states. It also acts naturally on the group algebra $\C  G$ and one might expect that this should play the role of a model spacetime for the theory, similarly to the way that $D(U(su_2))$ in  2+1 quantum gravity  acts on $U(su_2)$ as a `fuzzy $\R^3$' model quantum spacetime \cite{BatMa,FreLiv,Hoo}. Hence in this paper, we develop the role of $D(G)$ as such a Poincar\'e group action on $\C  G$ as a noncommutative space and show how quasi-particles as its irreps $V_{\CC,\pi}$ can be interpreted as a subset of fields on this noncommutative space. It is already known that irreps are labelled by conjugacy class $\CC\subset G$ and a representation $\pi$ of the isotropy group $C_G$ of a chosen representative $r\in \CC$, see \cite{Ma:cla,CowMa}, and our main result, in Section~\ref{secWigner},  is to understand their construction geometrically and  in line with the well-known Wigner construction of irreps of the Poincar\'e group of $\R^{1,3}$.

To do this, we will need to work with the `coordinate algebra' of the quantum group, which in our case means the dual Hopf algebra $D^\vee(G)=\C (G)\rcocross \C  G$ (in the role of functions on the Poincar\'e group). In this dual language, we will construct two different quantum homogeneous bundles given by Hopf algebra surjections $P\twoheadrightarrow H$, with base $A=P^H$ the invariant subalgebra. This uses  the language of quantum bundles \cite{BrzMa, BegMa} a.k.a. Hopf-Galois extensions. In the first bundle, $P=D^\vee(G)\twoheadrightarrow \C (C_G)\rcocross \C G$ and $A=\C (G/C_G)$, i.e. the base is a classical but discrete coset space $G/C_G$. In the second case, $P\twoheadrightarrow  \C (G)$ and the base is $\C G$. These correspond to two different bundles but we find that there is a natural `transfer map' from sections of corresponding associated bundles. This allows us to express the known $V_{\CC, \pi}$ first of all geometrically as sections of a bundle with base $G/C_G$ and then to transfer this to sections of the second bundle over `spacetime' $\C  G$.  A recap of the classical Wigner construction on spacetime, which serves as our model, is in Section~\ref{secW}, with more details to be found in physics texts or the review \cite{Fig}. As classically, the bundles are trivial so that the final result is that irreps appear as certain subspaces of vector-valued fields over $\C G$. 

Next, in Section~\ref{secPcalc}, we look at differential structures in the bundles that we used for the Wigner construction in Section~\ref{secWigner}. A rather special feature of `factorisable' quantum groups such as $D^\vee(G)$ is that nontrivial  irreps  or `quasiparticle types' themselves correspond bijectively with coirreducible bicovariant differential structures $\Omega^1(D^\vee(G))$. On a classical Lie group the differential structure is unique if we ask for a translation-invariant one, but this is not true for quantum groups. The quasi-particle $V_{\CC,\pi}$ viewed this way as a differential structure induces a tautological differential structure on $\Omega^1_{\tilde\pi}$ on $\C G$ as spacetime and we use this in Section~\ref{eq_free_field_equation} to present the quasiparticle fields as solutions of a $D^\vee(G)$-covariant wave equation which is a second order operator or `Laplacian' with respect to a suitable $D^\vee(G)$-covariant inner product metric on $\Omega^1_{\tilde \pi}$.  This is the analogue of irreps of the Poincar\'e group presented as solutions to the Klein-Gordon wave equation. We also give conditions for  and related to this, and elements of quantum Riemannian geometry on $\C  G$ such that $D(G)$ is indeed an `isometry quantum group' for a suitable inner product on $\Omega^1$. The first part of this entails a study of $D(G)$-covariant 2nd order `Laplacians' $L:\C  G\to \C G$ and their eigenvalues in terms of  conjugacy classes $\CC$ in the role of `mass'. The second part entails the construction of these from a suitable $D(G)$-covariant calculus $\Omega^1$ and `metric' on $\C  G$. Among the examples, we give full details for $G=S_3$, the group of permutations on 3 elements, and $\CC$ its order 2 conjugacy class. Then $C_G$ is $\Z_3$ and the role of `spin' is played by $j=0,1,2$ for its irreps. Our approach to these results leads us into a  remarkable duality, demonstrated further in Section~\ref{sec_duality} between differential structures and covariant wave equations  on $\C G$ and the same on the Fourier dual $\C(G)$. Whereas in Section~\ref{eq_free_field_equation}, differential structures  on $\C G$ correspond  to (sums of) irreducible representations of $G$ and eigenspaces of the $D^\vee(G)$ covariant wave operator to (unions of) conjugacy classes, we show that the role of this data is reversed for the  discrete geometry $\C(G)$ as spacetime. This generalises and understands abstractly a phenomenon first observed in \cite{MaTao}. 

In Section~\ref{secbrack}, we explore a corollary of $V_{\CC,\pi}$ providing the differential structure $\Omega^1_{\CC,\pi}$ on $D^\vee(G)$, namely, as for any coquasitriangular Hopf algebra, the dual of the space of left-invariant 1-forms is a braided-Lie algebra $\CL_{\CC,\pi}$, which we compute and study. The physical motivation here is more long-term but to fully understand the geometry and continuum limit of the quantum group gauge theory underlying the Kitaev and other models, one should work with gauge fields as `Lie-algebra valued' rather than, as normally done on a lattice, thinking of them as quantum group valued holonomy. Braided-Lie algebras \cite{Ma:lie,Ma:sol,GomMa} provide the relevant framework and feature a pentagonal Jacobi identity in order to make sense in a braided category. They were introduced  for quantum groups such as $U_q(su_2)$, where they solved the `Lie problem' (namely, a finite-dimensional Lie-algebra-like object that generates the quantum group). By contrast,  $\CL_{\CC,\pi}$ provide for the first time much simpler yet strictly braided examples that are  closer to super and symmetric-monoidal (colour) Lie algebras, with -1 factors in the super case when odd elements are transposed now  replaced by (\ref{GcrossPsi}). We describe these objects in a self-contained manner and compute their braided Killing forms and braided enveloping algebras $U(\mathcal{L}_{\CC, \pi})$ as well as the `covering' algebra map $U(\mathcal{L}_{\CC, \pi}) \to D(G)$. Here,  $U(\mathcal{L}_{\CC, \pi})$ is a transmutation of an ordinary FRT-type bialgebra $A_{\CC,\pi}$ which we also study and we show how, in some generality, both can be quotiented to braided and ordinary Hopf algebras, $B_{\CC,\pi}, H_{\CC,\pi}$ respectively. Everything is demonstrated for $G=S_3$ and the case $\CC$=2-cycles and $\pi=1$ (the trivial) provides a particularly nontrivial but amenable 9-dimensional strictly braided-Lie algebras. The analog of $\CL_{\CC,\pi}$ for $\C(G)$ (rather than for its quantum double) is not braided (the underlying braiding is trivial) but even these are of interest as `Lie algebras' of finite groups \cite{LMR,MaRie}. In this case the axioms of a braided-Lie algebra reduce to those of a rack, hence $\CL_{\CC,\pi}$ should be thought of as `braided racks'. 

The methods of the paper readily generalise to other quantum doubles $D(H)$ where $H$ is a finite-dimensional Hopf algebra and, with suitable duals, a more general Hopf algebra. Moreover, the braided racks that we construct can be seen as examples of a very general construction that associates a braided-Lie algebra to any rigid object of an Abelian braided category. This result and some formulae for the Hopf algebra case will be given elsewhere \cite{MaMc} and in principle cover the case of 2+1 quantum gravity as discussed above with $H=U(su_2)$ or $U_q(su_2)$.

\section{Algebraic preliminaries}\label{secpre}

We will use a small amount of Hopf algebra or quantum group methods, so we recall the basic definitions here. More details can be found in \cite{Ma}. General Hopf algebra and quantum differential constructions other than those involving $*$ work over any field  $k$. For the rest of the paper, however, we work with $k=\C$ both from the physical context, where compatibility with $*$ is needed for `unitarity' of the constructions, and to avoid issues when working with irreducible representations (irreps) or matrix blocks.  

\subsection{Hopf algebras} We recall that a bialgebra over a field $k$ is a unital algebra $H$ over $k$ which is also a counital coalgebra, i.e. there are linear maps $\Delta:H\to H\tens H$ and $\eps:H\to k$ obeying the arrow-reversal of the usual axioms of a unital algebra. We also require $\Delta,\eps$ to be algebra homomorphisms (later, for Hopf algebras in a braided category, there is a particular braided tensor product algebra $H\underline\tens H$ to use here where the two factors do not mutually commute). A Hopf algebra additionally has an `antipode' $S:H\to H$ obeying $h_1Sh_2=\eps(h)=(Sh_1)h_2$ for all $h\in H$ where we write $\Delta h=h_1\tens h_2$ as a `Sweedler-type' notation where the right hand side stands for a sum of such terms. Iterated coproducts are numbered similarly according to the position in the tensor product. For a tensor product of two arbitrary vector spaces $X \otimes Y$, a general element $P \in X \otimes Y$ will be written as $P^{(1)} \otimes P^{(2)}$ (sum of terms understood) when we want to isolate the tensor factors. 

An $H$-module just means a representation of $H$ as an algebra, but the coproduct provides a tensor product of left $H$-modules making the category ${}_H\CM$ of these monoidal. Similarly for $\CM_H$ the category of right $H$-modules. Dually, a right $H$-comodule means $(V,\Delta_R)$ where $\Delta_R:V\to V\tens H$ obeys $v^{11}\tens v^{12}\tens v^2=v^1\tens v^2{}_1\tens v^2{}_2$ using the notation $\Delta_R v=v^1\tens v^2$ (sum of such terms understood). This is merely the arrow-reversed notion of a (right) action and the category $\CM^H$ of all such is monoidal with tensor product provided by the algebra of $H$.  If $H$ is finite dimensional then $H^*$ is also a Hopf algebra (with invertible antipode) and a left action of $H^*$ can be canonically identified with evaluation against a right coaction of $H$. Similarly, the category ${}^H\CM$ of left $H$-comodules can be identified with $\CM_{H^*}$ for $H$ finite-dimensional. 

\subsection{Quantum bundles} A quantum principal bundle with universal calculus (a.k.a. Hopf-Galois extension) is an algebra $P$ which is covariant under a coaction of a Hopf algebra $H$ ($P$ is a right $H$-comodule algebra), such that the induced map $P\tens_A P\to P\tens H$ sending $p\tens q\mapsto p\Delta_R q$ is a bijection.  Here $A=P^H$ is the subalgebra of invariant elements in the coaction, i.e. of $a\in P$ such that $\Delta_Ra=a\tens 1$, and is the base of the quantum principal bundle. In algebra, one says that $P$ is a Hopf-Galois extension of $A$, while the geometric picture (with differential structures) is from \cite{BrzMa}. A bundle is {\em trivial} if there is a left $A$-module right $H$-comodule map isomorphism $\Phi:P\to A\tens H$. More details can be found in \cite[Chap. 6]{Ma}\cite[Chap~5]{BegMa}. The key example of a principal bundle of interest here is that obtained from a surjective Hopf algebra map $p:P\twoheadrightarrow H$ for which the induced coaction $\Delta_R=(\id\tens p)\Delta$ makes $P$ a quantum principal bundle. As such, $A$ is interpreted as functions on the resulting `coset manifold'. A key feature is that this quantum bundle is that it is additionally equipped with left $P$ coaction from the coproduct of $P$, which shall be called the `homogeneous structure' of the bundle.

For a principal bundle $P$ with structure group $H$, along with a right $H$-comodule $V$, we can define the $A$-bimodule of sections for its associated vector bundle as the vector space $E:= (P \otimes V)^{H}$ (and $A$-bimodule structure induced by the $P$-factor), in direct analogy of the fact that the space of sections of a classical associated vector bundle can be identified with the space of equivariant functions on the total space with values in $V$ (also known as the Mackey functions).  If $P$ is trivial as a bundle, then $\Phi$ induces a canonical isomorphism $\Theta_{H}^{V}: E \cong (A \otimes H \otimes V)^{H} \cong A \otimes V$ of bimodules (with the latter's bimodule structure simply induced by the isomorphism). In the case of the bundle induced by a Hopf algebra map $P \twoheadrightarrow H$, the space of sections $E$ retains a left coaction from the coproduct of $P$ as per the classical geometric construction or induced modules for a homogeneous space, now in an algebraic coaction setting. 

\subsection{Quantum differential structures} \label{section_quantum_diff_structures_intro}

On any unital algebra $A$, a (first order) differential calculus means  $(\Omega^1,\extd)$, where $\Omega^1$ is an $A$-bimodule in the role of `differential 1-forms' and $\extd: A\to \Omega^1$ a linear map obeying the Leibniz rule $\extd(ab)=(\extd a)b+a\extd b$ for all $a,b\in A$, in the role of the `exterior derivative'. We will additionally require $\Omega^1$ to be generated by $A,\extd A$. The calculus is then said to be \textit{connected} if $\textup{ker}(d) = k \cdot 1_{A}$ and \textit{inner} if there exists an element $\theta \in \Omega^{1}$ for which $\extd a = [\theta, a] := \theta \cdot a - a \cdot \theta$. Every algebra has the `universal calculus' $\Omega^1_u:=\ker(\cdot: A\tens A\to A)$ with $\extd_u a= 1\tens a-a\tens 1$, and is universal in the sense that any other calculus is a quotient of it. Next, a unital algebra map $\phi:A\to B$ is `differentiable' if $\phi_{*}(a\extd b):= \phi(a)\extd\phi(b)$ is a well-defined $A$-bimodule map $\Omega^{1}_{A} \to \Omega^{1}_{B}$. 

If $H$ is a Hopf algebra, a calculus $\Omega^{1}$ for it is `bicovariant' iff the following are well-defined left and right $H$-coactions on $\Omega^{1}$ 
\begin{equation} \label{eq_bicovariant_def} \Delta_L(h\extd g)=h_1 g_1\tens h_2\extd g_2,\quad \Delta_R(h\extd g)=h_1\extd  g_1\tens h_2g_2 \, , \end{equation}
i.e. they extend left and right translation (expressed in the coproduct) compatibly with $\extd$. By the Hopf-module lemma, it is well-known that we obtain the equivalent description $\Omega^1\isom H \otimes \Lambda^1$, where the right-hand side is free as left $H$-module with the rest of the structure maps given by 
\begin{equation} \label{recover_diffcalc}  v h=h_1 (v\ra h_2),\,  \, \Delta_L(hv)=h_1\tens h_2 v, \, \, \Delta_R (hv)=h_1 v^1\tens h_2 v^2, \, \, \extd h = h_{1} \varpi(\pi_{\epsilon}h_{2}), \end{equation} where $\pi_\eps h=h-1\eps h$ is the canonical projection whereby $H=k1\oplus H^+$. In the above, $\Lambda^{1}$ is quotient right $H$-crossed module of $H^{+}$, as described by the Maurer-Cartan form 
\begin{equation} \label{eq_def_maurer_cartan} \varpi: H^+ \twoheadrightarrow \Lambda^1,\quad \varpi(h)=Sh_1\extd h_2, \end{equation} with $H^{+}$ a crossed module by right multiplication and the right adjoint coaction $\Ad_R h=h_2\tens (Sh_1)h_3 $. By \eqref{recover_diffcalc}, $\Lambda^{1}$ is thus the space of left-invariant 1-forms, i.e. $\omega\in \Omega^1$ for which $\Delta_L(\omega)=1\tens \omega$ \cite{Wor}. Moreover, such $\Lambda^{1} \cong H^{+}/I$ are thus classified by subcrossed modules $I$ of $H^{+}$. We note $\Omega^{1}_{u} \cong H \otimes H^{+}$ is bicovariant. Also, a Hopf algebra map $\phi: H \to H'$ is differentiable as an algebra map iff $\phi(I) \subseteq J$ for the subcrossed modules $I, J$ defining the bicovariant calculi on $H$, $H'$, and moreover $\phi_{*}$ is then automatically a $B$-bicomodule map. More details of noncommutative differential geometry of the style that we need are in \cite{BegMa}.

For quantum principal bundles, we did not yet discuss differential structures and have so far described the theory with implicit use of the universal calculus. For a quantum $H$-principal bundle with general (non-universal) calculus \cite{BrzMa}, we need a calculus $\Omega^1_P$ on the algebra $P$ which is $H$-covariant in the sense that the right $H$-coaction $\Delta_R$ on $P$ extends to $\Omega^1_P$. We also need a bicovariant calculus $H \otimes \Lambda^{1}$ on $H$, and a calculus $\Omega^{1}$ on the algebra $A$ for which the inclusion $A \hookrightarrow P$ is differentiable. Finally, we require the following to be a short exact sequence of left $P$-modules \begin{equation} \label{eq_exact_bundle_def} 0\to \Omega^{1}_{hor}\to \Omega^1_P{\buildrel{\rm ver}\over\to} P\tens \Lambda^1 \to 0,\end{equation} with the `vertical map' to be well-defined by ${\rm ver}(p\extd q)=pq^{1} \otimes \varpi(q^{2})$ for $p,q\in P$. This is such that evaluation against its output sends an element of the Lie algebra $\Lambda^1{}^*$ to an associated `vertical vector field' $\Omega^1_P\to P$. The space $\Omega^1_{hor} :=P\Omega^1P$ are then the horizontal forms, pulled back from the base to be viewed as a subset of the differential forms $\Omega^{1}_{P}$ on $P$. A recent exposition of this theory is in \cite[Chap.~5]{BegMa}. 

\subsection{Quantum double $\boldsymbol{D(G)}$} \label{sec_DG_struc} Now let $G$ be a finite group. We denote by $k G$ its group Hopf algebra with product that of $G$ extended linearly, and $\Delta g=g\tens g$, $\eps g=1$, $Sg=g^{-1}$ for all $g\in G$. We denote by $k(G)$ its dual Hopf algebra of functions on $G$ with $\delta_g$ the Kronecker $\delta$-function a natural basis and 
\[ \delta_g\delta_h=\delta_{g,h}\delta_g,\quad 1=\sum_{g \in G} \delta_g,\quad \Delta\delta_g=\sum_{a, b \in G \, | \, ab=g}\delta_a\tens\delta_b,\quad \eps \delta_g=\delta_{g,e},\quad S\delta_g=\delta_{g^{-1}},\]
 where $e$ is the group identity. Then the quantum double $D(G)=k(G)\rtimes k G$ has  $k G$ and $k(G)$ as sub-Hopf algebras, with cross relations $h\delta_g =\delta_{hgh^{-1}}h$ and tensor product coalgebra. We will often prefer to refer to $D(G)$ explicitly on the tensor product vector space, then for example the cross relation appears explicitly as 
 \begin{equation}\label{DG} (1\tens h)(\delta_g\tens 1)= (\delta_{hgh^{-1}}\tens 1)(1\tens h)=\delta_{hgh^{-1}}\tens h, \end{equation}
 and antipode as $S(\delta_g\tens h)=\delta_{h^{-1}g^{-1}h}\tens h^{-1}$. $D(G)$ thus has semi-direct product algebra generated by the action of $kG$ on $k(G)$ given by the dual adjoint action $h \triangleright \delta_{g} = \delta_{hgh^{-1}}$. A left $D(G)$-module is exactly a Whitehead $G$-crossed module where $k(G)$ acts as $f\la v=f(|v|)v$ for $v$ of grade $|v|$. The resulting $D(G)$-action is then \begin{equation} \label{eq_DG_action} (\delta_{g} \otimes h)\la v = \delta_{g}\triangleright (h \triangleright v) = \delta_{g, h|v|h^{-1}} h \triangleright v , \end{equation} recalling that $|h \triangleright v| = h|v|h^{-1}$. The quantum group has a certain `quasitriangular structure' $\mathcal{R}:= \sum_{h \in G}\delta_{h} \otimes h\in D(G)\tens D(G)$ which induces the braiding (\ref{GcrossPsi}) for $G$-crossed modules.

\section{Wigner construction for irreps of $D(G)$}\label{secWigner}

Although the construction of irreps of $D(G)$ is known \cite{Ma:cla,CowMa}, these enter the Kiatev model as quasi-particles and in this section we exhibit them as constructed in an analogous manner to particles in Minkowski spacetime $M=\R^{1,3}$ as irreps of the Poincar\'e group. What plays the role of spacetime will be the group algebra $\C G$ regarded as a noncommutative coordinate algebra. 

\subsection{Recap of Wigner construction for the Poincar\'e group}\label{secW}

For the metric $\eta = \textup{diag}(+1, -1, -1, -1)$ on $\R^{1,3}$, the isometry group is $\mathbb{R}^{1, 3} \lcross O(1, 3)$ and particles are  unitary projective representations of this, and can be obtained from those of the identity component $\mathbb{R}^{1, 3} \lcross SO^{+}(1, 3)$. However, the latter are equivalent (by a theorem of Bargmann) to unitary representations of the universal cover 
\[ P^{1, 3} := \mathbb{R}^{1, 3} \lcross SL(2,\C).\]
Here $T$ is the group of translations but looking at translations from the origin, we an identify $T=M$ with Minkowski space. There is a standard identification of $M$ with hermitian $2\times 2$ matrices via the identity and Paul matrices and in this form the metric is the determinant and $L=SL(2,\C)$ acts by $x\mapsto g x g^\dagger$ for $g\in L$ and $x\in M$. 

For the standard construction of positive-energy irreps, e.g. \cite{Fig} for an introduction, we first fix an element $\alpha\in \hat T$, the set of unitary irreps of $T$. This space can be identified as another copy of $\R^{1,3}$ with we denote $M^*$, where $p_r\in M^*$ corresponds to a plane wave
\[ \alpha(x)= e^{\imath x^\mu p_r{}_\mu},\]
for all $x\in M$, and with the time component of $p_r$ non-negative. Since $L$ acts on $M$, it also acts on $M^*$ and we let $L_{p_r}$ denote the isotropy subgroup of Lorentz transformations that preserve $p_r$, more commonly known as the `little group'. (This is equivalent to $l\in L$ such that $\alpha(l\la t)=\alpha(t)$ for all $t\in T$). We let $H:= T \lcross L_{p_r}\subset P^{1,3}$ and consider the latter as a homogeneous principal bundle over
 \begin{equation} \label{eq_coset_manifold_sheet} P^{1, 3}/H \cong L/L_{p_{r}} \cong H^{1, 3}\subset M^*, \end{equation} 
 where $H^{1, 3}$ denotes the orbit of $p_{r}\in M^*$ under the Lorentz group. Such orbits are characterised by $||p||^2=||p_r||^2=m^2$ using the Lorentzian norm, i.e. by the rest mass $m$ of the base point $p_r$, and the condition $p_0\geq 0$.  There are two cases of interest: 
\[ m > 0: \quad  L_{p_{r}} \cong SU_{2},\]
\[ m =0:\quad  L_{p_{r}} \cong \R^2\lcross U(1), \]
(regarded as a double cover $\R^2 \lcross SO(2)$). Here $H^{1,3}$ for $m>0$ is the upper mass-hyperboloid in $M^*$ (with $p_{0} > 0$) and for $m=0$ this limits to the upper light cone. One can also consider $p_r=m=0$ with $L_{p_r}=L$ and $H^{1,3}$ gets replaced by a point. All bundles here are trivial. 
 
We now additionally choose a unitary irrep $\pi: L_{p_{r}} \to \textup{End}(V_{\pi})$ of  $L_{p_{r}}$ and extend this canonically to an irrep of $H$ by $\pi \otimes \alpha$ acting on the same vector space. Here,
\begin{equation}\label{eq_def_starting_irrep} (x \tens l).v=  e^{i x \cdot p_{r}}\pi(l)v,\end{equation}
for all $v\in V_\pi$. We denote $V_\pi$ with this action by $V_{p_r,\pi}$ and form the associated vector bundle 
\[ E:= P^{1, 3} \times_{H} V_{p_r,\pi}.\]
 We still have a left action by $P^{1,3}$ and the space of sections $\Gamma(E)$ inherits an action of $P^{1, 3}$, which turns out to be unitary and irreducible. Moreover, all positive energy irreps of $P^{1, 3}$ can be obtained in this way. Finally, choosing a different $p_r$ in the same orbit (i.e. the same mass) gives an isomorphic $V_{p_r,\pi}$ and hence isomorphic irrep of $P^{1,3}$. For $m>0$ the irreps of $L_{p_r}$ are labelled by a half-integer spin $j$. For $m=0$, the finite-dimensional irreps of $L_{p_r}$ are 1-dimensional and naturally labelled by a half-integer $j$, the helicity, where $(y\tens e^{\imath\theta}).v= e^{2j\theta}v$ for the action of $\R^2\lcross U(1)$ on $v\in \C$. 

These sections are fields on $H^{1,3}\subset M^*$. The final step in the Wigner construction is to map them over to fields in $M$ by some form of Fourier transform. To this end, we suppose that the $H$-module $V_{p_{r}, \pi}$ above embeds into a $P^{1, 3}$-module $(\rho,W_\rho)$, i.e.,  that the pull-back of $\rho:P^{1,3}\to \End(W_\rho)$ to an $H$-representation contains an isomorphic copy of $V_{p_r,
\pi}$ as a choice of summand. Associated to this data is a projection $P_{p_r,\pi}:W_\rho\to V_{p_r,\pi}$ which projects out the other components. Taking $W_\rho$ with its pulled back $H$-module structure defines a bundle and by the above inclusion and projection, bundle maps 
\[  E^W= P^{1,3}\times_H W_\rho {{\buildrel  P_{p_r,pi} \over \to }\atop \hookleftarrow} E.\]
These induce corresponding projection and inclusion $\Gamma(E^W){\rightarrow\atop \hookleftarrow}\Gamma(E)$.

We now consider the (trivial) homogeneous principal bundle associated to $L\subset P^{1,3}$ which has base
\[   P^{1,3}/L\isom M,\]
i.e. Minkowski space and the associated bundle
\[ E'=P^{1,3}\times_L W_\rho.\]
We still have an action of $P^{1,3}$ from the left and hence the sections $\Gamma(E')$ are a $P^{1,3}$ module. Sections $\Gamma(E')$ are $L$-equivariant functions $C^\infty_L(P^{1,3},W_\rho)$ and sections $\Gamma(E^W)$ are $H$-equivariant functions $C^\infty_H(P^{1,3},W_\rho)$. But any function in the latter space (in fact any element $f\in C^\infty(P^{1,3},W_\rho)$) can be averaged over $L$ so as to arrive in the former space:
 \begin{equation} \label{eq_averagin_integral}  \textup{av}: C^\infty(P^{1,3}, W_\rho) \to  C^\infty_L(P^{1,3}, W_\rho), \quad {\rm av}(f)(g) = \int_{L} d \mu(l)\rho(l)f(g l) \, ,  \end{equation} 
 where  $g \in P^{1,3}$, $l \in L$ and we use the Haar measure $\extd\mu$ on $L$.  Thus, we have maps
 \begin{equation} \label{eq_averaging_full} \Gamma(E) {\from\atop\hookrightarrow}\Gamma(E^W){\buildrel {\rm av}\over\to }\Gamma(E').\end{equation}
 by restricting \eqref{eq_averagin_integral} to $\Gamma(E^{W})$. Moreover, one can show that the ${\rm av}$ is injective on the image of $\Gamma(E)$  and hence we can by these maps realise our irrep $\Gamma(E)$ as a subspace of $\Gamma(E')$.
 
 So far, we have not used that our bundles are trivial. But since they are, we can identify for simplicity
 \[ \Gamma(E) \cong C^{\infty}(\sfrac{L}{L_{p_{r}}}, V_{p_{r}, \pi}),\quad \Gamma(E^{W})\cong C^{\infty}(\sfrac{L}{L_{p_{r}}}, W_\rho),\quad \Gamma(E')\cong C^\infty(M, W_\rho).\]
Then \eqref{eq_averaging_full} implies an embedding that sends a section $s\in \Gamma(E)$ to a section $s'\in \Gamma(E')$ given by \begin{equation} \label{eq_mass_shell_FT_pre}  s'(x) =  \int_{\sfrac{L}{L_{p_{r}}}} d\nu(p) e^{-ix \cdot p} s_{cov}(p),\end{equation} 
where $d\nu(p)$ is the induced measure on $L/L_{p_r}$ defined by pulling back a function on that to a function in $L$ and integrating there, $x\cdot p$ means the pairing with the corresponding $p\in H^{1,3}$, and $s_{cov}(p)=\rho(\sigma(p))s(p)$ for a choice $\sigma: p\mapsto \sigma(p) \in L$ of coset representatives for $\sfrac{L}{L_{p_{r}}}$ (the choice of this does not change the result).  The intermediate field $s_{cov}$ has the properties
\[ (g \cdot s_{cov})(p) = \rho(g)s_{cov}(g^{-1}\cdot p),\quad P^{cov}_{p_r,\pi}s_{cov}=s_{cov},\]
for $g\in P^{1,3}$, where 
\[ ( P^{cov}_{p_r,\pi}s_{cov})(p)=\rho(\sigma(p))P_{p_r,\pi}\rho(\sigma(p)^{-1})s_{cov}(p) \]
is how the projector looks in terms of $s_{cov}$. Moreover, the fields $s'(x)$ on Minkowski space that arise by this construction can be characterised as elements of $C^\infty(M,W)$ such that
\[ (\square + m^{2})s'=0,\quad P^{cov}_{p_r,\pi}s'=s',\]
where $\square=\del_t^2-\sum_i\del_i^2$ is the wave operator, i.e. waves on $M$ with values in $W$ and a constraint for the projection to $V_{p_r,\pi}$. Choosing a different embedding $W'$ of $V_{p_{r}, \pi}$ will give a different set of fields obeying the Klein-Gordon equation that describe the same physical particle as an element of $\Gamma(E)$. Note that $P^{cov}_{p_r,\pi}$ here is no longer given point wise by $W_\rho\to V_{p_r,\pi}$ (that was true only on $L/L_{p_r}$) and appears as a differential operator constraint on the field. The explanation here is to use the identification (\ref{eq_coset_manifold_sheet}) under which the measure $\nu$ transfers (up to a normalisation) to the measure $\frac{d^{3}p_{i}}{(2\pi)^{3}(2E(p))}$, where  $E(p) = p^{0}= \sqrt{m^{2} + \sum_{i}p_{i}^{2}}$ along with the spatial momenta $p_i$ are used as the coordinates of $H^{1,3}$. 
 \begin{equation} \label{eq_mass_shell_FT} s'(x) \propto  \int_{H^{1,3}} \frac{d^{3}p_{i}}{(2\pi)^{3}(2E(p))} e^{-ix \cdot p}s_{cov}(p) =
\int_{M^*} \extd^4 p\,  \delta(p^{2} - m^{2})|_{p^{0}>0}\, e^{-ix \cdot p}s_{cov}(p),\end{equation}
which we recognise as an inverse Fourier transform, from $M^*$ to $M$. Under this Fourier isomorphism, the condition that $p^2=m^2$ appears as the wave equation on $M$. More details can be found in works such as \cite{Fig}.

\subsection{Recap of irreps of $D(G)$} \label{sec_irreps_of_DG}

Working over $k=\C$, since $D(G)$ is semisimple, its representation theory is fully described by its irreducible modules. By explicitly constructing the central primitive idempotents of $D(G)$ (e.g. see \cite{Ma:cla}), then it is seen that these are classified by pairs $(\CC, \pi)$ with $\CC\subset G$ a conjugacy class and $\pi$ an irreducible representation of the centraliser subgroup $C_G$ of a fixed element $r\in\CC$, i.e. \[ C_{G}:= \{ n \in G \, | \, nr=r n \}\subseteq G\, . \] We denote by $V_\pi$ the carrier space of $\pi$. The choice of $r$ will not change the resulting $D(G)$ representation up to isomorphism (as the $C_{G}$ are isomorphic by conjugation, and so the irreps are equivalent) and nor does the following further choice of a map
\begin{equation}\label{qC} q:\CC\to G,\quad q_c r q_{c}^{-1}=c,\quad \forall c\in \CC,\end{equation}
where choose to fix in particular $q_{r}=e$.  This gives a group cocycle $\zeta:\CC\times G\to C_G$ respectively defined and characterised by \begin{equation}\label{def_cocycle} \zeta_c(g)=q_{gcg^{-1}}^{-1} g q_{c};\quad \zeta_c(gh)=\zeta_{hch^{-1}}(g)\zeta_c(h), \end{equation}
for all $c\in\CC$ and $g,h\in G$. Any irrep of $D(G)$ is then of the form $V_{\CC, \pi}= \C G \otimes_{\C C_{G}} V_{\pi}\cong \C  \CC\tens V_\pi$, since $\{ q_{c}\}_{c \in \CC}$ provides a canonical choice of a full set of representatives for the left coset space $G/C_G$, and so 
\begin{equation}\label{GCG} G/C_G\isom \CC,\end{equation}
 by sending $c$ to the equivalence class of $q_{c}$. By \eqref{qC}, the space $\C \CC\tens V_\pi$ has quantum double action 
\begin{equation}\label{Cpiirrep} h\la (c\tens w)=hch^{-1}\tens \pi(\zeta_c(h))w,\quad \delta_g\la (c\tens w)=\delta_{g,c}c\tens w.
\end{equation}
The action of $h\in G$ is simply the action of the induced representation $\C G \otimes_{\C C_{G}}V_{\pi}$ from the irrep $V_{\pi}$ of $C_{G}$, while the action of $\delta_g$ can be viewed as a grading $ |c\tens w|=c$ upgrading $V_{\CC, \pi}$ to a $G$-crossed module. In calculations it will also be useful to note that \begin{equation} \label{eq_partition}G= G/C_G\times C_G, \quad g = q_{grg^{-1}}\zeta_{r}(g) \mapsto q_{grg^{-1}} \times \zeta_{r}(g), \end{equation} as sets (i.e. just the partitioning of the group by the cosets). Equivalently, there is a unique factorisation
\[ G=q_\CC. C_G;\quad q_\CC=\{ q_{c}\}_{c \in \CC}.\]

\begin{example}\label{exSn} \rm  We take $\CC\subset G$ to be the 2-cycles conjugacy class in $G=S_n$ the group of permutations of $\{1,\cdots,n\}$. We let $r=(n-1,n)$, and as such $C_G=\Z_2\times S_{n-2}$ where $S_{n-2}$ permutes $\{1,\cdots,n-2\}$ and $\Z_2=\{e,(n-1,n)\}$. These elements centralise $r$ and then counting the elements tells us that this is all due to the unique factorisation $S_n=q_{\CC} C_G$. It is natural to decompute $\CC$ as
\[ \CC=\{r\}\amalg \{\bar i=(i,n-1), \underline{i}=(i,n)\}\amalg \{(ij)\in S_{n-2}\}.\]
Then a natural choice of $q_c$ is 
\[ q_r=e,\quad q_{\bar i}=\underline i,\quad q_{\underline i}=\bar i,\quad q_{(ij)}=\bar i \underline j;\quad i<j.\]
Note that all these elements are involutive as $\bar i$ and $\underline j$ commute. 

To describe the cocycle it is enough to give its restrictions to $\CC\times q_\CC$ and $\CC\times C_G$ and then use $\zeta_c(q_b n)=\zeta_{ncn^{-1}}(q_b)\zeta_c(n)$.

For $\zeta|_{\CC\times C_G}$, we have 
\[ \zeta_r(n)=n,\quad \zeta_{\bar i}(r)=\zeta_{\underline{i}}(r)=r,\quad \zeta_{(ij)}(r)=(ij),\]
hence we only need to know $\zeta_c(s)$ for $s\in S_{n-2}$ and can then use $\zeta_c(r s)=\zeta_{scs^{-1}}(r)\zeta_c(s)$. Similarly, as 2-cycles $(ij)$ generate $S_{n-2}$, it is enough to know the $\zeta_c((ij))$. These are
\[     \zeta_r((ij))=\zeta_{\bar k}((ij))=\zeta_{\underline k}((ij))=(ij),\quad \forall i,j,k,\]
where $k,l$ are assumed disjoint from $i,j$, while for the action of 2-cycles in $S_{n-2}$ we find $\zeta_{(kl)}((ij))=(ij)$  {\em except} for  three  cases computed as follows under the {\em convention that all 2-cycles are in taken standard form}, which here means that $k<l, i<j$. These are
\[ k=i,\ l<j;\quad \zeta_{(il)}((ij))=r (ilj),\quad l=j,\ i<k;\quad \zeta_{(kj)}((ij))=r (ijk),\]
\[k=i,\ l=j;\quad \zeta_{(ij)}((ij))=r.\]

 For $\zeta|_{\CC\times q_\CC}$, it is enough to know $\zeta_c(\bar i), \zeta_c(\underline{i})$ given that $\zeta_c(e)=e$ and $\zeta_c(\bar i \underline j)=\zeta_{\underline j c \underline j}(\bar i)\zeta_c(\underline{j})$. Now, $\zeta_r(q_b)=q^{-1}_{q_b r q_b^{-1}}q_b q_r=q^{-1}_b q_b=e$ for all $b\in\CC$, so in particular $\zeta_r(\bar i)=\zeta_r(\underline i)=e$. We also computed 
\[  \zeta_{\bar i}(\bar j)=\zeta_{\underline j}(\underline i)=\begin{cases}e & i>j\\ r & i=j\\ r (ij)& i<j\end{cases},\quad \zeta_{\bar i}(\underline j)=\zeta_{\underline i}(\bar j)=\begin{cases}e & i=j\\ (ij) &{\rm else}\end{cases}, \]
and in conventions where $(ij)$ is assumed in standard form with $i<j$, 
\[ \zeta_{(ij)}(\bar j)=\zeta_{(ij)}(\underline i)=r(ij),\quad \zeta_{(ij)}(\underline j)=\zeta_{(ij)}(\bar i)= e,\]
\[ \zeta_{(ij)}(\bar k)=(ik),\quad \zeta_{(ij)}(\underline k)=(jk), \]
for $k\ne i,j$. We see that in all cases the output of $\zeta$ is in $C_G$. 

This cocycle is then used to construct the representation $V_{\CC,\pi}$ associated to any irrep of $C_G$ according to (\ref{Cpiirrep}). These irreps are given by an irrep $\pi_0$ of $S_{n-2}$ and a choice of sign for the irrep of $\Z_2=\{e,r\}$, and so we have $\pi_\pm(s)=\pi_0(s)$,  $\pi_\pm(r s)=\pm\pi_0(s)$. For example, for $n=5$, $|\CC|=10$ and for this choice of $\CC$ there are 6 choices for $\pi$, leading to $V_{\CC,\pi}$ either 10-dimensional or 20-dimensional as $\pi_0$ has 3 choices: the trivial representation, the sign representation and a standard 2-dimensional one contained in the action of $S_3$ on $\C^3$ by permuting the bases.
\end{example}

\begin{example}\label{exS3}\rm For $G=S_3$, a full compilation of the resulting irreps of $D(S_3)$ over $\C$, equivalent to that in \cite{CowMa}, is as follows.

(i) For $\CC=\{e\}$ we have $r=e$ and $C_G=S_3$ and three irreps of this, one 2-dimensional and two of them 1-dimensional. Clearly $q_e=e$, $\zeta_e=\id$ and the action of $D(S_3)$ on the same vector spaces is by the given action of $S_3$ and $f\in \C(S_3)$ acting by $f(e)$. 

(ii) For $\CC=\{uv,vu\}$ where $u=(12), v=(23)$ and we also use $w=(13)$, we set $r=uv=(123)$ and have $C_G=\Z_3=\{e,r,r^2\}$ where $r^2=vu=(132)$. This has three irreps each dimension 1. We set $q_r=e$ and $q_{r^2}=u$ then $uru=vu=r^2$ as needed. Then 
\[ \zeta_{uv}(\{e,u,v,w,uv,vu\})=\{e,e,r,r^2,r,r^2\},\quad \zeta_{vu}(\{e,u,v,w,uv,vu\})=\{e,e,r^2,r, r^2,r \}.\]
This gives 3 of the irreps of $D(S_3)$ each of dimension 2. These correspond to the three irreps $\pi_j$ of $C_G$ for $j=0,1,2$ where  $\pi_j(r)=q^j$ and $q=e^{2\pi\imath\over 3}$ is a cube root of unity. Then the action of $D(S_3)$ is 
\[ u\la\{ uv,vu\}=\{vu,uv\},\quad  v\la \{uv,vu\}=\{ q^j vu, q^{-j} uv\},\quad f\la \{uv,vu\}=\{f(uv)uv,f(vu)vu\}.\]
The restriction to $\C S_3\subset D(S_3)$ is the 2-dimensional irrep of $S_3$ for $i=1,2$ but ${\rm sign}\oplus 1$  for $i=0$. This follows from the character table and the trace of the matrix for $u$ being 0, that for $uv$ being -1 when $i=1,2$ and $2$ when $i=0$. 

(iii) For $\CC=\{u,v,w\}$ we let $r=v=(23)$ as for $S_n$ above and have $C_G=\Z_2=\{e,r\}$ with two irreps, each 1-dimensional with $\pi_\pm(r)=\pm 1$. Here the general construction reduces to $q_v=e$, $q_u=w$, $q_w=u$. The cocycle is 
\[  \zeta_v(\{e,u,v,w,uv,vu\})=\{e,e,r,e,r,r\},\]
\[ \zeta_u(\{e,u,v,w,uv,vu\})=\{e, r,r,e,r,r\},\quad  \zeta_w(\{e,u,v,w,uv,vu\})=\{e,e,r,r,e,r \}.\]
This gives 2 of the irreps of $D(S_3)$, each 3-dimensional. The action of $D(S_3)$ is then
\[ u\la \{u,v,w\}=\{\pm u,w,v\},\quad v\la \{u,v,w\}=\{\pm w,\pm v,\pm u\},\quad f\la \{u,v,w\}=\{f(u)u,f(v)v,f(w)w\},\]
on the basis vectors for the two cases. The restriction to $\C S_3$ gives the $1\oplus 2$ representation of $S_3$ for $\pi_+$ and ${\rm sign}\oplus 2$ for $\pi_-$, again using the character table.
\end{example}

\subsection{The quantum homogeneous space $P\twoheadrightarrow H$} \label{sec_homogeneous_bundle_DG}

In the geometric picture we work in the dual `coordinate algebra' picture. So for the total space coordinate algebra, we take the dual Hopf algebra to $D(G) = k(G) \lcross kG$ but we do so in the slightly nonstandard form $P=D^{\vee}(G) := k(G) \rcocross kG$ defined by the non-standard duality pairing with $D(G)$
\begin{equation} \label{eq_new_pairing} \<\delta_{g} \otimes h,  \delta_{u} \otimes v\>= \delta_{g, v^{-1}} \delta_{h, u^{-1}} \, , \end{equation} where therefore the isomorphism with the linear dual Hopf algebra (Section~\ref{secpre}) is \begin{equation}\label{eq_iso_beta} \beta: kG \lcocross k(G)=(D(G))^{*} \cong D^{\vee}(G), \quad h \otimes \delta_{g} \mapsto \delta_{g} \otimes g^{-1}hg\, . \end{equation}
 $D^{\vee}(G)$ thus has the tensor product algebra and counit on $k(G) \otimes kG$, with 
 \begin{equation}\label{DvG} \Delta(\delta_{g} \otimes h) = \sum_{f \in G} \delta_{f} \otimes f^{-1}ghg^{-1}f \otimes \delta_{f^{-1}g} \otimes h\, ,  \quad S(\delta_{g} \otimes h) = \delta_{g^{-1}} \otimes gh^{-1}g^{-1} , \end{equation}  i.e. has semi-direct coproduct generated by the $k(G)$-coaction on $kG$ given by the coadjoint coaction \begin{equation} \label{eq_coadjoint_coaction} \Delta_{R} h = \sum_{g \in G} ghg^{-1} \otimes \delta_{g}, \quad h \in G.\end{equation} Its coquasitriangular structure is then given by \begin{equation} \label{eq_coquasi_DveeG}\mathcal{R}(\delta_{g} \otimes h \otimes \delta_{u} \otimes v) = \delta_{g, e}\delta_{h, u} .\end{equation} $D^\vee(G)$-comodules are again Whitehead $G$-crossed modules where, given one of the latter, the right $D^{\vee}(G)$ coaction is \begin{equation} \label{eq_comodule_DveeG} \Delta_R v= \sum\nolimits_{g \in G} g\la v\tens \delta_g \otimes |v|,\end{equation} dual to the action \eqref{eq_DG_action}, i.e. $\Delta_{R}v = \sum_{g, h \in G} (\delta_{g} \otimes h) \triangleright v \otimes \beta(g \otimes \delta_{h})$. We note that by  cosemisimplicity of the Hopf algebra, any comodule decomposes as a direct sum of comodules $V_{\CC, \pi}$ (with coaction `dualising' \eqref{Cpiirrep}). The constructions are general but we proceed with $k=\C$. 
 
We are now ready to give the geometric `Wigner construction' of coirreps of $D^{\vee}(G)$ following the template of Section~\ref{secW}. Given the semidirect product form, the conceptual dictionary that connects back to that is
\[ C^\infty(P^{1,3})\leftrightarrow P=D^{\vee}(G),\quad C^\infty(T)\isom C^\infty(M)\leftrightarrow \C  G,\quad C^\infty(L)\leftrightarrow \C (G).\]
In the role of $\alpha\in T^*$ labelled by $p_r\in H^{1,3}\subset M^*$, we now take co-irrep of $\C  G$, which is simply a 1-dim graded vector space, i.e. $ \C $ along with a choice of 
\[ r \in \CC\subseteq G,\]
and $\C  G$-coaction then given by $\Delta_{R}\lambda = \lambda \otimes r$, $\lambda \in \C $. Here $\CC$ is the orbit of $r$ in $G$ under the $G$-action (in our case the action is the adjoint one, which is indeed the dual action to the coaction \eqref{eq_coadjoint_coaction} defining the semi-direct product, so this is a conjugacy class) in the role of $H^{1,3}$. Then the isotropy subgroup $C_G$ plays the role of $L_{p_r}$ and for the fibre of one of our bundles, we will work with its algebra of functions and the subgroup property is now a Hopf algebra surjection
\begin{equation} \label{eq_surjection_kG_kCG} \C (G)\twoheadrightarrow \C (C_G),\end{equation}
given by restricting to a function to $C_G$. This is equivalent to multiplication by the characteristic function $\chi_{C_G}$ of $C_G$ as a projection in $\C(G)$, which in turn is given by
\[ \chi_{C_G}=(\delta_r\tens \id)\Delta_R r=\sum_{g\in G}\delta_r(g r g^{-1})\delta_g  = \sum_{g \in C_{G}} \delta_{g},\]
using \eqref{eq_coadjoint_coaction}. This is  1 on $C_G\subset G$ and zero elsewhere. 

Next, the analog of $C^\infty(T\lcross L_{p_r})$ is now the Hopf algebra $H:= \C (C_{G}) \rcocross \C G$ which comes with a surjection
$P\twoheadrightarrow H$ to give us a homogenous quantum $H$-principal bundle. The induced $H$-coaction on $P$ here is \[ \Delta_R(\delta_g\tens h)= \sum_{f^{-1}g\in C_G} \delta_{f} \otimes f^{-1}ghg^{-1}f \otimes \delta_{f^{-1}g} \otimes h\, , \]
so we get that the fixed point subalgebra $A$ is spanned by $\{c_g\delta_g\tens e\ | g\in G\}$ such that 
\[\sum_{a\in C_G}  \sum_g c_g  \delta_{g a^{-1}}\tens e\tens\delta_{a}\tens 1= \sum_g c_g \delta_{g}\tens e\tens 1\tens 1, \]
which happens precisely when $c_g \delta_g\in \C (G)$ is invariant under right translation by $C_G$, i.e in. $\C (G)^{C_G}= \C (G/C_G)$. Thus, (dropping $\tens e$) we have identified\begin{equation*} \label{eq_A_kC} A= P^H\isom \C (G/C_G) \end{equation*} as the base of the bundle. With respect to its natural basis given according to (\ref{GCG}) by $\{ \delta_{q_{c}}\}_{c \in \CC}$, the algebra $\C (G/C_G)$ thus sits in $P$ by its canonical inclusion in $\C (G)$ which in turn is included in $P$, and so is overall by
\begin{equation}\label{eq_AinDvee} \delta_{q_{c}} \mapsto \sum_{n \in C_{G}}\delta_{q_{c}n} \otimes e \in D^{\vee}(G).\end{equation}
 This bundle then has a $D^{\vee}(G)$ homogeneous structure induced by the left coaction $\Delta_{L} = \Delta$ on the total space $D^{\vee}(G)$. That we have a Hopf-Galois extension follows naturally from \cite[Thm.~5.9]{BegMa}. Finally, this bundle is in fact trivial by the following isomorphism of algebras \begin{equation} \label{eq_trivialisation_map} P\isom \C (G/C_G)\tens \C (C_G)\rcocross \C  G=A\tens H, \quad \delta_{g} \otimes h \mapsto \delta_{q_{grg^{-1}}} \otimes \delta_{\zeta_{r}(g)}\otimes h,\end{equation} with the tensor product algebra on the RHS (and inverse $\delta_{q_{c}g} \otimes h\mapsfrom \delta_{q_{c}} \otimes (\delta_{g} \otimes h)$). This used (\ref{eq_partition}) so that $\C (G)=\C (G/C_G)\tens \C (C_G)$ as algebras.
 
\subsection{Geometric construction of irreps $V_{r,\pi}$ of $H$} \label{sec_geometric_construct_irreps}

Continuing with the set-up  above, we now additionally choose a right coirrep $V_{\pi}$ of $\C (C_{G})$ and  set  $V_{r,\pi}=V_{\pi} \otimes \C  \cong V_{\pi}$ as vector spaces but now with the $H$-coaction
\begin{equation} \label{eq_comoduleVpi} \Delta_{R}v = \sum_{g \in C_{G}} g \triangleright v \otimes \delta_{g} \otimes r \, , \quad v \in V_{r, \pi}, \end{equation}where we wrote the right $\C (C_{G})$-coaction $v^{1} \otimes v^{2} = \sum_{g \in C_{G}}g \triangleright v \otimes \delta_{g}$ in terms of the equivalent left $C_{G}$-action. The $A$-bimodule of sections of its associated bundle is then \begin{equation} \label{eq_def_bundle1} E=(D^{\vee}(G)\tens V_{r, \pi})^{\C (C_{G}) \rcocross \C G}.\end{equation}  Moreover, $E$ has a left $D^{\vee}(G)$-comodule structure coming from the coproduct of $D^{\vee}(G)$. 

\begin{lemma} \label{lemma_Vcpi_geometrical} $E$ is coirreducible and all coirreducible $D^{\vee}(G)$-comodules are obtained in this way. As a corresponding left $D(G)$-module, $E\cong V_{\CC, \pi}$ (the irrep described in Section~\ref{sec_irreps_of_DG}).\end{lemma} \noindent \textit{Proof}: As mentioned in Section~\ref{secpre}, combining the isomorphism $((A \otimes H)\tens V_{r, \pi})^H\cong A \otimes V_{r, \pi}$ (by $\sum_{g\in C_{G}}\delta_{q_{c}} \otimes S^{-1}(\delta_{g} \otimes r) \otimes g \triangleright v \mapsfrom \delta_{q_{c}} \otimes v$) with the triviality of the principal bundle \eqref{eq_trivialisation_map} gives the following isomorphism of associated bundles \begin{equation} \label{eq_trivial_associated_bundle_H}\Theta_{H}^{V_{r, \pi}}: (D^{\vee}(G) \tens V_{r, \pi})^H \cong \C (G/C_G) \otimes V_{r, \pi},  \sum_{g \in C_{G}} \delta_{q_{c}g^{-1}} \otimes r^{-1} \otimes g \triangleright v \mapsfrom \delta_{q_{c}} \otimes v \, . \end{equation}Asking this to be a $D^{\vee}(G)$-comodule map induces the left-coaction on $\C (G/C_G) \otimes V_{r, \pi}$ 
\begin{align} \Delta_{L}&(\delta_{q_{c}} \otimes v)  = (\id \otimes \Theta_{H}^{V_{r, \pi}})\circ \sum_{g \in C_{G}} \Delta(\delta_{q_{c}g^{-1}} \otimes r^{-1}) \otimes g \triangleright v \nonumber\\ 
& = (\id \otimes \Theta_{H}^{V_{r, \pi}})\sum_{f \in G, g \in C_{G}} \delta_{f} \otimes f^{-1}q_{c}g^{-1}r^{-1}gq_{c}^{-1}f \otimes \delta_{f^{-1}q_{c}g^{-1}} \otimes r^{-1} \otimes g  \triangleright v)\nonumber\\
& = (\id \otimes \Theta_{H}^{V_{r, \pi}})\sum_{f \in G, g \in C_{G}} \delta_{f} \otimes f^{-1}c^{-1}f \otimes \delta_{q_{f^{-1}cf}(g\zeta_{c}(f^{-1})^{-1})^{-1}} \otimes r^{-1} \nonumber\\
&\qquad\qquad\qquad\qquad\qquad\qquad \otimes g \zeta_{c}(f^{-1})^{-1}\triangleright (\zeta_{c}(f^{-1}) \triangleright v) \nonumber\\
\label{eq_left_Dvee_coaction} & = \sum_{f \in G} \delta_{f} \otimes f^{-1}c^{-1}f \otimes \delta_{q_{f^{-1}cf}} \otimes \zeta_{c}(f^{-1}) \triangleright v . 
\end{align}
The corresponding left $D(G)$-action on $\C (G/C_G) \otimes V_{r, \pi}$ is then (recalling \eqref{eq_iso_beta}, which gives $\beta^{*}: (D^{\vee}(G))^{*} \cong D(G)$ by $g \otimes \delta_{h} \mapsto \delta_{ghg^{-1}} \otimes g$) \[ (\delta_{g} \otimes h) \triangleright (\delta_{q_{c}} \otimes v) = \langle (\beta^{-1})^{*}(\delta_{g} \otimes h), S(\delta_{q_{c}} \otimes v)^{0} \rangle \, (\delta_{q_{c}} \otimes v)^{1} = \delta_{g, hch^{-1}}\delta_{q_{hch^{-1}}} \otimes \zeta_{c}(h) \triangleright v.\] 
Lastly, we identify the vector space $\C (G/C_G)$ with $\C  \CC$ by identifying $\delta_{q_c}$ with $c\in \CC$.
\qed

The above results work over any field and are in line with a Clifford correspondence approach in \cite{Bur}. Over $\C$ we have further that $\C(G)$ and $\C G$ are dually paired as Hopf-$\star$-algebras by the canonical $\star$-structures (where $\star$ on $\C(G)$ is pointwise complex conjugation and on $\C G$ is group inversion on the basis $G$. As a result, $D(G)$ is a Hopf-$\star$-algebra given by 
\begin{equation} \label{eq_star_struc_DveeG}(\delta_{g} \otimes h)^{\star} = (1 \otimes h^{\star})(\delta_{g}^{\star} \otimes e) = \delta_{h^{-1}gh} \otimes h^{-1}.\end{equation}
Moreover, in a Hopf $\star$-algebra setting, a module $V$ of a Hopf-$\star$-algebra $H$ is `unitary' if there exists an sesquilinear inner product $(\ , \ )_{V}$ on $V$ for which 
\[ ( h \triangleright v, w )_{V} = ( v, h^{\star} \triangleright w)_{V}, \]
 for all $v, w \in V$ and $h \in H$. It is shown in  \cite[Lemma 4.1]{Gou} that any $D(G)$ module (or $D^{\vee}(G)$ comodule) is unitary with respect to some sesquilinear inner product (this is done by picking a sequilinear inner product and averaging over the action of the natural basis of $D(G)$). The same applies for $\C(G),\C G$. 

\subsection{Transfer to bundle $\boldsymbol{E'}$ over $\boldsymbol{\C  G}$}  \label{sec_second_bundle_construct}

It is well-known that both $\C  G$ and $\C (G)$ have unique normalised translation-invariant measures  \begin{equation}\label{eq_measure_kofG} \int_{\C  G}: \C  G \to \C  , \quad g \mapsto \delta_{g, e}, \quad \int_{\C  (G)}: \C  (G) \to \C  \quad f \mapsto \frac{1}{|G|}\sum_{g \in G}f(g),\end{equation} 
and that these establish a Fourier isomorphism $\mathcal{F}: \C G\to \C (G)$ which (in our conventions) sends 
\[ \CF(g)=(\int_{\C  G} g \sum_{h\in G} h)\delta_h= \delta_{g^{-1}}, \]
where $\sum_{h\in G} h\tens\delta_h$ plays the role of $e^{\imath x.p}$ for Fourier transform on $\R^{1,3}$, see \cite{Ma, BegMa}. So, the inverse FT is given by $\delta_{g} \mapsto g^{-1}$. We have explained how $G$ plays the role of momentum space, and that the orbit $\CC$ of $r$ under conjugation plays for $D(G)$ the role of $H^{1,3}$ in the Wigner construction. Under the inverse FT, $g^{-1}$ thus plays the role of a `plane wave $e^{-ix\cdot p}$' with fixed momentum $g \in G$. 

Moreover, the base algebra $A = \C (G/C_G)$ can be identified with $\C (\CC)$ since $c\in \CC$ can be identified with $[q_c]\in G/C_G$ in view of (\ref{GCG}). If we use the second form then the analogue of the  `mass-shell Fourier transform' from  functions on $H^{1, 3}$ to functions on position space is to extend a function $f$ on $\CC$ by zero to a function $\widetilde{f}$ of $G$ and apply the inverse finite group Fourier transform,
 \begin{equation} \label{eq_mass_shell_FT_hopf} \, \mathcal{F}^{-1}_{\CC}: \C (\CC) \to \C G\, , \quad f \mapsto  \mathcal{F}^{-1}(\widetilde{f}) = \sum_{c \in \CC} f(c) c^{-1},   \end{equation}   i.e. in analogy to a `sum of plane waves with momenta in $H^{1, 3}$', and coefficients $f(c)$. Here the role of exponential in the classical Fourier transform is given by $\sum_{c\in \CC} c^{-1}\tens \delta_c$ and the measure on $\CC$ is the counting one. If we revert back to $A=\C (G/C_G)$ then these become
 \[ \sum_{c \in \CC} c^{-1} \otimes \delta_{q_{c}} \in \C G \otimes A,\quad \int_{\C (G/C_G)} f= \frac{|C_{G}|}{|G|}\sum_{c \in \CC} f(q_{c}),\]
where  the latter is just the restriction of \eqref{eq_measure_kofG} to $\C (G/C_G)\subset \C (G)$. As previously mentioned, in this picture the analogue of the `spin' is then the parameters classifying the coirreps $V_{\pi}$ of $C_{G}$.

Next, $D^{\vee}(G)$ also forms a quantum principal bundle with base space given by the underlying algebra $\C G$ and structure group given by  $\C (G)$ in the role of functions on the Lorentz group. We have the usual surjection of Hopf algebras $D^{\vee}(G) \twoheadrightarrow \C (G)$ and the induced coaction is \begin{equation} \label{eq_coaction_DG_kG}\Delta_{R}(\delta_{g} \otimes h) = \sum_{f \in G} \delta_{f} \otimes f^{-1}ghg^{-1}f \otimes \delta_{f^{-1}g}\, , \end{equation} which indeed gives the base space (recalling the map $\beta$ of \eqref{eq_iso_beta}) as 
\begin{equation} \label{eq_inclusion_A_second_bundle} \Phi_{\C (G)}: D^{\vee}(G)^{\C (G)} \cong \C G\, , \quad \beta(g \otimes 1) \mapsfrom g \, .  \end{equation} 
Again, this bundle is $D^{\vee}(G)$-homogeneous by the left coregular coaction \eqref{DvG}, and is trivial by the algebra isomorphism $\beta^{-1}: D^{\vee}(G) \cong \C G \otimes \C (G)$. That we indeed have a Hopf-Galois extension again follows from \cite[Thm.~5.9]{BegMa}. For this bundle, the  left regular coaction of $D^\vee(G)$ on itself restricts to the base $\C  G$ 
\begin{equation} \label{eq_kG_DveeG_coaction_symmetry}  \Delta_{L}g = \sum_{f \in G} \delta_{f} \otimes f^{-1}gf \otimes f^{-1}gf \end{equation}
as a comodule algebra, in role of the action of the Poincar\'e group in spacetime. Next, given a right $\C (G)$-comodule $W$ written as $\Delta_R w= \sum_{f \in G} f \triangleright w \otimes \delta_{f}$ in terms of an equivalent left $G$-module structures, we can construct the $\C G$-bimodule of sections of the associated bundle as
\begin{equation} \label{eq_trivial_associated_bundle_Hprime} \Theta_{\C (G)}^{W}: E':=  (D^{\vee}(G) \tens W)^{\C (G)} \cong \C G \otimes W,  \quad \sum_{f \in G} \delta_{f^{-1}} \otimes fgf^{-1} \otimes f \triangleright w \mapsfrom g \otimes w,\end{equation} for the associated vector bundle, where for the latter form, the triviality of the principal bundle is again used. $E'$ will again inherit a left coaction of  $D^{\vee}(G)$ from its coproduct, which, as in the proof of Lemma~\ref{lemma_Vcpi_geometrical}, is
 \begin{equation} \label{eq_coact_Eprime} \Delta_{L}(g \otimes w) = \sum_{f \in G} \delta_{f} \otimes f^{-1}gf \otimes f^{-1}gf \otimes f^{-1} \triangleright w \end{equation}on $\C  G\tens W$ as an extension of (\ref{eq_kG_DveeG_coaction_symmetry}).

 We now want to follow the same steps as in Section~\ref{secWigner} for the Poincar\'e group to pass from sections of the associated bundle $E$ of \eqref{eq_def_bundle1} defined over $\C (G/C_G)$ to sections of the above bundle $E'$ defined over our `noncommutative spacetime' $\C G$, now using the Hopf algebraic analogue of the averaging map \eqref{eq_averagin_integral} given by \cite[Lemma 5.42]{BegMa}. Slightly adapted to our setting, this  states that given a Hopf algebra $H'$ with normalised left-integral $\int_{H'}$ and invertible antipode, then we have the following averaging map for any two right $H'$-comodules $X, Y$
 \begin{equation} \label{eq_general_passing_sections}\textup{av}: X \otimes Y \to (X \tens Y)^{H'}, \quad x \otimes y \mapsto x^{1'} \otimes y^{1'} \, \int_{H'}y^{2'}S^{2}(x^{2'}). \end{equation}  
Here, $x^{1'} \otimes x^{2'}$ denotes the $H'$-coaction on $x \in X$ and similarly for $Y$. In our context we recall that $H=\C (C_G)\rcocross \C  G$ has a corepresentation $V_{r, \pi}$ where $\pi$ is an irrep of $C_G$ the isotropy subgroup of an initial choice of $r\in \CC$. We use the averaging construction for $H'=\C (G)$ restricted to a subspace of a tensor product.

Following the geometric picture of 
 Section~\ref{secWigner}, we also define the left $D^{\vee}(G)$-comodule given by the sections of the associated bundle 
 \[ E^{W}:= (D^{\vee}(G) \tens W)^H\cong \C (G/C_G) \otimes W, \]
  again by the trivialisation \eqref{eq_trivial_associated_bundle_H}, which carries a left $D^\vee(G)$-coaction as associated to a quantum homogeneous bundle.  Finally, if $s\in  \C (G/C_G) \otimes W$, we define $s_{cov}$ in the same space by
\begin{equation} \label{eq_cov_map} (\delta_{q_{c}} \otimes w)_{cov}=\delta_{q_{c}} \otimes q_{c} \triangleright w. \end{equation}
  The corresponding left $D^\vee(G)$-coaction on $s_{cov}\in \C (G/C_G)\tens W$ is given now by 
 \begin{equation} \label{eq_tensor_product_coact} \Delta_{L}(\delta_{q_{c}} \otimes w) = (q_{c})^{1}S^{-1}w^2 \otimes (q_{c})^{2} \otimes w^{1} \end{equation} as the tensor product of the restriction of the left regular coaction of $D^{\vee}(G)$ to $\C (G/C_{G}) \hookrightarrow D^{\vee}(G)$ by \eqref{eq_AinDvee}, which is
 \[ \Delta_L(\delta_{q_c}) = \sum_{f \in G} (\delta_{f} \otimes e) \otimes \delta_{q_{f^{-1}cf}} \, , \] and the right coaction in $W$ converted to a left one 
\[ \Delta_{L}w = S^{-1}(w^{2}) \otimes w^{1}= \sum_{f \in G} \delta_{f^{-1}} \otimes f|w|^{-1}f^{-1} \otimes f \triangleright w.\] 
  
\begin{lemma} \label{lemma_transferance_of_sections_DG} Assume that  $V_{r, \pi}$ is contained in the pushout of a right $D^{\vee}(G)$-comodule $W$ with coaction written as $\Delta_{R}w = \sum_{g \in G} g \triangleright w \otimes \delta_{g} \otimes |w|$ as in \eqref{eq_comodule_DveeG}. The composition  
\[ E \hookrightarrow E^W {\buildrel {\rm av}\over \longrightarrow} E' \]
can be written explicitly as the following embedding of left $D^{\vee}(G)$-comodules
\[ f \otimes v \mapsto \sum_{c\in \CC} c^{-1} \otimes \Big( \int_{\C (G/C_G)} (f \otimes v)_{cov}^{(1)}\, \delta_{q_c} \Big) (f \otimes v)_{cov}^{(2)}= \frac{|C_{G}|}{|G|} \sum_{c \in \CC}f(q_{c})c^{-1} \otimes q_{c} \triangleright v, \]
for $f\tens v\in  \C (G/C_G)\tens  V_{r, \pi}$ and with result in $ \C G \otimes W$. 
\end{lemma}

\proof  Here, $X = D^{\vee}(G)$ a right $\C (G)$-comodule by \eqref{eq_coaction_DG_kG} and $Y = W$ a $\C (G)$-comodule by the pushforward along $D^{\vee}(G) \twoheadrightarrow \C (G)$ so that $\Delta_{R}w = \sum_{g \in G} g \triangleright w \otimes \delta_{g}$ for $w \in W$. The averaging map $\textup{av}: D^{\vee}(G) \otimes W \to E'$ from  (\ref{eq_general_passing_sections}) is
\begin{equation}\label{DGav} {\rm av}((\delta_{g} \otimes h) \otimes w)={1\over |G|}(\sum_{k \in G} \delta_{gk^{-1}} \otimes khk^{-1} \otimes k \triangleright w),\end{equation}
which combines with \eqref{eq_trivial_associated_bundle_Hprime} to land in $\C G\tens W$ to send 
\begin{align} (\delta_{g} \otimes h) \otimes w \mapsto 
&  \frac{1}{|G|}\Theta_{\C (G)}^{W}(\sum_{k \in G} \delta_{(kg^{-1})^{-1}} \otimes (kg^{-1})ghg^{-1}(kg^{-1})^{-1} \otimes (kg^{-1}) \triangleright (g\triangleright w))\nonumber\\
& = \frac{1}{|G|} ghg^{-1} \otimes g \triangleright w, \label{eq_averaging_map_kofG}\end{align}
and so the map $E^{W} \to E'$ is, using also \eqref{eq_trivial_associated_bundle_H}, given by
\begin{align*}\delta_{q_{c}} \otimes w  \mapsto & \textup{av}(\sum_{g \in C_{G}} \delta_{q_{c}g^{-1}} \otimes r^{-1} \otimes g \triangleright w)= \frac{1}{|G|}\sum_{g \in C_{G}}q_{c}g^{-1}r^{-1}gq_{c}^{-1} \otimes q_{c}g^{-1} \triangleright(g \triangleright w) \\
& = \frac{1}{|G|}\sum_{g \in C_{G}}c^{-1} \otimes q_{c} \triangleright w = \frac{|C_{G}|}{|G|} c^{-1} \otimes q_{c} \triangleright w  \end{align*}
on $\C (G/C_G) \otimes W$. Therefore, the transfer map is simply its restriction to $\C (G/C_G) \otimes V_{r, \pi}$. That this map intertwines coactions \eqref{eq_left_Dvee_coaction} and \eqref{eq_coact_Eprime} is immediate by the fact that the homogeneous structures of the two starting principal bundles are identical (i.e. the coregular left coaction), but it can be quickly verified explicitly. Moreover, $E^{W} \to E'$ is injective (i.e. supposing that $\delta_{q_{c}} \otimes w$ and $\delta_{q_{d}} \otimes w'$ have the same image, then $c = d$ by looking at the first tensor factor in the image, while the second tensor factor gives $q_{c} \triangleright w = q_{c} \triangleright w'$. Acting by $q_{c}^{-1}$ gives the result). So, the transfer map is also injective.  

It remains to give this the interpretation stated (to be in line with the Wigner construction) and to check various properties for this. First, we note that  $\Theta^{W}_{H}: E^{W}\cong \C (G/C_G) \otimes W$ is explicitly given by \begin{equation} \label{eq_ThetaWH} \sum_{n \in C_{G}} \delta_{q_{c}n^{-1}} \otimes n|w|^{-1}n^{-1} \otimes n \triangleright w \mapsfrom \delta_{q_{c}} \otimes w, \end{equation}
 and following the same procedure as for $E$, gives that the left $D^{\vee}(G)$-coaction on $\C (G/C_G) \otimes W$ is
  \begin{equation} \label{s_coact_coact} \Delta_{L}(\delta_{q_{c}} \otimes w) = \sum_{f \in G} \delta_{f} \otimes f^{-1}q_{c}|w|^{-1}q_{c}^{-1}f \otimes \delta_{q_{f^{-1}cf}} \otimes \zeta_{c}(f^{-1}) \triangleright w,  \end{equation} (for which $\Theta^{W}_{\C (G)}\textup{av}\circ (\Theta^{W}_{H})^{-1}$ is a left $D^{\vee}(G)$-comodule map). Meanwhile the coaction \eqref{eq_tensor_product_coact} on $s_{cov}$ comes out explicitly as
 \[ \Delta_{L}(\delta_{q_{c}} \otimes w)  = \sum_{f, k \in G} (\delta_{f} \otimes e)(\delta_{k^{-1}} \otimes k|w|^{-1}k^{-1}) \otimes \delta_{q_{f^{-1}cf}} \otimes k \triangleright w = \sum_{f \in G} \delta_{f} \otimes f^{-1}|w|^{-1}f \otimes \delta_{q_{f^{-1}cf}} \otimes f^{-1} \triangleright w, \] from which it is clear that $s\mapsto s_{cov}$  is a left $D^{\vee}(G)$-comodule map. It remains to check 
 \[ \sum_{d \in \CC} d^{-1} \otimes \Big(\int_{\C (G/C_{G})} (\delta_{q_{c}} \otimes v)^{(1)}_{cov}\delta_{q_{d}} \Big) (\delta_{q_{c}} \otimes v)^{(2)}_{cov} =\sum_{d \in \CC} d^{-1} \otimes \Big(\int_{\C (G/C_{G})} \delta_{q_{c}, q_{d}} \delta_{q_{c}} \Big) q_{c} \triangleright v =  c^{-1}\frac{|C_{G}|}{|G|}  \otimes q_{c} \triangleright v,\]
 so that the transfer map can be written as stated. \endproof
 
This also means that the transfer of sections in Lemma~\ref{lemma_transferance_of_sections_DG} coincides up to a normalisation with the `mass-shell Fourier transform' of vector-valued covariant sections:

\begin{corollary}\label{corollary_FT_description_VCpi} If we identify $A\cong \C (\CC)$ then in the latter terms the transfer map $E\to E'$ can be written as the left $D^{\vee}(G)$-comodule embedding
\[ \delta_{c} \otimes v\mapsto {|C_G|\over |G|}(\CF^{-1}_{\CC}\tens\id_W)((\delta_c\tens v)_{cov}),\]
where $(\delta_c\tens v)_{cov}=\delta_{c} \otimes q_{c} \triangleright v$. 
\end{corollary}
\proof This is immediate from the form of $\CF^{-1}_{\CC}$ in \eqref{eq_mass_shell_FT_hopf}. \endproof

It remains to understand the image of the transfer map in order to characterise coirreducible $D^{\vee}(G)$-comodules as certain $W$-valued sections. In describing projectors for this, it is convenient to write them as valued in the Hopf algebra $\C (G) \lcross \C  C_{G}$ that is dually paired as in \eqref{eq_new_pairing} with $H$, and acting by evaluation against a right $H$ coaction. We let $P_\pi\in \C  C_G$ be the projections associated with $V_\pi$. This is explicitly given by \begin{equation} \label{eq_projector_group}P_{\pi} = \frac{\textup{dim}(V_{\pi})}{|C_{G}|}\sum_{n \in C_{G}}\textup{Tr}_{\pi}(n^{-1}) n.\end{equation}
 
 \begin{lemma} \label{lemma_form_of_projectors} {\ }
 
 \noindent (i) Any coirreducible $H$-comodule is isomorphic to one of the form $V_{r, \pi}$ as in \eqref{eq_comoduleVpi}, for some $r \in G$ and $\pi$ an irreducible representation of $C_{G}$. 
 
\noindent (ii)  $P_{r, \pi}:= \delta_{r}\otimes P_{\pi} \in \C (G) \lcross \C C_{G}$ are a full set of centrally primitive idempotents satisfying $\sum_{(r, \pi)}P_{r, \pi} = 1$. 
 
\noindent (iii) The  $D^{\vee}(G)$-comodule $V_{\CC, \pi}$ defined in Section~\ref{sec_irreps_of_DG}, when pushed out to one of $H$, decomposes as $ \bigoplus_{(r, \pi) \, | r \in \CC} V_{r, \pi}$. In particular, $V_{r, \pi}$ is contained in the pushout of a  right  $D^{\vee}(G)$-comodule $W$ iff $W$ contains $V_{\CC, \pi}$ as a sub $D^\vee(G)$-comodule.  

\noindent (iv)  In the setting of Lemma~\ref{lemma_transferance_of_sections_DG} with  $V_{\CC, \pi}\subseteq W$, the  projection $E^{W} \to E$ given by the action of  $\id \otimes P_{r, \pi}$ corresponds to a projection $P_{r, \pi}^{cov}$ such that
\[ ((\id \otimes P_{r, \pi})s)_{cov}= P_{r, \pi}^{cov} s_{cov},\quad P_{r, \pi}^{cov}(\delta_{q_c}\tens w)=\delta_{q_{c}} \otimes (q_{c} P_{r, \pi}q_{c}^{-1}) \triangleright w.\]
\end{lemma}
\proof (i) $H$ is again cosemisimple (since $S^{2} = \id$) hence all its comodules can be decomposed as a direct sum of coirreducible ones. That $V_{r, \pi}$ is coirreducible is more easily seen in terms of the action of its `dual' Hopf algebra $H^{\vee}:= \C (G) \lcross \C C_{G}\subseteq D(G)$,  given by $(\delta_{g} \otimes h) \triangleright v = \langle (\beta^{*})^{-1}, v^{2} \rangle v^{1} = \delta_{g, r} h \triangleright v = \delta_{g, r} \pi(h) v$. Hence,  for an non-zero vector $v \in V_{r, \pi}= V_{\pi}$ as vector spaces,  by the latter's irreducibility we have that $\C C_{G} \triangleright v = V_{\pi} = V_{r, \pi}$. It is also clear that $(\C (G) \lcross \C C_{G})\la v = \C C_{G} \triangleright v = V_{r, \pi}$ as equalities of vector spaces, since $\delta_{r} \in \C (G)$. That all irreps of $\C (G) \lcross \C C_{G}$ are of this form is by simple counting arguments: $|H| = |G||C_{G}|$, and by irreducibility of the $V_{\pi}$ we know that $\sum_{\pi} \textup{dim}(V_{\pi}) = |C_{G}|$, and so indeed $\sum_{(r, \pi)} \textup{dim}(V_{r, \pi}) = \sum_{(r, \pi)} \textup{dim}(V_{\pi})  = \sum_{r \in G} |C_{G}| = |G||C_{G}|$. 

(ii) We note that $\{ P_{\pi}\}_{\pi}$ are a full set of centrally primitive orthogonal idempotents on $\C C_{G}$, and thus give the decomposition $\C G_{C}$-modules into its irreps $\pi$. That $P_{r, \pi}$ forms a full set for $\C (G) \lcross \C C_{G}$ is by direct computation. For $g \in G$, and $h \in C_{G}$
\begin{align*}P_{r, \pi}(\delta_{g} \otimes h) &= \frac{\textup{dim}(V_{\pi})}{|C_{G}|}\sum_{n\in C_{G}} \textup{Tr}_{\pi}(n^{-1})(\delta_{r} \otimes n)(\delta_{g} \otimes h) \\
&=  \frac{\textup{dim}(V_{\pi})}{|C_{G}|}\sum_{n\in C_{G}} \textup{Tr}_{\pi}(n^{-1}) \delta_{r, ngn^{-1}} \otimes nh = \delta_{g, r} \delta_{r} \otimes P_{\pi}h,\\
 (\delta_{g} \otimes h)P_{r, \pi} &=  \frac{\textup{dim}(V_{\pi})}{|C_{G}|}\sum_{n\in C_{G}} \textup{Tr}_{\pi}(n^{-1}) \delta_{g, hrh^{-1}}\delta_{g} \otimes hn = \delta_{g, r} \delta_{r} \otimes hP_{\pi} = \delta_{g, r} \delta_{r} \otimes P_{\pi}h,\end{align*} where we use that $n^{-1}rn = r$ for the first line and $hrh^{-1} = r$ for the second, as well as the fact that $P_{\pi}$ are central in $\C C_{G}$. Orthogonality is clear since 
 \begin{align*}P_{r, \pi}P_{r',\pi'} &=  \big(\frac{\textup{dim}(V_{\pi})}{|C_{G}|}\big)^{2}\sum_{n, m \in C_{G}}\textup{Tr}_{\pi}(n^{-1})\textup{Tr}_{\pi}(m^{-1}) \delta_{r, nr'r^{-1}} \delta_{r} \otimes nm = \delta_{r, r'} \delta_{r} \otimes P_{\pi}P_{\pi'} \\
 &= \delta_{r, r'} \delta_{\pi, \pi'}\delta_{r} \otimes P_{\pi} ,\end{align*} as the $\{ P_{\pi}\}$ are mutually orthogonal. Finally, $\sum_{r \in G}\sum_{\pi}(\delta_{r} \otimes P_{\pi}) = 1 \otimes \sum_{\pi } P_{\pi} =1 \otimes e$. 

(iii) From \eqref{Cpiirrep}, the $D(G)$-action on $V_{\CC, \pi}$ is $(\delta_{g} \otimes h) \triangleright (c \otimes v) = \delta_{g, hch^{-1}} hch^{-1} \otimes \pi(\zeta_{c}(h))v$, which restricts to give the action of $\C (G) \lcross \C C_{G}$ on $V_{\CC, \pi}$. In particular, the action of $P_{r, \pi}$ on $V_{\CC', \pi'}$ is thus given by 
\begin{align*} P_{r, \pi} \triangleright (c \otimes v) &=  \frac{\textup{dim}(V_{\pi})}{|C_{G}|}\sum_{n \in C_{G}} \textup{Tr}_{\pi}(n^{-1})\delta_{r, ncn^{-1}} ncn^{-1} \otimes \pi'(\zeta_{c}(n))v \\
&= \delta_{\CC, \CC'} \frac{\textup{dim}(V_{\pi})}{|C_{G}|}\sum_{n \in C_{G}} \textup{Tr}_{\pi}(n^{-1})\delta_{r, c} c \otimes \pi(\zeta_{r}(n))v = \delta_{r, c} \delta_{\CC, \CC'} c \otimes P_{\pi}\triangleright v \\
&= \delta_{r, c} \delta_{\CC, \CC'} \delta_{\pi, \pi'} c \otimes v , \end{align*} where $\delta_{r, ncn^{-1}}$ enforces $\CC = \CC'$, and as a consequence $\delta_{r, ncn^{-1}} = \delta_{r, c}$. The penultimate equality is since $\zeta_{r}(n) = n$, and using the definition \eqref{eq_projector_group}, as well as the Schur orthogonality relation to get $\delta_{\pi, \pi'}$ (where we note that this can only be applied since we have the support of $\delta_{\CC, \CC'}$, and thus we are now comparing irreps of the same centraliser subgroups). The direct sum decomposition follows. The above also gives the form of the action for any $W$ (by decomposing $W$ into its $D(G)$-irreps $V_{\CC, \pi}$). 

We also note that the $\C G$-module obtained by restriction of the $D(G)$ action on $V_{\CC, \pi}$ is simply the induced representation from $C_{G}$ (i.e. as previously noted below \eqref{Cpiirrep}), and thus this tells us how any $D(G)$-module $W$ restricts to $\C G$. 

(iv) This is immediate. \endproof

We can then map $P_{r, \pi}^{cov}$ over to $\C  (\CC)\tens W$ where it simply appears as \[P^{cov}_{r, \pi}(\delta_{c} \otimes w) = \delta_{c} \otimes (q_{c}P_{r, \pi}q_{c}^{-1}) \triangleright w, \] i.e. we obtain projectors whose action is dependent on the element $c \in \CC$, or in analogy with the Poincar\'e group case, the projector is dependent on the momentum of the plane wave. 

So, we can now describe the $D^{\vee}(G)$-comodule $E$ simply as a `mass-shell' Fourier transform of a subclass of $W$-valued sections on $\C (\CC)$ ($\cong \C (G/C_{G})$) to a subclass of $W$-valued sections on $\mathbb{C}G$, in analogy with Section~\ref{secW} with $\C (\CC)$ in the role of momentum space and $\C G$ in the role of position space. Explicitly, 

\begin{corollary} \label{cor_final_FT_form_of_E} In the setting of Lemma~\ref{lemma_transferance_of_sections_DG}, we have the following isomorphisms of left $D^{\vee}(G)$-comodules 
\[ E \cong \{ (\mathcal{F}^{-1}_{\CC} \otimes \id_{W})(s) \, | \, s \in \ker(\id-P_{r, \pi}^{cov})\subseteq \C (\CC)\tens W \} \subseteq \C G \otimes W, \]
where the RHS is a subcomodule of $E'$ with coaction obtained by restriction from \eqref{eq_coact_Eprime}.
\end{corollary}

\proof For completeness, we note that Lemma~\ref{lemma_form_of_projectors} we have the isomorphism of left $D^{\vee}(G)$-comodules $ E \cong \ker(\id-P_{r, \pi}^{cov})\subseteq \C (G/C_G)\tens W$ where the RHS has the coaction \eqref{eq_tensor_product_coact}. Combining this with the above mapping to $\C (\mathcal{C})$ and with Corollary~\ref{corollary_FT_description_VCpi} (up to a non-zero scalar $\frac{|C_{G}|}{|G|}$), gives the stated result.  \endproof

We also note that under $\mathcal{F}^{-1}_{\mathcal{C}}$, $P^{cov}_{r, \pi}$ in the above becomes defined on $\C \mathcal{C}^{-1} \otimes W$ \begin{equation} \label{eq_projector_on_position_space}P^{cov}_{r, \pi}(c^{-1} \otimes w) = c^{-1} \otimes (q_{c}P_{r, \pi}q_{c}^{-1}) \triangleright w \in \C G \otimes W , \end{equation} i.e. by simply asking $P^{cov}_{r, \pi}\circ (\mathcal{F}^{-1}_{\mathcal{C}} \otimes \id_{W}) = (\mathcal{F}^{-1}_{\mathcal{C}} \otimes \id_{W})\circ P^{cov}_{r, \pi}$. 

Overall, we have recovered the irrep $V_{\mathcal{C}, \pi}$ of $D(G)$ in Section~\ref{sec_irreps_of_DG} as a coirrep of $D^\vee(G)$ starting with sections $E$ of the bundle with base $\C (G/C_G)$ and transferring to a subspace of $E'$ as sections of a bundle over $\C G$, and charactised by (a) lying in a subspace $\C \CC^{-1}\tens W$ (analogous to the mass-shell FT of the Klein-Gordon equation) and (b) invariant under the projection $P_{r, \pi}^{cov}$ (analogous to the mass-shell FT of the additional differential constraints depending on the embedding $W$). The analogues of these free field equations will be discussed in Section~\ref{secPcalc}. 

\begin{example} \label{example_embedding} \rm  A canonical choice for $W$ is  $W = V_{\CC, \pi}$. Then the embedding of $H$-comodules $V_{r, \pi} \hookrightarrow V_{\CC,\pi}$ is just $v \mapsto r \otimes v$, and we have $E' = kG \otimes V_{\CC, \pi}$ with coaction \[ \Delta_{L}(g \otimes (d \otimes v)) = \sum_{f \in G} \delta_{f} \otimes f^{-1}gf \otimes f^{-1}gf \otimes f^{-1}df \otimes \pi(\zeta_{d}(f)) v,\] by \eqref{eq_coact_Eprime}. The covariantisation map \eqref{eq_cov_map} is  given by 
\[\textup{cov}: \delta_{c} \otimes (d \otimes v)  \mapsto \delta_{c} \otimes q_{c} \triangleright (d \otimes v) = \delta_{c} \otimes (q_{c}dq_{c}^{-1} \otimes \zeta_{d}(q_{c}) \triangleright v), \] 
when viewing $k(G/C_{G})$ as $k(\CC)$ in the usual way, and the map $E^{W} \to E'$ becomes 
\[ \delta_{q_{c}} \otimes (d \otimes v) \mapsto \frac{|C_{G}|}{|G|} c^{-1} \otimes q_{c}dq_{c}^{-1} \otimes \zeta_{d}(q_{c}) \triangleright v, \] and so the transference of sections map of Corollary~\ref{corollary_FT_description_VCpi} reduces to  
\begin{equation}\label{eximageE'} \delta_{c} \otimes v \mapsto \frac{|C_{G}|}{|G|} c^{-1} \otimes (c \otimes v) \in E'.\end{equation}
We now look to the content of Corollary~\ref{cor_final_FT_form_of_E} for this particular case, starting with  $q_{c}P_{r, \pi}q_{c}^{-1} = (1 \otimes q_{c})P_{r, \pi}(1 \otimes q_{c}^{-1})= \delta_{c} \otimes q_{c}P_{\pi}q_{c}^{-1}$ and
\[q_{c}P_{\pi}q_{c}^{-1} = \frac{\dim(V_{\pi})}{|C_{G}|} \sum_{n \in C_{G}}\textup{Tr}_{\pi}(n^{-1})q_{c}nq_{c}^{-1} =  \frac{\dim(V_{\pi})}{|C_{G}|} \sum_{m \in C_{G}^{q_{c}}}\textup{Tr}_{\pi^{q_{c}}}(m^{-1})m=:P_{\pi^{q_c}},\]
where $C_{G}^{q_{c}} = q_{c}C_{G}q_{c}^{-1}$ denotes the centraliser of $q_{c}rq_{c}^{-1} = c$ and $\pi^{q_{c}}$ the conjugate irrep to $\pi$. Then 
\[ q_{c}P_{r,\pi}q_{c}^{-1} = \delta_{c} \otimes P_{\pi^{q_{c}}} =: P_{c, \pi^{q_{c}}} \]
 is simply the projector $P_{r, \pi}$ for a different representative $c$ of the conjugacy class $\mathcal{C}$. By a similar calculation to  the proof of Lemma~\ref{lemma_form_of_projectors}(iii), we have
\[P_{c, \pi^{q_{c}}} \triangleright (d \otimes v) = \delta_{c, d}\delta_{\mathcal{C}, \mathcal{C}'}\delta_{\pi^{q_{c}}, \pi'^{q_{c}}} d \otimes v= \delta_{c, d}\delta_{\mathcal{C}, \mathcal{C}'}\delta_{\pi, \pi'} d \otimes v, \]
and hence 
\[P^{cov}_{r, \pi}(\delta_{c} \otimes (d \otimes v)) =  \delta_{c} \otimes P_{c, \pi^{q_{c}}} \triangleright (d \otimes v). \]
This is in the analogy with the Poincar\'e group case where the `covariantised projector' detects  the momentum $p \in H^{1, 3}$ in the first entry and accordingly applies the projector defined at $p_{r}$ shifted in momentum to $p$. Hence, in Corollary~\ref{cor_final_FT_form_of_E}, an element $s = \sum_{c, d\in \mathcal{C}} \lambda^{i}{}_{c, d}\delta_{c} \otimes d \otimes E_{i} \in \C (\mathcal{C}) \otimes (\C \mathcal{C} \otimes V_{\pi})$ lies in $\ker(\id-P_{r, \pi}^{cov})$ iff 
\[ \sum_{c, d\in \mathcal{C}} \lambda^{i}{}_{c, d}\delta_{c} \otimes d \otimes E_{i} = \sum_{c, d\in \mathcal{C}} \lambda^{i}{}_{c, d} \delta_{c} \otimes \delta_{c, d} d \otimes E_{i} = \sum_{c\in \mathcal{C}} \lambda^{i}{}_{c, c} \delta_{c} \otimes c \otimes E_{i},\] 
which implies $s$ is of the form $\sum_{c \in \mathcal{C}}f(c)\delta_{c} \otimes (c \otimes v) \in k(\mathcal{C}) \otimes V_{\CC, \pi}$ (where we suppress sums of tensors here) for some function $f: \CC \to k$. Applying $\mathcal{F}^{-1}_{\mathcal{C}} \otimes \id$ to the above gives that the right-hand side of the result in Corollary~\ref{cor_final_FT_form_of_E} consists of the subset of elements
\[ \sum_{c \in \mathcal{C}} f(c) c^{-1} \otimes c \otimes v  \in E',\]
in agreement with the image of the transfer map (\ref{eximageE'}). Applying $\epsilon$ to the first tensor factor gives a trivial isomorphism of this subset of $E' = kG \otimes V_{\CC, \pi}$ with  $V_{\CC, \pi}$, which is a $D^{\vee}(G)$-comodule map (with coaction on $E'$ above and the coaction on $V_{\CC, \pi}$ inferred from \eqref{eq_left_Dvee_coaction}).
\end{example}

By Lemma~\ref{lemma_form_of_projectors}, any other $W$ must contain $V_{\CC, \pi}$ and the above then agrees with the isomorphism in Corollary~\ref{cor_final_FT_form_of_E} for arbitrary $W$ by firstly decomposing it into its $D^{\vee}(G)$-coirreps.

\subsection{Transfer to bundle $\boldsymbol{E^\star}$ over $\C (G)$} We can also instead pass from sections of the original associated bundle $E \cong V_{\mathcal{C}, \pi}$ to sections of a bundle over $\C (G)$ associated to the principal bundle obtained from the alternative Hopf algebra surjection $D^{\vee}(G) \twoheadrightarrow \C G$ given by $\eps\tens\id$ to project out the left factor. The induced right $\C G$-coaction on the total space $D^{\vee}(G)$ of this principal bundle is simply \begin{equation} \label{eq_coaction_bundle_3}\Delta_R=\id\tens\Delta,\quad \textup{i.e. } \Delta_{R}(\delta_{g} \otimes h) = \delta_{g} \otimes h \otimes h, \end{equation} and indeed this is a trivial tensor product bundle with $D^{\vee}(G)^{\C G} \cong \C (G)$ imbedded as the first factor $\C (G)\tens 1$. Moreover, viewed as a quantum homogeneous bundle again by the coregular coaction on the total space, then  $D^{\vee}(G)$ left coacts on the base by \begin{equation} \label{eq_DveeG_coact_base_kofG} \Delta_{L}\delta_{g} = \sum_{f \in G} (\delta_{f} \otimes e) \otimes \delta_{f^{-1}g}. \end{equation} Now, given the same right $D^{\vee}(G)$-comodule $W$ as in Lemma~\ref{lemma_transferance_of_sections_DG}, we define the associated bundle \[E^{\star} := (D^{\vee}(G) \otimes W)^{\C G}  \] over $\C (G)$. Following the exact same procedure as in Section~\ref{sec_second_bundle_construct},  

\begin{lemma} \label{lem_transfer_to_Lorentz}  The transfer map
 \[ E \hookrightarrow E^{W}\xrightarrow{\textup{av}} E^\star\isom \C (G) \otimes W, \quad f \otimes v \mapsto \sum_{c \in \mathcal{C}, n \in C_{G}} f(q_{c}) \delta_{q_{c}n} \otimes n^{-1} \triangleright v \] 
 is an inclusion of $D^\vee(G)$-comodules. Its image can be described as 
\medskip
 \[ E \cong \textup{ker}(\id - P^{\C (G)}_{r, \pi}) \subseteq \C (G) \otimes W, \quad P^{\C (G)}_{r, \pi}(\delta_{g} \otimes w) =    \delta_{g}\otimes \zeta_{r}(g)^{-1}P_{r, \pi}\zeta_{r}(g) \triangleright w \] 
as $D^\vee(G)$-comodules, where  $g = q_{grg^{-1}}\zeta_{r}(g)$ under the factorisation of $G = \{ q_{c}\}_{c \in \mathcal{C}}\cdot C_{G}$ into its cosets. 
\end{lemma}

\proof For the first part, we note that for the associated vector bundle $E^{\star}$ we take $W$ a $\C G$-comodule by the pushforward coaction \begin{equation} \label{eq_pushforward_coact_3}\Delta_{R} w = w \otimes |w|, \end{equation} and so we have again using the triviality of the principal bundle  \begin{equation} \label{eq_Theta_3} \Theta^{W}_{\C G}: E^{\star} \cong \C (G) \otimes W, \quad \delta_{g} \otimes |w|^{-1} \otimes w \mapsfrom \delta_{g} \otimes w, \end{equation} with inherited left $D^{\vee}(G)$-coaction \begin{equation} \label{eq_coact_3} \Delta_{L}(\delta_{g} \otimes w) = \sum_{f \in G}\delta_{f} \otimes f^{-1}g|w|^{-1}g^{-1}f \otimes \delta_{f^{-1}g} \otimes w . \end{equation} Using the same average map \eqref{eq_general_passing_sections} but with $X = D^{\vee}(G)$ a right $\C G$-comodule by \eqref{eq_coaction_bundle_3} and $Y = W$ a $\C G$-comodule by \eqref{eq_pushforward_coact_3}, and this comes out as \begin{equation}\label{eq_averaging_map_kG}\textup{av}: D^{\vee}(G) \otimes W \to E^{\star}, \quad \delta_{g} \otimes h \otimes w \mapsto \delta_{g} \otimes h \otimes w \int_{\C G} |w|h= \delta_{h, |w|^{-1}} \delta_{g} \otimes h \otimes w .  \end{equation} When restricted to $E^{W}$, this gives a left $D^{\vee}(G)$-comodule map intertwining coactions \eqref{s_coact_coact} and \eqref{eq_coact_3}. By then similarly combining \eqref{eq_trivial_associated_bundle_H} with \eqref{eq_Theta_3}, the transfer map becomes \begin{align*}\delta_{q_{c}} \otimes v \mapsto & \Theta^{W}_{\C G}\circ \textup{av} \circ (\Theta^{V_{r, \pi}}_{H})^{-1}(\delta_{q_{c}} \otimes v) \\
& =\Theta^{W}_{\C G}(\sum_{g \in C_{G}} \delta_{r^{-1}, g|v|^{-1}g^{-1}}\delta_{q_{c}g^{-1}} \otimes r^{-1} \otimes g \triangleright v) = \sum_{g \in C_{G}} \delta_{q_{c}g^{-1}} \otimes g \triangleright v, \end{align*} hence the form in the lemma (recalling $|g \triangleright v| = g|v|g^{-1}$ and noting that $\delta_{|v|, r} = 1$ since $v \in V_{r, \pi}$) and indeed intertwines the coactions \eqref{eq_left_Dvee_coaction} and \eqref{eq_coact_3}, since applying the transfer map to \eqref{eq_trivial_associated_bundle_H} gives \[ \sum_{f \in G, g \in C_{G}} \delta_{f} \otimes f^{-1}c^{-1}f \otimes \delta_{q_{f^{-1}cf}g^{-1}} \otimes g\zeta_{c}(f^{-1}) \triangleright v, \] while coacting by \eqref{eq_coact_3} on the image of the transfer map gives 
\begin{align*} & \Delta_{L}( \sum_{g \in C_{G}} \delta_{q_{c}g^{-1}} \otimes g \triangleright v) = \sum_{f \in G, g \in C_{G}} \delta_{f} \otimes  f^{-1}q_{c}g^{-1}|v|^{-1}gq_{c}^{-1}f \otimes \delta_{f^{-1}q_{c}g^{-1}} \otimes g \triangleright v \\
& =  \sum_{f \in G, g \in C_{G}} \delta_{f} \otimes f^{-1}c^{-1}f \otimes  \delta_{q_{f^{-1}cf}\zeta_{c}(f^{-1})g^{-1}} \otimes g \triangleright v  \\
&= \sum_{f \in G, g' \in C_{G}} \delta_{f} \otimes f^{-1}c^{-1}f \otimes  \delta_{q_{f^{-1}cf}g'^{-1}} \otimes g' \zeta_{c}(f^{-1})\triangleright v, \end{align*} using again $|v| = r$, and relabelling $g' =g\zeta_{c}(f^{-1})^{-1} \in C_{G}$. It is also clear that the transfer map is injective (and split by the linear map $\delta_{g} \otimes w= \delta_{q_{c}n} \otimes w \mapsto \delta_{q_{c}} \otimes n^{-1} \triangleright$ using the factorisation $G = \{ q_{c}\}_{c \in \mathcal{C}}\cdot C_{G}$). 

For the second part, the form of the projector on $\C G \otimes W$ is obtained in the analogous same way from $P_{r, \pi}$, here by asking the projector to commute with the averaging map above in the following sense \[ \Big(\Theta^{W}_{\C G}\circ \textup{av}\circ (\Theta^{W}_{H})^{-1}\Big) \circ (\id \otimes P_{r, \pi})(\delta_{q_{c}} \otimes w) = P^{\C (G)}_{r, \pi} \circ \Big(\Theta^{W}_{\C G}\circ \textup{av}\circ (\Theta^{W}_{H})^{-1})\Big) , \] (where we recall that $\Theta^{W}_{H}$ is given by \eqref{eq_ThetaWH}). This comes out as the requirement that \[ \sum_{n \in C_{G}} \delta_{q_{c}n^{-1}} \otimes n P_{r, \pi} \triangleright w = P^{\C (G)}_{r, \pi}(\sum_{n \in C_{G}} \delta_{q_{c}n^{-1}} \otimes n \triangleright w) ,  \] hence the result (which is to be expected for a projector of $\C G$-comodules, or equivalenty $\C (G)$-modules). The identification with $E$ is immediate. \endproof

As one might expected, if $\pi=1$ so that the action on $V_{r,1}$ is trivial, then the transfer map in part (i) is simply the embedding of $\C (G/C_G)\subset \C (G)$ as in (\ref{eq_AinDvee}), but for a general $\pi$ we also act on $v$ in the sum as shown.  The Poincar\'e group analogue would be to transfer the irreps construction as sections over $L/L_{p_{r}}$ to sections of another bundle over the Lorentz group $L$ with fibre $T$, which is not something normally considered. The advantage of considering $E^*$ is that we now have one bundle over $\C (G)$ and the one we wanted before, $E'$ over $\C G$. The next lemma says that averaging from one to the other is just induced by Fourier transform. 

\begin{lemma}\label{lemFtransfer} If $X$ is a right $\C G$-comodule contained in the pushout of a right $D^{\vee}(G)$-comodule $Y$ then the transfer of sections map \[ \C (G) \otimes X \cong (D^{\vee}(G) \otimes X)^{\C G} \hookrightarrow (D^{\vee}(G) \otimes Y)^{\C G} \xrightarrow{\textup{av}}  (D^{\vee}(G) \otimes Y)^{\C (G)} \cong \C G \otimes Y\] is given by $f\otimes x \mapsto (\mathcal{F}^{-1}\otimes \id)(\sum_{g \in G}f(g)\delta_{|g \triangleright x|} \otimes g \triangleright x)$. 
\end{lemma}
\proof We note that the $\C G$-coaction on $X$ is just $\Delta_{R} x = x \otimes |x|$ for its grading $|x|$, and the $D^{\vee}(G)$-coaction on $Y$ is $\Delta_{R}y = \sum_{g \in G} g \triangleright y \otimes \delta_{g} \otimes |y|$. We therefore have that the pushout of $Y$ to a $\C G$-comodule contains $X$ iff $X$ and $Y$ have the same gradings. 

The averaging map is given by \eqref{eq_averaging_map_kofG} (using $Y$ a $\C (G)$-comodule by its pushout instead of $W$) and as such the transfer map is given by \begin{align*}
(\Theta^{Y}_{\C (G)}\circ &\textup{av}\circ (\Theta^{X}_{\C G})^{-1})\circ (\delta_{g} \otimes x) \\
& = \Theta^{Y}_{\C (G)}\circ \textup{av}(\delta_{g} 
\otimes |x|^{-1} \otimes x) = \Theta^{Y}_{\C (G)}(\sum_{k \in G} \delta_{gk^{-1}} \otimes k|x|^{-1}k^{-1} \otimes k \triangleright x) \\
& = \frac{1}{|G|} \Theta^{Y}_{\C (G)}(\sum_{f \in G} \delta_{f^{-1}} \otimes fg|w|^{-1}g^{-1}f^{-1} \otimes fg \triangleright x) =  \frac{1}{|G|} g|x|^{-1}g^{-1} \otimes g \triangleright x, 
\end{align*} which clearly intertwines the coactions \eqref{eq_coact_3} and \eqref{eq_coact_Eprime} and can be written as stated. \endproof
 
\begin{corollary} The transfer map $E\to E'$ constructed in Lemma~\ref{lemma_transferance_of_sections_DG} factors through $E^*$ as the composition of the transfer maps in Lemma~\ref{lem_transfer_to_Lorentz} and Lemma~\ref{lemFtransfer}.
\end{corollary}
\proof  Using the preceding lemma, the composition of the two relevant  transfer of sections map is then indeed given by \begin{align*} \delta_{q_{c}} \otimes v & \mapsto \sum_{n \in C_{G}} \delta_{q_{c}n^{-1}} \otimes n \triangleright v \mapsto \frac{1}{|G|}\sum_{n \in C_{G}}q_{c}n^{-1}(n|v|^{-1}n^{-1})nq_{c^{-1}} \otimes (q_{c}n^{-1})n \triangleright v\\
& = \frac{1}{|G|}\sum_{n \in C_{G}}c^{-1} \otimes q_{c} \triangleright v = \frac{|C_{G}|}{|G|} c^{-1} \otimes q_{c} \triangleright v,  \end{align*} using the fact $|v| = r$, as in Lemma~\ref{lemma_transferance_of_sections_DG}. \endproof

For completeness, we note that we can also construct a transfer of sections map $\C G \otimes Y \to \C (G) \otimes X$ under similar assumptions (and is explicitly given by $g \otimes y \mapsto \sum_{f \in G} \delta_{g, |y|^{-1}} \delta_{f^{-1}} \otimes f \triangleright y$). However, using $Y= X = W$ for the composition $E \hookrightarrow E' \to E^{\star}$ does not coincide with the direct transfer map of Lemma~\ref{lem_transfer_to_Lorentz}. Similarly, we can also build transfer of sections maps $E' \to E$ and $E^{\star} \to E$, but these do not have a particularly interesting form.

\begin{example} \label{example_transferance_of_sections}\rm We take $G=S_3$ and $\CC=\{uv,vu\}$ as Case (ii) from Example~\ref{exS3} in the role of choice of `mass', with  $C_G=\Z_3=\{e,r,r^2\}$ for basepoint $r=uv$ and $\pi_j(r)=q^j r$ for $j=0,1,2$ in the role of  the choice of `spin'. Here $q=e^{2\pi\imath\over 3}$.  In the geometric picture, $S_{3}/\mathbb{Z}_{3}$ has two equivalence classes, namely $[e]=\{e,uv,vu\}=C_G$ and $[u]=\{u,v,w\}$. The classes are bijective to $q_\CC=\{e,u\}$ which in turn is bijective to $\CC$. The algebra $\C(S_{3}/\mathbb{Z}_{3})$ can accordingly be viewed as $\C(S_{3})^{\Z_3}$ with basis 
\[ \delta_{[e]}\leftrightarrow \delta_e+\delta_{uv}+\delta_{vu}, \quad  \delta_{[u]}\leftrightarrow\delta_u+\delta_v+\delta_w, \]
or $\C(q_\CC)$ with basis $\delta_{q_{uv}}=\delta_e, \delta_{q_{vu}}=\delta_u$ (this was our preferred basis in proofs), or indeed with $\C(\CC)$ with basis $\delta_{uv},\delta_{vu}$.  

 The role of $H=\C (C_G)\rcocross \C  G$ when working with comodules is equivalently played by the Hopf algebra $U=\C (G)\lcross \C  C_G\subset D(G)$ when working with modules. Given an irrep $V_\pi$ of $C_G$, we extend this to a representation $V_{r,\pi}$ of $U$ by the same vector space and the action $f\la v= f(r) v$ of $f\in \C (G)$ on $v\in V_\pi$. For our present choices, $U=\C(S_3)\lcross \C \Z_3$ and $V_{r,\pi_j}=\C $ is 1-dimensional, and the action of $U$ is $f\la \lambda=f(uv)\lambda$, $(uv)^{\pm 1} \la \lambda= q^{\pm j}\lambda$ for $\lambda \in \C$. The sections $E=(D^{\vee}(S_{3}) \otimes V_{r,\pi_j})^{H}$ of the associated bundle can be identified according to \eqref{eq_trivial_associated_bundle_H} with $\C(S_{3}/\mathbb{Z}_{3})\tens V_{r,\pi_j}=\C(S_{3}/\mathbb{Z}_{3})$ in our case, sending the latter to
 \begin{equation} \label{eq_useful_identification} \delta_{[e]} \mapsto (\delta_{e} + q^{j}\delta_{vu} + q^{-j}\delta_{uv})\otimes r^{-1}, \quad \delta_{[u]} \mapsto (\delta_{u} + q^{j}\delta_{w} + q^{-j}\delta_{v}) \otimes r^{-1}, \end{equation}
using the actions given in Example~\ref{exS3}(ii). Here, $E$ itself is a $D(G)$-module with the same structure as in Example~\ref{exS3}(ii) on replacing $\{ uv, vu\}$ by $\{ \delta_{[e]}, \delta_{[u]}\}$.
 
 Following the prescription, we  next fix a representation $W$ of $D(S_{3})$ whose restriction to $U$ contains $V_{r, \pi_{j}}$. The simplest case it to just take $W=V_{\{ uv, vu\},\pi_j} = \C\{ uv, vu\}$. With this, the sections $E^{W}=(D^{\vee}(S_{3}) \otimes W)^{H}$ of the associated bundle can be identified according to the analogue  of \eqref{eq_trivial_associated_bundle_H} sending $\C(S_{3}/\mathbb{Z}_{3})\tens V_{\{ uv, vu\},\pi_j}$ to
 \[\delta_{[e]} \otimes r^{\pm 1} \mapsto (\delta_{e} + q^{\pm j}\delta_{vu} + q^{\mp j}\delta_{uv}) \otimes r^{-1} \otimes r^{\pm 1}, \]
 \[\delta_{[u]} \otimes r^{\pm 1} \mapsto (\delta_{u} + q^{\pm j}\delta_{w} + q^{\mp j}\delta_{v}) \otimes r^{-1} \otimes r^{\pm 1}. \] $E^{W}$ is a $D(G)$-module (corresponding to the coaction \eqref{s_coact_coact}) by  
 \[ (\delta_{g} \otimes h) \triangleright (\delta_{q_{c}} \otimes w) = \delta_{h^{-1}gh, q_{c}|w|q_{c}^{-1}}\delta_{q_{hch^{-1}}} \otimes \zeta_{c}(h) \triangleright w,\] which in our case using $\C (S_{3}/\mathbb{Z}_3)$ (rather than $\C(q_{\mathcal{C}})$) gives
  \[ u \triangleright \delta_{\{ [e], [u]\}} \otimes r^{\pm 1}= \delta_{[u], [e]} \otimes r^{\pm 1}, \ v \triangleright \delta_{\{[e], [u]\}} \otimes r^{\pm 1} =  \delta_{\{ [u], [e]\}} \otimes q^{\pm j}r^{\pm 1}, \] 
 \[ f \triangleright \{ \delta_{[e]} \otimes r^{\pm 1},  \delta_{[u]} \otimes r^{\pm 1}\} = \{ f(r^{\pm 1})\delta_{[e]} \otimes r^{\pm 1}, f(r^{\mp 1})\delta_{[u]} \otimes r^{\pm 1}\}, \] 
for the $D(S_3)$ action on $E^W$. The embedding $E \hookrightarrow E^W$ as a $D(S_3)$ module is simply 
 \[ \delta_{[e]} \mapsto \delta_{[e]} \otimes r, \quad \delta_{[u]} \otimes r.\]
  The sections $E'= (D^{\vee}(S_{3}) \otimes W)^{\C(S_{3})} $ of the bundle associated to $\C S_{3}\hookrightarrow D(S_{3}) \twoheadrightarrow \C(S_{3})$ can be identified according to \eqref{eq_trivial_associated_bundle_Hprime} with $\C S_{3} \otimes V_{\{ uv, vu\},\pi_j}$, sending its basis elements to
 \begin{align*} \{ e, uv, vu\} \otimes r^{\pm 1} \mapsto &(\delta_{e} + q^{\pm j}\delta_{vu} + q^{\mp j}\delta_{uv}) \otimes \{ e, uv, vu\} \otimes r^{\pm 1} \\
 & + (\delta_{u} + q^{\pm j}\delta_{v} + q^{\mp j}\delta_{w}) \otimes \{e, vu, uv\} \otimes r^{\mp 1} ,\\
  \{ u, v, w\} \otimes r^{\pm 1} \mapsto &(\delta_{e} \otimes \{u, v , w\} + q^{\pm j}\delta_{vu} \otimes \{ v, w, u\} +  q^{\mp j}\delta_{uv} \otimes \{ w, u, v\})\otimes r^{\pm 1} \\
 & + (\delta_{u} \otimes \{u, w, v \} + q^{\pm j}\delta_{v} \otimes \{ w, v, u\} + q^{\mp j}\delta_{w} \otimes \{ v, u, w\})\otimes r^{\mp 1}.\end{align*} 
Moreover, $E'$ is a $D(G)$-module (corresponding to \eqref{eq_coact_Eprime}) by
 \[ (\delta_{g} \otimes h) \triangleright (f \otimes w) = \langle (\beta^{-1})^{*}(\delta_{g} \otimes h), S(f \otimes w)^{0} \rangle (f \otimes w)^{1} = \delta_{g^{-1}, hfh^{-1}} hfh^{-1} \otimes h \triangleright w, \]
 for $f \in G$, which in our case is
 \[ u \triangleright \{e, u, v, w, uv, vu \} \otimes r^{\pm 1} = \{e, u, w, v, vu, uv \} \otimes r^{\mp 1},\]\[ v \triangleright \{e, u, v, w, uv, vu \} \otimes r^{\pm 1} = \{e, w, v, u, vu, uv \} \otimes q^{\pm j}r^{\mp 1},\]
 \[ f \triangleright \{e, u, v, w, uv, vu \} \otimes r^{\pm 1} = \{f(e)e, f(u)u, f(v)v, f(w)w, f(vu)uv, f(uv)vu \} \otimes r^{\pm 1}.\]
 Finally, the averaging map  (\ref{DGav}) upstairs (before restricting to the invariant subspaces) can then be explicitly computed in our case. For example,
 \[ {\rm av}(\delta_{u}\tens uv\tens r^{\pm 1})={1\over 6}\big((\delta_u+ q^{\pm j}\delta_w+ q^{\mp j}\delta_v)\tens uv\tens r^{\pm 1}+(\delta_e+q^{\pm j}\delta_{uv}+q^{\mp j}\delta_{vu})\tens vu\tens r^{\mp 1}\big), \]
 (and 35 similar).  Restricting to a map $E^W\to E'$ and using the above identifications of the associated bundles one can check this averaging map comes out as 
 \[ \delta_{[e]} \otimes r^{\pm 1} \mapsto \frac{1}{3} \, vu \otimes r^{\pm 1} , \quad \delta_{[u]} \otimes r^{\pm 1} \mapsto \frac{1}{3}\,uv \otimes r^{\mp 1},  \] 
in accordance with Example~\ref{example_embedding}. The $+$ cases  here are the images of $\delta_{[e]}$, $\delta_{[u]}$ for the transfer map $E\to E'$ in our description. 
\end{example}

\section{Tautological calculus $\Omega^1_{\mathcal{C}, \pi}(D^\vee(G))$ and  free field equations}\label{secPcalc}

We now study the quantum differential geometry of the bundles previously introduced in Section~\ref{secWigner} in order to recover a free field equations characterisation of the coirreps of $D^{\vee}(G)$. So far, these principal bundles use the universal calculus. For a quantum $H$-principal bundle with general calculus, we recall the definition around \eqref{eq_exact_bundle_def}. Since $A\subseteq P$ is a subalgebra, it inherits the natural calculus \[ \Omega^{1} = A(\extd A)A = A\extd A\] by restriction, which we use. With this, the algebra inclusion $A \hookrightarrow P$ is automatically differentiable. Moreover, as our bundles of interest are quantum homogeneous bundles obtained from a Hopf algebra surjection $\mathfrak{p}: P \twoheadrightarrow H$, we focus on bicovariant $\Omega^{1}_{P}$. Here the total space algebra $P$ has a left coregular $P$-coaction $\Delta_L=\Delta$ and bicovariance ensures that this  is differentiable as an algebra map when $P \otimes P$ is equipped with the canonical tensor product calculus $\Omega^{1}_{P \otimes P}:= (\Omega^{1}_{P} \otimes P)\bigoplus (P \otimes \Omega^{1}_{P})$, $\extd_{_{P \otimes P}} = (\extd_{_{P}} \otimes \id_{P})\oplus (\id_{P} \otimes \extd_{_{P}})$, i.e. the following map \begin{equation} \label{eq_def_differentiable} (\Delta_{L})_{*}: \Omega^{1}_{P} \to (\Omega^{1}_{P} \otimes P)\bigoplus (P \otimes \Omega^{1}_{P}), \quad p\extd q \mapsto p^{1}\extd q^{1} \otimes p^{2}q^{2} + p^{1}q^{1} \otimes p^{2}\extd q^{2} \end{equation} is well-defined automatically by \eqref{eq_bicovariant_def} (since here $p^{1} \otimes p^{2} = p_{1} \otimes p_{2}$). Moreover, restricting the coaction $\Delta_{L}$ to the base makes $A$ a differentiable $P$-comodule algebra, i.e. $(\Delta_{L})_{*}|_{\Omega^{1}}: \Omega^{1} \to (\Omega^{1}_{P} \otimes A) \bigoplus (P \otimes \Omega^{1})$ is automatically well-defined in an analogous way to \eqref{eq_def_differentiable}. By \cite[Lemma 5.46]{BegMa}, we can also say immediately that \begin{equation} \label{eq_calc_lambda1} \Lambda^{1} := \frac{H^{+}}{\mathfrak{p}(I)}, \quad \Omega^{1}_{P} \cong P \otimes \frac{P^{+}}{I}, \end{equation} for a subcrossed module $I \subseteq P^{+}$, and then the map $P \twoheadrightarrow H$ is also differentiable. Moreover, it is immediate by \cite[(5.19)]{BegMa} that the map $\textup{ver}$ in \eqref{eq_exact_bundle_def} is simply given in our case of interest by \begin{equation} \label{eq_map_ver}\textup{ver}: P \otimes (P^{+}/I) \to P \otimes (H^{+}/\pi(I)), \quad p \otimes \varpi_{_{P}}(\pi_{\epsilon}(q)) \mapsto p \otimes \varpi(\mathfrak{p}(\pi_{\epsilon}(q)))\end{equation} for the quotient maps $\varpi_{_{P}}: P^{+} \twoheadrightarrow P^{+}/I$ and $\varpi: H^{+} \twoheadrightarrow H^{+}/\mathfrak{p}(I)$ defining the bicovariant calculi (see Section~\ref{secpre}). The horizontal forms are thus given by the $P$-module 
\[ P \otimes \textup{ker}(P^{+}/I \twoheadrightarrow H^{+}/\mathfrak{p}(I)) = P \otimes \textup{ker}(P^{+} \twoheadrightarrow H^{+}/\mathfrak{p}(I))/I\subseteq P \otimes (P^{+}/I).\]
Thus, when $P$ has bicovariant calculus,  the differential geometry of the bundle is fully characterised by that of $P$. Our first task then is to recall the classification of the bicovariant calculi on $D^{\vee}(G)$ from \cite{Ma:cla}, but reworked in the conventions that we need. 

\subsection{Recap of bicovariant calculi $\Omega^1_{\mathcal{C}, \pi}$ on the total space $D^{\vee}(G)$} \label{sec_calc_on_DveeG}

A key property of the coquasitriangular Hopf algebra $D^{\vee}(G)$ is that it is factorisable, meaning its associated quantum Killing form \cite{Ma}, in our case $\CQ=\CR_{21}*\CR:D^{\vee}(G)^{\tens 2}\to k$,  is nondegenerate. Here we use the convolution product $*$ in $D^\vee(G)^{\tens 2}$  and explicitly, using \eqref{eq_coquasi_DveeG}, we therefore have have
\begin{equation} \label{eq_quantum_killing_form_DG} \CQ((\delta_{g} h)\tens (\delta_{u}  v)): =\CR((\delta_{u} v)_{1} \tens (\delta_{g} h)_{1})\CR((\delta_{g} h)_{2}\tens (\delta_{u} v)_{2}) = \delta_{g, uvu^{-1}}\delta_{h, u},\end{equation}
and this induces a linear isomorphism $\mathcal{Q}_{1}: D^{\vee}(G)\to (D^{\vee}(G))^{*}$ by $\delta_{g}h \mapsto \mathcal{Q}(\delta_{g}h \otimes \_)$.  In fact \[ \beta^{*}\circ \CQ_1 = \mathcal{Q}_{1} \circ \beta^{-1}:D^{\vee}(G)\to D(G)\] is simply the identity map as a linear map. Due to this factorisability, it follows that  bicovariant calculi on $D^{\vee}(G)$ are in 1-1 correspondence with two-sided ideals of $D(G)$\cite{Ma:cla}. By the latter's semisimplicity and finite dimensionality, coirreducible calculi (in the sense of having no proper quotients) correspond to simple two-sided ideals and hence to  irreps of $V_{\mathcal{C}, \pi}$ of $D(G)$. We already know how to construct these from Section~\ref{secWigner}.

Given a  fixed pair $(\mathcal{C}, \pi)$, we choose a basis $\{ v_{i}\}$ of $V_{\pi}$ with dual basis $\{ v^{i}\}$. This gives a basis $E_{ci}:= c \otimes v_{i}$ for $V_{\mathcal{C}, \pi}$ (where $c \in \mathcal{C}$) and $E^{ci} = \delta_{c} \otimes v^{i}$ for its linear dual $V_{\mathcal{C}, \pi}^*$. By \eqref{Cpiirrep}, the right $D^{\vee}(G)$-coaction on $V_{\mathcal{C}, \pi}$ and its categorical dual when written as a Whitehead $G$-crossed module \eqref{eq_comoduleVpi} is then given by \begin{equation} \label{eq_crossed_module_V} h \triangleright E_{ci} = \pi(\zeta_{c}(h))^{j}{}_{i}E_{(hch^{-1})j}, \ |E_{ci}| = c, \ h \triangleright E^{ci} = \pi(\zeta_{c}(h)^{-1})^{i}{}_{j}E^{(hch^{-1})j}, \ |E^{ci}| = c^{-1}. \end{equation} 
The latter are obtained by noting that an arbitrary $G$-crossed module $W$ with basis $\{ W_{i}\}$ has categorical dual realised on the linear dual $V^{*}$ with $G$-crossed structure 
\begin{equation} \label{eq_categorical_dual_crossed_module} \Delta_{R}\phi = \sum_{g \in G} \langle \phi, g \triangleright V_{i} \rangle V^{i} \otimes S(\delta_{g} \otimes |V_{i}|) = \sum_{g \in G} \langle \phi, g^{-1} \triangleright V_{i} \rangle \, V^{i} \otimes \delta_{g} \otimes g^{-1}|V_{i}|^{-1}g, \end{equation} 
for all  $\phi \in V^{\sharp}$. We take basis  $E_{ci}{}^{dj}:= E_{ci} \otimes E^{dj}$ of $\textup{End}(V_{\mathcal{C}, \pi})$ and with the convention of using left-invariant 1-forms on $D^{\vee}(G)$ (as opposed to right-invariant 1-forms  in \cite{Ma:cla}) the classification of calculi is as follows.

\begin{lemma} \label{lem_coirreducible_calc_DG}Coirreducible bicovariant calculi on $D^\vee(G)$ are of the form $\Omega^{1}_{\mathcal{C}, \pi} := D^{\vee}(G) \otimes \textup{End}(V_{\mathcal{C}, \pi})$ for a non-trivial pair $(\mathcal{C}, \pi)\neq (\{ e\}, 1) $). The  $D^{\vee}(G)$-bimodule structure of $\Omega^{1}_{\mathcal{C}, \pi}$ is by the free left action and the right action
\begin{equation}\label{eq_module_endV_struc} E_{ci}^{\hphantom{ci}dj} \cdot \delta_{g} = R_{c}(\delta_{g}) \otimes E_{ci}^{\hphantom{ci}dj},\quad E_{ci}^{\hphantom{ci}dj} \cdot h = chc^{-1} \otimes (E_{ci} \otimes h^{-1} \triangleright E^{dj}),  \end{equation}  where $R_{c}(\delta_{g}) = \delta_{gc^{-1}}$ is the right-translation operator. The  $D^{\vee}(G)$-bicomodule structure is by the free left coaction and the right $D^{\vee}(G)$-comodule structure induced by tensor product coaction on $V_{\mathcal{C}, \pi} \otimes V_{\mathcal{C}, \pi}^{\sharp} = \textup{End}(V_{\mathcal{C}, \pi})$, explicitly \begin{equation} \label{eq_right_comodule_EndVCpi} \, h \triangleright E_{ci}{}^{dj} = \pi(\zeta_{c}(h))^{k}{}_{i}\pi(\zeta_{d}(h)^{-1})^{j}{}_{l}E_{(hch^{-1})k}{}^{(hdh^{-1})l} , \quad |E_{ci}{}^{dj}| = cd^{-1}\, .\end{equation} The exterior derivative is given by \begin{equation} \label{eq_exterior_deriv_DveeG} \extd (\delta_{g}) =  \sum\limits_{c \in \mathcal{C}} \partial^{c}(\delta_{g}) \otimes E_{ci}^{\hphantom{ci}ci},\quad \extd h = \sum\limits_{c \in \mathcal{C}} chc^{-1} \otimes (E_{ci} \otimes h^{-1} \triangleright E^{ci}) - h \otimes E_{ci}{}^{ci}\, ,  \end{equation} and extended to products by the Leibniz rule. Here $\partial^{c}:= R_{c} - \id$. The calculus is inner by $\theta = 1 = \sum_{c \in \mathcal{C}}E_{ci}{}^{ci} \in \textup{End}(V_{\CC, \pi})$. 
\end{lemma} 
\proof  The result is in principle known \cite{Ma:cla} but we need the relevant formulae in our conventions. Viewing the $D(G)$-irrep $V_{\mathcal{C}, \pi}$ as a map of algebras $\rho_{_{\mathcal{C}, \pi}}: D(G) \to \End(V_{\mathcal{C}, \pi})$ the matrix $\rho_{_{\mathcal{C}, \pi}}(\delta_{g} \otimes h)$ has entries (using \eqref{Cpiirrep}) given by \begin{equation} \label{eq_matrix_coeff} \rho_{\mathcal{C}, \pi}(\delta_{g} \otimes h)^{ci}{}_{dj} = \delta_{c, hdh^{-1}}\delta_{g, hdh^{-1}}\pi(\zeta_{d}(h))^{i}{}_{j} = \delta_{g, c} \delta_{g, hdh^{-1}}\pi(\zeta_{d}(h))^{i}{}_{j} , \quad c, d \in \mathcal{C}\, . \end{equation} We then ask for the map \[\rho_{\mathcal{C}, \pi} \circ \beta^{*} \circ \mathcal{Q}_{1}: D^{\vee}(G) \twoheadrightarrow \textup{End}(V_{\mathcal{C}, \pi}), \, \,  \delta_{g} \otimes h \mapsto \sum_{c, d \in \mathcal{C}}\rho_{\mathcal{C}, \pi}(\delta_{g} \otimes h)^{ci}{}_{dj}E_{ci}{}^{dj} = \delta_{g, \mathcal{C}} \, E_{gi} \otimes h^{-1} \triangleright E^{gi}, \] to be a morphism of right $D^{\vee}(G)$-crossed modules, where   $D^{\vee}(G)$ is a right crossed module by $(\cdot, \textup{Ad}_{R})$ and $\textup{Ad}_{R}$ is the right adjoint coaction (for which an explicit formula is given later in \eqref{eq_adjoint_actions_DG_1}). This requirement induces the right $D^{\vee}(G)$-crossed comodule structure on $\textup{End}(V_{\mathcal{C}, \pi})$ given by the right $D^{\vee}(G)$-comodule structure of \eqref{eq_right_comodule_EndVCpi} (i.e. $\Delta_{R}E_{ci}{}^{dj} = \sum_{h \in G}h \triangleright E_{ci}{}^{dj} \otimes \delta_{h} \otimes |E_{ci}{}^{dj}|$) and right $D^{\vee}(G)$ module \[ E_{ci}{}^{dj} \triangleleft \delta_{g} = \delta_{g, c} E_{ci}{}^{dj}, \quad E_{ci}{}^{dj} \triangleleft h=  E_{ci} \otimes h^{-1} \triangleright F^{dj}. \] Subsituting 
\begin{equation} \label{eq_quantum_maurer_cartan_Cpi}\varpi:= \rho_{\mathcal{C}, \pi} \circ \beta^{*} \circ \mathcal{Q}_{1}|_{D^{\vee}(G)^{+}}, \end{equation}
 and the above data into \eqref{recover_diffcalc} recovers the stated structure maps for the calculus. \endproof

By previous discussion, a general bicovariant calculus on $D^{\vee}(G)$ is given by $\Omega^{1}_{D^{\vee}(G)} \cong D^{\vee}(G) \otimes (\bigoplus_{(\mathcal{C}, \pi)}\textup{End}(V_{\mathcal{C}, \pi})) \cong \bigoplus_{(\mathcal{C}, \pi)} \Omega^{1}_{\mathcal{C}, \pi}$ for a selection of non-trivial pairs $(\mathcal{C}, \pi)$ with the structure maps from Lemma~\ref{lem_coirreducible_calc_DG} on each block as a  $D^{\vee}(G)$-Hopf bimodule, and exterior derivative given by the direct sum of the above. We also know by the Artin-Wedderburn theorem  that $ \bigoplus_{(\mathcal{C}, \pi)}\rho_{\mathcal{C}, \pi}: D(G) \cong \bigoplus_{(\mathcal{C}, \pi)} \textup{End}(V_{\mathcal{C}, \pi})$. This has inverse $\bigoplus \phi_{\mathcal{C}, \pi}: \bigoplus_{(\mathcal{C}, \pi)} \textup{End}(V_{\mathcal{C}, \pi}) \to D(G)$, given by 
\begin{equation} \label{eq_artin_wedderburn_DG} \phi_{\mathcal{C}, \pi}: E_{ci}{}^{dj} \mapsto s_{ci}{}^{dj}:= \frac{\textup{dim}(V_{\pi})}{|C_{G}|} \sum_{n \in C_{G}} \pi(n^{-1})^{j}{}_{i} \, \delta_{c} \otimes q_{c}nq_{d}^{-1}, \end{equation}  
which has image $P_{\mathcal{C}, \pi}D(G)P_{\mathcal{C}, \pi} \subset D(G)$, i.e. $\textup{End}(V_{\mathcal{C}, \pi}) \cong P_{\mathcal{C}, \pi}D(G)P_{\mathcal{C}, \pi}$ (for detail, see \cite{CowMa}). We can use this to  obtain a description of the subcrossed modules $I \subseteq D^{\vee}(G)^{+}$ such that $\Omega^{1}_{D^{\vee}(G)} \cong D^{\vee}(G)^{+}/I$. Using $(\beta^{*} \circ \mathcal{Q}_{1})^{-1}  = \id$, to convert to $D^\vee(G)$, the relevant ideal for $\Omega_{\CC,\pi}$ is then
 \[ I_{\mathcal{C}, \pi} = \textup{Span}_{\mathbb{C}}\{ s_{ci}{}^{dj}\}.\]
Moreover, every subcrossed module $I$ of $D^{\vee}(G)^{+}$ is of the form \begin{equation} \label{eq_gen_calc_DveeG}I \cong \bigoplus_{(\mathcal{C}',  \pi') \notin \mathcal{S} }I_{\mathcal{C}', \pi'} \implies D^{\vee}(G)^{+}/I \cong \bigoplus_{(\mathcal{C}, \pi) \in \mathcal{S}} I_{\mathcal{C}, \pi}, \end{equation} for some subset $\mathcal{S} $ of the collection of all non-trivial pairs $(\mathcal{C}, \pi)\ne (\{ e\}, 1)$). In particular $\Lambda^{1}_{\mathcal{C}, \pi} \cong I_{\mathcal{C}, \pi}$, giving a useful basis of $\Lambda^{1}_{\mathcal{C}, \pi}$ for later calculations. The structure maps for $\Omega_{\CC,\pi}$ in the form $\Omega^{1}_{\mathcal{C}, \pi} \cong D^{\vee}(G) \otimes D^{\vee}(G)/I$ are just the same as in Lemma~\ref{lem_coirreducible_calc_DG}, but replacing $E_{ci}{}^{dj}$ with $s_{ci}{}^{dj}$. 
 
\subsection{Non-universal quantum homogeneous principal bundles $P \twoheadrightarrow H$} \label{sec_bundle_PL}

We now study the differential geometry of the principal bundles of Section~\ref{secWigner}. To this end, we first recall the well-known classifications for bicovariant calculi on $\C (G)$ and $\C G$.

\subsubsection{Recap of calculi on $\C (G)$ and $\C  G$} \label{sec_recap_calc_kG_kofG} By the form of $\textup{Ad}_{R}$ and right multiplication on $\C (G)$, it is immediate that bicovariant calculi on $\C (G)$ are all given by Ad-stable subsets $\CS\subseteq G \backslash \{ e \}$, with structure maps on the resulting calculi $\Omega^{1}_{\mathcal{S}} \cong \C (G) \otimes \Lambda^{1}_{\mathcal{\CS}}$ 
\begin{equation} \label{eq_kG_calc} \extd\delta_{g} = \sum_{c \in \mathcal{S}} \partial^{c}(\delta_{g}) \otimes e_{c},\quad e_{c} \cdot \delta_{g} = R_{c}(\delta_{g}) \otimes e_{c},\quad  \Delta_{R}e_{c} = \sum_{f \in G} e_{fcf^{-1}} \otimes \delta_{f},\end{equation}
where $\Lambda^{1}_{\mathcal{S}}$ is a vector space with basis denoted $\{ e_{s}\}_{s \in \mathcal{S}}$ (details in \cite[Proposition 1.52]{BegMa}). The cross-module map generating this calculi is $\varpi(\pi_{\epsilon}(\delta_{g})) = \delta_{g, \mathcal{S}}e_{g}$. It is clear that non-trivial conjugacy classes classify the coirreducible bicovariant calculi. Bicovariant calculi also correspond to Cayley graphs generated by $\mathcal{S}$ with arrows $\{ e_{x \to xc}\}_{x \in G, c \in \mathcal{S}}$ and these correspond to  $e_{x \to xc} = \delta_{x}\extd (\delta_{xc})$ with $\extd$ as above, or conversely $e_{c} = \sum_{x \in G}e_{x \to xc}$. The calculus is inner by $\theta = \sum_{c \in \mathcal{S}}e_{c}$ and connected iff $\mathcal{S}$ is a generating set for $G$ \cite{BegMa}. 

On the other hand, bicovariant calculi on $\C G$ are clearly just equivalent to right-ideals $I$ of $\C G$ (since $\textup{Ad}_{R}I \subseteq I \otimes \C G$ holds automatically as  $\textup{Ad}_{R}g = g \otimes e$). These can be conveniently described using (a type of) 1-cocycles $\zeta$ on $G$ (by taking $\varpi \circ \epsilon = \zeta$, with the latter linearised). From this, we define elements $e^{g} := \varpi(\pi_{\epsilon}(g)) = g^{-1}\extd g$ and note that $\{ e^{g}\}_{g \in G}$ spans $\Lambda^{1}$. By linearity of $\varpi$ and the Leibniz rule, we note the relations \begin{equation} \label{eq_linearity_e}e^{\lambda g} = \lambda e^{g}, \quad e^{g + h} = e^{g} + e^{h}, \quad e^{g}h = h(e^{gh} - e^{h}), \quad e^{e} =0 \quad \forall g, h \in G, \end{equation} since $\extd e = 0$. An example of particular interest, given a representation $\rho: G\to \End (V)$ we can take $\varpi = \rho|_{\C G^{+}}$, giving $\Omega^{1}_{\rho}=\C G \otimes \Lambda^{1}_{\rho} := \C  G\tens\rho(\C G^{+})$  with  
\begin{equation} \label{eq_calc_kG_rep} x.g = g \otimes x\rho(g),\quad \extd g=g\tens  \rho(g) - g \otimes 1,\quad \Delta_R x=x\tens 1,\end{equation} 
for all $x \in \textup{End}(V), g \in G$. Details are in \cite[Example 1.49]{BegMa}. By the Artin-Wedderburn theorem, such bicovariant calculi can be `broken down' its coirreducible calculi generated by the non-trivial irreps of $G$ contained in $\rho$. These calculi simply corrrespond to the case when the subcrossed module $I$ describing the calculus is a two-sided ideal of $\C G$. Such calculi are inner if $1 \in \textup{End}(V_{\rho})$ is contained in the image of $\varpi$, which is the case iff $V_{\rho}$ does not contain the trivial irrep in its decomposition into irreps and each irrep has multiplicity 1, and is connected iff $\rho$ is faithful. For the whole exterior algebra $(\Omega,\extd)$ a standard choice is the $\{e^g\}$ anticommuting so that $\Lambda$ is the Grassmann algebra and $\extd e^g=0$, but we do not assume this this. In the maximal prolongation, one always has $\extd e^g=- e^g\wedge e^g$ as in  \cite[Propn. 3.3]{MaTao}. For $\C (G)$, the standard exterior algebra is more complicated and uses the quantum double braiding to `skew symmetrize'. 

\subsubsection{The bundle $D^{\vee}(G) \twoheadrightarrow \C (G)$} \label{sec_diff_geom_first_bundle}

We now study the differential geometry of the $\C (G)$-principal bundle $D^{\vee}(G) \twoheadrightarrow \C (G)$ over the base quantum manifold $\C G$ discussed in Section~\ref{sec_second_bundle_construct} but with $D^{\vee}(G)$ equipped with non-universal calculi.  

\begin{proposition} \label{prop_bundle_over_kG}
For the bicovariant calculus $\Omega^{1}_{\mathcal{C}, \pi}$ on $D^{\vee}(G)$, we obtain a principal homogeneous bundle $\C G \hookrightarrow D^{\vee}(G) \twoheadrightarrow \C (G)$ with 
\[ \Omega^{1} \cong \Omega^{1}_{\widetilde{\pi}}\, , \quad \Lambda^{1} \cong \delta_{\pi, 1}\Lambda^{1}_{\mathcal{C}},\]
where $\Omega^{1}_{\widetilde{\pi}}$ is defined as above using the $G$-representation $\widetilde{\pi}$ obtained by restricting the action on $V_{\mathcal{C}, \pi}$ (i.e., with $G$-action from \eqref{Cpiirrep}), $\Lambda^{1}_{\mathcal{C}}$ is as above and $\Lambda^1$ is the calculus needed on $\C (G)$ as part of the bundle.
\end{proposition}

\proof Denoting $I:= \bigoplus_{(\mathcal{C}', \pi')|(\mathcal{C}', \pi') \neq (\mathcal{C}, \pi) }$ (for non-trivial pairs $(\mathcal{C}', \pi')$) the subcrossed module giving $\Lambda^{1}_{\mathcal{C}, \pi}$ (by \eqref{eq_gen_calc_DveeG}), then by \eqref{eq_calc_lambda1}, the resulting bicovariant calculus on $\C (G)$ is then obtained by quotienting with the subcrossed module $p_{1}(I)$, with $p_{1}: D^{\vee}(G) \twoheadrightarrow \C (G)$ is the Hopf algebra surjection onto the first factor. Under this map, \[ p_{1}(s_{ci}{}^{dj}) = \frac{\textup{dim}(V_{\pi'})}{|C_{G}|} \sum_{n \in C_{G}} \pi'(n^{-1})^{j}{}_{i} \delta_{c}= \delta_{\pi', 1}\delta_{c}, \quad s_{ci}{}^{dj} \in I_{(\mathcal{C}', \pi')},  \] by the Schur orthogonality relations. So $p_{1}: I \mapsto \bigoplus_{\mathcal{C}' \, | \, (\mathcal{C}', 1) \neq (\mathcal{C}, \pi), \, \mathcal{C}' \neq \{ e\}} I_{\mathcal{C}'}$, where $I_{\mathcal{C}'}$ denotes the subcrossed module of $\C (G)$ spanned by $\{ \delta_{c}\}_{c \in \mathcal{C}'}$. Noting that $\C (G)^{+}/I \cong \Lambda^{1}_{\mathcal{C}}$ by $\varpi$ (which is essentially just a `relabelling' $\delta_{c}$ with $e_{c}$), we thus have that \[ (p_{1})_{*}: \Lambda^{1}_{\mathcal{C}, \pi} \twoheadrightarrow \Lambda^{1} \cong \delta_{\pi, 1} \Lambda^{1}_{\mathcal{C}}, \quad I_{(\mathcal{C, \pi})}\ni s_{ci}{}^{dj} \mapsto \delta_{\pi, 1}e_{c}\, . \] The calculus on the base $\C G$ is then given under our conventions by $\C G(\extd \C G)\C G = \C G \extd \C G\subseteq \Omega^{1}_{\mathcal{C}, \pi}$, recalling how $\C G$ sits inside $D^{\vee}(G)$ by \eqref{eq_inclusion_A_second_bundle} in this case. Using this isomorphism of algebras $\C G \cong D^{\vee}(G)^{\C (G)}$, and asking it to be differentiable, then  \begin{align} \label{eq_hdg} h\extd g &= \sum_{k, f \in G} (\delta_{f} \otimes f^{-1}hf)\extd(\delta_{k} \otimes k^{-1}gk) \nonumber\\
&=  \sum_{f \in G, c \in \mathcal{C}}\delta_{f c^{-1}} \otimes c f^{-1}hgfc^{-1} \otimes E_{ci} \otimes f^{-1}g^{-1}f \triangleright F^{ci} - \delta_{f} \otimes f^{-1}hg f \otimes E_{ci}{}^{ci}\, , \end{align} 
for $h, g \in G$, using \eqref{eq_exterior_deriv_DveeG} for $\extd$. We thus get $\Omega^{1} = \textup{Span}_{\mathbb{C}}\{ h\extd g \, | h ,g \in G\} \subseteq \Omega^{1}_{\mathcal{C}, \pi}$. This subcalculus on the underlying algebra of $\C G$ in fact turns out to be left-covariant under the pushforward of the left $D^{\vee}(G)$-coaction on $\Omega^{1}_{\mathcal{C}, \pi}$ along the Hopf algebra surjection $p_{2}: D^{\vee}(G) \twoheadrightarrow \C G$, since indeed (it suffices by \cite[Lemma 2.25]{BegMa} to check that) \begin{align} \label{eq_left_kG_covariant_coact_check} \Delta_{L}(h\extd g) &= (p_{2} \otimes \id \otimes \id)\circ (\Delta \otimes \id)(h \extd g) \nonumber\\
&= \sum_{f, k \in G, c \in \mathcal{C}} \delta_{f, e} \,  f^{-1}hgf \otimes \delta_{f^{-1}kc^{-1}} \otimes ck^{-1}hgkc^{-1} \otimes E_{ci} \otimes k^{-1}g^{-1}k \triangleright F^{ci} \nonumber\\
& \quad  - \delta_{f, e} \,  f^{-1}hgf \otimes \delta_{f^{-1}k} \otimes k^{-1}hgk \otimes E_{ci}{}^{ci}= hg \otimes h \extd g = h_{1}g_{1} \otimes h_{2} \extd g_{2} , \end{align} and by cocommutativity of $\C G$, then this calculus is automatically bicovariant. By previous discussion, the calculus must be generated by a group 1-cocycle, given by (from \cite[Theorem 1.47]{BegMa})\[ \zeta: G \to \Lambda^{1}_{G}, \quad g \mapsto g^{-1} \extd g = \sum_{f \in G, c \in \mathcal{C}}\delta_{f c^{-1}} \otimes e \otimes E_{ci} \otimes f^{-1}g^{-1}f \triangleright F^{ci} - \delta_{f} \otimes e\otimes E_{ci}{}^{ci},\] and the space of left-invariant 1-forms $\Lambda^{1}_G$ of this calculus will be simply sums of the above. Its right $\C G$-crossed module structure is then given by the right action \begin{equation} \label{eq_right_action_base_calc_kG}g^{-1} \extd g \triangleleft h = h^{-1}g^{-1}( \extd g)h = \sum_{f \in G, c \in \mathcal{C}} \delta_{fc^{-1}} \otimes e \otimes E_{ci} \otimes f^{-1}(h^{-1}g^{-1} - h^{-1})f \triangleright F^{ci}, \end{equation} where we similarly use the right $\C G$-module structure on $\Omega^{1}_{\mathcal{C}, \pi}$ induced by the pullback along $\C G \cong D^{\vee}(G)^{\C (G)} \hookrightarrow D^{\vee}(G)$, and the trivial right $\C G$-coaction. This crossed-module $\Lambda^{1}_G$ is then isomorphic to $\Lambda^{1}_{\tilde\pi}$ as in \eqref{eq_calc_kG_rep} with $\tilde\pi$ the $G$-representation induced by the irrep $\pi$ of $C_{G}$ (or equivalently the restriction to $G$ of the $D(G)$ representation $V_{\CC,\pi}$), by 
 \[ \Lambda^{1}_G\cong \Lambda^{1}_{\tilde\pi}, \quad \sum_{f \in G, c \in \mathcal{C}}\delta_{f c^{-1}} \otimes e \otimes E_{ci} \otimes f^{-1}g^{-1}f \triangleright F^{ci} - \delta_{f} \otimes e\otimes E_{ci}{}^{ci} \mapsto \sum_{c \in \mathcal{C}}E_{ci} \otimes (g^{-1} - 1) \triangleright F^{ci}.\] This is clearly a $\C G$-crossed module map for the trivial coactions and the actions \eqref{eq_calc_kG_rep} and \eqref{eq_right_action_base_calc_kG} (and moreover intertwines the identity map $\C G \to \C G$ with the 1-cocycles).  $\Omega^{1}$ is isomorphic to $\Omega^{1}_{\tilde\pi}$ with structure maps \eqref{eq_calc_kG_rep} with the following matrix elements \[ \tilde\pi(g)^{ci}{}_{dj} = \delta_{c, gdg^{-1}}\pi(\zeta_{d}(g))^{i}{}_{j}.\] We note that this calculus on the base is often non-zero for pairs $(\mathcal{C}, \pi)$, since $\tilde\pi$ here is usually not the trivial representation, i.e. the representation is trivial iff the pair is of the form the pair $(\mathcal{C}, 1)$ for $\mathcal{C}$ a conjugacy that class contains only one-element, i.e. consists of an element in $Z(G)\backslash \{ e\}$ (neglecting the pair $(\{e \}, 1)$). \endproof

Note that if $(\mathcal{C}, \pi)$ has $\pi \neq 1$, then the map $\textup{ver}$ in \eqref{eq_exact_bundle_def} is the zero map, and the horizontal forms are simply given by $\Omega^{1}_{\mathcal{C}, \pi}$, while for pairs $(\mathcal{C}, 1)$, $\textup{ver}$ is given by \[ \textup{ver}: \Omega^{1}_{(\mathcal{C}, 1)} \to D^{\vee}(G) \otimes \Lambda^{1}_{\mathcal{C}}, \quad \delta_{g} \otimes h \otimes E_{c}{}^{d} \mapsto \delta_{g} \otimes h \otimes e_{c},\] by the above and \eqref{eq_map_ver} (where we drop the $i, j$ indices for $\textup{End}(V_{(\mathcal{C}, 1)})$ basis as $\pi = 1$), and the horizontal forms are given by $D^{\vee}(G) \otimes \textup{span}\{ E_{c}{}^{d} - E_{c}{}^{d'}\, | \, c,d ,d' \in \mathcal{C}\} \subseteq \Omega^{1}_{\mathcal{C}, \pi}$ (where as mentioned above even though the $C_{G}$-irrep is trivial, the induced $G$-rep is often non-trivial, and given by adjoint $G$-action on $\C \mathcal{C}$, and as such $\Omega^{1} \cong \C G \otimes \widetilde{\pi}(\C G^{+})$ is non-zero). 
 
\begin{remark}We  note that the above bundle is never `strong' in the sense of \cite[Section 5.4.2]{BegMa}, i.e. that it does not have the classical feature that $D^{\vee}(G) \Omega^{1} = \Omega^{1}D^{\vee}(G)$: using the fact that strongness is equivalent to the fact that $(\Omega^{1}_{hor})^{\C (G)} = \Omega^{1}$ by \cite[Corollary 5.53]{BegMa} (where $\Omega^{1} \subseteq (\Omega^{1}_{hor})^{\C (G)}$ by definition of the base $A$ and $\Omega^{1}_{hor}$), then we firstly note the following isomorphism $(\Omega^{1}_{\mathcal{C}, \pi})^{\C (G)} \cong (D^{\vee}(G) \otimes V_{\mathcal{C}, \pi})^{\C (G)} \cong ((\C G \otimes \C (G))\otimes V_{\mathcal{C}, \pi})^{\C (G)} \cong \C G \otimes V_{\mathcal{C}, \pi}$, given on the basis by 
\[ \sum\nolimits_{f \in G} \delta_{f} \otimes f^{-1}gf \otimes f^{-1} \triangleright E_{ci}{}^{dj}\mapsfrom g \otimes E_{ci}{}^{dj}, \] 
with $f \triangleright E_{ci}{}^{dj}$ from \eqref{eq_right_comodule_EndVCpi}. By the form of $\Omega^{1}_{hor}$ further above, we thus have that $(\Omega^{1}_{hor})^{\C (G)}  \cong \C G \otimes V_{\mathcal{C}, \pi}$ if $\pi \neq 1$, and $(\Omega^{1}_{hor})^{\C (G)}  \cong \C G \otimes \textup{span}\{ E_{c}{}^{d} - E_{c}{}^{d'}\, | \, c,d ,d' \in \mathcal{C}\}$ otherwise (with basis by the above). We also have that $\Omega^{1} \subseteq \C (G) \otimes \textup{End}(V_{\widetilde{\pi}})$ (where $V_{\widetilde{\pi}} = V_{\mathcal{C}, \pi}$ as a vector space) given inside $ \Omega^{1}_{\mathcal{C}, \pi}$ by the span of the elements 
\begin{align*}  g \otimes e^{h} =  gh^{-1}\extd h & \mapsto  \sum_{f \in G, c\in \mathcal{C}} \delta_{f} \otimes f^{-1}gf \otimes E_{ci} \otimes (c^{-1}f^{-1}h^{-1}fc - e) \triangleright F^{ci} \\
&= \sum_{f \in G, d \in \mathcal{C}} \delta_{f} \otimes f^{-1}gf \otimes f^{-1} \triangleright (E_{di} \otimes (d^{-1}h^{-1}d - e) \triangleright F^{di}), \end{align*} by \eqref{eq_hdg}, and where in particular $\Omega^{1} \cong \C (G) \otimes \textup{End}(V_{\widetilde{\pi}})$ iff $\widetilde{\pi}$ is an irrep of $G$. So, for $\pi \neq 1$ and $\widetilde{\pi}$ not an irrep, the bundle is not strong as seen immediately by looking at the dimensions of the two spaces. One can also show that if $\CC=\{e\}$ then the bundle is strong for all $\pi$. 
\end{remark}

For a general bicovariant calculus $\Omega^{1}_{D^{\vee}(G)}$ on $D^{\vee}(G)$, the induced calculi on the structure group and base is obtained by breaking $\Omega^{1}_{D^{\vee}(G)}$ down into its coirreducible calculi $\Omega^{1}_{\mathcal{C}, \pi}$ and using the above. Remarkably, just as in the case of the Poincar\'e group $\R^{1,3}$ which gets its standard translation invariant calculus (even though $P/L$ is not a group), the fixed subalgebra $\C  G\hookrightarrow D^{\vee}(G)^{\C (G)}$ gets an induced calculus which is bicovariant. The origin of this `translation-invariance' on the base is by the `larger symmetry' it satisfies, arising from the differentiable homogeneous structure of the bundle, i.e. we saw that the base $\C G$ has an induced $D^{\vee}(G)$-coaction by \eqref{eq_kG_DveeG_coaction_symmetry}, and thus the analogue of \eqref{eq_def_differentiable} is well-defined, and as part of this $\Omega^{1}$ thus has the following well-defined induced $D^{\vee}(G)$-coaction 
\begin{equation} \label{eq_DvvG_coaction_base_calc}\Delta_{L}(h \extd g) = h^{1}g^{1} \otimes h^{2}\extd g^{2} = \sum_{f \in G}\delta_{f} \otimes f^{-1}hgf \otimes f^{-1}hf \extd (f^{-1}gf), \quad g, h \in G\, , \end{equation} (where in the above $g^{1} \otimes g^{2}$ denotes the coaction \eqref{eq_kG_DveeG_coaction_symmetry}). The pushforward of these two coactions to $\C G$ is precisely the left coregular coaction of $\C G$ on itself and the left-covariance condition \eqref{eq_left_kG_covariant_coact_check} respectively. More generally, we say that a calculus $\Omega^1$ on $\C G$ is  $D^\vee(G)$-covariant calculi if \eqref{eq_DvvG_coaction_base_calc} is well-defined.

\begin{lemma} \label{lem_DG_covariant_calc_kG}
$\Omega^{1}$ on $\C G$ is $D^\vee(G)$-covariant iff $\Omega^{1} = \Omega^{1}_{\rho}$ , i.e. $\Omega^{1}$ is of the form (\ref{eq_calc_kG_rep}) given by a representation $\rho$ of $G$. In this case, if $\Omega^{1}_{D^{\vee}(G)}$ is a bicovariant calculus on $D^{\vee}(G)$ then the coaction \eqref{eq_kG_DveeG_coaction_symmetry} is differentiable with respect to $\Omega^{1}_{D^{\vee}(G)}$ iff the algebra map $\Phi:\C G \hookrightarrow D^{\vee}(G)$ in \eqref{eq_inclusion_A_second_bundle} is differentiable.
\end{lemma}

\proof Given a calculus $\Omega^{1}$ on $\C G$ for which the $D^{\vee}(G)$-coaction \eqref{eq_DvvG_coaction_base_calc} is well-defined, then we have already explained in the above why it is then automatically bicovariant (i.e. by pushforwarding the coaction to a $\C G$-coaction and noting that the resulting condition is precisely that of \eqref{eq_bicovariant_def}). So, by Section~\ref{section_quantum_diff_structures_intro} we have that $\Omega^{1} \cong \C G \otimes \C G^{+}/I$ for a subcrossed module $I$ of $\C G^{+}$. Next, we note that on the presentation $\C G \otimes \C G^{+}/I$ of the calculus, then the coaction \eqref{eq_DvvG_coaction_base_calc} takes the form $\Delta_{L}(h \otimes e^{g})) = \Delta_{L} (hg^{-1}dg) = \sum_{k \in G} \delta_{k} \otimes k^{-1}hk \otimes k^{-1}hg^{-1}k d(k^{-1}gk)  = \sum_{k \in G} \delta_{k} \otimes k^{-1}hk \otimes k^{-1}hk \otimes e^{k^{-1}gk}$ for the quotient crossed-module map $\varpi: \C G^{+} \twoheadrightarrow \C G^{+}/I$ defining the bicovariant calculus. As such, $\C G^{+}/I$ is a $D^{\vee}(G)$-comodule by \begin{equation} \label{eq_second_form_of_coaction}\Delta_{L} \varpi(\pi_{\epsilon}(g)) = \sum_{f \in G}(\delta_{f} \otimes e) \otimes \varpi(\pi_{\epsilon}(f^{-1}gf)),\end{equation} (and the $D^{\vee}(G)$-comodule structure on $\C G \otimes \C G^{+}/I$ above is thus simply the tensor product one). So, we see that \eqref{eq_DvvG_coaction_base_calc} is well-defined iff \eqref{eq_second_form_of_coaction} is well-defined. This latter $D^{\vee}(G)$-coaction is equivalent to the right $\C G$-crossed module structure
 \[e^{g} \triangleleft h = e^{h^{-1}gh}, \quad \Delta_{R}e^{g} =  e^{g} \otimes e. \] 
The trivial coaction is clearly always well-defined, while the right-action being well-defined is equivalent to asking $I$ to be submodule of the adjoint action. Since $I$ is also a right ideal, then it is equivalently a two-sided ideal of $\C G$. This was previously seen in Section~\ref{sec_recap_calc_kG_kofG} to correspond to $\Omega^{1}$ being of the form $\Omega^{1}_{\rho}$. 

By \eqref{eq_def_differentiable}, the second part of the differentiability criteria is that the map $ \Omega^{1} \to \Omega^{1}_{D^{\vee}(G)} \otimes \C G$ give by $h\extd g \mapsto h^{1}\extd g^{1} \otimes h^{2}g^{2}$ is also well-defined. Noting that the coaction $h^{1} \otimes h^{2}$ of \eqref{eq_kG_DveeG_coaction_symmetry} can be written as $\Phi_{\C (G)}^{-1}(h)_{1} \otimes \Phi_{\C (G)}(\Phi_{\C (G)}^{-1}(h)_{2})$ by definition, we require that  $h\extd g \mapsto \Phi_{\C (G)}^{-1}(h)_{1} \extd \Phi_{\C (G)}^{-1}(g)_{1}\otimes \Phi_{\C (G)}(\Phi_{\C (G)}^{-1}(h)_{2})\Phi_{\C (G)}(\Phi_{\C (G)}^{-1}(g)_{2})$, or equivalently  $h \extd g \mapsto \Phi_{\C (G)}^{-1}(h)_{1} \extd \Phi_{\C (G)}^{-1}(g)_{1}\otimes \Phi_{\C (G)}^{-1}(h)_{2}\Phi_{\C (G)}^{-1}(g)_{2}$ is well-defined. If the calculus on $D^{\vee}(G)$ is bicovariant then 
\[ \Phi_{\C (G)}^{-1}(h) \extd \Phi_{\C (G)}^{-1}(g) \mapsto \Phi_{\C (G)}^{-1}(h)_{1} \extd \Phi_{\C (G)}^{-1}(g)_{1}\otimes \Phi_{\C (G)}^{-1}(h)_{2}\Phi_{\C (G)}^{-1}(g)_{2}\]
 is well-defined, and so we simply require the algebra map $h\extd g \mapsto \Phi_{\C (G)}^{-1}(h) \extd \Phi_{\C (G)}^{-1}(g)$ to be well-defined. By definition,  this is just asking for differentiability of $\Phi$. \endproof

Having studied the differential geometry of the above principal bundle, we now discuss the Riemannian geometry of the base $\C G$ such that $D^{\vee}(G)$ becomes an `isometry group' analogously to the Poincar\'e group. As for any differential algebra, a  quantum metric on $\C G$ means a bimodule map  $(\ , \ ): \Omega^{1} \otimes_{\C G} \Omega^{1} \to \C G$ or `bimodule inner product', usually with quantum symmetry and nondegeneracy conditions which we do not impose at this stage. In our case, for bicovariant $\Omega^{1}$ on $\C G$ we have seen that $\Omega^1$ is  $D^{\vee}(G)$-coaction \eqref{eq_DvvG_coaction_base_calc} and we ask for $(\ ,\ )$ to be likewise $D^{\vee}(G)$-covariant, so that  $D^{\vee}(G)$ can be seen as in the role of isometry group in $\C G$ for $(\ , \ )$. 

\begin{lemma} \label{lemma_metric_crossed_module}
Let  $\Omega^{1} \cong \C G \otimes \Lambda^{1}$ be bicovariant and equip $\Omega^{1}\otimes_{\C G} \Omega^{1}$ with the descended tensor product left $D^{\vee}(G)$-comodule structure. The following are equivalent

\noindent (i)  $(\ , \ ): \Omega^{1} \otimes_{\C G} \Omega^{1} \to \C G$ is a left $D^{\vee}(G)$-covariant bimodule inner product.

\noindent (ii) $(\ , \ ): \Lambda^{1} \otimes \Lambda^{1} \to \C $ obeys
 \begin{equation}\label{eq_conditions_on_metric}(e^{g}, e^{h}) = (e^{gu}, e^{hu}) - (e^{gu}, e^{u}) - (e^{u}, e^{hu}) + (e^{u}, e^{u}), (e^{g}, e^{h}) = (e^{u^{-1}gu}, e^{u^{-1}hu}),\end{equation} 
 for all $g, h , u \in G$. 
\end{lemma}

\proof We firstly note that the coaction \eqref{eq_DvvG_coaction_base_calc} is equivalent to the following right Whitehead $G$-crossed module \begin{equation} \label{eq_DvvG_coaction_base_calc_alternative}(h\extd g) \triangleleft f^{-1}hf \, (\extd f^{-1}gf), \quad |h\extd g| = hg. \end{equation}

We then note that since $\Omega^{1}$ is a $\C G$-Hopf bimodule, it is immediate that the $\C G$-coaction on $\Omega^{1} \otimes \Omega^{1}$ descends to $\C G$-coactions on $\Omega^{1}\otimes_{\C G} \Omega^{1}$ and that $\Omega^{1}\otimes_{\C G}\Omega^{1}$ is also then a $\C G$-Hopf bimodule. The same applies for the $D^{\vee}(G)$-coaction, i.e. we simply need to check the $G$-action in \eqref{eq_DvvG_coaction_base_calc_alternative} (since the grading is included in $\C G$-Hopf bimodule structure) is compatible with $\otimes_{\C G}$, which is true for the same reasons. Explicitly, we  check that 
 \[l^{-1}hl(\extd( l^{-1}gl))l^{-1}f l \otimes_{\C G} l^{-1}ul \extd (l^{-1}vl) =  l^{-1}hl \extd(l^{-1}gl)\otimes_{\C G} l^{-1}f ul \extd( l^{-1}vl), \]
for all $h, g, f, u, v, l \in G$. Then, using the fact that $\Omega^{1}\otimes_{\C G}\Omega^{1} \cong (\C G \otimes \Lambda^{1}) \otimes_{\C G}(\C G \otimes \Lambda^{1}) \cong \C G\otimes \Lambda^{1} \otimes \Lambda^{1}$ as Hopf bimodules,  clearly $(\Omega^{1}\otimes_{\C G}\Omega^{1})^{\C G} \cong \Lambda^{1} \otimes \Lambda^{1}$. Hence, by the Hopf module lemma, we have that the bimodule inner product $(\ , \ ): \Omega^{1}\otimes_{\C G} \Omega^{1} \to \C G$ is equivalent to the right $\C G$-crossed module map $( \ , \ )|_{\Lambda^{1} \otimes \Lambda^{1}}$ on  $\Lambda^{1} \otimes \Lambda^{1}$  with the tensor product crossed module structure. Here, $\C $ has the trivial crossed module structure and $\Lambda^{1}$ is a crossed-module by \eqref{eq_calc_kG_rep}, or in terms of the $\{e^{g}\}$ span by \begin{equation} \label{eq_crossed_module_eg} e^{g}\xtriangleleft u = e^{gu} - e^{u}, \quad (e^{g})^{1} \otimes (e^{g})^{2} = e^{g} \otimes 1. \end{equation} We then note that the $\C G$-action on $\Omega^{1}$ in \eqref{eq_DvvG_coaction_base_calc_alternative} restricts to a $\C G$-action on $\Lambda^{1}$, as verified most simply by looking at elements $e^{g}$ of the span, i.e. \begin{equation} \label{eq_kG_ad_action} e^{g} \triangleleft l = (g^{-1}\extd g) \triangleleft l = l^{-1}g^{-1}l(\extd l^{-1}gl) = e^{l^{-1}gl} \in \Lambda^{1}, \quad \forall g, l \in G . \end{equation} So, $(\ , \ )|_{\Lambda^{1} \otimes \Lambda^{1}}$ is then additionally a right $G$-module map under the tensor product action on $\Lambda^{1} \otimes \Lambda^{1}$ and the trivial action on $\C $. Writing out these conditions on the span $\{ e^{g}\}_{g \in G}$ gives \[ ((e^{g})^{1}, (e^{h})^{1}) \otimes (e^{g})^{2}(e^{h})^{2} = (e^{g}, e^{h}) \otimes 1, \quad (e^{gu} -e^{u}, e^{hu} - e^{h}) = (e^{g}, e^{h}) \epsilon(u),\]\[ (e^{u^{-1}gu}, e^{u^{-1}hu}) = (e^{g}, e^{h}) \epsilon(u), \] The first condition empty and the rest give \eqref{eq_conditions_on_metric}.  \endproof

We will be interested in such $D^\vee(G)$-covariant bimodule inner products. These are automatically $\C G$-bicomodules for the same reasons as before. Given such data, we denote the `length' of the differential form $e^{g}$ by $l_{\mathcal{C}}:= (e^{g},e^{g}) \in \mathbb{C}$ for $g \in \mathcal{C}$, since this quantity only depends on the conjugacy class of $g$ in $G$ by the second relation \eqref{eq_conditions_on_metric}. We note that the first relation in \eqref{eq_conditions_on_metric} gives that \[ (e^{g}, e^{g}) = (e^{gg^{-1}}, e^{gg^{-1}}) - (e^{gg^{-1}}, e^{g^{-1}}) - (e^{g^{-1}}, e^{gg^{-1}}) + (e^{g^{-1}}, e^{g^{-1}}) = (e^{g^{-1}}, e^{g^{-1}}),  \]\begin{equation} \label{eq_length_C_Cinv} \implies l_{\mathcal{C}} = l_{\mathcal{C}^{-1}},\end{equation} since $e^{e} = 0$ (which also gives $l_{\{ e\}} = 0$). We shall say that a $D^\vee(G)$-covariant bimodule inner product is {\em regular} if all the $l_{\mathcal{C}}$'s differ for distinct conjugacy classes that are not related by inversion. 

Next we note that since $\{ e^{g} \}_{g \in G\backslash \{ e\}}$ spans $\Lambda^{1}$,  we can pick a subcollection $\{ e^{i} \}_{i \in \mathcal{B}}$, $\mathcal{B} \subseteq G\backslash \{ e\}$ to form a basis of $\Lambda^{1}$ and denote by $\{e_i\}$ the dual basis. Writing $N_{\mathcal{B}}$ for the number of elements in $\mathcal{B}$ which live in distinct conjugacy classes unrelated by inversion,  we note that a necessary condition for regularity of $(\ , \ )$ is then
\begin{equation} \label{eq_regularity_inner_prod}
N_{\mathcal{B}} \geq N_{\rm conj}-1,
\end{equation}
i.e. the calculus needs to be large enough. \color{black} We denote $(e^{i}, e^{j}) =: g^{ij} \in \C $ as the components of the $|\mathcal{B}| \times |\mathcal{B}|$-`bimodule inner product matrix' $\underline{\underline{g}}$. From \eqref{eq_crossed_module_eg} and the $\C G$-action \eqref{eq_kG_ad_action}, we can define the following matrices \begin{equation} \label{eq_rho_gamma_def}\gamma(g)^{i}{}_{j} = \langle e^{ig} - e^{g}, e_{j} \rangle, \quad \rho(g)^{i}{}_{j} = \langle e^{g^{-1}ig}, e_{j} \rangle,\end{equation} 
which satisfy the commutation relation $\gamma(g)\rho(k) = \rho(k)\gamma(k^{-1}gk)$ by the fact that $(e^{g} \xtriangleleft h) \triangleleft k = e^{k^{-1}ghk} - e^{k^{-1}hk} = e^{k^{-1}gk} \xtriangleleft k^{-1}hk = (e^{g}\triangleleft k) \xtriangleleft k^{-1}hk$. With these, we can write the $D^{\vee}(G)$-covariant condition on the basis in the more compact index form \begin{equation} \label{DveeG_covariant_bimodule_inner} \gamma(u)^{i}{}_{k} \gamma(u)^{j}{}_{l}g^{kl} = g^{ij}, \quad \rho(u)^{i}{}_{k}\rho(u)^{i}{}_{l}g^{kl} = g^{ij}, \end{equation} or more simply $\gamma(u) \underline{\underline{g}} \gamma(u)^{T} = \underline{\underline{g}}$ and similarly for $\rho$. Nondegeneracy for the bimodule inner product $(\  , \ ): \Omega^{1} \otimes_{\C G} \Omega^{1} \to \C G$ is imposed by requiring an  inverse   `metric element' $\metric \in \Omega^{1} \otimes_{\C G} \Omega^{1}$ if \[ ((\omega, \ ) \otimes \id)\metric = \omega = (\id \otimes (\ , \omega))\metric, \] 
for all $\omega\in \Omega^1$. In our case with  $\Omega^{1}$ and $(\ , \ )$ being $D^{\vee}(G)$-covariant, this exists iff the matrix $\underline{\underline{g}}$ above is invertible. In this case, we denote the matrix elements of $(\underline{\underline{g}}^{-1})_{i, j} =: g_{ij}$ so that $\metric = g_{ij} e^{i} \otimes e^{j} \in \Lambda^{1} \otimes \Lambda^{1}$. The $D^{\vee}(G)$-covariance is also clear from this by inverting \eqref{DveeG_covariant_bimodule_inner} to get $(\gamma(u)^{-1})^{T}\underline{\underline{g}}^{-1}\gamma(u)^{-1} = \underline{\underline{g}}^{-1}$ for all $u \in G$ and similarly for $\rho$. These conditions can be written in an analogous index form to \eqref{DveeG_covariant_bimodule_inner} as $\gamma(v)^{k}{}_{i}\gamma(v)^{l}{}_{j}g_{kl} = g_{ij}$ and $\rho(v)^{k}{}_{i}\rho(v)^{l}{}_{j}g_{kl} = g_{ij}$ for all $v \in G$. 

Finally, we briefly discuss connections on the bundle of Proposition~\ref{prop_bundle_over_kG}: `Ehresmann connections' on quantum principal bundles are defined in the obvious way as a choice of equivariant projector $\Pi: \Omega^{1}_{P} \to \Omega^{1}_{P}$ that defines a complement to $\Omega^{1}_{hor}$, i.e. the vertical forms \cite{BegMa}. By \cite[Proposition 5.41]{BegMa}, these are in correspondence with the dual analogue of `principal $H$-connections', i.e. right $H$-comodule maps $\omega: \Lambda^{1} \to \Omega^{1}_{P}$ such that $\textup{ver} \circ \omega = 1 \otimes \id$. So, from the form of $\textup{ver}$ Proposition~\ref{prop_bundle_over_kG}, we observe that a bundle $\C G \hookrightarrow D^{\vee}(G) \twoheadrightarrow \C (G)$ does not admit a connection if $\Omega^{1}_{D^{\vee}(G)}$ does not contain a pair $(\mathcal{C}, 1)$ in its decomposition. Hence, to construct covariant derivatives on associated quantum vector bundle of interest using arbitrary $\Omega^{1}_{\mathcal{C}, \pi}$ on $D^{\vee}(G)$, we proceed to these directly  rather than as associated to a connection on the principal bundle. 

Here, as for any differential algebra and bimodule $F$ over it, in our case a left covariant derivative on a $\C G$ bimodule $F$ is a $\C $-linear map $\nabla_{F}: F \to \Omega^{1} \otimes_{\C G} F$ obeying the left Leibniz rule \cite[Chapter 3]{BegMa}\[ \nabla_{E}(g \cdot f) =  \extd g \otimes_{\C G}f + g \cdot \nabla_{F}f, \quad g \in G, f \in F ,  \] and its curvature is the left $\C G$-module map $R_{F}: F \to \Omega^{2} \otimes_{\C G} F$ defined by \begin{equation} \label{eq_curvature}R_{F} = (\extd \otimes \id - \id \wedge \nabla_{F}) \circ \nabla_{F}. \end{equation} Of particular interest is when $F$ is a $\C G$-bimodule  equipped with a covariant derivative $\nabla_{F}$ for which the map $\sigma_{F}: F \otimes_{\C G} \Omega^{1} \to \Omega^{1} \otimes_{\C G} F$ given by \[\sigma_{F}(f \otimes \extd g) = \nabla_{F}(f \cdot g) - (\nabla_{F}f)\cdot g\] is well-defined (allowing for covariant derivatives on tensor product bundles), and for which the curvature is additionally a bimodule map for identical reasons (and this feature is also known as `Riemann compatibility' \cite{BegMa}). 

In our cases of interest, $F$ will be additionally $D^\vee(G)$-covariant  in a manner compatible with the bimodule structure, i.e. $F \in \, _{\C G}^{D^{\vee}(G)}\mathcal{M}_{\C G}$ in the notation of \cite{BegMa}. Noting that a left Whitehead $G$-crossed module $V$ as in \eqref{eq_comodule_DveeG} is equivalent to a right one by \begin{equation} \label{eq_left_white_headcrossed_module}\Delta_{L}v = \sum_{g \in G} \delta_{g} \otimes |g^{-1}\la v|^{-1} \otimes g^{-1}\la v = \sum_{g \in G} \delta_{g} \otimes g^{-1}|v|^{-1}g \otimes g^{-1}\la v, \quad v \in V,\end{equation} where here the corresponding right $G$-crossed module has  right action is $v\triangleleft g:= g^{-1} \triangleright v$ and grading given by the inverse of that in \eqref{eq_comodule_DveeG}, then the compatibility $F$ with the bimodule structure comes down to \begin{equation} \label{eq_hopf_bimodule_F} \sum_{l \in G} \delta_{l} \otimes l^{-1}|(g \cdot f \cdot h)|l \otimes (g \cdot f \cdot h) \triangleleft l = \sum_{l \in G}\delta_{l} \otimes l^{-1}g|f|hl \otimes l^{-1}gl \cdot (f \triangleleft l) \cdot l^{-1}hl,  \end{equation} 
i.e. $|g \cdot f \cdot h| = g|f|h$, and $(g \cdot f \cdot h) \triangleleft l = l^{-1}gl \cdot (f \triangleleft l) \cdot l^{-1}hl$, where $\cdot$ highlights the bimodule actions. When this holds, it is natural to ask that  $\nabla_{F}$ be $D^\vee(G)$-covariant with respect to tensor product coaction on $\Omega^{1} \otimes_{\C G}F$. In this scenario, since $F$ is a left $\C G$-Hopf module by pushforward of the $D^\vee(G)$ coaction, then $F \cong \C G \otimes \, {}^{\C G}F$ as a free left $\C G$-module by the Hopf module-Lemma. By this freeness, it follows as in the general argument\cite[(3.14)]{BegMa} that connections on $F$ can all be expressed as \[ \nabla_{F}(g\cdot f) = \extd g \otimes_{\C G} f + g \cdot \nabla_{F}^{L}f;\quad \nabla_{F}^{L}: {}^{\C G}F \to \Lambda^{1} \otimes \, {}^{\C G}F, \]
for all $g \in \C G, f \in \, {}^{\C G}F$. Here $\nabla_{F}^{L}$ is the right $\C G$-module map given by restriction of $\nabla$, so that $ (\nabla_{F}^{L} f) \triangleleft g = \nabla_{F}^{L}(f \triangleleft g)$ for all $g \in G$. If ${}^{\C G}F$ is finite dimensional with basis $\{ f^{i}\}$ and dual $\{ f_{i}\}$, writing $\nabla_{F}^{L}$ in terms of its Christoffel symbols  \[ \nabla^{L}_{F} f^{i}= - \Gamma^{i}{}_{jk}e^{j} \otimes f^{k} , \quad \Gamma^{i}{}_{jk} \in \C ,\]  gives that $D^\vee(G)$-covariant connections are all of the following form on the basis \[ \nabla_{F}(g \otimes  f^{i}) = -\big( -\delta_{k}^{i}\partial_{j}g  + g\Gamma^{i}{}_{jk}\big) \, e^{j} \otimes f^{k}, \]
for all $g\in G$. Here, the partial derivatives $\{ \partial_{i}: \C  G \to \C  G \}_{i \in \mathcal{B}}$ for the calculus on $\C  G$ are defined by $\extd g = \partial_{i}g \otimes e^{i}$, and since $\extd g=g\tens e^g$, are  given by \begin{equation} \label{eq_form_of_partial} \partial_{i}g = g \langle e_{i}, e^{g} \rangle \in \C G.\end{equation}The $\Gamma$ coefficients are required to satisfy \begin{equation} \label{eq_Gamma_covariant} \rho^{F}(g)^{i}{}_{l}\Gamma^{l}{}_{jk} = \Gamma^{i}{}_{ml} \rho(g)^{m}{}_{j}\rho^{F}(g)^{l}{}_{k}, \end{equation}
for all $g \in G$, where we similarly defined $f^{j}\rho^{F}(g)^{i}{}_{j}:= f^{i}\triangleleft g$. 

\begin{lemma}\label{lembicnablaF} A left $\C G$-covariant left covariant derivative $\nabla_{F}$ is a bimodule covariant derivative with respect to the bicovariant calculus $\Omega^{1}_{\rho}$ on $\C G$ for some representation $\rho$ iff  $f^{i} \otimes e^{g}\mapsto \gamma^{F}(g^{-1})^{i}{}_{l}\Gamma^{l}{}_{pm}\gamma^{F}(g)^{m}{}_{k}f^{k}$ is well-defined for all $g$, where $f^{j}\gamma^{F}(g)^{i}{}_{j}:= f^{i} \xtriangleleft g$ defines $\gamma^{F}$ as the matrix for $\xtriangleleft$ the action from the $\C G$-crossed module structure on $F$. 
 \end{lemma}
\proof  This follows from the discussion around \cite[Prop. 3.73]{BegMa} that the condition we need is that
\[ \sigma^{L}(f^{i} \otimes e^{g}) - e^{g} \otimes f^{i}=  \nabla^{L}_{F}f^{i} \epsilon(g) -  \nabla^{L}_{F}(f^{i} \xtriangleleft g^{-1}) \xtriangleleft g = -\Big(\Gamma^{i}{}_{jk} \epsilon(g)- \gamma^{F}(g^{-1})^{i}{}_{l}\Gamma^{l}{}_{pm}\gamma(g)^{p}{}_{j}\gamma^{F}(g)^{m}{}_{k}\Big)\, e^{j} \otimes f^{k}\] 
is a well-defined map for all $g \in G$. So for $\Omega^{1}_{\rho}$ generated by a representation $\rho$, if $e^{g} = e^{h}$ (and as such $\rho(g) - \id = \rho(h) - \id$ implying $\rho(g) = \rho(h)$), then $e^{lg} = \rho(lg) - \id = \rho(l)\rho(g) - \id = \rho(l)\rho(h) - \id = \rho(lh) - \id=  e^{lh}$ for all $l \in G$ and as such the association $e^{g} \mapsto \gamma(g)$ is well-defined, and so the above map is well-defined iff the condition stated holds. \endproof

We also note that $\sigma_{F}$ is automatically a $D^{\vee}(G)$-comodule map by its definition, and similarly to before, it can be quickly shown explicitly by noting that by \eqref{eq_hopf_bimodule_F} we have the compatibility $(f \xtriangleleft g) \triangleleft k = (f \triangleleft k) \xtriangleleft k^{-1}gk$, which in terms of the matrices defined above gives the same commutation relation $\gamma^{F}(g)\rho^{F}(k) = \rho^{F}(k)\gamma^{F}(k^{-1}gk)$. 

Next, the curvature of $\nabla_{F}$ is given in index-form by \begin{equation} \label{eq_index_form_curvature} R_{F}(g \otimes f^{i}) = g \cdot R_{F}^{L}(f^{i}) =  g\Big( -\Gamma^{i}{}_{jk}\extd e^{j} - \Gamma^{i}{}_{jm}\Gamma^{m}{}_{lk} \, e^{j} \wedge e^{l}\Big) \otimes_{\C G}(e \otimes f^{k}),\end{equation} 
for all $g\in G$, where we recall that $\extd e^{g} = - e^{g} \wedge e^{g}$ for any choice of $\Lambda^{2}$. $R_{F}$ is also automatically a $D^{\vee}(G)$-comodule map by its definition.

\begin{proposition} $R_{F}$ is left $\C G$-module (and thus bimodule) map iff 
\begin{align*} &\Gamma^{i}{}_{kj} (e^{k} \wedge e^{k}) - \Gamma^{i}{}_{kl}\Gamma^{l}{}_{mj} e^{k} \wedge e^{m} \\
&= \gamma^{F}(h^{-1})^{i}{}_{k}\Big(\Gamma^{k}{}_{lm}\gamma(h)^{l}{}_{p}(e^{p} \wedge e^{p}- e^{h} \wedge e^{p} - e^{p} \wedge e^{h}) - \Gamma^{k}{}_{ln}\Gamma^{n}{}_{pm} e^{q} \wedge e^{s} \gamma(h)^{l}{}_{q} \gamma(h)^{p}{}_{s} \Big) \gamma^{F}(h)^{m}{}_{j} , \end{align*}
for all $h \in G$.
\end{proposition}
\proof We need
\[ R_{F}((g \otimes f^{i})\cdot h) = g \cdot R^{L}_{F}(f^{i})\cdot h = -g\Big(\Gamma^{i}{}_{jk}\extd e^{j} + \Gamma^{i}{}_{jm}\Gamma^{m}{}_{lk} \, e^{j} \wedge e^{l}\Big) \otimes_{\C G}(h \otimes f^{n}\gamma^{F}(h)^{k}{}_{n}), \]
which means
\[ gh R_{F}^{L}(f^{j})\gamma^{F}(h)^{i}{}_{j} = -gh\Big(\Gamma^{i}{}_{jk}(\extd e^{j})\xtriangleleft h + \Gamma^{i}{}_{jm}\Gamma^{m}{}_{lk} \, (e^{j} \wedge e^{l})\xtriangleleft h\Big) \otimes_{\C G}(e \otimes f^{n}\gamma^{F}(h)^{k}{}_{n}).\]
Using \cite[Prop. 3.3]{MaTao} again to note that $(\extd e^{i}) \xtriangleleft h= \extd (e^{i}\xtriangleleft h) + e^{h} \wedge (e^{i} \xtriangleleft h) + (e^{i} \xtriangleleft h) \wedge e^{h}$ and $(e^{i} \wedge e^{j}) \xtriangleleft h = (e^{i} \xtriangleleft h) \wedge (e^{j} \xtriangleleft h)$ by \cite[Prop. 3.4]{MaTao}, then we simplify to the stated condition. \endproof

Of particular interest in our case arecovariant derivatives on $E'$ from Section~\ref{sec_second_bundle_construct} that are $D^\vee(G)$-covariant. As seen immediately from \eqref{eq_coact_Eprime} pushforwarded to $\C G$, we have that ${}^{\C G}E' \cong W$. As noted in Section~\ref{secpre}, $(D^{\vee}(G) \otimes W)^{\C (G)} \cong E'$ has $\C G$-bimodule structure induced by the first tensor factor (and the inclusion \eqref{eq_inclusion_A_second_bundle}), which under the isomorphism \eqref{eq_trivial_associated_bundle_Hprime} gives the obvious bimodule structure on $\C G \otimes W$, i.e.  \[ g \cdot (h \otimes w) = gh \otimes v, \quad (h \otimes w) \cdot g = hg \otimes w , \]  Thus, we see that $w \xtriangleleft g = w$, implying $\gamma^{E'}(g) = 1$. Thus by the above Lemma, we have that $\sigma^{L}$ is always well-defined, and hence all $\C G$-covariant left covariant derivates are automatically bimodule covariant derivatives. 

The second quantum vector bundle of interest is $\Omega^{1}$ itself, where we take the basis $f^{i} = e^{i}$ in the above Lemma, and as such $\gamma^{F}(G) = \gamma(g)$ and $\rho^{F}(g) = \rho(g)$, and so the relation \eqref{eq_Gamma_covariant} comes out more cleanly as just $\Gamma^{i}{}_{jk} = \rho(g)^{i}{}_{l}\Gamma^{l}{}_{mn}\rho(g^{-1})^{m}{}_{j}\rho(g^{-1})^{n}{}_{k}$. Since we saw that the association $e^{g} \mapsto \gamma(g)$ was well-defined, then all $\C G$-covariant left covariant derivatives are again bimodule covariant derivatives. 

Next, if $\Omega^{1}$ has a $D^\vee(G)$-covariant metric $\metric$ then of particular interest are $D^\vee(G)$-covariant quantum Levi-Civita connections (QLC's), defined as those for which \[ \nabla_{\Omega^{1} \otimes \Omega^{1}}\metric := (\nabla_{\Omega^{1}} \otimes \id + (\sigma_{\Omega^{1}} \otimes \id)\circ (\id \otimes \nabla_{\Omega^{1}})) \metric = 0, \quad T_{\nabla} := \wedge \nabla_{\Omega^{1}} - \extd = 0, \] where the left $\C G$-module map $T_{\nabla}: \Omega^{1} \to \Omega^{2}$ plays the role of torsion \cite{BegMa}, i.e. $\nabla_{\Omega^{1}}$ is metric compatible and torsion free. Since we are dealing with $D^\vee(G)$-covariant $\Omega^{1}$ and $\nabla_{F}$, then we note that on the basis \[ T_{\nabla}(g \otimes e^{i}) = g\cdot (\wedge \nabla_{\Omega^{1}}^{L} - d)(e^{i}) = - \Gamma^{i}{}_{jk} g e^{j} \wedge  e^{k} - g\extd e^{i},  \] and so torsion freeness means $\extd e^{i} (= - e^{i} \wedge e^{i})= - \Gamma^{i}{}_{jk} e^{j} \wedge  e^{k} \in \Lambda^{2}$. Similarly,  metric compatibility is equivalent to \begin{equation} \label{eq_metric_compatible_index} g_{lk}\Gamma^{l}{}_{ij} + g_{jl}\Gamma^{l}{}_{ik} - g_{lm}\Gamma^{l}{}_{ij}\Gamma^{m}{}_{pk}\epsilon(p) + g_{lm}\Gamma^{m}{}_{pk}\gamma(p^{-1})^{l}{}_{a}\Gamma^{a}{}_{bc}\gamma(p)^{b}{}_{i}\gamma(p)^{c}{}_{j}  = 0 , \end{equation} 
where we note that $p \in \mathcal{B}$ here. There are also weaker notions of QLCs in \cite[Chapter 8]{BegMa} which could be similarly analysed.

We finalise the discussion by bringing the $\star$-structure described around \eqref{eq_star_struc_DveeG} to $\star$-structures on the above. $\Omega^{1}$ can be made into a $\star$-calculus by simply defining $(h\extd g)^{\star}:= g^{\star}\extd h^{\star}$, which gives $(e^{g})^{\star} = -e^{g}$. We also note that with these $\star$-structures, the `homogeneous structure' arising from the $D^{\vee}(G)$-coaction on $\C G$ by \eqref{eq_kG_DveeG_coaction_symmetry} is unitary as a comodule algebra in the sense that $\Delta_{L} \circ \star = (\star \otimes \star)\Delta_{L}$. Then, defining a bimodule inner product $(\ , \ )$ on $\C G$ as $\star$-compatible if $(\ , \ )\circ \star = \star\circ (\ , \ )$, where $\star$ extends to tensor products by $\star$ on each factor and an additional flip (as denoted by $\dagger$ in \cite{BegMa}). Explicitly,  \begin{equation} \label{eq_star_compatible_inner} (h\extd g, u\extd v)^{\star} = ((u\extd v)^{\star}, (h\extd g)^{\star}), \quad h, g, u, v \in G, \end{equation} then a biinvariant bimodule inner product is $\star$-compatible iff its restriction $(\ , \ )|_{\Lambda^{1} \otimes \Lambda^{1}}$ is $\star$-compatible, i.e. satisfies $(e^{g}, e^{h})^{*} = (e^{h}, e^{g})$: the $\impliedby$ direction is clear since $(e^{g}, e^{h})^{\star} = (e^{g}, e^{h})^{*}$, while the $\implies$ direction is by the following, \begin{align*} ((&u\extd v)^{\star}, (h\extd g)^{\star}) = (e^{v}v^{-1}u^{-1}, e^{g}g^{-1}h^{-1}) = v^{-1}u^{-1}g^{-1}h^{-1}(e^{u^{-1}} - e^{v^{-1}u^{-1}}, e^{g}) \\
& =  v^{-1}u^{-1}g^{-1}h^{-1}[\big( -(e^{u}, e^{gu} ) + (e^{u}, e^{u}) \big) - \big( -(e^{uv}, e^{guv} ) + (e^{uv}, e^{uv}) \big)]\\
& = v^{-1}u^{-1}g^{-1}h^{-1}[-(e^{u}, e^{gu} ) + (e^{u}, e^{u}) + \big( (e^{u}, e^{gu} )  + (e^{v}, e^{guv}) + (e^{uv}, e^{v}) - (e^{uv}, e^{v})\big) \\
& \hphantom{xxxxxxxx}- \big( (e^{u}, e^{u}) +(e^{uv}, e^{v}) + (e^{v}, e^{uv}) - (e^{v}, e^{v}) \big)]\\
& = v^{-1}u^{-1}g^{-1}h^{-1} [(e^{v}, e^{guv}) - (e^{v}, e^{uv})] = \Big(hguv [(e^{guv}, e^{v}) - (e^{uv}, e^{v})]\Big)^{\star}\\
& = \Big( hg(e^{g}uv, e^{v})\Big)^{\star}= (hge^{g}, uve^{v})^{\star} = (h\extd g, u\extd v)^{\star}, \end{align*} using  \eqref{eq_conditions_on_metric} accordingly for the second and third equalities. We note that in particular $(e^{g}, e^{g}) \in \mathbb{R}$ for all $g \in G$, and so all lengths $l_{\mathcal{C}}$ are real. We also see that $(\ , \ )$ takes values in $\mathbb{R}$ iff $(\ , \ )$ is symmetric, or in index form $g^{ij} = g^{ji}$ for all $i, j \in \mathcal{B}$, or equivalently the inverse matrix $g_{ij}$ for $\metric$ is symmetric.  An interesting feature in this latter case is that specifying the lengths fully specifies the metric as follows.

\begin{lemma} \label{eq_lemma_squared_lengths_metric}
A $D^{\vee}(G)$-covariant $\star$-compatible bimodule inner product on $\Omega^{1}$ is fully determined by real squared lengths $\{l_{\mathcal{C}} \}_{\mathcal{C} \neq \{ e\}}$ by
 \[ (e^{h}, e^{u}) = -\frac{1}{2}(l_{\mathcal{C}_{hu^{-1}}} - l_{\mathcal{C}_{h}} - l_{\mathcal{C}_{u}}), \] (and with $l_{\{e \}} = 0$ as previously seen). 
\end{lemma}
\proof \eqref{eq_conditions_on_metric} gives that $l_{\mathcal{C}_{g}} = l_{\mathcal{C}_{gu}} - (e^{gu}, e^{u}) - (e^{u}, e^{gu}) + l_{\mathcal{C}_{u}}$, and as such rearranging this and taking $g = hu^{-1}$ gives that $(e^{h}, e^{u}) + (e^{u}, e^{h})  = (e^{(hu^{-1})u}, e^{u}) + (e^{u}, e^{(hu^{-1})u}) = -(l_{\mathcal{C}_{hu^{-1}}} - l_{\mathcal{C}_{h}} - l_{\mathcal{C}_{u}})$, $\forall h, u \in G $. So, if $(\ , \ )$ is real (and thus symmetric), then the equation stated follows. We note that in fact we need only specify the squared lengths $l_{\mathcal{C}_{i}}$ for $i \in \mathcal{B}$, i.e. $N_{\mathcal{B}}$-numbers. \endproof

Next, given a star structure on $F$, then a left-bimodule covariant derivative $\nabla_{F}$ is said to be `$\star$-preserving' if \cite{BegMa} \[ \nabla_{F}(f^{\star}) = \sigma_{F} \circ \star \circ \nabla_{F}(f), \quad \forall f \in F. \] We note that for left-invariant connections, \begin{align*} \nabla_{F}&((g \otimes f)^{\star}) - \sigma_{F}\nabla_{F}(g \otimes f)^{\star} \\
& = \big( \nabla_{F}^{L}(f^{\star}) \cdot g^{-1} + \sigma_{F}(f^{\star} \otimes \extd g^{-1})\big) - \sigma_{F}\big( f^{\star} \otimes \extd g^{-1} + \nabla^{L}_{F}(f)^{\star} \cdot g^{-1}\big)\\
& = \nabla^{L}_{F}(f^{\star})\cdot g^{-1} - \sigma_{F}(\nabla_{F}(f^{\star})\cdot g^{-1}) =  \nabla_{F}^{L}(f^{\star})\cdot g^{-1} - \sigma^{L}_{F}(\nabla^{L}_{F}(f^{\star}))\cdot g^{-1}. \end{align*} Hence, similarly, asking for star compatibility of $(\nabla_{F}, \sigma_{F})$ is equivalent to asking for $\star$-compatibility of the restriction $(\nabla_{F}^{L}, \sigma_{F}^{L})$. The latter then happens iff 
\begin{equation*}  -\Gamma^{l}{}_{jk} \, \langle f_{l}, (f^{i})^{\star}\rangle - (\Gamma^{i}{}_{jl})^{*} \, \langle f_{k},  (f^{l})^{\star} \rangle = (\Gamma^{i}{}_{uv})^{*}\big( -\Gamma^{q}{}_{jk} \epsilon(u) + \gamma^{F}(u^{-1})^{q}{}_{l}\Gamma^{l}{}_{pm} \gamma^{F}(u)^{p}{}_{j}\gamma(u)^{m}{}_{k}\big) \, \langle f_{q} , (f^{v})^{\star} \rangle , \end{equation*} which for $F = \Omega^{1}$ thus simplifies to 
\begin{equation}\label{eq_star_compatibility_Gamma} 2\textup{Re}(\Gamma^{i}{}_{jk})  = (\Gamma^{i}{}_{uv})^{*}\big( \Gamma^{v}{}_{jk} \epsilon(u) - \gamma(u^{-1})^{v}{}_{l}\Gamma^{l}{}_{pm} \gamma(u)^{p}{}_{j}\gamma(u)^{m}{}_{k}\big). \end{equation}
For additional or weaker compatibility conditions between the calculi and the above, see \cite[Chapter 8, p574]{BegMa}.

The most trivial example of a bimodule covariant derivatives on $(\Omega^{1}, (\ , \ ))$ satisfying all the above is by simply setting $\Gamma^{i}{}_{jk} = 0$, which gives the Maurer-Cartan bimodule covariant derivative $\nabla_{F}(g \otimes \omega) = \extd g \otimes_{\C G}\omega $ for $\omega \in \Lambda^{1}$, and $\sigma^{L} = \textup{flip}$. This covariant derivative is flat $\star$- and metric- compatible. It is is additionally torsion free, and thus a QLC, iff $\extd e^{g} = -e^{g}\wedge e^{g}= 0$ for all $g \in G$, i.e. iff $\Lambda^{2}$ is some quotient of the usual Grassmannian one.

\begin{example} \label{example_covariant_deriv_on_calc} \rm For $G=S_3$, we know the possible $(\CC,\pi)$ and the induced representation $\tilde\pi$ from  Example~\ref{exS3} $\Omega^{1}_{\widetilde{\pi}} \cong kG \otimes \Lambda^{1}_{\widetilde{\pi}}$. For $\Lambda^1_{\tilde\pi}$, we take $\textup{End}(\ )$ of each component of $\tilde\pi$ except the trivial one. Hence:

(i) For $\mathcal{C} = \{ e\}$, $\Lambda^{1}_{\widetilde{\pi}} = \Lambda^{1}_{\pi} $ is given by $\textup{End}(\textup{sign})$ or $\textup{End}(2)$ (with $\pi = 2$ denoting the two-dimensional irrep of $S_{3}$), 

(ii) For $\mathcal{C} = \{ uv, vu\}$, we have $\Lambda^{1}_{\widetilde{\pi_{1, 2}}} = \textup{End}(2)$ or $\Lambda^{1}_{\widetilde{\pi_{0}}} = \textup{End}(\textup{sign})$.

(iii) For $\mathcal{C} = \{ u, v, w\}$ gives $\Lambda^{1}_{\widetilde{\pi_{+}}} = \textup{End}(2)$ and $\Lambda^{1}_{\widetilde{\pi_{-}}} = \textup{End}(\textup{sign}) \oplus \textup{End}(2)$.

From this list, we see that there are only three possibilities for $\Lambda^{1}_{\widetilde{\pi}}$, namely $\textup{End}(\textup{sign})$, $\textup{End}(2)$ and their sum. 

(a) We begin with the more interesting case of the 4-dimensional calculus $\Lambda^{1}_{\widetilde{\pi}} = \textup{End}(2)$, with the following an extension of \cite[Proposition 4.5]{MaTao}. To be concrete we take the $\End(2)$ arising from Case (ii) and $\widetilde{\pi_2}$ there. There is a similar story for $\widetilde{\pi_1}$ with $q = e^{\frac{2\pi i}{3}}$ replaced by $q^{-1}$. From the actions given in Example~\ref{exS3} applied to the construction of $\Omega(S_3)$, we obtain 
\begin{equation} \label{eq_matrices_S3} e^{u} = \begin{pmatrix}
    -1      & 1\\
    1 & -1
\end{pmatrix}, e^{v} = \begin{pmatrix}
    -1      & q^{-1}\\
    q & -1
\end{pmatrix}, e^{uv} = \begin{pmatrix}
    q-1     & 0\\
    0 & q^{-1}-1
\end{pmatrix}, e^{vu} = \begin{pmatrix}
    q^{-1}-1     & 0\\
    0 & q-1
\end{pmatrix}.
\end{equation}
This still obeys $e^{u} + e^{v} + e^{w} = e^{uv} + e^{vu}=-3\theta$ as in \cite[Proposition 4.5]{MaTao}, as it must be as the representations are equivalent. Working in the basis order $\{e^u,e^v,e^{uv},e^{vu}\}$ for the calculus, we have for $\rho$ in (\ref{eq_rho_gamma_def}) as
\[ \rho(u) = \begin{pmatrix}
        1  &0  & 0&0 \\
 -1   & -1&1 &1 \\
     0   &0 & 0&1\\
         0   &0 & 1&0
\end{pmatrix}, \rho(v) = \begin{pmatrix}
        -1  & -1 & 1& 1\\
   0 & 1& 0&0 \\
     0   &0 & 0&1\\
         0   &0 & 1&0
\end{pmatrix}, \gamma(u) = \begin{pmatrix}
         -1 &  0& 0&0 \\
    -1& 0& 0&1 \\
      -2  &-1 & 1&1\\
          -1  & 1&0 & 0
\end{pmatrix}, \gamma(v) = \begin{pmatrix}
         0 &-1  & 1& 0\\
   0 &-1 & 0& 0\\
        1& -1& 0&0\\
          -1  & -2& 1& 1
\end{pmatrix}.\] 
We also fix a $D^{\vee}(S_3)$-$\star$-compatible bimodule inner product $(\ , \ )$ on this $\Omega^1$ as follows. As two conjugacy classes are present in the basis $\mathcal{B}$,  by Lemma~\ref{eq_lemma_squared_lengths_metric} the metric is determined by two arbitrary real numbers $l_{\mathcal{C}_{u}}=: l_{1}$, $l_{\mathcal{C}_{uv}}=: l_{2}$ as 
 \[ \underline{\underline{g}} = \frac{1}{2}\begin{pmatrix}
       2l_{1}  &2l_{1} - l_{2} &l_{2}&  l_{2}\\
    2l_{1} - l_{2} &2l_{1} & l_{2}& l_{2}\\
       l_{2} & l_{2}& 2l_{2}& l_{2}\\
             l_{2}& l_{2}& l_{2}& 2l_{2}
\end{pmatrix}. \] This has determinant $\frac{1}{16}l_{2}^{3}(12l_{1} - 7l_{2})$, and so a metric $\metric$ exists with respect to this bimodule inner product iff $l_{2} \neq 0$, and $l_{1} \neq \frac{7}{12}l_{2} $. Here,  \cite[Proposition 4.5]{MaTao} is recovered for $l_{2} = 1$, $l_{1} = \frac{2}{3}$.

Denoting the $4$ matrices $\{ \Gamma^{i}\}$ by $(\Gamma^{i})_{j, k}:= \Gamma^{i}{}_{jk}$, where $i,j,k$ run over $u,v,uv,vu$,   the $D^{\vee}(S_3)$-covariance condition $ \rho(g)^{i}{}_{l}\Gamma^{l} = \rho(g)^{T}\cdot
    \Gamma^{i}\cdot \rho(g)$ from \eqref{eq_Gamma_covariant} written on the generators $u,v$ can be simplified to the following $6$ equalities
  \[ \rho(u)^{T} \Gamma^{vu} \rho(u) = \rho(v)^{T} \Gamma^{vu} \rho(v) = \Gamma^{uv},  \quad \Gamma^{u} = \rho(u)^{T} \Gamma^{u} \rho(u), \quad \Gamma^{v} = \rho(v)^{T} \Gamma^{v} \rho(v),\]\begin{equation} \label{eq_Gamma_just_covariance} \rho(u)^{T}\Gamma^{v}\rho(u) = \rho(v)^{T} \Gamma^{u}\rho(v) = - (\Gamma^{u} + \Gamma^{v}) + (\Gamma^{uv} + \Gamma^{vu}) . \end{equation}
  Solving these (most simply done in the order given) results in a complex 10-parameter moduli space \[  \hspace{-20mm} \Gamma^{u} ={\small \begin{pmatrix}
       r  & -(a+b)-r& \hphantom{-}-(a+b)-(c+f) + 2s +p+t & s\\
 -(a+b)-r & \hphantom{-}-2(a+b)-2r &-(c+f) +t +r +s+p & 2(a+b) +(c+f) -t+r-s-p  \\
       (a+b)+r-t+p &  \hphantom{-}2(a+b) -t+2r&  (c+f) -r-s-p& \hphantom{-}-(a+b) - \frac{(c+f)}{3} +\frac{(d+e)}{3} + t-r\\
         p  & t& \hphantom{-}(a+b) + \frac{2(c+f)}{3} + \frac{(d+e)}{3}-s-t-p&0   \end{pmatrix}},\] \[
         \Gamma^{vu} = \begin{pmatrix}
        -(a+b)  &b - 2a  & a&a \\
 a- 2b   & -(a+b)&b &b \\
     b   &a & c&d\\
         b   &a & e&f
\end{pmatrix}.
         \] 
 The matrix         $\Gamma^{v}$ has the same entries as $\Gamma^{u}$ but swapping 
 \[ 11 \leftrightarrow 22,\ 12 \leftrightarrow 21,\ 13 \leftrightarrow 24,\ 14 \leftrightarrow 23,\  31 \leftrightarrow 42, \ 32 \leftrightarrow 41,\ 33 \leftrightarrow 44,\ 34 \leftrightarrow 43.\]
 The same rearrangement of entries gives $\Gamma^{uv}$ from $\Gamma^{vu}$. To solve the equations for a QLC it is simplest to first solve for a  weak QLC (WQLC) in the sense \cite{BegMa} of a  torsion and cotorsion free covariant derivative, where cotorsion is defined by $\textup{coT}_{\nabla} = (\extd \otimes \id - \id \wedge \nabla_{\Omega^{1}})\metric$. To fix $\Omega^2$, we note that for any $G$, a natural $D^{\vee}(G)$-covariant and quantum symmetric choice of $\Lambda^{2}$ is that given by the Grassmann relations that the $\{e^i\}$ basis elements anticommute. Then torsion freeness requires each matrix $\Gamma^{i}$ to be symmetric. In our case, it is sufficient to impose this just for $\Gamma^{u}$ and $\Gamma^{vu}$, which gives the further relations
         \[ b = a, \quad e = d, \quad p = s, \quad 2t = 4a + (c+f) +r - 2s. \] 
Cotorsion freeness requires the matrices $(\Theta^{k})_{i, j}:= g_{im}\Gamma^{m}{}_{jk}$ to also be symmetric (and for the same reasons this need only be checked for $\Theta^{u}$ and $\Theta^{vu}$), which in our case gives additionally
          \[ a = 3s, \quad (c+f) = r + 2s, \quad (d+3f) = 2r + 4s;\quad \begin{cases}  \textup{ $r$, $s$ free parameters, } & \textup{if $l_{1} = \frac{2}{3}l_{2}$}\\
        \,  r = s = 0, & \textup{otherwise}  \end{cases} .\] 
 This gives us a moduli space or WQLCs with 3 parameters $r,s,f$ or 1 parameter $f$ for the two cases, with Christoffel symbols 
        \begin{equation} \label{eq_WQLC}\hspace{-10mm}\Gamma^{u} ={\small  \begin{pmatrix}
          r& -(6s+r) & s& s\\
  -(6s+r) & -2(6s+r)& 6s+r& 6s+r \\
        s&6s+r &0 &r+2s - 2f\\
          s  &6s+r & r+2s - 2f &0
\end{pmatrix}, \quad \Gamma^{vu} = \begin{pmatrix}
          -6s& -3s & 3s& 3s\\
  -3s &\hphantom{-}-6s & 3s& 3s\\
       3s & 3s& -f+r+2s& \hphantom{-}2r+4s-3f\\
          3s  &3s & 2r+4s-3f& f
\end{pmatrix}}, \end{equation} (with $s = r = 0$ if $l_{1} \neq \frac{2}{3}l_{2}$).  Asking for a QLC, i.e., full metric compatibility, forces $r=s=0$ independently of $l_1,l_2$.  Next, we impose $\star$-compatibility as in \eqref{eq_star_compatibility_Gamma}, which  gives 
\[ s = 0, \quad r =f , \quad \textup{Re}(f) = 0,\] 
and so there is only one imaginary degree of freedom $f$ if $l_{1} = \frac{2}{3}l_{2}$, and all parameters zero otherwise. We note that either case makes $\sigma$ the flip map on the basis (since this happens iff $s = 0$, $f = r$). Going back to a WQLC, one can compute the Riemann curvature 
\[ R_{\nabla}^{L}(e^{i}) =- {1\over 2}R^{i}{}_{jk l} \, e^{j} \wedge e^{k} \otimes e^{l}; \quad  R^{i}{}_{jk l}=- \Gamma^{i}{}_{jm}\Gamma^{m}{}_{kl} + \Gamma^{i}{}_{km}\Gamma^{m}{}_{jl}, \]
and find that Riemann compatibility imposes  \[ s = 0, \quad f = \frac{r}{4} \textup{ or } r  .\] Thus, $\star$-compatible among WQLCs implies Riemann compatible case. We also note that the case $s=0, f=r$ imaginary is the one in \cite{MaTao}, which has formulae for the Riemann curvature.  The Ricci tensor as in \cite[Chap.~8]{BegMa} for the antisymmetric lift is  $\textup{Ricci}_{ij} = \frac{1}{2}(R^{k}{}_{kij} - R^{k}{}_{ikj})$ in the $\star$-compatible case becomes \[(\textup{Ricci})_{ij} =  \frac{f^{2}}{2}\begin{pmatrix}
          6& 3 &-3 &-3 \\
   3& 6& -3& -3\\
        -3& -3& 4& 1\\
          -3 &-3 & 1& 4
\end{pmatrix}. \] 
So, the only non-zero $\nabla^L$ which is  $D^{\vee}(G)$-covariant, Riemann and $\star$-compatible is in the case $l_{1} = \frac{2}{3}l_{2}$ and is completely determined by a single imaginary parameter $f \in i \mathbb{R}$ and has this Ricci curvature and with only $f=0$ the QLC case.  The Ricci tensor can similarly be computed for any WQLC from \eqref{eq_WQLC}. Here we just give the Ricci scalar given by evaluating the metric against this, namely
\[ S = l_{2}(12f^{2} - 15fr - 54fs + 6r^{2} + 48rs + 105s^{2}).\]

(b) For the case of $\Lambda^{1}_{\widetilde{\pi}} = \textup{End}(\textup{sign})$, we again take case (ii) and now $\tilde{\pi_0}$ for our calculations. We take $e^u$ as a basis of $\Lambda^1$ which as 1-1 matrix has value -2. For other group elements we have $e^g=\sign(g)-1$. For the matrices defined by \eqref{eq_rho_gamma_def}, we therefore get $ \rho(g)=  1,$ $ \gamma(g) = \textup{sign}(g)$ for all $g \in G$. The bimodule inner product is described by a single real parameter $\underline{\underline{g}} = g^{uu} = l \in \mathbb{R}$ (and is invertible if $l \neq 0$), and the Christoffel symbols by $\Gamma^{u}{}_{uu} =: \Gamma \in \mathbb{C}$. If we use the Grassmann relations then $\Omega^2$ is zero and all connections are WQLCs, while  $\star$-compatibility forces $\textup{Re}(\Gamma) = |\Gamma|^{2}$ or $\Gamma = x \pm i \sqrt{x(1-x)}$ for real parameter $ 0\leq x \leq 1$. The curvature here is zero as $\Omega^2=0$. The case of a QLC is $\Gamma = 1$ or $\Gamma = 0$. \end{example}

\subsubsection{The bundles $D^{\vee}(G) \twoheadrightarrow \C (C_{G}) \rcocross \C G$ and $D^{\vee}(G) \twoheadrightarrow  \C G$}

Next, we study the quantum differential geometry of the original quantum principal bundle used in the Wigner construction, i.e. which recovered geometrically the coirrep $V_{\mathcal{C},\pi}$ through  associated vector bundle to the universal quantum principal bundle. Since we have seen that $V_{\mathcal{C}, \pi}$ itself corresponds to the $\Omega^1_{\mathcal{C}, \pi}$ coirreducible calculi on the total space, the non-universal bundle below is a tautological quantum homogeneous bundle $P \twoheadrightarrow H$ associated to the irrep $V_{\mathcal{C}, \pi}$.

\begin{proposition} \label{prop_induced_calc_H13} For the bicovariant calculus $\Omega^{1}_{\mathcal{C}, \pi}$ on $D^{\vee}(G)$:

\noindent (i) We obtain a principal homogeneous bundle $\C (G) \hookrightarrow D^{\vee}(G) \twoheadrightarrow \C G$ with \[ \Omega^{1} \cong \Omega^{1}_{\mathcal{C}}\, , \quad \Lambda^{1} \cong \delta_{\mathcal{C}, \{ e\}}\Lambda^{1}_{\pi},\]
where $\Omega^{1}_{\mathcal{C}}$ is defined as above using the conjugacy class $\mathcal{C} \subseteq G\backslash \{ e\}$ and $\pi$ a $G$-representation.

\noindent (ii) We obtain a principal homogeneous bundle $\C (G/C_{G}) \hookrightarrow D^{\vee}(G) \twoheadrightarrow \C (C_{G}) \rcocross \C G$ with \[ \Omega^{1} = \textup{span}_{\C }\{ e_{c \to dcd^{-1}} , \, c, d \in \mathcal{C}\, | \, , dcd^{-1} \neq c\}, \quad \Lambda^{1} \cong \frac{(\C (C_{G})\rcocross \C G)^{+}}{\textup{span}_{\mathbb{C}}\{ s_{ci}{}^{dj} \, | \, c \in C_{G}\}},\] 
where $\Omega^{1}$ is the calculus on the algebra of functions $\C (G/C_{G})$ corresponding to the graph with vertices $\mathcal{C}$ and the subcollection of the arrows $\{ e_{a \to b}\}_{a, b \in \mathcal{C}}$.

\noindent (iii) The pair $\id: D^{\vee}(G) \to D^{\vee}(G)$ and $p_{2}|_{H}: \C (C_{G}) \rcocross \C G \twoheadrightarrow \C G$ form a morphism from the principal bundle in (ii) to the principal bundle in (i).
\end{proposition}

\proof (i) As in the proof of Lemma~\ref{prop_bundle_over_kG}, we note that under the Hopf algebra surjection $p_{2}: D^{\vee}(G) \twoheadrightarrow \C G $ onto the second factor,  \[ p_{2}(s_{ci}{}^{dj}) = \frac{\textup{dim}(V_{\pi})}{|C_{G}|} \sum_{n \in C_{G}} \pi(n^{-1})^{j}{}_{i} \, \delta_{c, e} q_{c}nq_{d}^{-1} = \delta_{\mathcal{C}, \{ e\}} \frac{\textup{dim}(V_{\pi})}{|C_{G}|} \sum_{n \in C_{G}} \pi(n^{-1})^{j}{}_{i}n, \quad s_{ci}{}^{dj} \in I_{(\mathcal{C}', \pi')}, \] (where in the case that $\mathcal{C} = \{ e\}$,  we note that $\pi$ is an irrep of $C_{G}= G$). Therefore, for the same ideal $I$ defining $\Lambda^{1}_{\mathcal{C}, \pi}$, then $p_{2}: I \mapsto \bigoplus_{\pi' \in \textup{Irr}(G) \, | (\{ e\}, \pi') \neq (\mathcal{C}, \pi), \pi' \neq 1} I_{\pi}$, and so  \[ (p_{2})_{*}: \Lambda^{1}_{\mathcal{C}, \pi} \twoheadrightarrow \delta_{\mathcal{C}, \{ e\}}\Lambda^{1}_{\pi}, \quad s_{ci}{}^{dj} \mapsto \delta_{\mathcal{C}, \{ e\}} E_{i}{}^{j}, \] for the basis $\{ E_{i}{}^{j}\}$ of $\textup{End}(V_{\pi})$, since $\varpi( \frac{\textup{dim}(V_{\pi})}{|C_{G}|} \sum_{n \in C_{G}} \pi(n^{-1})^{i}{}_{j}n) = E_{i}{}^{j}$ by the Artin-Wedderburn Theorem. Here ${\rm Irr}(G)$ denotes the irreps of $G$. 

The calculus on the base $\C (G)$ is given by $\C (G)(\extd \C (G))\C (G) = \C (G)\extd \C (G) \subseteq \Omega^{1}_{\mathcal{C}, \pi}$, where,  as $\C (G)$ sits trivially inside $D^{\vee}(G)$,  we simply use Lemma~\ref{lem_coirreducible_calc_DG} directly to get $\Omega^{1}$ is given by the span of the following elements \[ \delta_{h}\extd (\delta_{g}) = \sum_{c \in \mathcal{C}} \delta_{h}(R_{c} - \id)(\delta_{g}) \otimes E_{ci}{}^{ci}, \] and this can be quickly verified to indeed be bicovariant under the pushforward coactions. Then it is clear that $\Lambda^{1} = (\Omega^{1})^{\C (G)}$ is then  spanned only by $\{ E_{ci}{}^{ci}\}_{c \in \mathcal{C}}$, and so we then have the obvious isomorphism on the bases by \[ \Lambda^{1} \cong \Lambda^{1}_{\mathcal{C}}, \quad E_{ci}{}^{ci} \mapsto e_{c} \, ,  \] if $\mathcal{C} \neq \{ e\}$, and $\Lambda^{1} = 0$ otherwise. 

If $(\mathcal{C},\pi)$ has $\mathcal{C} \neq \{ e\}$ then $\textup{ver}$ is the zero map (and the horizontal forms by $\Omega^{1}_{\mathcal{C}, \pi}$) and otherwise is given by the identity map $\delta_{g} \otimes h \otimes E_{i}{}^{j} \mapsto \delta_{g} \otimes h \otimes E_{i}{}^{j}$. This would give that horizontal forms are given by the zero space, which agrees with the fact that from the above $\Omega^{1} = \C (G) \otimes 0 = 0$ in this case. 

(ii) Denoting the quotient Hopf algebra map $q: \C (G) \twoheadrightarrow \C (C_{G})$ of \eqref{eq_surjection_kG_kCG},  
\[ (q \otimes \id)(s_{ci}{}^{dj}) = \frac{\textup{dim}(V_{\pi})}{|C_{G}|} \sum_{n \in C_{G}} \pi(n^{-1})^{j}{}_{i} \delta_{c, C_{G}} \delta_{c} \otimes q_{c}nq_{d}^{-1} = \delta_{cr, rc} s_{ci}{}^{dj}, \quad s_{ci}{}^{dj} \in I_{(\mathcal{C}', \pi')}. \] 
This will induce a subcrossed module $(q \otimes \id)(I)$ of $\C (C_{G}) \rcocross \C G$ generating a bicovariant calculus on $\C (C_{G}) \rcocross \C G$. (Note that $C_{G}$ is never trivial, since if $r = e$, then $C_{G} = G$, and if $r \neq e$ then $r, e \in C_{G}$, and so this will never be a calculus just on $\C G$.) 

For the base algebra $\C (G/C_{G})$, the calculus is $\Omega^{1} := \C (G/C_{G})(\extd \C (G/C_{G}))\C (G/C_{G}) \subseteq \Omega^{1}_{\mathcal{C}, \pi}$, where here we recall that $\C (G/C_{G})$ sits inside $D^{\vee}(G)$ via \eqref{eq_AinDvee}. By the form of this latter inclusion map, it is clear that $\Omega^{1}$ is a subcalculus of $\C (G) (\extd \C (G))\C (G)$, where the latter was seen (in part (i) above) to be isomorphic to $\Omega^{1}_{\mathcal{C}}$ if $\mathcal{C} \neq \{ e\}$ and $0$ otherwise. In other words, the inclusion $\C (G/C_{G}) \hookrightarrow \C (G)$ is differentiable with respect to the calculus inclusion $\Omega^{1} \subseteq \Omega^{1}_{\mathcal{C}}$. So, in the trivial case of $\mathcal{C} = \{ e\}$, then we get that $\Omega^{1} = 0$ also. In the non-trivial case that $\mathcal{C} \neq \{ e\}$,  since $\C (G/C_{G})$ and $\C (G)$ are algebras of functions over finite sets, we can use the (contravariant) functional correspondence  between directed graphs and calculi on algebras of functions over finite sets, to conclude that $\Omega^{1}$ corresponds to a directed graph on the set $G/C_{G}$. Moreover, by the functoriality, it must correspond more strongly to a quotient graph of the Cayley graph on $G$ (since this latter graph was seen to correspond to $\Omega^{1}_{\mathcal{C}}$), i.e. all arrows on the graph corresponding to $\Omega^{1}$ must originate from arrows on the Cayley graph. This quotient map is then given on the vertices by \[ G \twoheadrightarrow G/C_{G} \cong \{q_{c}\}_{c \in \mathcal{C}} \cong \mathcal{C}, \quad g \mapsto grg^{-1}, \] for all $g\in G$ (where the final isomorphisms are by $[g] \mapsto q_{grg^{-1}}\mapsto grg^{-1}$ as explained in \eqref{eq_partition}), and thus on the arrows by \[ (g \to gc) \mapsto (grg^{-1} \to (gc)r(gc)^{-1}), \] for all $g \in G, c \in \mathcal{C}$. Writing an arbitrary element $d \in \mathcal{C}$ as $grg^{-1}$ for some $g$,  the end-point of the arrow is given by $gcrc^{-1}g^{-1} = gcg^{-1}dgcg^{-1}$,  where $gcg^{-1} \in \mathcal{C}$. Hence the result  follows. Here, if the head and tail of the arrow are equal under the above map, then the arrow goes to zero under the map. 

(iii) For the pair of maps $(\id, p_{2}|_{H})$, we note that map on the total spaces is clearly a $\C G$-comodule differentiable map (with the $\C G$-comodule structure on the domain by pushforward along $p_{2}|_{H}$), and indeed induces (the obvious) differentiable map on the base spaces $\C (G/C_{G}) \hookrightarrow \C (G)$ (since we saw that $\Lambda^{1}$ was induced by the calculus on $\C (G)$). Also, the map $p_{2}|_{H}$ is clearly a differentiable map of Hopf algebras, since $p_{2}|_{H}(J) = J'$ for the two subcrossed modules $J$ and $J'$ defining the bicovariant calculi, since $e \in C_{G}$. Moreover, this pair intertwines the two sequence of maps $\C (G/C_{G}) \hookrightarrow D^{\vee}(G) \twoheadrightarrow \C (C_{G}) \rcocross \C G$ and $\C (G) \hookrightarrow D^{\vee}(G) \twoheadrightarrow \C G$, as well as the relevant exact sequence of $D^{\vee}(G)$-modules of \eqref{eq_exact_bundle_def}. \endproof

Again, the above gives the form of the calculi on the base space and structure group when starting with an arbitrary calculi $\Omega^{1}_{D^{\vee}(G)}$ on the total space, using the decomposition of the latter's into its coirreducible blocks as explained before. That the calculi on the base space $\C (G)$ in the above is bicovariant is unsurprising since this is the coinvariant space of a normal quotient Hopf algebra (or in the Poincar\'e group analogy, it is the quotient by a normal subgroup).  

\begin{example} \rm Induced calculi $\Omega^{1}$ on $\mathbb{C}(G/C_{G})$ from Proposition~\ref{prop_induced_calc_H13}(ii) for the case of $G = S_{3}$ is simply $\Omega^{1} = 0$ when $\mathcal{C} = \{ e\}$ and $\mathcal{C} = \{ uv, vu\}$, and is the complete graph when $\mathcal{C} = \{ u, v, w\}$.
\end{example}

\subsection{The free field equation description of coirreps of $D^{\vee}(G)$} \label{eq_free_field_equation}Having understood the quantum differential geometry of the principal bundles studied in Section~\ref{secWigner}, we now return to the description of the coirreps $V_{\mathcal{C}, \pi}$ of Corollary~\ref{corollary_FT_description_VCpi}. The aim of this section is to obtain a description of these coirreps as spaces of solutions of some free field equations on  $\C G$ as a quantum Riemmanian geometry in the role of spacetime. 

In particular, we would  like a description of the output of the $W$-valued mass-shell Fourier transform $\C (\mathcal{C}) \otimes W \to \C \mathcal{C}^{-1} \otimes W \subseteq \C G \otimes W$ of Corollary~\ref{corollary_FT_description_VCpi} as solutions to a `Klein-Gordan equation', i.e. as eigenspaces of a $W$-valued $D^{\vee}(G)$-coinvariant Laplacian on $\C G$ for which the eigenvalues play the role of  `mass'. By \eqref{eq_left_white_headcrossed_module}, the $D^{\vee}(G)$-coaction \eqref{eq_kG_DveeG_coaction_symmetry} that we need to be covariant under can be equivalently written as \begin{equation} \label{eq_kG_DveeG_coaction_symmetry_alternative}\Delta_{L}g = \sum_{f \in G} \delta_{f} \otimes f^{-1}|g|f \otimes g \triangleleft f, \quad |g| = g, \quad g \triangleleft f = f^{-1}gf. \end{equation} 

\begin{lemma} \label{lemma_eigenspaces_linear_op}
A linear map  $L:\C G \to \C G$ is a left $D^{\vee}(G)$-comodule map with respect to \eqref{eq_kG_DveeG_coaction_symmetry_alternative} iff the eigenspaces are direct sums of spaces of the form $\C \CC$ spanned by the conjugacy classes $\mathcal{C}$ of $G$. \end{lemma}

\proof This can be shown concretely by directly testing the comodule condition: denoting the ouput $L(g) = \sum_{k \in G}c^{g}{}_{k}k$ for arbitrary $g \in G$ and asking for invariance under the grading of \eqref{eq_kG_DveeG_coaction_symmetry_alternative} gives that $L(g \triangleleft \delta_{h}) = \delta_{g, h}L(g) = \delta_{g, h}\sum_{k \in G}c^{g}{}_{k}k$ must equal $L(g) \triangleleft \delta_{h} =\sum_{k \in G}c^{g}{}_{k}k \triangleleft \delta_{h} = \sum_{k \in G}c^{g}{}_{k}\delta_{k, h} k =c^{g}{}_{h}h$ for all $h \in G$, and so $L(g) = \lambda_{g}g $ for some $\lambda_{g} \in \mathbb{C}$. Thus, each $g \in G$ is an eigenfunction of $L$ with eigenvalue $\lambda_{g}$. Moreover, asking for invariance under the right $G$-action in \eqref{eq_kG_DveeG_coaction_symmetry_alternative} gives that $L(g\triangleleft k) = L(k^{-1}gk) = \lambda_{k^{-1}gk}k^{-1}gk$ must equal $L(g) \triangleleft k = \lambda_{g}k^{-1}gk$, and so $\lambda_{g} = \lambda_{k^{-1}gk}$ for all $k \in G$. Hence eigenspaces consist of the spans of ad-invariant subsets of $G$. \endproof

It follows that $D^\vee(G)$-covariant operator $L: \C G \to \C G$ is fully characterised by a set of not necessarily distinct  eigenvalues $\{ \lambda_{\mathcal{C}}\}$, and is of the form \begin{equation} \label{eq_numbers_to_linear_operators} L = \sum_{\mathcal{C}} \lambda_{\mathcal{C}}\sum_{c \in \mathcal{C}} \delta_{c} \otimes c  \in (\C (G) \otimes \C G)^{D^{\vee}(G)} \cong \textup{Mor}_{^{D^{\vee}(G)}\mathcal{M} }(\C G, \C G) \cong  \mathbb{C}^{N_{\rm conj}},  \end{equation} where $N_{\rm conj}$ is the number of conjugacy classes in $G$.  Also note that the eigenfunctions corresponding to the conjugacy class $\mathcal{C}^{-1}$ of the linear operator $L \otimes \id_{W}: \C G \otimes W \to \C G \otimes W$ can be written as \[ \sum_{c \in \mathcal{C}^{-1}}\alpha_{c^{-1}}c^{-1} \otimes w = \sum_{c \in \mathcal{C}} f(c)c^{-1} \otimes w, \] for some coefficients $\alpha_{c^{-1}} \in \mathbb{C}$ which we choose to rewrite as a function $f: \mathcal{C} \to \mathbb{C}$ to match our transfer map in Corollary~\ref{corollary_FT_description_VCpi}  via (\ref{eq_mass_shell_FT_hopf}). We also note that the $D^{\vee}(G)$-comodule condition on $L$ contains the natural left $\C G$-coaction on itself as a Hopf algebra (i.e. the coregular coaction), and by cocommutivity of $\C G$, $L$ is then automatically also a right $\C G$-comodule map under the right coregular coaction. 

We could have started more generally with a $D^{\vee}(G)$-comodule map $L: \C G \otimes W \to \C G \otimes W$, however, repeating the steps in the proof of Lemma~\ref{lemma_eigenspaces_linear_op} gives with respect to a basis $\{ W_{i}\}$ of $W$ that $L(g \otimes W_{i}) = (\lambda_{g})^{j}{}_{i} g \otimes W_{j}$ for some matrix $\lambda_{g}$ satisfying $\rho(k)\lambda_{g}\rho(k^{-1}) = \lambda_{k^{-1}gk}$ for all $k \in G$. So, we will not necessarily have that eigenspaces are of the form above. Moreover, in our classical case of interest, the Laplacian on vector-valued functions is simply taken to be the scalar Laplacian on each orthogonal component. In the next section we will be interested in the particular case where $L$ is an invariant 2nd order differential operator defined over $\C G$ as a quantum Riemannian geometry.

\subsubsection{2nd order Laplacians on $\boldsymbol{\C G}$} \label{sec_Laplacians_on_CG}

 Given an arbitrary calculus $\Omega^{1}$ on $\C G$ and a bimodule inner product $(\ , \ ): \Omega^{1}\otimes_{\C G} \Omega^{1} \to \C G$, we recall that a 2nd order Laplacian is defined in quantum differential geometry as a linear map $\square: \C G \to \C G$ for which 
\begin{equation} \label{eq_def_laplacian} \square(gh) = \square(g)h + g\square(h) + 2(\extd g, \extd h) , \quad \forall g, h \in G. \end{equation}

This can then be extended to $W$-valued Laplacians simply by $\square \otimes \id_{W}: \C G \otimes W \to \C G \otimes W$. Laplacians can be defined more generally on any $\C G$-bimodules $F$ that is equipped with a covariant derivative simply as operators $F \to F$ satisfying a reasonable 2nd order Leibniz rule (i.e. similarly to \cite[Lemma 8.6]{BegMa}). However, as in previous comments, we only are interested in Laplacians on $\C G \otimes W$ arising from Laplacians on $\C G$ as above. 

We next observe that by \eqref{eq_def_laplacian} and the fact that $\extd e = 0$,  \[ \square(ee) = \square(e) = \square(e)e + e\square(e)+ 2(\extd e, \extd e) = 2\square(e), \] so that  $\square(e) = 0$. Hence, we always have $\lambda_{\{ e\}} = 0$ for a second order operator in this sense. 

\begin{lemma} \label{lemma_metric_laplacian_relation}
Let $\Omega^1$ be a $D^\vee(G)$-covariant calculus on $\C G$ and $L: \C G \to \C G$ a $D^\vee(G)$-covariant linear map with not necessarily distinct eigenvalues $\{ \lambda_{\mathcal{C}}\}_{\mathcal{C}}$. Then  $\lambda_{\{ e\}} = 0$ iff \begin{equation} \label{eq_metric_from_laplacian} (e^{g}, e^{h}) := -\frac{1}{2}g^{-1}\big[ L(gh^{-1}) - L(g)h^{-1}  - gL(h^{-1})\big]h\end{equation} defines a $D^\vee(G)$-covariant bimodule inner product on $\C G$. Moreover, $\square=L$ is then a second order differential operator with respect to this and has eigenvalues  obeying \begin{equation} \label{eq_relation_mass_length} \lambda_{\mathcal{C}} + \lambda_{\mathcal{C}^{-1}} =2  l_{\mathcal{C}}, \end{equation} with $l_{\mathcal{C}}=l_{\mathcal{C}^{-1}}$ defined with respect to \eqref{eq_metric_from_laplacian}. 
\end{lemma}

\proof The expression \eqref{eq_metric_from_laplacian} arises from supposing that the expression $( \extd g, \extd h):= L(gh) - L(g)h - gL(h)$ defines a bimodule map $\Omega^{1}_{\rho} \otimes_{\C G} \Omega^{1}_{\rho} \to \C G$  and then noting that under this assumption, $(e^{g}, e^{h}) = (g^{-1}\extd g, h^{-1} \extd h) = g^{-1}( \extd g, \extd (h^{-1}h) - (\extd h^{-1})h) = -g^{-1}(\extd g, \extd h^{-1})h$ by the Leibniz rule. We then rearrange this to obtain  \eqref{eq_def_laplacian} for $(\extd g, \extd h^{-1})$. Now starting with the definition \eqref{eq_metric_from_laplacian},  firstly it is clear that it is bilinear using \eqref{eq_linearity_e}. Since $L(g) = \lambda_{\mathcal{C}_{g}}g$ for $\mathcal{C}_{g}$ the conjugancy class of $G$ that contains $g \in G$ (by covariance of $L$),  \eqref{eq_metric_from_laplacian} indeed defines a map $(\ , \ ): \Lambda^{1}_{\rho} \otimes \Lambda^{1}_{\rho} \to k$, and is more simply given by \begin{equation} \label{eq_relation_metric_lambda} (e^{g}, e^{h}) = -\frac{1}{2} [\lambda_{\mathcal{C}_{gh^{-1}}} - \lambda_{\mathcal{C}_{g}} - \lambda_{\mathcal{C}_{h^{-1}}}].\end{equation} By Lemma~\ref{lemma_metric_crossed_module}, we check that this expression satisfies \eqref{eq_conditions_on_metric}. Here, the adjoint action invariance is immediate by the full covariance of $L$, i.e. \[ (e^{u^{-1}gu}, e^{u^{-1}hu}) = -\frac{1}{2}[\lambda_{\mathcal{C}_{u^{-1}gh^{-1}u}} - \lambda_{\mathcal{C}_{u^{-1}gu}} - \lambda_{\mathcal{C}_{u^{-1}h^{-1}u}}] = (e^{g}, e^{h}), \quad \forall u \in G, \] while the first condition of \eqref{eq_conditions_on_metric} (i.e. the right $\C G$-crossed module condition) comes out as the requirement that, starting off from the right-hand side of \eqref{eq_conditions_on_metric}\begin{align*}
 -\frac{1}{2} &[\lambda_{\mathcal{C}_{gu (hu)^{-1}}} - \lambda_{\mathcal{C}_{gu}} - \lambda_{\mathcal{C}_{(hu)^{-1}}}] + \frac{1}{2} [\lambda_{\mathcal{C}_{gu (u)^{-1}}} - \lambda_{\mathcal{C}_{gu}} - \lambda_{\mathcal{C}_{u^{-1}}}] \\
& + \frac{1}{2} [\lambda_{\mathcal{C}_{u (hu)^{-1}}} - \lambda_{\mathcal{C}_{u}} - \lambda_{\mathcal{C}_{(hu)^{-1}}}] - \frac{1}{2} [\lambda_{\mathcal{C}_{u (u)^{-1}}} - \lambda_{\mathcal{C}_{u}} - \lambda_{\mathcal{C}_{u^{-1}}}]\\
& = \frac{1}{2}[-\lambda_{\mathcal{C}_{gh^{-1}}} + \lambda_{\mathcal{C}_{g}} + \lambda_{\mathcal{C}_{h}}] - \frac{1}{2}\lambda_{\mathcal{C}_{e}} = (e^{g}, e^{h}) - \frac{1}{2}\lambda_{\mathcal{C}_{e}}. 
\end{align*} Thus, the condition is satisfied iff $\lambda_{\mathcal{C}_{e}}$ = 0. 

Moreover, from \eqref{eq_relation_metric_lambda} with $h = g$, we  have $(e^{g}, e^{g}) = -\frac{1}{2}[\lambda_{\{ e\}} - \lambda_{\mathcal{C}} - \lambda_{\mathcal{C}^{-1}}] = \frac{1}{2}(\lambda_{\mathcal{C}} + \lambda_{\mathcal{C}^{-1}})$ for $g \in \mathcal{C}$, i.e. that the eigenvalues satisfy the relation $ \lambda_{\mathcal{C}} + \lambda_{\mathcal{C}^{-1}} = 2(e^{g}, e^{g})$, which gives \eqref{eq_relation_mass_length} in view of \eqref{eq_length_C_Cinv}. \endproof

Thus, for any $D^{\vee}(G)$-covariant calculus $\Omega^{1}$ on $\C G$, we can choose $N_{\rm conj}-1$ not necessarily distinct $\{\lambda_\CC\}_{\CC\ne\{e\}}$ and set $\lambda_{\{e\}}=0$ to  generate a $D^\vee(G)$-covariant 2nd order Laplacian of the form \eqref{eq_numbers_to_linear_operators} with respect to a canonically defined $D^\vee(G)$-covariant bimodule inner product \eqref{eq_relation_metric_lambda}. Conversely, given any $D^\vee(G)$-covariant linear operator $L$, $\square=L - \lambda_{\{e \}}$ has eigenvalues $\{\lambda_{\mathcal{C}} -  \lambda_{\{e \}}\}$ and hence defines a bimodule inner product with respect to which it is a second order Laplacian.

The eigenvalues $\lambda_{\mathcal{C}}$ of $\square$ play the role of `squared masses' in the classical Poincar\'e group case associated to the eigenspace $\C \mathcal{C}^{-1}$ under the transfer/Fourier transform map. For this to match up better, we should work over $\C$ and require  $\lambda_{\mathcal{C}} = l_{\mathcal{C}} \in \mathbb{R}$. This is because  (from the discussion in Section~\ref{sec_second_bundle_construct} around the Fourier transform) $g^{-1} \in G $ is analogous to a plane wave $e^{-ix\cdot p}$ of  momentum $p$, hence $e^{g}$ is analogous to $ e^{-ix\cdot p} \partial_{\mu} e^{ix\cdot p}\extd x^{\mu} = ip_{\mu}\extd x^{\mu}$, which indeed has `length' $p^{2} = m^{2} \in \mathbb{R}$. We indeed recover the required relation with the lengths in the $\star$-compatible setting

\begin{corollary} \label{cor_star_compatibility_laplacian} If $\Omega^1$ is a $D^\vee(G)$-covariant $\star$-calculus and $\square$ is a $D^\vee(G)$-covariant operator second order Laplacian that commutes with $\star$ then the metric $(\ ,\ )$ defined by Lemma~\ref{lemma_metric_laplacian_relation} is $\star$-compatible. Moreover, if $\square$ has real-valued eigenvalues then $(\ , \ )$ is real-valued. Conversely, if $(\ ,\ )$ is $D^\vee(G)$, $\star$-compatible and real-valued then we have a unique 2nd order Laplacian that commutes with $*$ and has real eigenvalues, defined by $\lambda_{\mathcal{C}}= l_{\mathcal{C}} \in \mathbb{R}$
\end{corollary}

\proof $\star$ commuting with a $D^\vee(G)$-covariant operator $L$ is equivalent to simply asking $\lambda_{\mathcal{C}} = \lambda_{\mathcal{C}^{-1}}^{*}$. Its eigenvalues are then real iff $\lambda_{\mathcal{C}} = \lambda_{\mathcal{C}^{-1}}$. Hence, if  $\square$ is $D^\vee(G)$-covariant and commutes with $*$  then by \eqref{eq_relation_metric_lambda}, we must have that \[(e^{g}, e^{h})= -\frac{1}{2} [\lambda_{\mathcal{C}_{gh^{-1}}} - \lambda_{\mathcal{C}_{g}} - \lambda_{\mathcal{C}_{h^{-1}}}] =-\frac{1}{2} [(\lambda_{\mathcal{C}_{hg^{-1}}})^{*} - (\lambda_{\mathcal{C}_{g^{-1}}})^{*} - (\lambda_{\mathcal{C}_{h}})^{*}] =(e^{h}, e^{g})^{*}, \] so that  $(\ , \ )$ obeys the correct `reality' or $\star$-compatibility condition as in \cite{BegMa}. We also get by Lemma~\ref{lemma_metric_laplacian_relation} that $l_{\mathcal{C}} = \textup{Re}(\lambda_{\mathcal{C}}) \in \mathbb{R}$. As a partial converse,  if $(\ , \ )$ is $\star$-compatible and there exists a Laplacian $\square$ with respect to it for which $\textup{Re}(\lambda_{\mathcal{C}}) = \textup{Re}(\lambda_{\mathcal{C}^{-1}})$, then $\square$ must commute with $\star$. Here, by  \eqref{eq_relation_mass_length}, we get that $\textup{Re}(\lambda_{\mathcal{C}} + \lambda_{\mathcal{C}^{-1}}) = 2l_{\mathcal{C}}$ and $\textup{Im}(\lambda_{\mathcal{C}} + \lambda_{\mathcal{C}}^{-1}) = 0$, implying that $\lambda_{\mathcal{C}} = (\lambda_{\mathcal{C}^{-1}})^{*}$. 

For a $D^\vee(G)$-covariant real-valued 2nd order Laplacian that commutes with $\star$ (and so $\lambda_{\mathcal{C}} = \lambda_{\mathcal{C}^{-1}}$), we see by the same approach  that $(\ , \ )|_{\Lambda^{1}\otimes \Lambda^{1}}$ is symmetric and thus real by the above, and $\lambda_{\mathcal{C}} = l_{\mathcal{C}} \in \mathbb{R}$. Conversely, by Lemma~\ref{eq_lemma_squared_lengths_metric} we can construct such an operator $\square$ on this bimodule inner product by  taking $\lambda_{\mathcal{C}} = l_{\mathcal{C}} \in \mathbb{R}$. In fact, by our results, this is the unique $\square$ with the stated properties. \endproof

From these results, we see that in the setting of  real-valued $\star$-compatible bimodule inner products on $\C G$, then there exists a unique 2nd order Laplacian real-valued $\square$ on it that commutes with $\star$, defined by $\lambda_{\mathcal{C}} = l_{\mathcal{C}}$. Moreover, these eigenvalues necessarily have  $\lambda_{\mathcal{C}}= \lambda_{\mathcal{C}^{-1}}$, i.e. the eigenspaces associated to $\mathcal{C}$ and $\mathcal{C}^{-1}$ cannot be distinguished solely by the eigenvalue. For the interpretation that we want as `squared mass' we do want the eigenvalues to be otherwise distinct. This leads us to the notion of a {\em regular spectrum} $\{\lambda_\CC\in \R\}$ which obey $\lambda_\CC=\lambda_{\CC^{-1}}$ and are otherwise distinct, and  obey $\lambda_{\{e\}}=0$.   The corresponding operator will be called {\em regular} also and the canonical bimodule inner product with respect to which this is 2nd order will be regular in the sense defined previously (recalling here the necessary condition \eqref{eq_regularity_inner_prod}). Intrinsic to this approach is that the value of $\lambda_\CC$ then does not determine whether we refer to the $\CC$ or $\CC^{-1}$ eigenspace and we would need to designate one of each pair of non-involutive conjugacy classes for the irrep of $D(G)$ corresponding to a particular `mass'. This is exactly analogous to what happens for the group $\mathbb{R}^{1, 3} \lcross SO^{+}(1, 3)$ of Section \ref{secW}, where mass and spin only classify the irreps if we restrict to those of positive energy $p_0>0$.

\subsubsection{Geometric Laplacians on $\C G$} \label{sec_geometric_laplacians}

Having determined what are the 2nd order Laplacians $\square$ of interest on $\C G$, it remains to see if the associated bimodule inner product $(\ , \ )$ can be equipped with bimodule covariant derivative $(\nabla, \sigma)$ on $\Omega^1$ such that the Laplacian takes the geometric form 
 \begin{equation}\label{geoL} \square = (\ , \ ) \circ \nabla\circ \extd. \end{equation} 
 We note that this will satisfy the second order Leibniz relation $\square(gh) = \square(g)h + g \square(h) +2 (\extd g, \extd h)' $ where the bimodule inner product $(\ , \ )'$ is the linear combination \cite[Lemma 8.6]{BegMa} \begin{equation} \label{eq_linear_combo_innerproduct} (\ , \ )' = \frac{1}{2}(\ , \ )\circ (\id + \sigma). \end{equation} 
 Therefore, we take care to distinguish these two, denoting the bimodule inner product of Section~\ref{eq_free_field_equation} and its lengths by a prime. It is immediate that if $\Omega^{1}$, $\nabla$ and $(\ , \ )$ are all $D^\vee(G)$-covariant, then so is  $\square$ in (\ref{geoL}), and we assume this.  We use the left invariant basis $\{e^i\}$ of $\Omega^1$ as done previously,  in which case (\ref{geoL}) becomes
 \begin{equation} \label{eq_geometric_Laplacian_index_form} \square g = g^{ij}\partial_{i} \partial_{j}g - (\partial_{i}g)\Gamma^{i}{}_{jk}g^{jk}. \end{equation} 
Finally, working over $\C$, the requirements of a quantum Riemannian geometry include that $(\ , )$ is  $\star$-compatible sense or `real' in sense of \cite{BegMa} and that $\nabla$ is $\star$-preserving.  This ensures contact with real metric and connection coefficients in the classical limit. Here, we further limit attention to $(\ ,\ )$ real-valued on the left-invariant 1-forms, so that $g^{ij}$ are real and symmetric, and the corresponding $l_\CC$ are real. We do not necessarily assume that $\nabla$ is torsion free and metric compatible (a QLC), though this would be a natural choice. 
 \begin{proposition} \label{prop_relationship_mass_lenghts} A choice of 2nd order Laplacian dictated as in the preceding section by a choice of eigenvalues $\{ \lambda_{\mathcal{C}} \in \mathbb{R} \, | \, \lambda_{\mathcal{C}} = \lambda_{\mathcal{C}^{-1}}\}$  can be obtained from  $(\Omega^{1}, (\ , \ ),\nabla)$ as above if and only if
\[ - (\lambda_{\mathcal{C}_{hg^{-1}}} - \lambda_{\mathcal{C}_{g}} ) + (l_{\mathcal{C}_{hg^{-1}}} - l_{\mathcal{C}_{g}} ) = g^{ak}(h^{-1}\partial_{i}h)\gamma(g^{-1})^{i}{}_{l} \Gamma^{l}{}_{ak}. \]
 Here, $\{ l_{\mathcal{C}}\}$ are the square-lengths associated with $(\ ,\ )$.
\end{proposition}
\proof Here $\{\lambda_\CC\}$ are also used to define the square-lengths $l'_\CC$ associated with $(\ , \ )'$ in Corollary~\ref{cor_star_compatibility_laplacian}  so that $\square$ is 2nd order with respect to this. We evaluate \eqref{eq_linear_combo_innerproduct} on $e^{h} \otimes e^{g}$ and use \eqref{eq_relation_metric_lambda} with $\lambda_{\mathcal{C}} = \lambda_{\mathcal{C}^{-1}}$,  Lemma \ref{eq_lemma_squared_lengths_metric}, and the formula for $\sigma$ in the proof of Lemma~\ref{lembicnablaF} to obtain this as equivalent to
 \begin{align*} - (\lambda_{\mathcal{C}_{hg^{-1}}} - \lambda_{\mathcal{C}_{g}} - \lambda_{\mathcal{C}_{h}}) + (l_{\mathcal{C}_{hg^{-1}}} - l_{\mathcal{C}_{g}} - l_{\mathcal{C}_{h}}) &= - g^{ak}(h^{-1}\partial_{i}h)\big( \Gamma^{i}{}_{ak}\epsilon(g) - \gamma(g^{-1})^{i}{}_{l} \Gamma^{l}{}_{pm} \gamma(g)^{p}{}_{a}\gamma(g)^{m}{}_{k}\big) , \end{align*} where  $\langle e_{i}, e^{h} \rangle  = h^{-1}\partial_{i}h$ from the definition of $\partial_{i}$ in \eqref{eq_form_of_partial}. Moreover, 
 \begin{align*} g^{ij}\partial_{i}\partial_{j}h=(e^{i}, e^{j}) \partial_{i}\partial_{j}h= (e^{i}, e^{j}) \partial_{i}(h) \langle e_{j}, e^{h} \rangle = (e^{i}, e^{j}) h \langle e_{i}, e^{h} \rangle \langle e_{j}, e^{h} \rangle = (e^{h}, e^{h})h = l_{\mathcal{C}_{h}}h.\end{align*}
 so that \eqref{eq_geometric_Laplacian_index_form} is equivalent to \[ \lambda_{\mathcal{C}_{h}} - l_{\mathcal{C}_{h}}=  -(h^{-1}\partial_{i}h)\Gamma^{i}{}_{jk}g^{jk}, \]
 for all $h \in G$. Inserting this into the first equation, and using \eqref{DveeG_covariant_bimodule_inner} then gives the stated condition. Here, to go backwards, we take $g = e$ in the stated condition.\endproof
 
Note that $\square$ by assumption here commutes with $\star$ as we assumed real eigenvalues. Also, if $\nabla$ is $\star$-preserving then $\sigma$ is invertible and  $(\ ,\ )$ $\star$-compatible or `real' in the sense \cite{BegMa} then $\star \circ \square \circ \star =  (\ , \ ) \circ \star \circ \sigma\circ \star \circ \nabla \circ \extd=(\ , \ ) \circ \sigma^{-1}\circ \nabla \circ \extd=\square$, which  is automatic if  $(\ ,\ )\circ\sigma=(\ ,\ )$, which is equivalent to $(\ ,\ )'=(\ ,\ )$. For an example, we can take  $\nabla$ to be the Maurer-Cartan connection on $\C G$ given by $\nabla^L=0$ and $\square h = g^{ij}\partial_{i} \partial_{j}h $ which indeed has eigenvalues $\lambda_{\mathcal{C}} = l_{\mathcal{C}}$ by the calculation in the proof above. In this case $\sigma$ is the flip map on the basis of $\Lambda^1$ and the condition in the proposition holds automatically. We note also the partial converse that if $\lambda_{\mathcal{C}} = l_{\mathcal{C}}$ then $\square h = g^{ij}\partial_{i} \partial_{j}h $.

More generally, a geometric Laplacian on $E'$ would take the form $\square = (\ , \ ) \circ \nabla_{\Omega^{1}}\circ \nabla_{E'}$, and so the above Laplacian is simply that using $\nabla_{E'} = \extd \otimes \id$, i.e.  the Maurer-Cartan connection. We note that there are other methods of geometrically constructing quantum Laplacians with the same starting data, namely by the `divergence of the gradient' \cite[page 491]{BegMa}, however this construction is not further investigated here.

\begin{example} \rm Continuing the setting of Example~\ref{example_transferance_of_sections}, we fix a regular quantum Laplacian defined by the set of real-numbers $\{0,  \lambda_{1}, \lambda_{2} \}$ with $\lambda_{1}$ corresponding to the conjugacy classes $\{ u, v, w\}$ and $\lambda_{2}$ to $\{ uv, vu\}$ (and $0$ to $\{ e\}$). 

(a) For the case of $\Lambda^{1}_{\widetilde{\pi}} = \textup{End}(2)$, the vectors with components $\alpha(g)_{i}:= (g^{-1}\partial_{i}g) = \langle e^{g}, e_{i} \rangle$ take the form $ \alpha(e) = 0$, $\alpha(u, v, uv, vu)_{i} = \delta_{i}^{u, v, uv, vu}$, $ \alpha(w)^{T} = \begin{pmatrix} -1 & -1& 1 & 1\end{pmatrix}. $ Going through the pairs $(h, g) \in S_{3} \times S_{3}$ in Proposition~\ref{prop_relationship_mass_lenghts} gives the  relations \[ \lambda_{1} - l_{1} = \frac{2}{3}(\lambda_{2} - l_{2}), \quad \lambda_{1} - l_{1} = - \textup{Tr}(\Gamma^{u}\cdot \underline{\underline{g}}) = - \textup{Tr}(\Gamma^{v}\cdot \underline{\underline{g}}) = -\frac{2}{3}\textup{Tr}(\Gamma^{uv} \cdot \underline{\underline{g}}) = -\frac{2}{3}\textup{Tr}(\Gamma^{vu} \cdot \underline{\underline{g}}). \] 
If the Christoffel symbols are $D^{\vee}(S_{3})$-covariant (i.e. satisfy \eqref{eq_Gamma_just_covariance}) then the last three of the conditions above are automatic. So, for a $D^{\vee}(S_{3})$-covariant WQLC's (i.e. \eqref{eq_WQLC}), choosing say the equality $\lambda_{2} - l_{2} = - \textup{Tr}(\Gamma^{vu}\cdot \metric)$ from the above,  gives only one additional relation \[ \lambda_{2} = l_{2}[1 + 3(f - r-3s)], \] for any value of $l_{1}$, which gives the constraint $f - r - 3s \in \mathbb{R}$ among the three complex parameters.  In particular, for case of the $\star$-compatible case with $l_{1} = \frac{2}{3}l_{2}$, this gives $\lambda_{1} = l_{1} = \frac{2}{3}\lambda_{2} = \frac{2}{3}l_{2}$.

(b) For the case of $\Lambda^{1}_{\widetilde{\pi}} = \textup{End}(\textup{sign})$, we have $h^{-1}\partial_{u}h = \langle e^{h}, e_{u} \rangle= 1$ for $h \in \{ u, v ,w \}$ and $0$ otherwise. Proposition~\ref{prop_relationship_mass_lenghts} with the pair $(h, g) = (uv, e)$  then gives $\lambda_{2} = 0$. So, we cannot realise a regular quantum Laplacian geometrically by \eqref{geoL} in this case (as we would need $\lambda_{2} \neq 0$). This is because the necessary condition of \eqref{eq_regularity_inner_prod} does not hold.
\end{example}

\subsubsection{The additional free field equations} By extending the projector $P^{cov}_{r, \pi}$ of \eqref{eq_projector_on_position_space} defined on $\C \mathcal{C}^{-1} \otimes W$ to one on the whole of $\C G \otimes W$ by \[ P^{cov}_{r, \pi}(g \otimes w) = g \otimes (q_{g}P_{r, \pi}q_{g}^{-1}) \triangleright w, \] where $q_{g}$ denotes the cocycle $q$ defined with respect to a representative $r'$ of the conjugacy class $\mathcal{C}_{g}$, then we get 

\begin{corollary} \label{eq_corollary_klein_gordan_charac_ofVcpi} In the setting of Lemma~\ref{lemma_transferance_of_sections_DG} with $E \cong V_{\mathcal{C}, \pi}$ and a choice of regular $D^{\vee}(G)$-covariant $\star$-compatible Laplacian specified by the real-eigenvalues $\{ \lambda_{\mathcal{C}}\}$, we have the following isomorphism of left $D^{\vee}(G)$-comodules \[ E \cong \{ s \in \C G \otimes W \, | \, \big((\square + \lambda_{\mathcal{C}^{-1}}) \otimes \id_{W}\big)s = 0 \textup{ and } P^{cov}_{r, \pi}s = s\} \subseteq \C G \otimes W.  \] 
\end{corollary}

In other words, the analogue of the Klein Gordan equation detects the conjugacy class, imposing $s \in \C \mathcal{C}^{-1} \otimes W$, and the projector detects the restriction of $W$ to $C_{G}$. Overall, the solutions to these two equations together are those of the sections of the form $\sum_{c \in \mathcal{C}}f(c)c^{-1} \otimes q_{c} \triangleright w$ (from Lemma~\ref{lemma_transferance_of_sections_DG}).  

In the Poincar\'e group case, it is known for say spin-1 particles, with the 1-representation of $SU_{2}$ embedded inside the $(\frac{1}{2}, \frac{1}{2})$ representation of $SL_{2}\mathbb{C}$, then the projector condition can again be written as an additional free field equation. Whether similar additional free field equations can be obtained here using now covariant derivatives on $E'$ to play the role of the projector condition is less clear in our finite group case.

\subsection{Duality between Laplacians on $\C G$ and on $\C (G)$}\label{sec_duality}

The quantum double of a finite-dimensional Hopf algebra involves both the original Hopf algebra and its dual in a symmetric way. In this section, we explore a resulting duality between Laplacians and their eigenspaces on $\C G$ as studied above and the same on $\C(G)$. The aim is to generalise and put in context an interesting duality first observed in \cite{MaTao}. 

We recall that we saw in Lemma~\ref{lemma_eigenspaces_linear_op} that eigenspaces of $D^{\vee}(G)$-covariant maps $L: \C G \to \C G$ correspond to $\textup{Ad}$-stable subsets of $G$, of which the building blocks are conjugacy classes of $G$. We also recall from Section~\ref{sec_recap_calc_kG_kofG} that (non-trivial) conjugacy classes are building blocks for bicovariant calculi on $\C (G)$. We can phrase this relationship between such maps on $\C G$ and bicovariant calculi on $\C (G)$ as arising from the bijection 
\[ \{\textup{subobjects of $\C G \in {}^{D^{\vee}(G)}\mathcal{M}$ } \} \leftrightarrow  \{ \textup{subobjects of $\C(G) \in \XM{}^{\C (G)}_{\C (G)}$ }\}, \]
 given by taking annihilators (and viewing $D^\vee(G)$-comodules on the left as equivalently $\C(G)$-crossed modules on the right (these can also be identified with $G$-crossed modules).  Here $\C(G)$ has the right adjoint coaction and the right regular action, which is the standard way that bicovariant calculi on Hopf algebras are  constructed, as in Section~\ref{section_quantum_diff_structures_intro}. 

Dually, we have seen in Lemma~\ref{lem_DG_covariant_calc_kG} that $D^{\vee}(G)$-covariant calculi on $\C G$ correspond to two-sided ideals of $(\C G)^{+}$, whose building blocks are (non-trivial) irreps of $G$. Noting the duality between irreps and conjugacy classes, we consider the bijection 
\[ \{ \textup{subobjects of $\C G \in {}_{\C G}\mathcal{M}_{\C G}$}\} \leftrightarrow \{ \textup{subobjects of $\C (G) \in {}^{\C (G)}\mathcal{M}^{\C (G)}$} \}, \] given by again taking annihilators. The left-hand side is by the regular actions by multiplication, so the right-hand side is by the regular coactions by the coproduct. Analogously to the first bijection, the right-hand side of this second bijection can indeed be viewed as eigenspaces of certain linear operators $L^*:\C (G) \to \C (G)$ corresponding to irreps of $G$, as we now show.

\begin{lemma} \label{lem_L_decompose_2}
A linear map  $L^*:\C (G) \to \C (G)$ is a $\C( G)$-bicomodule map with respect to the left and right coregular actions iff the eigenspaces are direct sums of spaces of the form $\textup{span}_{\mathbb{C}}\{ \rho^{i}{}_{j}\}$ for irreps $\rho$ of $G$, where $\rho^{i}{}_{j}\in \C(G)$ by $\rho^i{}_j(g)=\rho(g)^{i}{}_{j}$. 
\end{lemma}
 \proof The Peter-Weyl decomposition can be  viewed as providing an isomorphism of $\C G$-bimodules (or $\C (G)$-bicomodules),
 \[ \C (G) \cong \bigoplus_{\rho \in \textup{Irr}(G)} \textup{End}(V_{\rho}),\]
where $\textup{Irr}(G)$ denotes the sets of irreps of $G$,  given by a direct sum of maps $\C (G) \to \textup{End}(V_{\rho})$, $\delta_{g}\mapsto \langle \delta_{g}, s_{j}{}^{i} \rangle E_{i}{}^{j}$ and  $s_{i}{}^{j} = \frac{\textup{dim}(V_{\rho})}{|G|}\sum_{g \in G} \rho(g^{-1})^{j}{}_{i}g$ (c.f. \eqref{eq_artin_wedderburn_DG}), so that $\{ s_{i}{}^{j}\}$ gives the dual basis to the space with basis $\{ \rho^{i}{}_{j}\}$. $\C (G)$ is a $\C G$-bimodule by the coregular actions (i.e. $h \triangleright \delta_{g} = \delta_{hg}$, $\delta_{g} \triangleleft h =\delta_{gh}$), while $\textup{End}(V_{\rho})$ is an indecomposable $\C G$-bimodule by 
  \[ h \triangleright E_{i}{}^{j} = E_{i}{}^{k}\rho(h^{-1})^{j}{}_{k} , \quad E_{i}{}^{j} \triangleleft h = \rho(h^{-1})^{k}{}_{i}E_{k}{}^{j}, \]
  where $\{E_{i}{}^{j} \}$ denotes as before the basis for $\textup{End}(V_{\rho})$ built out of a basis $\{ E_{i}\}$ of $V_{\rho}$ and its dual $\{ E^{i}\}$ as before. The inverse of the Peter-Weyl decomposition is given the direct sum of the following map on each block
  \[ \rho^{j}{}_{i}\mapsfrom E_{i}{}^{j}, \quad E_{i}{}^{j} \in \textup{End}(V_{\rho}). \] 
Combining this with Schur's Lemma gives the statement.  \endproof

Just as for \eqref{eq_numbers_to_linear_operators}, it then follows that any such linear operator $L^*: \C (G) \to \C (G) $ is fully characterised by a set of not necessarily distinct eigenvalues $\{ \lambda^*_{\rho}\}_{\rho \in \textup{Irr}(G)}$, and is of the form
\[ L^* (\delta_h)= \sum_{\rho \in \textup{Irr}(G)} \frac{\textup{dim}(V_{\rho})}{|G|}\lambda^*_{\rho} \sum_{ f \in G}  \textup{Tr}_{\rho}(h^{-1}f) \delta_{f}, \]
where we used
\begin{equation} \label{eq_deltag_decomp}\delta_{g} = \sum_{\rho \in \textup{Irr}(G)} \frac{\textup{dim}(V_{\rho})}{|G|} \rho(g^{-1})^{j}{}_{i} \rho^{i}{}_{j}, \end{equation}
from the Peter-Weyl decomposition in the proof of Lemma~\ref{lem_L_decompose_2}. 

We can put our two bijections on a more parallel footing. For the first bijection, $D^{\vee}(G)$-covariant linear operators $L:\C G\to \C G$ can be viewed as $D^\vee(G)$-$\C G$-bicomodule maps with the right regular coaction of $\C G$. Similarly, bicovariant calculi on $\C(G)$ can be viewed as $D^\vee(G)$-$\C (G)$-bicovariant with trivial left coaction $\tens e$ of the $\C G\subset D^\vee(G)$. Dually for the second bijection, we can equivalently regard $L^{*}$ as a $D^\vee(G)$-$\C(G)$-bicomodule map where the left coaction extends trivially the coaction of $\C(G)$ by $\tens e$. Similarly, $D^\vee(G)$-covariant calculi on $\C G$ can be viewed as calculi which are $D^\vee(G)$-$\C G$-bicovariant (because $\C G$ is cocommutative,  a left-covariant calculus is automatically bicovariant).

For the rest of the section, we now fix a bicovariant $\Omega^{1}_{\mathcal{S}} \cong \C(G) \otimes \Lambda^{1}_{\mathcal{S}}$ with $\mathcal{S}$ an Ad-stable subset of $G \backslash \{ e\}$, i.e. $\mathcal{S} = \cup_{\mathcal{C}} \mathcal{C}$ a union of non-trivial conjugacy classes $\CC\ne\{e\}$. Its structure maps were given by \eqref{eq_kG_calc}. We next study bicovariant bimodule inner products $(\ , \ ): \Omega^{1}_{\mathcal{S}} \otimes_{\C(G)} \Omega^{1}_{\mathcal{S}} \to \C (G)$ in the usual sense (or to match with the above, equivalently as a $D^\vee(G)$-$\C (G)$-bicovariant $(\ , \ )$). As a (possibly degenerate) bimodule inner product on the Cayley graph corresponding to the calculus, we know that this has to have the form 
\[ (\omega_{x \to xc}, \omega_{y \to yd}) = l_{x \to xc} \delta_{x, yd} \delta_{xc, y} \delta_{x}, \] 
where $\omega_{x\to xc}=\delta_x\extd\delta_{xc}$ is the 1-form corresponding to an arrow generated by $c\in\CS$, for some edge weights $\{l_{x\to xc}\in \C\}$. We also have a  left-invariant basis $\{e_{c}:= \sum_{x \in G} \omega_{x \to xc}\}$ of 1-forms over $\C(G)$, where 
\[ (e_{c}, e_{d}) = \delta_{c, d^{-1}} \sum_{x \in G} l_{x \to xc} \delta_{x}\]
is zero if $c\in \CS$ has no inverse in $\CS$. These generate the 1-way arrows. We let  the complement $\mathcal{S}^{+} \subseteq \mathcal{S}$ denote the union  of conjugacy classes $\mathcal{C}$ contained in $\mathcal{S}$ for which $\mathcal{C}^{-1} \subseteq \mathcal{S}$, which generate the bidirectional arrows of the Cayley graph. 

\begin{lemma} \label{lem_metric_k_of_G} For  a bicovariant calculus $\Omega^{1}_{\mathcal{S}}(\C(G))$, a bicovariant bimodule inner product $(\ , \ ): \Omega^{1}_{\mathcal{S}} \otimes_{\C(G)} \Omega^{1}_{\mathcal{S}} \to \C(G)$ corresponds to data $\{l^*_{\CC}\in \C\}$ for each $\CC\subseteq \CS_+$ and has  values 
\[ (e_c,e_d)=\begin{cases}\delta_{c,d^{-1}}l^*_{\CC} & {\rm if\ } c\in \CC\subseteq\CS_+, d\in \CC^{-1}\subseteq\CS_+\\ 0 & {\rm otherwise}\end{cases}.\]
\end{lemma}
\proof For bicovariance in the $\C G$ bimodule form the actions are $h\la\delta_{x}=\delta_{hx}$, $h\la e_c=e_c$, $\delta_x\ra h=\delta_{xh}$ and $e_c\ra h=e_{h^{-1}ch}$ for all $h\in G$. It is then clear that left-invariance of the inner product for $c,d\in \CS^+$ holds if an only if $l_{hx\to hx c}=l_{x\to xc}$ for all $x\in G$, so given by $l_{e\to c}$ at the group identity. It is also clear that right invariance then needs $l_{e\to c}=l_{e\to h^{-1}ch}$, hence a constant, which we denote $l^*_\CC$, on the conjugacy class containing $c$. This is also well-known in other contexts for the $\CS^+=\CS$ case. \endproof

We are now ready to discuss $D^{\vee}(G)$-covariant 2nd order Laplacians on $\C(G)$. As before, the $\C G\subset D^\vee(G)$ part oof this is trivial, so we just mean bicovariant ones, i.e. as in \eqref{eq_def_laplacian} now for elements in $\C(G)$, and $(\ ,\ )$ bicovariant as above. As in Section~\ref{sec_Laplacians_on_CG}, we must have $\square(1) = 0$. Hence, denoting the eigenvalues of $\square$ by $\lambda_{\rho}^*$, we must have $\lambda^*_1= 0$ for the trivial representation (where $\rho^i{}_j=1$ as a function on $G$). 

\begin{proposition} \label{prop_relationship_l_lambda_kofG}
For $\Omega^{1}_{\mathcal{S}}(\C(G))$ bicovariant, $L^*$ in Lemma~\ref{lem_L_decompose_2}  with not necessarily distinct eigenvalues $\{\lambda^*_\rho\}$ such that $\lambda^*_1=0$ is a bicovariant 2nd-order Laplacian $\square$ with respect to a bicovariant inner product given by $\{l^*_\CC\}$ as in Lemma~\ref{lem_metric_k_of_G}, iff 
\begin{equation*} \sum_{\rho \in \textup{Irr}(G)} \frac{\textup{dim}(V_{\rho})^2}{|G|}  \lambda^*_{\rho} \chi_\rho(\CC)  = 0 \end{equation*}
for all conjugacy classes $\CC\subset G\setminus\{e\}$  {\em not} contained in $\CS^+$ and
\begin{equation*} \sum_{\CC\subseteq \CS^+}\sum_{\rho \in \textup{Irr}(G)} \frac{\textup{dim}(V_{\rho})^2}{|G|}  \lambda^*_{\rho} (1+\chi_\rho(\CC)  |\CC|)=0,\quad l^*_{\mathcal{C}}= -\frac{1}{2}\sum_{\rho \in \textup{Irr}(G)} \frac{\textup{dim}(V_{\rho})^{2}}{|G|}\lambda^*_{\rho}\chi_{\rho}(\mathcal{C})^{*}.\end{equation*}
In this case, 
\[ \lambda^*_{\rho} = 2\sum_{\mathcal{C} \subseteq \mathcal{S}^{+}} l^*_{\mathcal{C}}(1 - \chi_{\rho}(\mathcal{C}))|\mathcal{C}|.  \]
\end{proposition}
\proof We note that by \eqref{eq_kG_calc}, 
\begin{align*} (\extd \delta_{g}, \extd \delta_{h}) & = \sum_{c, d \in \mathcal{S}} (\delta_{gc^{-1}} - \delta_{g}) \, (e_c, (\delta_{hd^{-1}} - \delta_{h}) e_d)=\sum_{c,d\in \CS}(\delta_{gc^{-1}}-\delta_g)(\delta_{hd^{-1}c^{-1}}-\delta_{hc^{-1}})\delta_{c,d^{-1}}l^*_{[c]}\\
&= \sum_{c\in \CS^+}l^*_{[c]}(\delta_{gc^{-1}}\delta_h + \delta_g\delta_{hc^{-1}}-\delta_{g,h}\delta_g - \delta_{g,h}\delta_{gc^{-1}}),
\end{align*} 
using Lemma~\ref{lem_metric_k_of_G}, where we set $d=c^{-1}$ if $c\in \CS^+$ and $[c]$ is the conjugacy class containing $c$. Using $L^*(\rho^{i}{}_{j}) = \lambda^*_{\rho}\rho^{i}{}_{j}$  and \eqref{eq_deltag_decomp},  we have that the 2nd order Leibniz rule evaluated on the product $\delta_{g}\delta_{h}$ becomes the equality 
\begin{align*}
 \sum_{\rho} \frac{\textup{dim}(V_{\rho})}{|G|}  \lambda^*_{\rho} \delta_{g, h}  \textup{Tr}_{\rho}(g^{-1}f)  = & \sum_{\rho}  \frac{\textup{dim}(V_{\rho})}{|G|} \lambda^*_{\rho}\Big(\delta_{h, f} \textup{Tr}_{\rho}(g^{-1}h) + \delta_{g, f}\textup{Tr}_{\rho}(h^{-1}g)\Big)\\
 & + 2\sum_{\mathcal{C} \subseteq \mathcal{S}^{+}} l^*_\CC \sum_{c \in \mathcal{C}}  \big( \delta_{g, hc} \delta_{h, f}+ \delta_{h, gc} \delta_{g, f} - \delta_{g, h}\delta_{g, f}- \delta_{g, h}\delta_{gc^{-1}, f}  \big)
\end{align*} 
for all $g, h, f \in G$. We now cover the different possibilities of $g, h, f $. If $g,h,f$ are distinct then both sides as zero and we have no condition. If  $g = h\neq f$ then we need
\[ \sum_{\rho} \frac{\textup{dim}(V_{\rho})}{|G|}  \lambda^*_{\rho} \textup{Tr}_{\rho}(g^{-1}f)  = - 2\sum_{\mathcal{C} \subseteq \mathcal{S}^{+}} l^*_\CC \sum_{c \in \mathcal{C}} \delta_{c, gf^{-1}}. \] So if $gf^{-1} \in \mathcal{C}$ for some conjugacy class $\mathcal{C} \subseteq \mathcal{S}^{+}$, then this gives the stated expression for $l^*_\CC$ (using also that $\textup{Tr}_{\rho}(c^{-1}) = \textup{Tr}_{\rho}(c)^{*}$).  If $m=gf^{-1}\notin\CS^+$ then we get the first constraint stated, since $m\ne e$ under our assumptions. The cases $h=f\ne g$ and $g = f \neq h$ give nothing new, leaving the case $g = h = f$. Substituting this into the original equality gives 
\begin{equation}\label{3rdconstraint}\sum_{\rho} \frac{\textup{dim}(V_{\rho})^{2}}{|G|}  \lambda^*_{\rho}  =  2\sum_{\mathcal{C} \subseteq \mathcal{S}^{+}} l_{\mathcal{C}}^*|\CC|, \end{equation}
since $\{e\}$ is not an allowed conjugacy class in $\CS$. Given the formula for $l^*_\CC$ already obtained, this is equivalent to the second constraint constraint stated on $\{\lambda_\rho^*\}$. 

Finally, for an irrep of interest $\rho'$, we multiply our original condition by $\rho'(f^{-1})^{m}{}_{n}$ and sum over $f \in G$. This gives 
 \begin{align*} &\lambda^*_{\rho'} \delta_{g, h}\rho'(g^{-1})^{m}{}_{n} \\
 &=  \sum_{\rho} \frac{\textup{dim}(V_{\rho})}{|G|}  \lambda^*_{\rho} \big( \textup{Tr}_{\rho}(g^{-1}h) \rho'(h^{-1})^{m}{}_{n} + \textup{Tr}_{\rho}(h^{-1}g)\rho'(g^{-1})^{m}{}_{n}\big)\\
 &\quad + 2\sum_{\mathcal{C} \in \mathcal{S}^{+}}l^*_{\mathcal{C}}\sum_{c \in \mathcal{C}}\big(\delta_{g, hc} \rho'(h^{-1})^{m}{}_{n} + \delta_{h, gc} \rho'(g^{-1})^{m}{}_{n}  - \delta_{g, h}\rho'(h^{-1})^{m}{}_{n} - \delta_{g, h}\rho'(ch^{-1})^{m}{}_{n}\big), 
 \end{align*}
where the left-hand side is by Schur orthogonality. Setting $g = h = e$ and $m = n$, and summing over $m$ then gives
\begin{align*} \lambda^*_{\rho'}\textup{dim}(V_{\rho}') & = 2\sum_{\rho}\frac{\textup{dim}(V_{\rho})^{2}}{|G|}  \lambda^*_{\rho} \textup{dim}(V_{\rho'}) + 2 \sum_{\mathcal{C} \subseteq \mathcal{S}^{+}}\sum_{c \in \mathcal{C}} l^*_{\mathcal{C}} \big(  - \textup{dim}(V_{\rho'})  - \textup{Tr}_{\rho'}(c))\big)\\
& = 4 \textup{dim}(V_{\rho'}) \sum_{\mathcal{C} \subseteq \mathcal{S}^{+}}|\mathcal{C}| \, l^*_{\mathcal{C}}  - 2 \sum_{\mathcal{C} \subseteq \mathcal{S}^{+}} l^*_{\mathcal{C}}( |\mathcal{C}|\textup{dim}(V_{\rho'})+|\mathcal{C}|\chi_{\rho'}(\mathcal{C})\textup{dim}(V_{\rho'})), 
\end{align*} 
using (\ref{3rdconstraint}).  This  gives the stated equation for $\lambda_{\rho}^*$. 
\endproof

Parallel to Section~\ref{sec_Laplacians_on_CG}, we define a bicovariant $L^*$ to be regular if $\lambda^*_{\rho} \neq \lambda^*_{\rho'}$ for all $\rho\ne \rho' \in \textup{Irr}(G)$. From our formula for $\lambda^*_\rho$ in Proposition~\ref{prop_relationship_l_lambda_kofG}, we see that this needs $\CS^+$ to be non-empty. For example, if  $\mathcal{S}^{+} = \mathcal{C}=\CC^{-1}$ is a single conjugacy class then $\square$ is regular iff $l^*_{\mathcal{C}} \neq 0$ and $\chi_{\rho}(\mathcal{C}) \neq \chi_{\rho'}(\mathcal{C})$ for all $\rho\ne \rho' \in \textup{Irr}(G)$. 

Finally, we discuss $\star$-compatibility. We note that the $\star$-structure on $\C(G)$ described in \eqref{eq_star_struc_DveeG} extends to give a $\star$-structure on $\Omega^{1}_{\mathcal{S}}$ iff $\mathcal{S} = \mathcal{S}^{+}$, so we proceed in this case. Here, on the basis 1-forms we have  $e_{c}^{\star} = -e_{c^{-1}}$ \cite{BegMa}.  Next, by the general form of $(e_{c}, e_{d})$ it is clear that  $(\ , \ )$ is $\star$-compatible (in the analogous way to \eqref{eq_star_compatible_inner}) if and only if the $l_{x \to xc}$ (i.e. the weights $l^*_{\mathcal{C}}$ in our context) are all real. If the inner product is `real' in this sense and for example $\CS=\CC=\CC^{-1}$ is a single conjugacy class, our formula for $\lambda^*_\rho$ in Proposition~\ref{prop_relationship_l_lambda_kofG} tells us that $\square$ has real eigenvalues $\lambda^*_\rho$ (because $\chi_{\rho}(\mathcal{C})^* = \chi_{\rho}(\mathcal{C}^{-1})= \chi_{\rho}(\mathcal{C})$ in this case). 

To further manifest the duality, one can in principle explicitly decompose the $D^\vee(G)$-covariant calculi on $\C G$  in Lemma~\ref{lem_DG_covariant_calc_kG} in the form $\Omega^1_T=\oplus_{\rho\in T}\Omega^1_\rho$ on $\C G$ where $T\subseteq {\rm Irr}(G)\setminus\{1\}$ and express the analysis in Lemma~\ref{lemma_metric_laplacian_relation} for a $D^\vee(G)$-covariant $L:\C G\to \C G$ to be a 2nd order Laplacian with respect to a $D^\vee(G)$-covariant inner product $(\ ,\ )$ in these terms. Returning to the $\C(G)$ side, one can next ask for which $\square$ can be realised geometrically by a QRG and WQLCs on $\C(G)$. For example, \cite[Ex. 8.21]{BegMa} explicitly computes the moduli of the latter for $\C (S_{3})$ with respect to the calculus given by $\CC=
\{u,v,w\}$. The matter of a free field description of irreps of $D^\vee(G)$ (or coirreps of $D(G)$) as a dual Wigner construction on $\C(G)$ in the role of `spacetime' appears less interesting due to the simple tensor algebra structure on $D^\vee(G)$.

\begin{example} \label{example_kofS3}\rm We continue our study for the case of $G = S_{3}$, and its three irreps $\rho_i$ with $\rho_0$ the trivial representation, $\rho_{1}$ the sign irrep and $\rho_{2}$ the two-dimensional irrep. We denote the corresponding eigenvalues of a $D^{\vee}(S_{3})$-Laplacian $L^*$ by $\lambda^*_i$, where we previously noted that we must have $\lambda_{0} = 0$. We already made use of the character table 
\[\chi_{\rho_{0}}(g) = 1, \  \chi_{\rho_{1}}(g)= \textup{sgn}(g),\ \chi_{\rho_{2}}(\{  e\}) = {1\over 2},\  \chi_{\rho_{2}}(\{  uv, vu\}) =  - {1\over 2},\  \chi_{\rho_{2}}(\{  u, v, w\}) =  0, \]
for the decomposition into  irreps in Example~\ref{exS3}. There are three possible bicovariant calculi on $\C(S_{3})$ given by taking $\mathcal{S}$ as being $\{uv, vu\}$, $\{ u, v, w\}$ and their union (with this latter calculus being the universal calculus.  We now verify the constraints of Proposition~\ref{prop_relationship_l_lambda_kofG}, noting that $\mathcal{S}^{+} = \mathcal{S}$ for all three choices here since $\mathcal{C}^{-1} = \mathcal{C}$ for all conjugacy classes of $S_{3}$. 

(i) The case $\mathcal{S} = \{ uv, vu\}$. By Proposition~\ref{prop_relationship_l_lambda_kofG}, we have a $D^{\vee}(S_{3})$-covariant bimodule inner product specified by a single number $l^*_{1}$ and we require
\[ \ \frac{1}{6}( - \lambda^*_{1}) = 0,\quad l^*_{1} =  -\frac{1}{12}(\lambda^*_{1} - 2\lambda^*_{2}), \] where the first and second constraints on $\lambda_\rho^*$ both give the first condition above. The unique solution is thus
 \[\lambda^*_{0} =  \lambda^*_{1} = 0 , \quad \lambda^*_{2} = 6l^*_{1}. \] Although there is a one parameter moduli of $D^{\vee}(S_{3})$-covariant 2nd order Laplacians on this calculus, they are not regular. 
 
(ii) The case $\mathcal{S} = \{ u, v, w\}$. By Proposition~\ref{prop_relationship_l_lambda_kofG}, we have a $D^{\vee}(S_{3})$-covariant bimodule inner product specified by a single number $l^*_{2}$ and we require
\[ \frac{1}{6}(\lambda^*_{1} - 2\lambda^*_{2}) = 0,\quad l^*_{2} = \frac{1}{12}\lambda^*_{1}, \] 
where the first and second constraints give the first condition above. This has unique solution
\[ \lambda^*_{0} = 0, \quad \lambda^*_{1} = 12l^*_{2}, \quad \lambda^*_{2} = 6l^*_{2}. \] 
There is a one parameter moduli of $D^{\vee}(S_{3})$-covariant 2nd order Laplacians on this calculus, which are regular iff $l_{2} \neq 0$. 

(iii) The case $\mathcal{S} =\{ uv, vu\} \cup \{ u, v, w\}$. By Proposition~\ref{prop_relationship_l_lambda_kofG}, we have a $D^{\vee}(S_{3})$-covariant bimodule inner product specified by a two numbers $l^*_{1}, l^*_{2}$ corresponding to $\{ uv, vu\}, \{ u,v, w\} \subseteq \mathcal{S}$ respectively and we require 
\[{1\over 12}\lambda^*_1+{1\over 3}\lambda^*_2=0,\quad  l^*_{1} =  -\frac{1}{12}(\lambda^*_{1} - 2\lambda^*_{2}), \quad l^*_{2} = \frac{1}{12}\lambda^*_{1},\]
where the first constraint of Proposition~\ref{prop_relationship_l_lambda_kofG} does not apply since $G\backslash (\mathcal{S} \cup \{ e\})$ is empty. This has solution 
\[ \lambda^*_{0} = 0, \quad \lambda^*_{1} = -8 l^*_{1}, \quad \lambda^*_2=2 l^*_1, \quad l^*_2=-{2\over 3}l^*_1. \] There is thus a 1-parameter moduli of $D^{\vee}(S_{3})$-covariant 2nd order Laplacians on this calculus. The Laplacian is regular iff $l^*_1\ne 0$. 

For all three cases, we have $\star$-compatibility iff $l^*_{1}, l^*_{2} \in \mathbb{R}$, which would then imply $\lambda^*_{1}, \lambda^*_{2} \in \mathbb{R}$ by the above equations.
\end{example}

\section{Braided-Lie algebras $\CL_{\mathcal{C}, \pi}$}\label{secbrack}

Again since $D^\vee(G)$ is coquasitriangular, its bicovariant differential calculi correspond to braided-Lie algebra structures on $\Lambda^1{}^*$, just as classically for Lie groups, where the dual of the space of left-invariant 1-forms recovers the Lie algebra. In this section, we describe these for our calculi $\Omega^1_{\mathcal{C}, \pi}$. Since we won't be considering $\star$-structures, we work in this section over a general field $k$ for the general theory and over $\C$ for examples. 

\subsection{Recap of braided-Lie algebras and braided-Hopf algebras} \label{sec_braided_Lie_recap}

\begin{figure}
\[\includegraphics[scale=0.9]{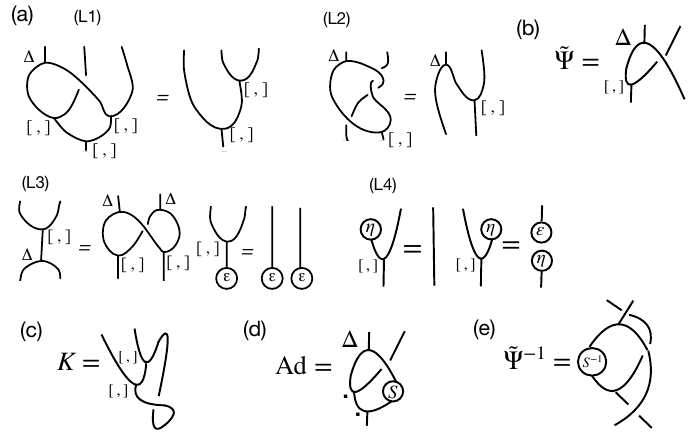}\] 
 \vspace{-5mm}
\caption{\label{figLie} (a) Axioms of a braided-Lie algebra $\CL$ in terms of string diagrams read down the page and with inputs and outputs tensor powers of $\CL$. (b) Fundamental braiding (c) Braided-Killing form (d) Braided-adjoint action of a braided-Hopf algebra with inputs and outputs $B$. Diagrams from \cite{Ma:lie,Ma:alg}. (e) Inverse of fundamental braiding for a left braided-cocommutative Hopf algebra $B$ with respect to $\textup{Ad}$.}  \end{figure}

Informally, we recall that a braided category is a category with a tensor product functor and an `associator' natural isomorphism obeying Mac Lane's `pentagon' identity (i.e. a monoidal category), together with a natural `braiding' natural isomorphism $\Psi: \tens\Rightarrow \tens^{op}$, denoted above as isomorphisms $\Psi_{V,W}:V\tens W\to W\tens V$ for all objects $V,W$, and subject to `hexagon' identities. There is a natural diagrammatic notation for `braided algebra' in such categories \cite{Ma:alg,Ma}, where $\tens$ is indicated by juxtaposition, the unit object $\underline{1}$ for the tensor product is denoted by omission and most morphisms are denoted by strings flowing down the page. The $\tens$ can w.l.o.g. be taken strictly associative (inserting brackets as needed to make sense of expressions), the braiding is denoted $\Psi=\includegraphics{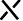}$ and $\Psi^{-1}$ is denoted  by the reverse braid crossing. For a coalgebra in the braided category we write $\Delta=\includegraphics{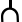}$, while for an algebra in the braided category (as well as for the braided-Lie bracket) we use $\includegraphics{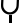}$. With this notation, a {\em left braided-Lie algebra} $\CL$ \cite{Ma:lie} is a coalgebra in a braided category together with a left `bracket' $[\ ,\ ]:\CL\tens\CL\to \CL$ subject to the axioms (L1)-(L3) in Figure~\ref{figLie}. $\mathcal{L}$ is {\em unital} if there exists a map of coalgebras $\eta: \underline{1} \to \mathcal{L}$ (with $\underline{1}$ grouplike) such that (L4) is additionally satisfied. The map $\tilde\Psi: \CL\tens\CL\to \CL\tens\CL$ shown in part (b) obeys the braid relations and is called the {\em fundamental braiding} and the braided-Lie algebra is called {\em regular} if this is invertible. It is also possible to write the axioms of a braided-Lie algebra in terms of this morphism.

To construct braided-Lie algebras, it will be useful to also have the notion of a {\em braided-Hopf algebra} $B$ in a braided category \cite{Ma,Ma:alg}. Like a braided-Lie algebra, this is a coalgebra in the category but now equipped with an associative algebra product $\cdot:B\tens B\to B$ in the category obeying a `bialgebra' axiom the same as (L3) but for the product in place of the bracket. The unit of the algebra is now a map $\eta:\underline{1}\to B$. There is also an antipode $S:B\to B$ with axioms of the same form as for a usual Hopf algebra but in the category. As for ordinary Hopf algebras, every braided-Hopf algebra acts on itself by a braided-adjoint action $\Ad$, shown for the left one in part (d). The axioms of a left braided-Lie algebra are such that $B$ itself with $[\ ,\ ]$ the braided left adjoint action obeys (L1),(L3), (L4) while for (L2) $B$ would need to be {\em left braided-cocommutative} in a certain sense with respect to $\Ad$. In this case the braided-Lie algebra is regular if the braided-antipode is invertible, with $\widetilde{\Psi}^{-1}$ given by Figure~\ref{figLie}(e) as can be verified by composition. More details are in \cite{Ma:alg,BegMa}.

Now let $H$ be a coquasitriangular Hopf algebra. Then its category $\CM^H$ of right comodules, say, is braided. Braided Hopf algebras in such categories can be canonically obtained through transmutation \cite{Ma}, which turns any coquasitriangular Hopf algebra $H$ into its braided version $B(H)$ in the category $\CM^H$  by the right adjoint coaction $\Ad_R$ and a modified product. A special property of these is that the associated quantum Killing form $\mathcal{Q}: H \otimes H \to k$ upgrades to a self-duality pairing $\mathcal{Q}: B(H) \otimes B(H) \to k$ of braided Hopf algebras \cite{Ma}, which in the finite dimensional case gives a braided Hopf algebra map \begin{equation} \label{eq_morphism_braided_hopf_algebras} \mathcal{Q}_{1}: B(H) \to B(H^{*}), \quad h \mapsto \mathcal{Q}(h \otimes (\ ))\end{equation} to the left module transmutation $B(H^{*})$ of the quasitriangular Hopf algebra $H^{*}$ (linearly) dual to $H$. We identify left $H^*$ modules and right $H$-comodules by evaluating the latter. The module transmutation theory  is dual and modifies instead the coproduct of $H^*$. In fact $B(H^*)$ is the left categorical dual object to $B(H)$ as an $H$-comodule, via the evaluation and coevaluation maps in $\mathcal{M}^{H}$ \begin{equation} \label{eq_coev_ev} \textup{ev}(\phi \otimes h) = \langle S\phi, h\rangle , \textup{coev}(1) = e_{a} \otimes S^{-1}(f^{a}) , \quad \phi \in H^{*}, h \in H, \end{equation} for a basis $\{ e_{a}\}$ of $H$ with dual $\{ f^{a}\}$. In the factorisable case, the module and comodule transmutations are thus isomorphic, and so these are self-dual braided Hopf algebras. 

Another key feature of transmutation theory is that $B(H^*)$ is always braided cocommutative with respect to its left braided adjoint action, and therefore also forms a unital regular left braided Lie algebra by previous discussion. Also, the left braided adjoint action here coincides on the vector space with the usual left adjoint action of $H^*$ on itself. Less obvious is that $B(H)$ is also braided cocommutative in the required sense and hence is also a unital regular braided-Lie algebra. Moreover,  $\mathcal{Q}_{1}$ intertwines the braided-Lie brackets and is thus a map of unital left braided Lie algebras. More details are in \cite{BegMa,MaMc}.

\subsection{Braided-Lie algebra $BD(G)$}

In this section we now specialise these constructions to $H=D^\vee(G)$. A version of the transmutation theory for the quantum double first appeared in \cite{Ma:skl} but we give a much simpler and more direct realisation. 
   
 \begin{lemma} \label{lem_braided_Hopf_DG}  $B(D^{\vee}(G)) = B(D(G))$ can be identified as braided Hopf algebras, and we will refer to them together as $BD(G)$. This is a braided-Hopf algebra in the category of $G$-crossed modules and has algebra (\ref{DG}) of $D(G)$ and the coalgebra (\ref{DvG}) of $D^\vee(G)$, and invertible braided antipode \[ \underline{S}(\delta_{g} \otimes h)  = \delta_{\{ g^{-1}, h\}g^{-1}} \otimes gh^{-1}g^{-1}. \] The $G$-crossed-module structure is 
 \begin{equation} \label{eq_crossed_DG_transmute} f\la (\delta_g\tens h)=\delta_{fgf^{-1}}\tens fhf^{-1},\quad  |\delta_g\tens h|=\{g^{-1},h\}, \end{equation}
 for all $f,g,h\in G$ and $\{ h, g\}:= h^{-1}g^{-1}hg$  the group commutator. 
 \end{lemma}
\proof A key feature of comodule transmutation is that a map of coquasitriangular Hopf algebras $f: H \to H'$ becomes a map $f: B(H) \to B(H^{'})$ of the comodule transmutations when viewing $B(H)$ inside $\mathcal{M}^{H'}$ by pushforward. Thus, the isomorphism of coquasitriangular Hopf algebras $\beta^{-1}: D^{\vee}(G) \cong D(G)^{*}$ of \eqref{eq_iso_beta} gives an isomorphism $B(D^{\vee}(G)) \cong B(D(G)^{*})$ of braided Hopf algebras. Composing this with \eqref{eq_morphism_braided_hopf_algebras} (with $\mathcal{Q}_{1}$ here the isomorphism \eqref{eq_quantum_killing_form_DG}), and recalling that $\beta^{-1} \circ \mathcal{Q}_{1} = \id$ gives that $B(D^{\vee}(G)) = B(D(G))$. 

The structure maps of $BD(G)$ are now computed using comodule transmutation of $D^{\vee}(G)$: the braided Hopf algebra $B(D^{\vee}(G))$ lives in $\mathcal{M}^{D^{\vee}(G)}$ by the right adjoint coaction of $D^{\vee}(G)$, which can be computed as 
 \begin{equation}\label{eq_adjoint_actions_DG_1} \Ad_{R}(\delta_{g} \otimes h)  = \sum_{f \in G} \delta_{fgf^{-1}} \otimes fhf^{-1} \otimes \delta_{f} \otimes \{ g^{-1}, h\},\end{equation} which by comparing with \eqref{eq_comodule_DveeG} gives the stated $G$-crossed module structure. Moreover, by general theory of comodule transmutation the coalgebra structure of $B(D^{\vee}(G))$ coincides with that of $D^{\vee}(G)$. Similarly, by module transmutation we would get that $B(D(G))$'s algebra structure is that of $D(G)$, which is thus also that of $B(D^{\vee}(G))$ as $B(D^{\vee}(G)) = B(D(G))$. The braided antipode is then recovered from general transmutation theory \cite{Ma}, and it can be quickly verified that it is indeed the antipode (and that $\underline{S}= \underline{S}^{-1}$ since $S = S^{-1}$ for $D(G)$ and $D^{\vee}(G)$). \endproof
 
 We remark that $k(G)$ is trivially coquasitriangular and so its braided version is the same but with structure maps regarded in the category of $G$-modules (or right $k(G)$ comodules), where $G$ acts by the dual adjoint action on $k(G)$. We can regard it in the category of $G$-crossed modules with trivial grading. Likewise $k G$ is trivially quasitriangular and its braided version has the same structure maps and can also be regarded in the $G$-crossed module category with the adjoint action on $k G$ by conjugation, and trivial grading. Then 
 \begin{equation} \label{eq_exact_sequence_braided_hopf} k(G) \hookrightarrow BD(G) \twoheadrightarrow k G, \end{equation}
  in the category of $G$-crossed modules is `exact'  in the sense that $k(G)$ is the braided Hopf algebra kernel of $BD(G) \twoheadrightarrow kG$ and $kG$ is the braided Hopf algebra cokernel of $k(G) \hookrightarrow BD(G)$. We also know from the general theory that $BD(G)$ is a braided-Lie algebra.

\begin{corollary} \label{cor_braided_lie_DG}
The underlying braided coalgebra of $BD(G)$ is a unital regular left braided Lie algebra in the category of $G$-crossed modules, with Lie bracket and fundamental braiding given by \[  [\delta_{u}  v, \delta_{g}  h] = v \triangleright \delta_{g}h\, \delta_{u,|v\triangleright \delta_{g}h|},\quad \widetilde{\Psi}(\delta_{u} v \otimes \delta_{g}  h) = v \triangleright \delta_{g}h \otimes \delta_{|v \triangleright \delta_{g}h|^{-1}u} \otimes v. \] We denote this braided Lie algebra by $\mathcal{L}_{BD(G)}$. \end{corollary}

\proof By left module transmutation theory for a quasitriangular Hopf algebra $H$, the left braided adjoint action of $B(H)$ on itself coincides with the ordinary braided adjoint action. In our case of interest, noting that $B(D(G))$ sits in the category $G$-crossed modules by its adjoint action, then comparing the $G$-crossed module structure of Lemma~\ref{lem_braided_Hopf_DG} with \eqref{eq_DG_action} gives that $\Ad_{\delta_u\tens v}(\delta_g\tens h)$ is given by the above, which gives the bracket on $BD(G)$. The fundamental braiding is computed from Figure~\ref{figLie}(b) as
\begin{align*} 
\widetilde{\Psi}(\delta_{u}v \otimes \delta_{g}h) & = \sum_{f \in G}([\_, \_] \otimes \id) \circ (\id \otimes \Psi_{\mathcal{L}_{BD(G)}, \mathcal{L}_{BD(G)}})(\delta_{f} f^{-1}uvu^{-1}f \otimes \delta_{f^{-1}u} v \otimes \delta_{g}h) \\
& =\sum_{f \in G} [\delta_{f}v \{uf^{-1}, v\}^{-1} \delta_{\{ uf^{-1}, v\}g\{uf^{-1}, v\}^{-1}}\{ uf^{-1}, v\}h\{ uf^{-1}, v\}^{-1}] \otimes \delta_{f^{-1}u} v \\
& = \sum_{f \in G}  \delta_{f, |\delta_{vgv^{-1}}vhv^{-1}|} \delta_{vgv^{-1}}vhv^{-1} \otimes \delta_{f^{-1}u}v = \delta_{vgv^{-1}}vhv^{-1} \delta_{v\{g^{-1}, h\}v^{-1}u}v,
 \end{align*} 
 where we recall to use the coproduct of $D^{\vee}(G)$ for $BD(G)$, and the braid crossing $\Psi_{\mathcal{L}_{BD(G)}, \mathcal{L}_{BD}(G)}$ is given on the basis by 
 \[\Psi_{\mathcal{L}_{BD(G)}, \mathcal{L}_{BD}(G)}(\delta_{u}v \otimes \delta_{g}h) = |\delta_{u}v| \triangleright \delta_{g}h \otimes \delta_{u}v = \delta_{\{ u^{-1}, v\} g\{ u^{-1}, v\}^{-1}} \otimes \{ u^{-1}, v\} h\{ u^{-1}, v\}^{-1} \otimes \delta_{u}v,\]
  using \eqref{GcrossPsi} with the $G$-crossed module structure of Lemma~\ref{lem_braided_Hopf_DG}.  We note that in the second line, we rewrote $f^{-1}uvu^{-1}f = vv^{-1}f^{-1}uvu^{-1}f = v \{v, uf^{-1}\} = v\{ uf^{-1}, v\}^{-1}$.  From Lemma~\ref{lem_braided_Hopf_DG} we know that the braided antipode is invertible, hence so is $\widetilde{\Psi}$. \endproof

The braided-Lie bracket in $k(G)$ regarded as a braided-Hopf algebra is trivial $[f,g]=\eps(f)g$ with trivial braiding. Its fundamental braiding is also the trivial flip map. The braided-Lie bracket on $kG$ is $[u,v]=uvu^{-1}$ and again with the trivial braiding and fundamental braiding given by conjugation as for racks. It is easy to see that 
\begin{equation} \label{eq_sequence_braided_lie} \mathcal{L}_{k(G)} \hookrightarrow \mathcal{L}_{BD(G)} \twoheadrightarrow \mathcal{L}_{kG}, \end{equation} 
when all three are viewed in the category of $G$-crossed modules. One can also check that  the fundamental braiding of $\CL_{BD(G)}$ as a braided-Lie algebra when viewed functorially in the category of $D(G)$-crossed modules \cite{Ma} coincides with the left $D(G)$-crossed module braiding obtained from viewing $D(G)$ as a $D(G)$-crossed module by the left adjoint action and the regular action, namely 
\begin{equation} \label{eq_fund_equals_crossed} \xPsi_{D(G), D(G)}(\delta_{u} v \otimes \delta_{g} h) = \Ad_{(\delta_{u} v)_{1}}(\delta_{g}  h) \otimes (\delta_{u} v)_{2} \, . \end{equation} 
Dualising this and working equivalently with right  $D^\vee(G)$-crossed modules, one finds that  $D^{\vee}(G)$ has the right adjoint coaction and right regular action. This is exactly the  $D^{\vee}(G)$-crossed module structure on $D^{\vee}(G)$ used to construct bicovariant calculi on $D^\vee(G)$ as seen in Section~\ref{secPcalc}. More precisely, $\CL_{BD(G)}^+$ (the kernel of the counit) is a braided-Lie algebra dual to $\Lambda^1=D^\vee(G)^+$ for the universal calculus on $D^\vee(G)$. 
 
\subsection{Braided-Lie algebra associated to $\mathcal{C},\pi$} 

The braided-Lie algebra $\CL_{BD(G)}$ restricts to a braided Lie algebra structure on each block $\textup{End}(V_{\mathcal{C},\pi})$. Slightly more generally, we consider
\begin{equation}\label{Lsum} \CL=\bigoplus_{(\CC,\pi)\in S} \textup{End}(V_{\mathcal{C},\pi}),\end{equation}
where each block is a matrix braided-Lie algebra $\CL_{\CC,\pi}$ but we sum over some collection $S$ of non-trivial blocks. It turns out that there are non-trivial brackets and fundamental braidings between these blocks which we also provide.

\begin{proposition} \label{prop_braided_Lie_end}Given $V_{\mathcal{C},\pi}$,  $\End(V_{\mathcal{C}, \pi})$ is a $G$-crossed module structure by \eqref{eq_right_comodule_EndVCpi} and becomes a unital regular left braided-Lie algebra $\CL_{\CC,\pi}$ in the category of $G$-crossed modules, with the matrix coalgebra structure,   \[ \Delta E_{ai}{}^{bj} = \sum_{c \in \mathcal{C}}E_{ai}{}^{ck} \otimes E_{ck}{}^{bj}, \quad \epsilon(E_{ai}{}^{bj}) = \delta_{a}{}^{b}\delta_{i}{}^{j}, \] 
and  braided Lie bracket given by  
\[  [E_{ai}{}^{bj}, E_{ck}{}^{dl}]=\delta_{|E_{ai}{}^{bj}||E_{ck}{}^{dj}|, |b^{-1} \triangleright E_{ck}{}^{dl} |} \, \pi(\zeta_{a}(|b^{-1} \triangleright E_{ck}{}^{dl}|^{-1}))^{j}{}_{i}\,  b^{-1} \triangleright E_{ck}{}^{dl} , \] and fundamental braiding by \[ \widetilde{\Psi}(E_{ai}{}^{bj} \otimes E_{ck}{}^{dl}) =  b^{-1} \triangleright E_{ck}{}^{dl} \otimes |b^{-1} \triangleright E_{ck}{}^{dl}|^{-1} \triangleright E_{a i}{}\otimes E^{bj}. \] 
There are similar formulae between blocks in the case of (\ref{Lsum}). 
\end{proposition}
\proof Recall from Section~\ref{secPcalc} that for the construction of bicovariant calculi on $D^{\vee}(G)$ with $\Lambda^{1} \cong \bigoplus_{(\mathcal{C},\pi) \in \mathcal{S}} \textup{End}(V_{\mathcal{C},\pi})$ for some subcollection $\mathcal{S}$ of the set of non-trivial pairs $(\mathcal{C},\pi)$, we used the fact that $\varpi = \bigoplus_{(\mathcal{C},\pi) \in \mathcal{S}}\rho_{\mathcal{C},\pi}\circ \mathcal{Q}_{1}\circ \beta^{-1}: D^{\vee}(G) \twoheadrightarrow \Lambda^{1}$ was a $D^{\vee}(G)$-crossed module quotient map. The given formulae provide calculi over any field, only the result that all calculi were of this form needed the Artin-Wedderburn theorem which we took over $\C$.  Forgetting the action part, this is thus a right $D^{\vee}(G)$-comodule map, or equivalently $G$-crossed module map with $D^{\vee}(G)$ a $G$-crossed module by \eqref{eq_crossed_DG_transmute}, via $D(G)$ with the same $G$-crossed module structure. The space $ \Lambda^{1}$ is a $G$-crossed module by \eqref{eq_right_comodule_EndVCpi}. We also saw that $\mathcal{Q}_{1} \circ \beta^{-1}: B(D^{\vee}(G)) \to BD(G)$ was (trivially) a morphism of braided Hopf algebras with respect to these $G$-crossed module structures. Also, as usual we get that $\bigoplus_{(\mathcal{C},\pi) \in \mathcal{S}}\rho_{\mathcal{C},\pi}: D(G)\to \bigoplus_{(\mathcal{C},\pi) \in \mathcal{S}} \textup{End}(V_{\mathcal{C},\pi})$ is an algebra map with $\bigoplus_{(\mathcal{C},\pi) \in \mathcal{S}} \textup{End}(V_{\mathcal{C},\pi})$ an algebra by $E_{ai}{}^{bj}\cdot E_{c'k'}{}^{d'l'} = \delta_{b, c'}\delta_{\pi, \pi'}\delta^{j}_{k'} E_{ai}{}^{d'l'}$ with respect to the basis $\{ E_{ai}{}^{bj}\}$  of $\textup{End}(V_{\mathcal{C},\pi})$ and $\{ E_{c'k'}{}^{d'l'} \}$ of $\textup{End}(V_{\mathcal{C}', \pi'})$ (and the result is indeed again an element of $\bigoplus_{(\mathcal{C},\pi) \in \mathcal{S}} \textup{End}(V_{\mathcal{C},\pi})$ due to $\delta_{b, c'}$ which imposes $\delta_{\mathcal{C}, \mathcal{C}'}$) and $1 = \bigoplus_{(\mathcal{C},\pi) \in \mathcal{S}} \sum_{c \in \mathcal{C}} E_{ci}{}^{ci} $. It can be verified that this algebra map is also compatible with the $G$-crossed module structures,
 \begin{align*} (h \triangleright & E_{ai}{}^{bj}) \cdot (h \triangleright E_{c'k'}{}^{d'l'})  \\
 &= \pi_{\zeta_{a}}(h)^{x}{}_{i}\pi(\zeta_{b}(h)^{-1})^{j}{}_{y}\pi'_{\zeta_{c'}}(h)^{u'}{}_{k'}\pi'(\zeta_{d'}(h)^{-1})^{l'}{}_{v'} E_{hah^{-1}x}{}^{hbh^{-1}y} \cdot E_{hc'h^{-1}u'}{}^{hd'h^{-1}v'} \\
& =\pi_{\zeta_{a}}(h)^{x}{}_{i}\pi(\zeta_{b}(h)^{-1})^{j}{}_{y}\pi'_{\zeta_{c'}}(h)^{u'}{}_{k'}\pi'(\zeta_{d'}(h)^{-1})^{l'}{}_{v'} \delta_{hbh^{-1}, hc'h^{-1}}\delta_{\pi, \pi'} \delta^{y}_{u'} E_{hah^{-1}x}{}^{hd'h^{-1}v'}\\
& = \delta_{\mathcal{C}, \mathcal{C}'} \delta_{\pi, \pi'}\delta_{b, c}\pi_{\zeta_{a}}(h)^{x}{}_{i}\pi(\zeta_{b}(h)^{-1}\zeta_{b}(h))^{j}{}_{k}\pi(\zeta_{d}(h)^{-1})^{l}{}_{v}  E_{hah^{-1}x}{}^{hdh^{-1}v} \\
& = \delta_{\mathcal{C}, \mathcal{C}'} \delta_{\pi, \pi'} \delta_{b, c}\delta^{j}_{k} \pi_{\zeta_{a}}(h)^{x}{}_{i}\pi(\zeta_{d}(h)^{-1})^{l}{}_{v}  E_{hah^{-1}x}{}^{hdh^{-1}v} = \delta_{\mathcal{C}, \mathcal{C}'} \delta_{\pi, \pi'} \delta_{b, c} \, h \triangleright (E_{ai}{}^{bj}\cdot E_{c'k'}{}^{d'l'}), \end{align*} 
\begin{align*}
\sum_{f \in G}(\delta_{f} \triangleright E_{ai}{}^{bj}) \cdot  (\delta_{f^{-1}g} &\triangleright E_{c'k'}{}^{d'l'})  = \sum_{f \in G} \delta_{f, |E_{ai}{}^{bj}|} \delta_{f^{-1}g, |E_{c'k'}{}^{d'l'}|} E_{ai}{}^{bj}\cdot E_{c'k'}{}^{d'l'} \\
&= \delta_{g, ab^{-1}c'd'^{-1}} \delta_{b, c'} \delta_{\pi, \pi'} \delta^{j}_{k} E_{ai}{}^{d'l'} = \delta_{g, ad'^{-1}} \delta_{b, c'}\delta_{\pi, \pi'} \delta^{j}_{k} E_{ai}{}^{d'l'}\\
& =\delta_{b, c'}\delta_{\pi, \pi'} \delta^{j}_{k}  \delta_{g} \triangleright  E_{ai}{}^{d'l'}  = \delta_{g} \triangleright (E_{ai}{}^{bj} \cdot E_{c'k'}{}^{d'l'}).\end{align*} 
Hence, as the algebra structure of $D(G)$ is that of $BD(G)$, when we dualise $\varpi$ in the category of $G$-crossed modules, we  obtain an inclusion of braided coalgebras \[ (\bigoplus_{(\mathcal{C},\pi) \in \mathcal{S}}\rho_{\mathcal{C},\pi}\circ \mathcal{Q}_{1}\circ \beta^{-1})^{\sharp}: \bigoplus_{(\mathcal{C},\pi) \in \mathcal{S}} \textup{End}(V_{\mathcal{C},\pi})^{\sharp} \hookrightarrow B(D^{\vee}(G))^{\sharp}. \] We saw in Lemma~\ref{lem_braided_Hopf_DG} that we can take $B(D^{\vee}(G))^{\sharp} = BD(G)^{\sharp}\cong BD(G)$ as a $G$-crossed module via the evaluation and coevaluation maps 
\begin{equation} \label{eq_coev_ev_DG}\textup{ev}(\delta_{g}h \otimes \delta_{u}v) = \langle S(\delta_{g}h), \beta^{-1}(\delta_{u}v) \rangle , \quad \textup{coev}(1) = \sum_{u, v \in G} \beta(u \delta_{v})\otimes S^{-1}(\delta_{u}v), \end{equation} 
by \eqref{eq_coev_ev}, and the isomorphism $\beta$. We also have that $\textup{End}(V_{\mathcal{C},\pi})^{\sharp} = (V_{\mathcal{C},\pi}^{\sharp})^{\sharp} \otimes V_{\mathcal{C},\pi}^{\sharp}$, where we previously have seen that we can realise $V_{\mathcal{C},\pi}^{\sharp}$ on the linear dual $V_{\mathcal{C},\pi}^{*}$ with $G$-crossed module structure given by the second part of \eqref{eq_crossed_module_V}. Using the categorical dual crossed module \eqref{eq_categorical_dual_crossed_module} again, it can be shown that $(V_{\mathcal{C},\pi}^{\sharp})^{\sharp}$ realised on $(V_{\mathcal{C},\pi}^{*})^{*} \cong V_{\mathcal{C},\pi}$ has $G$-crossed module given by the first part of \eqref{eq_crossed_module_V}. Hence, we obtain that $\textup{End}(V_{\mathcal{C},\pi})^{\sharp} = \textup{End}(V_{\mathcal{C},\pi})$ as a $G$-crossed module (with respect to the standard linear pairing, but extended to tensor products in the categorical way). With respect to this pairing, the coalgebra structure on the block $\textup{End}(V_{\mathcal{C},\pi})$ can be computed as in \cite{BegMa} as
 \begin{align*} &\Delta E_{ai}{}^{bj} = (\textup{ev} \otimes \id^{\otimes 2})\circ (\id \otimes \cdot \otimes \id^{\otimes 2}) \circ (\id^{\otimes 2} \otimes \textup{coev}(1) \otimes \id)\circ (\id \otimes \textup{coev}(1))E_{ai}{}^{bj}\\
& = \sum_{(\mathcal{C}', \pi'), (\mathcal{C}'', \pi'')} \sum_{c', d' \in \mathcal{C}', f'', g'' \in \mathcal{C}''}(\textup{ev}\otimes \id^{\otimes 2})(E_{ai}{}^{bj} \otimes E_{c'k'}{}^{d'l'} \cdot E_{f''p''}{}^{g''q''} \otimes E_{g''q''}{}^{f''p''} \otimes E_{d'l'}{}^{c'k'})\\
& =\sum_{(\mathcal{C}', \pi'), (\mathcal{C}'', \pi'')} \sum_{c', d' \in \mathcal{C}', f'', g'' \in \mathcal{C}''} \textup{ev}(E_{ai}{}^{bj}, E_{c'k'}{}^{g''q''}) \delta_{d', f''} \delta_{\pi', \pi''} \delta^{l'}_{p''}  \otimes  E_{g''q''}{}^{f''p''} \otimes E_{d'l'}{}^{c'k'}\\
& = \sum_{(\mathcal{C}', \pi')} \sum_{c', d', g' \in \mathcal{C}'} \delta_{b, c'}\delta_{a, g'}\delta_{\pi, \pi'}\delta^{j}_{k'}\delta_{i}^{q'} \otimes  E_{g'q'}{}^{d'l'} \otimes E_{d'l'}{}^{c'k'} = \sum_{d \in \mathcal{C}} E_{ai}{}^{dl} \otimes E_{dl}{}^{bj}, \end{align*} as stated. The coproduct on $\bigoplus_{(\mathcal{C},\pi) \in \mathcal{S}}\textup{End}(V_{\mathcal{C},\pi})$ is then simply the direct sum of the above. Here, we note that we must use the categorical pairings $\textup{ev}(E_{ai}{}^{bj} \otimes E_{c'k'}{}^{d'l'}) = \delta_{a, d'} \delta_{b, c'} \delta_{\pi, \pi'}\delta_{i}^{l'} \delta_{k'}^{j}$ and $\textup{coev}(1) = \sum_{\mathcal{C},\pi} \sum_{a, b \in \mathcal{C}} E_{ai}{}^{bj} \otimes E_{bj}{}^{ai}$. 

Next, as $\mathcal{Q}_{1} \circ \beta^{-1} = \id $, we can compute the map $(\bigoplus_{(\mathcal{C},\pi) \in \mathcal{S}}\rho_{\mathcal{C},\pi})^{\sharp} = \bigoplus_{(\mathcal{C},\pi) \in \mathcal{S}}\rho_{\mathcal{C},\pi}^{\sharp}: \textup{End}(V_{\mathcal{C},\pi}) \hookrightarrow BD(G)$ as 
 \begin{align*}
\rho_{\mathcal{C},\pi}^{\sharp}(E_{ai}{}^{bj}) & = (\textup{ev} \otimes \id)\circ  (\id \otimes \rho_{\mathcal{C},\pi} \otimes \id)\circ (E_{ai}{}^{bj} \otimes \textup{coev}(1)) \\
&= \sum_{u, v \in G} \textup{ev}(E_{ai}{}^{bj} \otimes \rho_{\mathcal{C},\pi}(\beta(u \delta_{v}))) \, S^{-1}(\delta_{u}v)\\
& = \sum_{u, v \in G} \textup{ev}(E_{ai}{}^{bj} \otimes \rho_{\mathcal{C},\pi}(\delta_{v} \otimes v^{-1}uv) \, \delta_{v^{-1}u^{-1}v} \otimes v^{-1} \\
& = \sum_{u, v \in G, c, d \in \mathcal{C}} \delta_{v, c} \delta_{v, v^{-1}uvdv^{-1}u^{-1}v}\pi(\zeta_{d}(v^{-1}uv))^{k}{}_{l}\textup{ev}(E_{ai}{}^{bj} \otimes E_{ck}{}^{dl}) \, \delta_{v^{-1}u^{-1}v} \otimes v^{-1}\\
& = \sum_{u \in G} \pi(\zeta_{a}(b^{-1}ub))^{j}{}_{i} \delta_{b, b^{-1}ubab^{-1}u^{-1}b} \delta_{b^{-1}u^{-1}b} \otimes b^{-1} = \sum_{n \in C_{G}} \pi(n)^{j}{}_{i} \delta_{q_{a}n^{-1}q_{b}^{-1}} \otimes b^{-1},
\end{align*}
where $S^{-1}$ is the antipode of $D(G)$. In the last equality, we defined $n:= \zeta_{a}(b^{-1}ub)$ and note that 
\[ n \delta_{b, b^{-1}ubab^{-1}u^{-1}b} = q_{b^{-1}ubab^{-1}u^{-1}b}^{-1}b^{-1}ubq_{a}\delta_{b, b^{-1}ubab^{-1}u^{-1}b} = q_{b}^{-1}b^{-1}ubq_{a} \delta_{b, b^{-1}ubab^{-1}u^{-1}b},\] 
implying $b^{-1}ub\delta_{b, b^{-1}ubab^{-1}u^{-1}b} = q_{b}nq_{a}^{-1} \delta_{b, b^{-1}ubab^{-1}u^{-1}b}.$  We emphasise that this is not the same map as \eqref{eq_artin_wedderburn_DG}. 

Finally, we note that since the fundamental braiding in Corollary~\ref{cor_braided_lie_DG} is described completely in terms of the $G$-crossed module structure of $BD(G)$, and since $\rho^{\sharp}_{\mathcal{C},\pi}$ is $G$-crossed module map then it is unsurprising that we get the following result 
 \begin{align*} 
& \widetilde{\Psi}(\rho^{\sharp}_{\mathcal{C},\pi}(E_{ai}{}^{bj}) \otimes \rho^{\sharp}_{\mathcal{C}', \pi'}(E_{c'k'}{}^{d'l'})) \\
&=  \sum_{n \in C_{G}, m \in C'_{G}}\pi(n)^{j}{}_{i}\pi'(m)^{l'}{}_{k'}\widetilde{\Psi}( \delta_{q_{a}n^{-1}q_{b}^{-1}} b^{-1} \otimes \delta_{q_{c'}m^{-1}q_{d'}^{-1}} d'^{-1})\\
& = \sum_{n \in C_{G}, m \in C'_{G}}\pi(n)^{j}{}_{i}\pi'(m)^{l'}{}_{k'} b^{-1} \triangleright \delta_{q_{c'}m^{-1}q_{d'}^{-1}} d'^{-1} \otimes \delta_{| b^{-1} \triangleright \delta_{q_{c'}m^{-1}q_{d'}^{-1}} d'^{-1}|^{-1}q_{a}n^{-1}q_{b}^{-1}}b^{-1}\\
& =\sum_{n \in C_{G}, m \in C'_{G}}\pi(n)^{j}{}_{x}\pi(\zeta_{a}(b^{-1}d'c'^{-1}b))^{x}{}_{i}\pi'(m)^{l'}{}_{k'} b^{-1} \triangleright \delta_{q_{c'}m^{-1}q_{d'}^{-1}} d'^{-1} \otimes \delta_{q_{s}n^{-1}q_{b}^{-1}}b^{-1}\\
&  = \sum_{m \in C'_{G}}\pi(\zeta_{a}(b^{-1}d'c'^{-1}b))^{x}{}_{i}\pi'(m)^{l'}{}_{k'} b^{-1} \triangleright \delta_{q_{c'}m^{-1}q_{d'}^{-1}} d'^{-1} \otimes \rho^{\sharp}_{\mathcal{C},\pi}(E_{sx}{}^{bj})\\
& = b^{-1} \triangleright \rho^{\sharp}_{\mathcal{C}', \pi'}(E_{c'k'}{}^{d'l'}) \otimes \rho^{\sharp}_{\mathcal{C},\pi}(b^{-1}d'c'^{-1}b \triangleright E_{ai} \otimes E^{bj}) \\
& = \rho^{\sharp}_{\mathcal{C}', \pi'}(b^{-1} \triangleright E_{c'k'}{}^{d'l'}) \otimes \rho^{\sharp}_{\mathcal{C},\pi}(|b^{-1}\triangleright E_{c'k'}{}^{d'l'}|^{-1} \triangleright E_{ai} \otimes E^{bj}),\end{align*}
where $s := (b^{-1}d'c'^{-1}b)a(b^{-1}d'c'^{-1}b)^{-1} \in \mathcal{C}$. For the third equality, we used that 
\[ | b^{-1} \triangleright \delta_{q_{c'}m^{-1}q_{d'}^{-1}} d'^{-1}| = b^{-1}\{ q_{d'}mq_{c'}^{-1}, d'^{-1} \}b = b^{-1}q_{c'}m^{-1}q_{d'}^{-1}d'q_{d'}mq_{c'}^{-1}d'^{-1}b = b^{-1}c'd'^{-1}b.\]
We  then wrote 
\[ | b^{-1} \triangleright \delta_{q_{c'}m^{-1}q_{d'}^{-1}} d'^{-1}|^{-1} q_{a}n^{-1}q_{b}^{-1} =b^{-1}d'c'^{-1}b q_{a}n^{-1}q_{b}^{-1} = q_{s}\zeta_{a}(b^{-1}d'c'^{-1}b)n^{-1}q_{b}^{-1},\]
 (hence the definition of $s$), and then a change of summed variable $n \mapsto n\zeta_{a}(b^{-1}d'c'^{-1}b)^{-1} \in C_{G}$. Thus, since the $\rho^{\sharp}_{\mathcal{C},\pi}$ are all injective, then we get that the fundamental braiding on the subcoalgebra $\bigoplus_{(\mathcal{C},\pi) \in \mathcal{S}} \textup{End}(V_{\mathcal{C},\pi})$ of $BD(G)$ is well-defined. Invertibility of $\widetilde{\Psi}$ of $\mathcal{L}_{\CC, \pi}$ follows by restricting the inverse of  $\widetilde{\Psi}$ for $\mathcal{L}_{BD(G)}$. Taking $\mathcal{S} = \{(\mathcal{C},\pi)\}$ gives the fundamental braiding in the statement. The Lie bracket follows from this in the usual way, giving \[ [E_{ai}{}^{bj}, E_{c'k'}{}^{d'l'}] = \delta_{|E_{ai}{}^{bj}||E_{c'k'}{}^{d'j'}|, |b^{-1} \triangleright E_{c'k'}{}^{d'l'} |}\pi(\zeta_{a}(|b^{-1} \triangleright E_{c'k'}{}^{d'l'}|^{-1})^{j}{}_{i} b^{-1} \triangleright E_{c'k'}{}^{d'l'}. \]The above gives the braided Lie algebra structure on $\bigoplus_{(\mathcal{C},\pi) \in \mathcal{S}} \textup{End}(V_{\mathcal{C},\pi})$ for any subcollection $\mathcal{S}$. Again, focusing on just the one block gives the stated Lie algebra. 

Lastly, we note that if we adjoined the pair $(\{ e\}, 1)$ to $\mathcal{S}$, then $\bigoplus_{(\mathcal{C},\pi) \in \mathcal{S}} \textup{End}(V_{\mathcal{C},\pi})$ is  trivially a unital braided Lie algebra by $\eta(1) = 1 \in \textup{End}(V_{(\{ e\}, 1)}) \cong k$. \endproof

The merit of this approach is that, by construction, we immediately have:

\begin{corollary} \label{cor_form_of_inclusion_braided_lie}
There is an inclusion of left regular braided Lie algebras \[ \rho_{\mathcal{C},\pi}^{\sharp}: \mathcal{L}_{\mathcal{C},\pi} \hookrightarrow BD(G), \quad E_{ai}{}^{bj} \mapsto r_{ai}{}^{bj}:= \sum_{n \in C_{G}} \pi(n)^{j}{}_{i} \delta_{q_{a}n^{-1}q_{b}^{-1}} \otimes b^{-1} \] in the category of $G$-crossed modules. Likewise $\rho^\#:\CL\hookrightarrow BD(G)$ for direct sum $\CL$ in (\ref{Lsum}) defined by $\rho_{\CC,\pi}^\#$ on each component. 
\end{corollary}

\begin{remark} \label{rem_braided_matrix_lie}
$\mathcal{L}_{\mathcal{C},\pi}$ turns out to be a `matrix braided Lie algebra' in the sense of \cite[Corollary 4.3]{Ma:sol}. We recall that  for a finite dimensional vector space $W$ of dimension $n$ with basis $\{ E_{i}\}$ and a bi-invertible matrix solution $R \in M_{n}(k) \otimes M_{n}(k)$ of the YBE (i.e. both $R^{-1}$ and the `second inverse' defined by $\widehat{R}^{i}{}_{u}{}^{v}{}_{l}R^{u}{}_{j}{}^{k}{}_{v} = \delta^{i}_{j}\delta^{k}_{l} = R^{i}{}_{u}{}^{v}{}_{l}\widehat{R}^{u}{}_{j}{}^{k}{}_{v}$ exist), one has that
\[  \Delta(E^{i}{}_{j}) = E^{i}{}_{k} \otimes E^{k}{}_{j}, \quad \epsilon(E^{i}{}_{j}) = \delta^{i}{}_{j},\]
 \begin{equation} \label{eq_matrix_braided_lie}  \widetilde{\Psi}(E^{I} \otimes E^{K}) = E^{M} \otimes E^{P} \, \,  \widetilde{R}^{I}{}_{P}{}^{K}{}_{M}, \quad \widetilde{R}^{I}{}_{P}{}^{K}{}_{M} :=  \widehat{R}^{v}{}_{j}{}^{k}{}_{s}  \mathcal{R}^{q}{}_{v}{}^{r}{}_{l}\mathcal{R}^{n}{}_{r}{}^{u}{}_{p}(\mathcal{R}^{-1})^{s}{}_{m}{}^{i}{}_{u}, \end{equation} 
 is a braided Lie algebra on $W^{\sharp} \otimes W$ in the braided category generated by $W$ and $R$ with this fundamental braiding. Here, $E^{I}:= E^{i}{}_{j}$ as a multi-index shorthand used such that $\widetilde{R}$ is a YBE solution clearly valued in $M_{n}(k)^{\otimes 2} \otimes M_{n}(k)^{\otimes 2}$. The braided-Lie bracket is given by applying $\epsilon(E^{p}{}_{q}) = \delta^{p}_{q}$ to the right factor of the output of the fundamental braiding, so that
 \[ [E^{i}{}_{j}, E^{k}{}_{l}]= \mathfrak{d}^{i}{}_{j}{}^{k}{}_{l}{}^{n}{}_{m} E^{m}{}_{n},\quad \mathfrak{d}^{i}{}_{j}{}^{k}{}_{l}{}^{n}{}_{m} := \widehat{R}^{v}{}_{j}{}^{k}{}_{s}  \mathcal{R}^{p}{}_{v}{}^{r}{}_{l}\mathcal{R}^{n}{}_{r}{}^{u}{}_{p}(\mathcal{R}^{-1})^{s}{}_{m}{}^{i}{}_{u}.\]
 
 To explain these formulae in an equivalent form that we prefer, we recall that if  $H$ is a Hopf agebra and $V$ is a finite dimensional right $H$-comodule with basis $\{ E_{i}\}$, and denoting the corresponding algebra map $\rho: H^{*} \to \textup{End}(V)$ with matrix elements $\rho^{i}{}_{j} = \rho(f^{a})^{i}{}_{j}e_{a} \in H$ with respect to a basis $\{ e_{a}\}$ of $H$ with dual $\{ f^{a}\}$, the coaction on $V$ can be written with respect to its basis as $\Delta_{R}(E_{i}) = E_{j}\otimes \rho^{j}{}_{i}$. Moreover, if $H$ is coquasitriangular,   we obtain \cite{Ma} a bi-invertible matrix solution of the YBE  associated to $V$ in the sense 
 \begin{equation} \label{eq_R_matrices} \Psi^{\mathcal{R}}(E_{i} \otimes E_{k}) = E_{l} \otimes E_{j} \, \mathcal{R}^{j}{}_{i}{}^{l}{}_{k}, \quad \mathcal{R}^{j}{}_{i}{}^{l}{}_{k}:= \mathcal{R}(\rho^{j}{}_{i} \otimes \rho^{l}{}_{k}), \end{equation} 
 with  $\widehat{\mathcal{R}}^{i}{}_{j}{}^{k}{}_{l}= \mathcal{R}(\rho^{i}{}_{j} \otimes S\rho^{k}{}_{l})$ and $(\mathcal{R}^{-1})^{i}{}_{j}{}^{k}{}_{l}:= \mathcal{R}(S\rho^{i}{}_{j} \otimes \rho^{k}{}_{l})$. 
Now in our case, we are interested in a braided Lie algebra structure on the $H$-comodule $\textup{End}(V)^{\sharp} = (V^{\sharp})^{\sharp} \otimes V^{\sharp}$, and we take $W = V^{\sharp}$. The coaction on $V^{\sharp}$ is then $\Delta_{R}(E^{i}) = E^{j} \otimes S\rho^{i}{}_{j}$, and by properties of  $\mathcal{R}: H \otimes H \to k$,  its associated $R$-matrix is again \eqref{eq_R_matrices}, but with $i \leftrightarrow j$ and $k \leftrightarrow l$ (in the second equation). Hence, the matrix braided Lie algebra for $\textup{End}(V)^{\sharp}$ realised on its basis $\{ E_{i}{}^{j}\}$ is given by \eqref{eq_matrix_braided_lie}, but with $E^{i}{}_{j} \leftrightarrow E_{i}{}^{j}$, namely
\[ \widetilde{\Psi}(E_{I} \otimes E_{K}):=  \widetilde{\Psi}(E_{i}{}^{j} \otimes E_{k}{}^{l}) = E_{m}{}^{n} \otimes E_{p}{}^{q} \widehat{\mathcal{R}}^{j}{}_{v}{}^{s}{}_{k}  \mathcal{R}^{v}{}_{q}{}^{l}{}_{r}\mathcal{R}^{r}{}_{n}{}^{p}{}_{u}(\mathcal{R}^{-1})^{m}{}_{s}{}^{u}{}_{i} =: E_{M} \otimes E_{P}\widetilde{\mathcal{R}}^{P}{}_{I}{}^{M}{}_{K}.\]
 Applying $\id \otimes \epsilon$ then gives the bracket. 

These are general remarks, but we now check that the choice $H=D^\vee(G)$ indeed lands us on the structure maps in Proposition~\ref{prop_braided_Lie_end} for $\CL_{\CC,\pi}$ as well as more generally for $\CL$ given by direct sums. The former computation for one block has been done in \cite[Eqn.~(33)]{Ma:cla} but with respect to different conventions, requiring fresh calculations. The coquasitriangular structure is given by \eqref{eq_coquasi_DveeG} and we take matrix coefficients $\rho_{\mathcal{C},\pi}: D(G) \to \textup{End}(V_{\mathcal{C},\pi})$ with respect to the natural basis of $D^{\vee}(G)$ given by \eqref{eq_matrix_coeff} and  $\beta$. Thus, 
\begin{align} \label{eq_matrix_coeff_dual} (\rho_{\mathcal{C},\pi})^{ci}_{\hphantom{ci}dj} &= \sum\limits_{g, h \in G} \rho_{\mathcal{C},\pi}(\delta_{g} \otimes h)^{ci}_{\hphantom{ci}dj} \beta(g \otimes \delta_{h}) = \sum_{h \in G} \pi(\zeta_{d}(h))^{i}_{\hphantom{i}j} \delta_{h} \otimes h^{-1}ch\delta_{c, hdh^{-1}} \nonumber \\
& = \sum_{n \in C_{G}} \pi(n)^{i}_{\hphantom{i}j} \delta_{q_{c}nq_{d}^{-1}} \otimes d,\end{align} 
where at the end we set $n:= \zeta_{d}(h)$ and using $n\delta_{c, hdh^{-1}} = q_{hdh^{-1}}^{-1}hq_{d}\delta_{c, hdh^{-1}}  = q_{c}^{-1}hq_{d}\delta_{d, hch^{-1}}$. We must then use \eqref{eq_coquasi_DveeG} in \eqref{eq_R_matrices}, but to cover the case where $\CL$ is a direct sum, we will do the calculation between components $(\CC,\pi)$ and $(\CC',\pi')$. This gives the matrices
\begin{align} \label{eq_R_matrix_start} \mathcal{R}^{a'i'}{}_{b'j'}{}^{ck}{}_{dl} &= \sum_{n, m  \in C_{G}} \pi'(n')^{i'}{}_{j'}\pi(m)^{k}{}_{l}\mathcal{R}( \delta_{q_{a'}n'q_{b'}^{-1}} \otimes b' \otimes  \delta_{q_{c}mq_{d}^{-1}} \otimes d) \nonumber \\
& =  \sum_{n' \in C'_{G}, m \in C_{G}} \pi'(n')^{i'}{}_{j'}\pi(m)^{k}{}_{l}\delta_{q_{a'}n'q_{b'}^{-1}, e}\delta_{b', q_{c}mq_{d}^{-1}} \nonumber\\
&= \delta_{a', b'} \delta_{c, a'da'^{-1}} \delta^{i'}_{j'} \pi(\zeta_{d}(a'))^{k}{}_{l}, \end{align} 
where, for the final equality, we used that $n = q_{a'}^{-1}q_{b'} \in C'_{G}$ iff $a' = b'$, and $m = q_{c}^{-1}b'q_{d} \in C_{G}$ iff $c = b'db'^{-1}$ by Section~\ref{sec_irreps_of_DG}. By similar calculations, we also get \[ ( \mathcal{R}^{-1})^{a'i'}{}_{b'j'}{}^{ck}{}_{dl} = \delta_{a', b'}\delta_{c, a'^{-1}da'} \delta^{i'}_{j'} \pi(\zeta_{d}(a'^{-1}))^{k}{}_{l} = \widetilde{\mathcal{R}}^{a'i'}{}_{b'j'}{}^{ck}{}_{dl}. \] The fundamental braiding is thus indeed given by 
\begin{align*} & \widetilde{\Psi}(E_{ai}{}^{bj} \otimes E_{c'k'}{}^{d'l'}) \\
&= \sum_{g', f', z'_{1}, z'_{2}\in \mathcal{C}', \, x, h, y_{1}, y_{2} \in \mathcal{C}}E_{g'm'}{}^{f'n'} \otimes E_{xp}{}^{hq}\\
&\qquad\qquad  \big(\delta_{b, y_{1}}\delta_{z_{1}', b^{-1}c'b}\delta^{j}_{v}\pi'(\zeta_{c'}(b^{-1}))^{s'}{}_{k'} \big) \big(  \delta_{y_{1}, h}\delta_{d', y_{1}z_{2}'y_{1}^{-1}}\delta^{v}_{q}\pi'(\zeta_{z_{2}'}(y_{1}))^{l'}{}_{r'}\big) \\
 & \qquad\qquad \big(\delta_{z_{2}', f'} \delta_{x, z_{2}'y_{2}z_{2}'^{-1}} \delta^{r'}_{n'} \pi(\zeta_{y_{2}}(z_{2}'))^{p}{}_{u} \big) \big(\delta_{g', z_{1}'}\delta_{y_{2}, g'^{-1}ag'} \delta^{m'}_{s'} \pi(\zeta_{a}(g'^{-1}))^{u}{}_{i} \big)\\
& = E_{b^{-1}c'bm'}{}^{b^{-1}d'bn'} \otimes E_{b^{-1}d'c'^{-1}bab^{-1}c'd'^{-1}b p}{}^{bj} \\
 &\qquad\qquad  \pi'(\zeta_{c'}(b^{-1}))^{m'}{}_{k'}  \pi'(\zeta_{b^{-1}d'b}(b))^{l'}{}_{n'} \pi(\zeta_{b^{-1}c'^{-1}bab^{-1}c'b}(b^{-1}d'b))^{p}{}_{u}  \pi(\zeta_{a}(b^{-1}c'^{-1}b))^{u}{}_{i} \\
 &  =  E_{b^{-1}c'bm'}{}^{b^{-1}d'bn'} \otimes E_{b^{-1}d'c'^{-1}bab^{-1}c'd'^{-1}b p}{}^{bj}  \\
 &\qquad\qquad  \pi'(\zeta_{c'}(b^{-1}))^{m'}{}_{k'} \pi'(\zeta_{d'}(b^{-1})^{-1})^{l'}{}_{n'} \pi(\zeta_{a}(b^{-1}d'c'^{-1}b))^{p}{}_{i}\\
& = b^{-1} \triangleright E_{c'k'}{}^{d'l'} \otimes (|b^{-1} \triangleright E_{c'k'}{}^{d'l'}|^{-1} \triangleright E_{a i}{}\otimes E^{bj}), \end{align*}using \eqref{def_cocycle} to note that $\zeta_{b^{-1}c'^{-1}bab^{-1}c'b}(b^{-1}d'b)\zeta_{a}(b^{-1}c'^{-1}b) = \zeta_{a}(b^{-1}d'c'^{-1}b)$, and that $\zeta_{b^{-1}d'b}(b) = \zeta_{d'}(b^{-1})^{-1}$. Reading off from the above also gives that the matrix elements 
\begin{align}\label{eq_R_widetilde}  \widetilde{\mathcal{R}}^{P}{}_{I}{}^{M'}{}_{K'} & := \widetilde{\mathcal{R}}^{(^{xp}{}_{hq})}{}_{{}{}(^{bj}{}_{ai})}{}^{{}{}(^{g'm'}{}_{f'n'})}{}_{{}{}(^{d'l'}{}_{c'k'})} \nonumber \\
& = \delta_{g', b^{-1}c'b} \delta_{f', b^{-1}d'b} \delta_{x, b^{-1}d'c'^{-1}bab^{-1}c'd'^{-1}b}\delta_{h, b}\,  \delta^{j}{}_{q} \nonumber\\
&\qquad \pi'(\zeta_{c'}(b^{-1}))^{m'}{}_{k'} \pi'(\zeta_{d'}(b^{-1}))^{l'}{}_{n'}\pi(\zeta_{a}(b^{-1}d'c'^{-1}b))^{p}{}_{i}\end{align} 
between the blocks $\textup{End}(V_{\mathcal{C},\pi})$ and $\textup{End}(V_{(\mathcal{C}', \pi')})$. Moreover, applying $( \id \otimes \epsilon)$ to $\widetilde{\Psi}(E_{ai}{}^{bj} \otimes E_{ck}{}^{dl})$ above the structure constants of the braided-Lie bracket as
\begin{align*} \mathfrak{c}_{ai}{}^{bj}{}_{c'k'}{}^{k'l'}{}_{g'n'}{}^{f'm'}&= \delta_{g', b^{-1}c'b} \delta_{f', b^{-1}d'b} \delta_{b, b^{-1}d'c'^{-1}bab^{-1}c'd'^{-1}b}\\
&\qquad  \pi'(\zeta_{c'}(b^{-1}))^{m'}{}_{k'} \pi'(\zeta_{d'}(b^{-1}))^{l'}{}_{n'}\pi(\zeta_{a}(b^{-1}d'c'^{-1}b))^{j}{}_{i},\end{align*}
 thus indeed recovering the braided Lie algebra of Proposition~\ref{prop_braided_Lie_end}. \end{remark}

Continuing the relationship with bicovariant calculi on $D^{\vee}(G)$ discussed around \eqref{eq_fund_equals_crossed}, the fundamental braiding for $\mathcal{L}_{\mathcal{C},\pi}$ similarly viewed functorially in the right $D^\vee(G)$-crossed-module category dualises to the braiding on the space of left-invariant 1-forms $\Lambda^{1} = \textup{End}(V_{\mathcal{C},\pi})$ (generated by \eqref{eq_module_endV_struc} and \eqref{eq_right_comodule_EndVCpi}) for the bicovariant calculus $\Omega^{1}_{\mathcal{C},\pi}$. This 1-1 correspondence between bicovariant calculi on arbitrary coquasitriangular Hopf algebras $H$ and subbraided Lie algebras of $\mathcal{L}_{B(H^{*})}$ is a result from \cite{GomMa} (with different conventions). Further results at this general level, including a generalisation of matrix braided Lie algebras to arbitrary braided abelian categories will appear in \cite{MaMc}.

Moreover, given this relationship with calculi, we indeed recover a sequence of braided Lie algebras corresponding to the bundle of Proposition~\ref{prop_induced_calc_H13}(i) characterising $D^{\vee}(G)$ as an analogue of `extension of groups',

\begin{corollary} \label{cor_sequence_braided_lie2}
There is a sequence of regular braided Lie algebras \[ \delta_{\CC, \{e \}}\mathcal{L}_{\pi} \hookrightarrow \mathcal{L}_{\CC, \pi} \twoheadrightarrow \mathcal{L}_{\CC} \] in the category of $G$-crossed modules, where $\mathcal{L}_{\CC}$ is the braided Lie algebra associated to $\Lambda^{1}_{\CC}$ and $\delta_{\mathcal{C}, \{e\}}\mathcal{L}_{\pi}$ is the restriction of $\mathcal{L}_{\CC, \pi} \subseteq BD(G)$ to $k(G)$. This sequence commutes with the sequence of \eqref{eq_sequence_braided_lie} with vertical inclusion maps. 
\end{corollary}
\proof The categorical dual to the inclusion $(i_{1})_{*}: \Lambda^{1}_{\mathcal{C}} \hookrightarrow \Lambda^{1}_{\mathcal{C},\pi} $, $e_{c} \mapsto E_{ci}{}^{ci}$ from Proposition~\ref{prop_induced_calc_H13} (where $\{ e_{c}\}_{c \in \mathcal{C}}$ was a basis for $\Lambda^{1}_{\CC}$) is can be shown to be given by the $k(G)$-comodule surjection \[ i_{1}{}_{*}^{\sharp}: E_{ai}{}^{bj} \mapsto \delta^{j}_{i} \delta_{a, b}e^{a^{-1}},\] (with $\{e^{c}\}_{c \in \mathcal{C}}$ the dual basis to $\{ e_{c}\}_{c \in \mathcal{C}}$. Applying this to the structure maps of $\mathcal{L}_{\CC, \pi}$ from Proposition~\ref{prop_braided_Lie_end} induces the following structure maps for $\mathcal{L}_{\CC}$ \[ \Delta e^{c^{-1}} = e^{c^{-1}} \otimes e^{c^{-1}}, \quad \epsilon(e^{c^{-1}}) = 1,\] \[ [e^{c^{-1}}, e^{d^{-1}}] = e^{c^{-1}d^{-1}c}, \quad \widetilde{\Psi}(e^{c^{-1}} \otimes e^{d^{-1}}) = e^{c^{-1}d^{-1}c^{-1}} \otimes e^{c^{-1}}, \] by simply asking $i_{1}{}_{*}^{\sharp}$ to be a braided Lie algebra map. For example, this indeed intertwines the brackets 
\begin{align*} 
i_{1}&{}_{*}^{\sharp} ([E_{ai}^{\hphantom{ai}bj}, E_{ck}^{\hphantom{ck}dl}] )\\
& = \pi(\zeta_{c}(b^{-1}))^{n}_{\hphantom{n}k}\pi(\zeta_{b^{-1}db}(b))^{l}_{\hphantom{l}m} \pi(\zeta_{a}(b^{-1}dc^{-1}b))^{j}_{\hphantom{j}i} \delta_{ab^{-1}cd^{-1},b^{-1}cd^{-1}b} \, i_{1}{}_{*}^{\sharp}  (E_{(b^{-1}cb)n}^{\hphantom{(b^{-1}cb)n} (b^{-1}db)m)} )\\
& = \pi(\zeta_{c}(b^{-1}))^{n}_{\hphantom{n}k}\pi(\zeta_{b^{-1}db}(b))^{l}_{\hphantom{l}m} \pi(\zeta_{a}(b^{-1}dc^{-1}b))^{j}_{\hphantom{j}i} \delta_{ab^{-1}cd^{-1},b^{-1}cd^{-1}b}\delta^{n}_{m}\delta_{b^{-1}cb, b^{-1}db} e^{b^{-1}d^{-1}b}\\
&   = \delta_{c, d} \pi(\zeta_{b^{-1}db}(b)\zeta_{d}(b^{-1}))^{l}_{\hphantom{k}k} \pi(\zeta_{a}(e))^{j}_{\hphantom{j}i} \delta_{a, b} e^{b^{-1}d^{-1}b} = \delta^{j}_{i}\delta_{c, d}\delta_{a, b} \delta^{l}_{k} e^{b^{-1}d^{-1}b} \\
& = \delta^{j}_{i} \delta^{l}_{k}[\delta_{a, b} e^{b^{-1}}, \delta_{c, d} e^{d^{-1}}] = [(i_{*})^{\sharp}(E_{ai}^{\hphantom{ai}bj}), (i_{*})^{\sharp}(E_{ck}^{\hphantom{ck}dl})] .
\end{align*} 
Next, the inclusion $k(G) \hookrightarrow BD(G)$ restricts to an isomorphism $\textup{End}(V_{\pi}) \cong \textup{End}(V_{\{ e\}, \pi})$ of braided Lie algebras, for an irrep $\pi$ of $G$. So, $\mathcal{L}_{\CC, \pi}$ viewed in $BD(G)$ restricts to $k(G)$ iff $(\mathcal{C}, \pi) = (\{ e\}, \pi)$, and gives the braided Lie algebra $\delta_{\mathcal{C}, \{e \}}\mathcal{L}_{\pi} := \delta_{\mathcal{C}, \{e \}}\mathcal{L}_{\{e \}, \pi}$. In a more formal language, this is the pullback of $k(G) \to BD(G) \leftarrow \mathcal{L}_{\CC, \pi}$ in the category of left braided Lie algebras in the $G$-crossed module category. We can verify that the exact sequence stated commutes with \eqref{eq_sequence_braided_lie}.
\endproof

The formulae for the our braided-Lie algebras simplify in two cases. (a)   $\CC=\{e\}$ with any  irrep $\pi$ of $G=C_G$, where we just have $\{E_i{}^j\}$ as the basis with matrix coalgebra and
\[ [E_i{}^j,E_k{}^l]=\delta_{i}^{j} E_k{}^l,\quad \tilde \Psi(E_i{}^j\tens E_k{}^l)=E_k{}^l\tens E_i{}^j,\]
as the trivial braided-algebra of the relevant dimension. (b)  A conjugacy class in $\CC\subset G\setminus\{e\}$ and $\pi$ a  1-dimensional irrep of $C_G$, where we just have $\{E_c{}^d\}$ as the basis with matrix coalgebra and
\[ [E_a{}^b,E_c{}^d]=\delta_{bab^{-1}cd^{-1}b^{-1},cd^{-1}}\pi(a,b,c,d)\ E_{b^{-1}cb}{}^{b^{-1}db},\]
\begin{equation}\label{Psi1d}\tilde\Psi(E_a{}^b\tens E_c{}^d)=\pi(a,b,c,d)\ E_{b^{-1}cb}{}^{b^{-1}db}\tens E_{b^{-1}dc^{-1}bab^{-1}cd^{-1}b}{}^b,\end{equation}
\[ \pi(a,b,c,d)={\pi(\zeta_c(b^{-1}))\pi(\zeta_a(b^{-1}d c^{-1}b))\over \pi(\zeta_d(b^{-1}))}.\]
This simplifies considerably when $c=d$ to 
\begin{equation}\label{c=d}  [E_a{}^b,E_c{}^c]=\delta_{a,b} E_{b^{-1}cb}{}^{b^{-1}cb},\quad \tilde\Psi(E_a{}^b\tens E_c{}^c)=E_{b^{-1}cb}{}^{b^{-1}cb}\tens E_a{}^b.\end{equation}
The bracket and fundamental braiding on the diagonal restrict to those of the standard braided Lie algebra $k\CC$ but note that the coproduct has off-diagonal terms.

\begin{example} \label{example_braided_Lie}\rm For $G=S_3$, we look at the 3 conjugacy classes. 

(i) For $\CC=\{e\}$ we have 3 trivial braided-Lie algebras of the 1st kind above. For the other two classes we have $\pi$ 1-dimensional, hence of the 2nd kind above. 

(ii) For $\CC=\{uv,vu\}$ and $\pi_j$, $j=0,1,2$, we have the elements $c\in \CC$ commute and using $\zeta$ given before, one can compute that
\[ [E_a{}^b,E_c{}^d]=\delta_{a,b}\pi(a,a,c,d)\ E_{c}{}^{d},\quad \tilde\Psi(E_a{}^b\tens E_c{}^d)=\pi(a,b,c,d)\ E_{c}{}^{d}\tens E_a{}^b,\]
\[ \pi(a,b,c,d)=q^{j(\delta_{a,c}+\delta_{b,c}-\delta_{a,d}-\delta_{b,d})},\]
using $\zeta_a(b)=r^{(2\delta_{a,b})-1}$ from calculations in Example~\ref{exS3}.

(iii) For $\CC=\{u,v,w\}$ and $\pi_\pm$, we have the $c=d$ case (\ref{c=d}) and for $c\ne d$:
\[ [E_a{}^a, E_c{}^d]=0,\quad \tilde\Psi(E_a{}^a\tens E_c{}^d)=\begin{cases}\pi(a,a,c,d)\ E_{d}{}^{c}\tens E_{c}{}^a  & {\rm if\ }a,c,d\ {\rm distinct}\\
 \pi(a,a,a,d)\ E_{a}{}^{\eps(a,d)}\tens E_{d}{}^a  & {\rm if\ }a=c\ne d\\
 \pi(a,a,c,a)\ E_{\eps(a,c)}{}^a\tens E_{\eps(a,c)}{}^a  & {\rm if\ }a=d\ne c\end{cases}, \]
where $\eps(a,c)$ for distinct $a,c\in\CC$  means the other element of $\CC$. This leaves us the cases $a\ne b$, $c\ne d$ with
\[ [E_a{}^b,E_c{}^d]=\begin{cases} \pi(a,b,b,d)\ E_b{}^a&{\rm if\ }a, b=c, d\ {\rm distinct}\\
\pi(a,b,a,b)\ E_{\eps(a,b)}{}^b &{\rm if\ }a=c \ne b=d\\
\pi(a,b,c,a)\ E_a{}^c &{\rm if\ }a=d, b, c\ {\rm distinct}\\ 0 &{\rm otherwise}.\end{cases}
\]
For the fundamental braiding in the cases $a\ne b$, $c\ne d$,
\[ \tilde\Psi(E_a{}^b\tens E_c{}^d)=\begin{cases} \pi(a,b,b,d)\ E_b{}^a\tens E_b{}^b&{\rm if\ }a, b=c, d\ {\rm distinct}\\
\pi(a,b,c,b)\ E_a{}^b\tens E_c{}^b&{\rm if\ }a, b=d, c\ {\rm distinct}\\
\pi(a,b,a,d)\ E_d{}^a\tens E_d{}^b&{\rm if\ }a=c, b, d\ {\rm distinct}\\
\pi(a,b,a,b)\ E_{\eps(a,b)}{}^b\tens E_b{}^b &{\rm if\ }a=c \ne b=d\\
\pi(a,b,c,a)\ E_a{}^c \tens E_b{}^b&{\rm if\ }a=d, b, c\ {\rm distinct}\\
\pi(a,b,b,d)\ E_b{}^{\eps(a,b)} \tens E_{\eps(a,b)}{}^b&{\rm if\ }a=d\ne b=c\end{cases}, \]
where the nine relevant $\pi$ (for distinct values of the labels) are only needed in the case $\pi=\pi_-$, in which case
\[ \pi(a,a,c,d)=-\pi_{ca}\pi_{ad}\pi_{da}     ,\quad \pi(a,a,a,d)= -\pi_{\eps(a,b)a}\pi_{ad}\pi_{da}  ,\quad  \pi(a,a,c,a)=-\pi_{ac}\pi_{ca},    \]
\[ \pi(a,b,b,d)= \pi_{ad}\pi_{db}    ,\quad \pi(a,b,c,b)=-\pi_{cb}\pi_{ab}   ,\quad  \pi(a,b,a,d) = \pi_{ad}\pi_{ab}\pi_{ba},   \]
\[\pi(a,b,a,b)= -\pi_{\eps(a,b)a}    ,\quad \pi(a,b,c,a)=-\pi_{cb}\pi_{ac}\pi_{ab},\quad  \pi(a,b,b,d)= \pi_{ad}\pi_{db},   \]
where $\pi_{ab}:=\pi_-(\zeta_a(b))$. These can be read off from
\[  \pi_{\pm}(\zeta_{a}(b)) = (\pm 1)^{\delta_{a, b}+ (1 - \delta_{a, b})\delta_{b, v} }. \]
 \end{example}

\subsection{The universal enveloping algebra $U(\mathcal{L}_{\mathcal{C},\pi})$} 

For every braided-Lie algebra $\CL$ in an abelian braided category, one has a universal enveloping algebra $U(\CL)$ defined as the tensor algebra $T(\mathcal{L}):= \bigoplus_{n \geq 0} \mathcal{L}^{\otimes n}$ in $\CL$ modulo braided-commutativity with respect to the fundamental braiding, i.e. the relations
\begin{equation} \label{eq_defining_relation_UL} \cdot \tilde\Psi_{\CL,\CL}=\cdot \ .\end{equation}
The coproduct of $\CL$ extends to make $T(\mathcal{L})$ and $U(\CL)$ bialgebras in the braided category \cite{Ma:lie} and $\CL\subset U(\CL)$ since the relations are in degree 2. Note that given a classical Lie algebra $\cg$, it is $\CL=k c\oplus\cg$ which is a braided-Lie algebra, with result that as algebras $U(\CL)\twoheadrightarrow U(\cg)$ by setting the generator $c$ of $k$ here to 1. Indeed, in this case $U(\CL)$ is a quadratic homogenisation of $U(\cg)$ and is an ordinary bialgebra.  On the other hand $U(\mathcal{L})$ never has a (braided) antipode \cite{GomMa}. For matrix braided-lie algebras \cite{Ma:lie} one has that $U(\CL)=B(R)$ the braided matrix bialgebra (related by transmutation to the FRT bialgebra $A(R)$). Hence, by  Remark~\ref{rem_braided_matrix_lie} we have:

\begin{corollary} \label{cor_UL_is_braided_matrices}
$U(\mathcal{L}_{\mathcal{C},\pi}) \cong  B(\mathcal{R})$ the braided-matrix bialgebra for $\mathcal{R}$ given by \eqref{eq_R_matrix_start} applied in the case of one block, and likewise for direct sum $\CL$ as in (\ref{Lsum}). 
\end{corollary}

\proof The case of one block is from \cite[Proposition 5.2]{Ma:cla}. More generally for the braided Lie algebra $\mathcal{L}:= \bigoplus_{(\mathcal{C},\pi) \in \mathcal{S}} \mathcal{L}_{\mathcal{C},\pi} = \textup{End}(V_{(\mathcal{C}_{1}, \pi_{1})}) \oplus ... \oplus \textup{End}(V_{(\mathcal{C}_{s}, \pi_{s})})$, we get that $U(\mathcal{L})$ is given by the braided tensor algebra 
\[\bigoplus_{n \geq 0} \mathcal{L}^{\otimes n} =  \bigoplus_{n \geq 0} \big(\bigoplus_{(\mathcal{C},\pi) \in \mathcal{S}} \textup{End}(V_{\mathcal{C},\pi})\big)^{\otimes n} = \bigoplus_{n_{1}, ..., n_{s}}  \textup{End}(V_{(\mathcal{C}_{1}, \pi_{1})})^{\otimes n_{1}} \otimes ... \otimes \textup{End}(V_{(\mathcal{C}_{s}, \pi_{s})})^{\otimes n_{s}},\]
 modulo the  relations $E_{I}E_{K'} = E_{M'}E_{P}\widetilde{\mathcal{R}}^{P}{}_{I}{}^{M'}{}_{K'}$ for $E_{I} \in \textup{End}(V_{\mathcal{C}_{i}, \pi_{i}})$ and $E_{K'} \in \textup{End}(V_{(\mathcal{C}_{k}, \pi_{k})})$ with $\widetilde{\mathcal{R}}^{P}{}_{I}{}^{M'}{}_{K'}$ by \eqref{eq_R_widetilde}. The braided coalgebra structure is  just the braided tensor product coalgebra structure with the matrix form on each block. \endproof

 Also of interest is the corresponding FRT bialgebra $A(\CR)$ (an ordinary $k$-bialgebra), related to $U(\CL)$ by transmutation. This is  generated 1 and $\sum_{(\mathcal{C},\pi )\in \mathcal{S}} |\mathcal{C}|^{2}\textup{dim}(V_{\pi})^{2}$ indeterminates $\{t_{ai}{}^{bj}  \}$ for each pair $\mathcal{C},\pi$ with coalgebra structure given by the direct sum of the matrix coalgebra structure on each block, i.e. $\Delta t_{ai}{}^{bj} = \sum_{c \in \mathcal{C}}t_{ai}{}^{ck} \otimes t_{ck}{}^{bj}$, $\epsilon(t_{ai}{}^{bj}) = \delta_{a, b} \delta_{i}^{j}$, and algebra structure \begin{equation} \label{eq_relations_AR} \sum_{ f \in \mathcal{C}, g' \in \mathcal{C}' }\mathcal{R}^{fp}{}_{ai}{}^{g'q'}{}_{c'k'}t_{fp}{}^{bj}t_{g'q'}{}^{d'l'} = \sum_{ f \in \mathcal{C}, g' \in \mathcal{C}' } t_{c'k'}{}^{g'q'}t_{ai}{}^{fp}\mathcal{R}^{bj}{}_{fp}{}^{d'l'}{}_{g'q'}, \end{equation} with $ \mathcal{R}^{ai}{}_{bj}{}^{c'k'}{}_{d'l'}$ from \eqref{eq_R_matrix_start} applied between $(\CC',\pi')$ and $(\CC,\pi)$. Moreover,  $A(\mathcal{R})$ has coquasitriangular structure $\mathcal{R}^{A}: A(\mathcal{R}) \otimes A(\mathcal{R}) \to k$ defined on degree 1 elements by $\mathcal{R}^{A}(t_{ai}{}^{bj} \otimes t_{c'k'}{}^{d'j'}) := \mathcal{R}^{bj}{}_{ai}{}^{d'k'}{}_{c'k'}$ and extended to products as a bialgebra bicharacter. 
 
\begin{corollary} For $\CL_{\CC,\pi}$, we denote by $A_{\CC,\pi}$ the  standard the FRT bialgebra given by the $R$-matrix from \eqref{eq_R_matrix_start} for one block, with coquasitriangular structure from \cite{Ma}. Its relations  simplify to
\[ \sum_{m} \pi(\zeta_{c}(a))^{m}{}_{k} \, t_{ai}{}^{bj}t_{a^{-1}cam}{}^{dl} = \sum_{m} \,  \pi(\zeta_{d}(b))^{l}{}_{m} t_{ck}{}^{bdb^{-1}m}t_{ai}{}^{bj}. \] 
where $i,j,k,l,m\in \{1,\cdots,\dim(V_\pi)\}$ and $a,b,c,d\in \CC$ as per our standing conventions. 
\end{corollary}
\proof This is $A(\CR)$ for one block, with relations given by substituting \eqref{eq_R_matrix_start} for this case into \eqref{eq_relations_AR}. \endproof

In the extreme case where $\CC=\{e\}$, these simplify to the generators being commutative, so $A_{\{e\},\pi}=k[t_i{}^j]$ with matrix coalgebra. In the other extreme case where $\CC$ is general and $\pi$ is 1-dimensional, we have 
 \begin{equation}\label{A1dpi}  \pi(\zeta_{c}(a)) \, t_{a}{}^{b} t_{a^{-1}ca}{}^{d} = \pi(\zeta_{d}(b))\, t_{c}{}^{bdb^{-1}} t_{a}{}^{b},\end{equation}
 for all $a,b,c,d\in \CC$.

 \begin{example} \label{example_UL} \rm
 For the case of $G = S_3$ and  the various pairs $(\CC, \pi)$, the braided bialgebras  $U(\mathcal{L}_{\CC, \pi})$ for are as follows. We use $\tilde\Psi$ already computed in Example~\ref{example_braided_Lie} and the matrix form of coalgebra on the generators. 
 
 (i) For $\CC=\{e\}$ we have $\tilde\Psi$ the flip map, hence $U(\CL)=\C[E_i{}^j]$ with its matrix coproduct. The braiding as a braided-bialgebra is trivial.
  
 (ii) For $\mathcal{C} = \{ uv, vu\}$ and $\pi_{j}$, $j = 0, 1, 2$, the relations are
 \[ E_a{}^b E_c{}^d=q^{j(\delta_{a,c}+\delta_{b,c}-\delta_{a,d}-\delta_{b,d})}E_c{}^d E_a{}^b,\]
 but these $q$-commutation relations amount to $\C[E_a{}^b]$ for all $j$ modulo the further relations
 \[ E_1{}^1 E_1{}^2=E_1{}^1 E_2{}^1= E_2{}^2 E_1{}^2=E_2{}^2 E_2{}^1=0,\]
 for $j=1,2$. We used basis order $uv,vu$ in labelling the generators. 
 
 (iii) For $\CC=\{u,v,w\}$ and $\pi_\pm$, we let $e_a:=E_a{}^a$ and from the $c=d$ case of the $\tilde\Psi$ we get
\begin{equation}\label{S3iiirel}  e_a e_b=e_{\eps(a,b)}e_a,\quad E_a{}^b e_{\eps(a,b)}=e_a E_a{}^b,\quad E_a{}^be_a=e_{\eps(a,b)}E_a{}^b,\quad E_a{}^b e_b=e_b E_a{}^b,\end{equation}
 for distinct $a,b$. We see that the $\{e_a\}$ obey  the usual relations of $U(\CL_\CC)$ and in addition we see the exchange relations between the diagonal and the off-diagonal generators.  
 
 Further relations coming from (\ref{Psi1d})  involving the diagonals are 
\[ e_a E_b{}^c=\pi(a,a,b,c) E_c{}^b E_b{}^a,\  e_a E_a{}^b=\pi(a,a,a,b) E_a{}^{\eps(a,b)} E_b{}^a,\   e_a E_b{}^a=\pi(a,a,b,a)(E_{\eps(a,b)}{}^a)^2,\]
for $a,b,c$ distinct. This leaves the relations for $E_a{}^bE_c{}^d$ for $a\ne b$, and $c\ne d$ with various cases from the listed $\tilde\Psi$ giving 6 relations
\[ E_a{}^bE_b{}^c=\pi(a,b,b,c)E_b{}^ae_b,\  E_a{}^b E_c{}^b(1-\pi(a,b,c,b))=0,\  E_a{}^bE_a{}^c=\pi(a,b,a,c)E_c{}^a E_c{}^b,\]
\[ (E_a{}^b)^2=\pi(a,b,a,b)E_{\eps(a,b)}{}^b e_b,\  E_a{}^bE_c{}^a=\pi(a,b,c,a)E_a{}^ce_b,\ E_a{}^bE_b{}^a=\pi(a,b,b,a)E_b{}^{\eps(a,b)} E_{\eps(a,b)}{}^a, \]
 for distinct $a,b,c$. The second relation is empty for $\pi_+$ but for $\pi_-$ we obtain $E_a{}^b E_c{}^b=0$ in all cases of distinct $a,b,c$.  
 
 For more details we focus on the $\pi_+$ case where we drop the $\pi$. Then our full set of relations after the above observations become (\ref{S3iiirel}) and
 \[  E_a{}^bE_a{}^{\eps(a,b)}=E_{\eps(a,b)}{}^a E_{\eps(a,b)}{}^b,\quad E_a{}^bE_b{}^a=E_b{}^{\eps(a,b)} E_{\eps(a,b)}{}^b,\]
 \[  E_{\eps(a,b)}{}^b E_b{}^a=e_a E_b{}^{\eps(a,b)},\quad  E_a{}^{\eps(a,b)} E_b{}^a=e_a E_a{}^b,\quad   (E_{\eps(a,b)}{}^a)^2=e_a E_b{}^a,\]
 for distinct $a,b$. Writing out the relations explicitly, we find that $U(\CL_{\CC,\pi_+})$ is a  reasonable quadratic algebra with dimension 33 in degree 2. 
 \end{example}
 
 \begin{example} \label{example_rel_AR} \rm
 For the case of $G = S_3$ and  the various pairs $(\CC, \pi)$,  the ordinary bialgebras  $A_{\CC, \pi}$ for the various pairs $(\CC, \pi)$ are as follows. We use the two extreme cases, notably (\ref{A1dpi}) with the matrix coalgebra. 
 
 (i) For $\mathcal{C} = \{ e\}$, we have $A_{\{e\},
\pi}=\C[t_i{}^j]$ as in general,   the same  as $U(\CL_{\CC,\pi_j})$.
  
 (ii) For $\mathcal{C} = \{ uv, vu\}$ and $\pi_{j}$, $j = 0, 1, 2$ we use  $\zeta_{a}(b) = r^{(2\delta_{a, b}) - 1}$ from Example~\ref{example_braided_Lie}(ii) to find the relations as
 \[ t_{a}{}^{b}t_{c}{}^{d} = q^{j(\delta_{a, c} - \delta_{b, d})}t_{c}{}^{d}t_{a}{}^{b},  \]
but these $q$-commutations amount to $\C[t_a{}^b]$ for $j=0$ and the further relations
\[ t_1{}^1 t_1{}^2=t_1{}^1 t_2{}^1= t_2{}^2 t_1{}^2=t_2{}^2 t_2{}^1=0,\]
 for $j=1,2$,   the same  as $U(\CL_{\CC,\pi_j})$.

(iii) For $\mathcal{C} = \{ u, v, w \}$ and $\pi_{\pm}$, with $\pi_\pm(\zeta_{a}(b))$ from Example~\ref{example_braided_Lie}(iii), we have relations
\[ t_{a}{}^{b}t_{c}{}^{d} = (\pm 1)^{\delta_{a, c}(1 - \delta_{a, v}) + \delta_{b, d}(1- \delta_{b, v})+ \delta_{a, v} + \delta_{b, v} }t_{\eps(a,c)}{}^{\eps(b,d)}t_{a}{}^{b},  \] provided we understand $\eps(a,a)=a$ when its arguments coincide. This expands out more explicitly to
\[  t_{a}{}^{b}t_{c}{}^{d} = \begin{cases} (\pm 1)^{1 + \delta_{b, v}}t_{c}{}^{\eps(b,d)}t_{a}{}^{b} & {\rm if\ }a = c {\rm \ and \ } b \neq d\\
(\pm 1)^{1 + \delta_{a, v}}t_{\eps(a,c)}{}^{d}t_{a}{}^{b} & {\rm if\ }a \neq c {\rm \ and \ } b=d\\ 
(\pm 1)^{\delta_{a, v} + \delta_{b, v}}t_{\eps(a,c)}{}^{\eps(b,d)}t_{a}{}^{b} & {\rm if\ } a \neq c {\rm \ and \ }b \neq d.\end{cases}\] 
The first case amounts to
\[t_{a}{}^{u}t_{a}{}^{v} = t_{a}{}^{v}t_{a}{}^{w} = \pm t_{a}{}^{w}t_{a}{}^{u}, \quad t_{a}{}^{u}t_{a}{}^{w} = \pm t_{a}{}^{v}t_{a}{}^{u} = \pm t_{a}{}^{w}t_{a}{}^{v},  \]
for any $a \in \mathcal{C}$ and the second case the same swapping upper and lower indices. The third case has a minus in the $\pi_-$ case precisely when either $a=v$ or $b=v$ but not both. 

For more details, we focus on the $\pi_+$ case without sign factors and we write $t_a=t_a{}^a$. Then the relations can be written as
\[ t_a t_b=t_{\eps(a,b)}t_a,\quad  t_a t_b{}^{\eps(a,b)}=t_{\eps(a,b)}{}^{b}t_a,\]
for distinct $a,b$ and 6 cases for products of off-diagonals. One of these  is empty and the others are
\[ t_a{}^b t_b{}^{\eps(a,b)}=t_{\eps(a,b)}{}^a t_a{}^b,\quad  t_a{}^bt_{\eps(a,b)}{}^a=t_b{}^{\eps(a,b)} t_a{}^b,\]
\[t_a{}^b t_{\eps(a,b)}{}^b=t_b t_a{}^b,\quad t_a{}^b t_a{}^{\eps(a,b)}=t_a t_a{}^b, \quad t_a{}^b t_b{}^a=t_{\eps(a,b)}t_a{}^b,\]
for distinct $a,b$. Writing  out the relations explicitly, we find that $A_{\CC,\pi_+}$ is a  reasonable quadratic algebra with dimension 33 in degree 2 as for $U(\CL_{\CC,\pi_+})$, as required since the two are related by transmutation, and also has a similar general structure. 
\end{example}

\subsection{Surjectivity of the map $U(\CL)\to BD(G)$}

For the case of a rack given by a conjugacy class, one has a map $U(\CL_\CC)\to k G$ which is surjective when $\CC$ generates $G$ (the case of a connected calculus, which in turn corresponds to a connected Cayley graph on $G$). In this section we want to show analogous results for $U(\CL_{\CC,\pi})$. As before, we also allow $\CL$ which are direct sums of such blocks. 

\begin{lemma} \label{lem_map_from_enveloping_to_og}
There is a map of braided bialgebras $U(\mathcal{L}) \to BD(G)$ given on the generators by $E_{ai}{}^{bj} \mapsto r_{ai}{}^{bj}$ (with $r_{ai}{}^{bj}$ defined in Corollary~\ref{cor_form_of_inclusion_braided_lie}), which commutes with the inclusion of $\mathcal{L}$. Moreover, this map is the transmutation of the map of coquasitriangular bialgebras $A(\mathcal{R}) \to D^{\vee}(G)$ given on the generators by $t_{ai}{}^{bj} \mapsto r_{ai}{}^{bj}$.
\end{lemma}

\proof We start with the simpler map $\Upsilon: A(\mathcal{R}) \to D^{\vee}(G)$ and again work more generally with $\mathcal{L}:= \bigoplus_{(\mathcal{C},\pi) \in \mathcal{S}} \mathcal{L}_{\mathcal{C},\pi}$, for which the stated map generalises in the obvious way to $t_{ai}{}^{bj}...t_{c'k'}{}^{d'l'} \to r_{ai}{}^{bj}...r_{c'k'}{}^{d'l'}$ for indeterminates from any pairs $(\mathcal{C},\pi), (\mathcal{C}', \pi') \in \mathcal{S}$ (and sends the unit of $A(\mathcal{R})$ to the unit of $D^{\vee}(G)$). 

Since we saw in Corollary~\ref{cor_form_of_inclusion_braided_lie} that $\bigoplus_{(\mathcal{C},\pi) \in \mathcal{C}}\rho^{\sharp}_{\mathcal{C},\pi}$ was a morphism of coalgebras and the coalgebra of  $\mathcal{L}$ coincides with that of $A(\mathcal{R})$,  it is immediate that the map $t_{ai}{}^{bj} \mapsto r_{ai}{}^{bj}$ is a morphism of $k$-coalgebras. Extending this coalgebra map as a unital algebra map then gives the stated map. So, we need only check that the map is indeed compatible with the defining relations \eqref{eq_relations_AR} of $A(\mathcal{R})$, i.e. that we have \[ \Upsilon\Big(\sum_{ f \in \mathcal{C}, g' \in \mathcal{C}' }\mathcal{R}^{fp}{}_{ai}{}^{g'q'}{}_{c'k'}t_{fp}{}^{bj}t_{g'q'}{}^{d'l'} \Big)= \Upsilon\Big( \sum_{ f \in \mathcal{C}, g' \in \mathcal{C}' } t_{c'k'}{}^{g'q'}t_{ai}{}^{fp}\mathcal{R}^{bj}{}_{fp}{}^{d'l'}{}_{g'q'}\Big). \] This is shown firstly by noting that using \eqref{DvG}, then
\begin{align} \label{eq_relation_r_rho} S(r_{ai}{}^{bj}) &= \sum_{n \in C_{G}} \pi(n)^{j}{}_{i} S(\delta_{q_{a}n^{-1}q_{b}^{-1}} \otimes b^{-1}) =\sum_{n \in C_{G}} \pi(n)^{j}{}_{i} \delta_{q_{b}nq_{a}^{-1}} \otimes (q_{a}n^{-1}q_{b}^{-1})b(q_{b}^{-1}nq_{a}^{-1}) \nonumber\\
& = \sum_{n \in C_{G}} \pi(n)^{j}{}_{i} \delta_{q_{b}nq_{a}^{-1}} \otimes a= (\rho_{\mathcal{C},\pi})^{bj}{}_{ai},\end{align} 
 i.e. the matrix coefficients from \eqref{eq_matrix_coeff_dual}. On $r_{ai}{}^{bj}, r_{c'k'}{}^{d'l'} \in D^{\vee}(G)$, we know by the quasicommutivity relation given by the coquasitriangular structure of $D^{\vee}(G)$ that the following holds
 \[  \sum_{f \in \mathcal{C}, g' \in \mathcal{C}'}\mathcal{R}(r_{ai}{}^{fp} \otimes r_{c'k'}{}^{g'q'}) r_{fp}{}^{bj}r_{g'q'}{}^{d'l'}  = \sum_{f \in \mathcal{C}, g' \in \mathcal{C}'}r_{c'k'}{}^{g'q'}r_{ai}{}^{fp} \mathcal{R}(r_{fp}{}^{bj} \otimes r_{g'q'}{}^{d'l'}), \] 
and so using $\mathcal{R} = \mathcal{R} \circ (S \otimes S)$ for any coquasitriangular Hopf algebra we can use  \eqref{eq_relation_r_rho} and \eqref{eq_matrix_coeff_dual} to obtain \[ \mathcal{R}(r_{ai}{}^{fp} \otimes r_{c'k'}{}^{g'q'}) = \mathcal{R}((\rho_{\mathcal{C},\pi})^{fp}{}_{ai} \otimes ((\rho_{\mathcal{C},\pi})^{g'q'}{}_{c'k'}) = \mathcal{R}^{fp}{}_{ai}{}^{g'q'}{}_{c'k'}.\] 
Doing the same to the right-hand side of the above equation then indeed gives compatibility of $\Upsilon$ with the required relation. By $A(\mathcal{R})$'s previously given coquasitriangular structure, it is then clear that the map $\Upsilon$ is a map of coquasitriangular bialgebras by the above.

Next, since $B(\mathcal{R}) = U(\mathcal{L})$ is the comodule transmutation of $A(\mathcal{R})$ by Corollary~\ref{cor_UL_is_braided_matrices}, it is immediate from transmutation theory that the map $U(\mathcal{\mathcal{L}}) \to BD(G)$ stated is a well-defined map of braided bialgebras (and can be verified explicitly in an analogous way to the above).  \endproof

Our first characterisation of when this map is surjective is motivated from representation theory.

\begin{lemma} \label{corollary_surjective_UL_inner_faith} For a direct sum $\CL$ as above, we let $\rho:D(G)\to \End(V)$ be given by (\ref{eq_matrix_coeff}) is applied to each block of $V=\oplus_{(\CC,\pi)\in S}$. Then the map $U(\mathcal{L}) \to BD(G)$ of Lemma~\ref{lem_map_from_enveloping_to_og} is surjective iff $\rho$ is `inner faithful' in the sense of \cite{paper_inner_faith}, i.e. does not factor through any proper bialgebra quotients of $D(G)$. 
\end{lemma}
\proof Since $D(G) = D^{\vee}(G)$ as vector spaces, we note that subspaces can be viewed as subspaces of either.  By definition, the map $U(\mathcal{L}) \to BD(G)$ is clearly surjective iff $1$ and $\CL$ generate the underlying unital algebra of $BD(G)$, which is that of $D(G)$ by Lemma~\ref{lem_braided_Hopf_DG} (and where here $\mathcal{L}$ is a viewed as a subspace via $\rho^{\sharp}$ in Corollary~\ref{cor_form_of_inclusion_braided_lie}). This can be  phrased more conveniently as the statement that if $A\subseteq D(G)$ is a unital algebra which contains the subvector space $\CL$ then $A = D(G)$. Dualising this statement gives equivalently that  $T= \textup{ker}(\rho)$ is a {\em cogenerating subvector space} of the underlying coalgebra of $D^{\vee}(G)$, in the sense that if $J\subseteq T$ is a coideal  of $D^{\vee}(G)$ then $J = \{ 0\}$. Recall that a coideal $J$ of a counital coalgebra, here $D^{\vee}(G)$, means $\Delta(J)\subseteq J\tens D^{\vee}(G)+ D^{\vee}(G)\tens J$ and $\eps(J)=0$. We also used in the dualizing process that $\textup{ker}(D^{\vee}(G)\twoheadrightarrow L^{*}) = \textup{ker}(\rho)$ since the inclusion of $L \hookrightarrow D(G)$ was by $\rho^{\sharp}$ and the adjoint $D^{\vee}(G)\twoheadrightarrow L^{*}$ of this with respect to our (slightly non-standard) duality pairing \eqref{eq_coev_ev_DG} between $D(G), D^\vee(G)$ returns $\rho$.

We next recall that $I$ being a proper biideal of $D(G)$ here  means is it both a two-sided ideal of $D(G)$ and a coideal of $D(G)$. We also note the bijection 
\begin{equation} \label{eq_correspondence_biideals}
\{ \textup{biideals of $D(G)$}\} \leftrightarrow \{\textup{ biideals of $BD(G)$}\}, \end{equation} where $I$ being a biideal of the braided Hopf algebra $BD(G)$ additionally requires $I$ to be a sub $G$-crossed module of $BD(G)$ (such that the quotient is again in the category). This is a general fact for any transmutation of a quasitriangular Hopf algebra to a braided one in the module category. Since this has the same algebra, an ideal in one is the same as an ideal on the other. If it is also a biideal, since the transmutation changes the coproduct by multiplication by certain elements from either side, it also remains a biideal. Conversely, the transmutation is reversible by similar formulae. 

Given these facts, we now show that $\textup{ker}(\rho)$ being cogenerating for $D^{\vee}(G)$ is equivalent to inner faithfulness of $\rho$. In one direction, assume that $\textup{ker}(\rho)$ is cogenerating for $D^{\vee}(G)$, and suppose for contradiction that there exists a non-zero biideal $I \subseteq \textup{ker}(\rho)$ of $D(G)$. By \eqref{eq_correspondence_biideals}, $I$ is a biideal of $BD(G)$ hence a coideal of $D^{\vee}(G)$ (since $BD(G)$ has the coalgebra of the latter). $I$ is thus zero by assumption of $\textup{ker}(\rho)$ being cogenerating. Conversely, suppose $J \subseteq \textup{ker}(\rho)$ is a coideal of $D^{\vee}(G)$. Consider the $D(G)$-module $J'$ generated by $J$ under the adjoint action (with which $BD(G)$ sits in the $G$-crossed module category, as given in Corollary~\ref{cor_braided_lie_DG}). Since the coproduct of $D^{\vee}(G)$ is that of $BD(G)$ and the latter is a $D(G)$-module map under the adjoint action, it follows that $J'$ is again a coideal of $D^{\vee}(G)$ and hence of $BD(G)$.  This implies a biideal $J''=BD(G)\cdot J'\cdot BD(G)$ of $BD(G)$ (using the braided bialgebra axioms) and we still have $J''\subseteq \textup{ker}(\rho)$ since $\rho$ is an algebra map for the algebra of $BD(G)$. Moreover,  $J''$ is a biideal of $D(G)$ by the correspondence and hence must be zero by inner faithfulness of $\rho$.  But $J \subseteq J''$ so that $J = \{ 0\}$ as required.\endproof

We then recover a result expected by analogy with classical Lie theory. 

\begin{corollary}
If $U(\mathcal{L}) \to BD(G)$ is surjective, then the corresponding calculus $\Omega^{1}$ as discussed after  Lemma~\ref{lem_coirreducible_calc_DG} is connected.
\end{corollary}
\proof We note that the set $J:= \textup{ker}\big( (\id \otimes (\rho_{_{\CC, \pi}}) - 1\circ \epsilon) \circ \Delta \big)^{+}$ is a coideal of $D^{\vee}(G)$. If $j  \in J$ then by definition $j_{1} \otimes \rho(j_{2}) - j \otimes 1 = 0$. Applying $\Delta \otimes \id$ to this then shows that $\Delta J \subseteq D^{\vee}(G) \otimes J \subseteq D^{\vee}(G) \otimes J + J \otimes D^\vee(G)$. Moreover, applying instead $\epsilon \otimes \id$ shows that $J \subseteq \textup{ker}(\rho)$. So $J = \{ 0\}$ by Lemma~\ref{corollary_surjective_UL_inner_faith}. 
Next, by \eqref{eq_quantum_maurer_cartan_Cpi} and \eqref{recover_diffcalc}, we have  \[ \extd (\delta_{g} h) = (\delta_{g} h)_{1} \otimes (\rho \circ \beta^{*} \circ \mathcal{Q}_{1} \circ \pi_{\epsilon})(\delta_{g}h)_{2} =  (\delta_{g} h)_{1} \otimes \rho ((\delta_{g}h)_{2}) - \delta_{g}h \otimes 1,\] (since we have previously seen that $\beta^{*} \circ \mathcal{Q}_{1}  = \id$). So, $\textup{ker}(\extd)^{+}$ can be identified with $J$, and as such we have connectedness. \endproof

We note that starting with the penultimate equality of \eqref{eq_matrix_coeff_dual},  we get from \eqref{eq_relation_r_rho} that equivalently $r_{ai}{}^{bj} = \sum_{h \in G} \pi(\zeta_{a}(h^{-1}))^{j}{}_{i} \delta_{h} \otimes b^{-1} \, \delta_{a, hbh^{-1}}$, and as such 
\[ r_{a_{1}i_{1}}{}^{b_{1}j_{1}}...r_{a_{n}i_{n}}{}^{b_{n}j_{n}} = \sum_{h \in G} \pi_{1}(\zeta_{a_{1}}(h^{-1}))^{j_{1}}{}_{i_{1}}...\pi_{n}(\zeta_{a_{n}}(h^{-1}))^{j_{n}}{}_{i_{n}} \delta_{a_{1}, hb_{1}h^{-1}}....\delta_{a_{n}, hb_{n}h^{-1}} \, \delta_{h} \otimes b_{1}^{-1}...b_{n}^{-1} , \]
for any such collection of $r$'s, regardless of which pairs $(\mathcal{C},\pi)$ from $\mathcal{S}$ they come from. With this it is clear that a necessary condition for $U(\mathcal{L}) \to BD(G)$ to be surjective is that $G_{\mathcal{C}^{\mathcal{S}}} = G$ for the subgroup $G_{\mathcal{C}^{\mathcal{S}}}$ of $G$ generated by the Ad-stable subset $\mathcal{C}^{\mathcal{S}}:= \cup_{(\mathcal{C},\pi) \in \mathcal{S}} \mathcal{C}$.

\begin{remark} \label{remark_kG_kofG} The  relationship between $I$ being a cogenerating subspace and connectedness of the calculus is also holds for $k G$ and for $k(G)$ as follows.  For $k(G)$, a bicovariant calculus corresponds to $I = \{ \delta_{c}\}_{c \in \mathcal{C}}$ for an Ad-stable subset $\mathcal{C}$ of $G\backslash \{ e\}$ as recapped in Section~\ref{sec_recap_calc_kG_kofG}. This subset being cogenerating is equivalent to $(k(G)/I)^{*} \cong k\mathcal{C} \subset kG$ being a generating set, which is equivalent to $\mathcal{C}$ generating the group $G$. This is the condition for the connectedness of the calculus.

For  $kG$, we follow a similar strategy to the proof of Lemma~\ref{corollary_surjective_UL_inner_faith}. Firstly note that inner faithfulness of the algebra map $\rho: kG \to \textup{End}(V)$ is equivalent to faithfulness in the usual group-theoretic sense  for group homomorphism $\rho$ \cite{paper_inner_faith}, which is equivalent to connectedness of the calculi in Section~\ref{sec_recap_calc_kG_kofG}. In this case,  if  $J \subseteq \textup{ker}(\rho)$ (where here we mean the algebra kernel) is a coideal of $kG$ then $kG \cdot J \cdot kG$ is a biideal of $kG$ contained again in $\textup{ker}(\rho)$. By inner faithfulness of $\rho$, $kG \cdot J \cdot kG = \{ 0 \}$ implying $J = \{ 0 \}$. Hence $\textup{ker}(\rho)$ is cogenerating. Conversely, assume that $\textup{ker}(\rho)$ is cogenerating.  If $J\subseteq \ker(\rho)$ is a biideal of $kG$, then it is in particular a coideal of $k G$ and hence $J = \{ 0 \}$. Hence $\rho$ is faithful.  \end{remark}

 Another result of interest is an extension of Corollary~\ref{cor_sequence_braided_lie2}.

\begin{corollary} \label{cor_sequence_braided_enveloping}
There is a sequence of braided bialgebras $U(\mathcal{L}_{\pi}) \hookrightarrow U(\mathcal{L}_{\CC, \pi}) \twoheadrightarrow U(\mathcal{L}_{\CC})$ for which 
\begin{equation*} 
          \begin{tikzcd}
         \delta_{\CC, \{ e\}}\mathcal{L}_{\pi}  \arrow[hookrightarrow]{r}\arrow[hookrightarrow]{d}& \mathcal{L}_{\CC, \pi}  \arrow[twoheadrightarrow, swap]{r} \arrow[hookrightarrow]{d}& \mathcal{L}_{\CC}\arrow[hookrightarrow]{d}\\[-10pt]
         U(\delta_{\CC, \{ e\}}\mathcal{L}_{\pi})  \arrow[hookrightarrow]{r} \arrow[d] & U(\mathcal{L}_{\CC, \pi} ) \arrow[twoheadrightarrow, swap]{r} \arrow[d]& U(\mathcal{L}_{\CC}) \arrow[d]\\[-10pt]
                  k(G) \arrow[hookrightarrow]{r}  & BD(G) \arrow[twoheadrightarrow, swap]{r} & kG \\[-10pt]
\end{tikzcd}
\end{equation*}  is a commutative diagram.
\end{corollary}
\proof The braided Hopf algebra map in the category of $G$-crossed modules $U(\mathcal{L}_{\CC, \pi}) \twoheadrightarrow U(\mathcal{L}_{\CC})$ is given on the generators by $i_{1}{}_{*}^{\sharp}$ from Corollary~\ref{cor_sequence_braided_lie2} which was seen to be a braided Lie algebra map, and hence will intertwine the fundamental braidings. So, this map is indeed compatible with the relations defining the enveloping algebras. Similarly for the map $U(\delta_{\CC, \{ e\}}\mathcal{L}_{\pi}) \hookrightarrow U(\delta_{\CC, \{ e\}}\mathcal{L}_{\CC, \pi})$ given by $E_{i}{}^{j} \mapsto E_{i}{}^{j}$ if $\mathcal{C} = \{ e\}$ and $1 \mapsto 1$ otherwise (where $U(\delta_{\CC, \{ e\}}\mathcal{L}_{\pi}) = k$ here since $\mathcal{L}_{\pi} = 0$). 

We note that the map $U(\mathcal{L}_{\CC}) \to kG$ and $U(\mathcal{L}_{\pi}) \to k(G)$ are then given on the generators by \begin{equation} \label{eq_maps_on_gen} e^{c^{-1}} \mapsto c^{-1}, \quad  E_{i}{}^{j} \mapsto  \sum_{n \in G}\pi(n)^{j}{}_{i} \, \delta_{n^{-1}}, \end{equation} i.e. are given simply by the inclusion $\mathcal{L}_{\CC} \hookrightarrow kG$ and $\mathcal{L}_{\pi} \hookrightarrow k(G)$. 

For commutivity of the lower right cell, going from  $U(\mathcal{L}_{\CC, \pi}) \to U(\mathcal{L}_{\CC}) \to kG$
\[ E_{ai}{}^{bj} \mapsto \delta_{i}^{j} \delta_{a, b} e^{a^{-1}} \mapsto \delta_{i}^{j} \delta_{a, b} a^{-1} , \] 
while the map going $U(\mathcal{L}_{\CC, \pi}) \to BD(G) \to kG$ gives the same result
\[ E_{ai}{}^{bj} \mapsto r_{ai}{}^{bj} \mapsto \sum_{n \in C_{G}} \pi(n)^{j}{}_{i} \delta_{q_{a}n^{-1}q_{b}^{-1}, e} b^{-1} = \sum_{n \in C_{G}} \pi(n)^{j}{}_{i} \delta_{a, b} \delta_{n, e} b^{-1} = \delta^{j}_{i} \delta_{a, b}b^{-1}. \] Similarly for the lower left cell.  \endproof

We then also obtain 

\begin{corollary} \label{corollary_surjectivity_other}

(i) The map $U(\mathcal{L}_{\CC}) \to kG$ is surjective iff $\mathcal{C}$ is a generating set. 

(ii)  $U(\mathcal{L}_{\pi}) \to k(G)$ is surjective iff $\pi$ is faithful as a representation of $G$. 
\end{corollary}
\proof (i) is clear from \eqref{eq_maps_on_gen}.  For (ii),  we follow the same ideas as in the proof of Lemma~\ref{corollary_surjective_UL_inner_faith}. Surjectivity of the map is equivalent to the statement that $S= \{1\}\cup \CL_\pi$  being a generating set for the algebra $k(G)$, or $\ker\pi$  cogenerating $kG$. This follows by dualising the map  \eqref{eq_maps_on_gen} to the map
\[ g \mapsto \langle g \sum_{n \in G} \pi(n)^{j}{}_{i} \delta_{n^{-1}} \rangle \, E_{j}{}^{i} = \pi(g^{-1})^{j}{}_{i}E_{j}{}^{i} = \pi(g^{-1}), \quad 1 \mapsto 1. \] We then use the second part of Remark~\ref{remark_kG_kofG}. \qed

Remark~\ref{remark_kG_kofG} already tells us that these surjectivity properties are equivalent to  connectedness of the calculi on $k(G)$ and $k G$ respectively. The surjectivity in case (ii) also implies that the map $U(\mathcal{L}_{\{ e, \pi\}}) \to BD(G)$ has image  $k(G) \otimes e \subset BD(G)$. 

\begin{example} \label{example_image_map}\rm For $G = S_{3}$ and the pairs $(\mathcal{C}, \pi)$ in Example~\ref{exS3}, we study the image of the map $U(\mathcal{L}_{\CC, \pi}) \to BD(G)$. 

(i) For $\mathcal{C} = \{ e\}$, and $\pi=1$ or $\sign$, the $n$-th power of the generator of $\CL_\pi$ maps to 
\[  (\delta_{e} + \delta_{uv} + \delta_{vu}) + (\pm 1)^{n} (\delta_{u} + \delta_{v} + \delta_{w})\big) \otimes e .\] 
 Hence the image of  $U(\mathcal{L}_{\{ e\}, \pi})$ is  $\C1\subset  BD(S_{3})$ for $\pi = 1$ and $ \C(S_{3})^{\mathbb{Z}_{3}} \otimes e \subset BD(S_{3})$ for $\pi = \textup{sign}$.  For $\pi$ the two-dimensional irrep, we already know by Corollary~\ref{corollary_surjectivity_other}, since $\pi$ is faithful,  that the image of $U(\mathcal{L}_{\{ e\}, 2})$ is $C(S_{3}) \otimes e \subset BD(S_{3})$. One can verify this explicitly using the matrix representation  \eqref{eq_matrices_S3}, from which one has
 \[ r_{1}{}^{1} = (\delta_{e} + q\delta_{uv} + q^{-1}\delta_{vu}) \otimes e, \quad r_{1}{}^{2} = (\delta_{u} + q\delta_{v} + q^{-1}\delta_{w})\otimes e, \]
 \[ r_{2}{}^{1} =  (\delta_{u} + q^{-1}\delta_{v} + q\delta_{w}) \otimes e, \quad\quad r_{2}{}^{2} = (\delta_{e} + q^{-1}\delta_{uv} + q\delta_{vu}) \otimes e, \] 
 with $q = e^{2\pi i /3}$, with all elements of $\C(S_3)\tens e$ reachable by these and their 2nd order products.

(ii) For $\mathcal{C} = \{ uv, vu\}$ and $\pi_j$,  the  $\{ r_{a}{}^{b}\}$ are all given by 
\[ r_1{}^1 = (\delta_{e}  + q^{-j} \delta_{uv} + q^{j}\delta_{vu}) \otimes vu, \quad r_1{}^2= (\delta_{u} + q^{j}\delta_{v} + q^{-j}\delta_{w}) \otimes uv, \]
\[ r_2{}^1 = (\delta_{u} + q^{-j} \delta_{v} +  q^{j}\delta_{w}) \otimes vu, \quad r_2{}^2 = (\delta_{e} + q^{j}\delta_{uv} + q^{-j} \delta_{vu}) \otimes uv, \]
numbered in the basis order, with $j$ labelling the choice of irrep as in Example~\ref{exS3}. The only non-zero products of order 2 are the following 
\[ (r_1{}^1)^{2} = r_2{}^2, \quad (r_2{}^2)^{2} = r_1{}^1, \quad r_1{}^1r_2{}^2 = (\delta_{e} + \delta_{vu} + \delta_{uv}) \otimes e, \]
\[ (r_1{}^2)^{2} = q^{-j} r_2{}^1, \quad ( r_2{}^1)^{2} = q^{-j} r_1{}^2, \quad r_1{}^2r_2{}^1 = q^{j}(\delta_{u} + \delta_{v} + \delta_{w}) \otimes e, \] 
given the second order relations in Example~\ref{example_UL}. Similarly to part (i), one can then verify that for $j = 0$  the image of $U(\mathcal{L}_{\mathcal{C}, \pi_{0}})$ is $\C(S_{3})^{\mathbb{Z}_{3}} \otimes \C \mathbb{Z}_{3} \subset D^{\vee}(S_{3})$. The image of $U(\mathcal{L}_{\mathcal{C}, \pi_{j}})$ for $j=1,2$ can also identified with 3 copies of $\C(S_{3})^{\mathbb{Z}_{3}}$, namely $(\delta_{e} + \delta_{vu} + \delta_{uv}) \otimes e, (\delta_{u} + \delta_{v} + \delta_{w}) \otimes e$ as one copy, then by \eqref{eq_useful_identification}, the span of $\{ r_{1}{}^{1}, r_{2}{}^{1}\}$ is $\C(S_{3})^{\mathbb{Z}_{3}} \otimes V_{\{uv,vu\}, \pi_{j}} \cong (D^{\vee}(S_{3}) \otimes V_{\{uv, vu\}\pi_{j}} )^{H} \subseteq D^{\vee}(S_{3})$ (where $V_{\{uv,vu\}, \pi_{j}}$ is one-dimensional). The same can be done for the span of $\{ r_{1}{}^{2}, r_{2}{}^{2}\}$. 

(iii) for $\mathcal{C} = \{ u, v, w\}$ (which generates the group $S_{3}$) and $\pi_{\pm}$, we obtain 
\[ r_{u}{}^{u} = (\delta_{e} \pm \delta_{u}) \otimes u, \quad r_{u}{}^{v} = (\delta_{w} \pm \delta_{vu}) \otimes v, \quad r_{u}{}^{w} = (\delta_{uv} \pm \delta_{v}) \otimes w, \]
\[ r_{v}{}^{u} = (\delta_{w} \pm \delta_{uv}) \otimes u, \quad r_{v}{}^{v} = (\delta_{e} \pm \delta_{v}) \otimes v , \quad r_{v}{}^{w} = (\delta_{u} \pm \delta_{vu}) \otimes w, \]\[r_{w}{}^{u} = (\delta_{vu} \pm \delta_{v}) \otimes u, \quad r_{w}{}^{v} = (\delta_{u} \pm \delta_{uv}) \otimes v, \quad r_{w}{}^{w} = (\delta_{e} \pm \delta_{w}) \otimes w. \]
Working out their nonzero products and relations among them, as well as the relations in Example~\ref{example_UL}, one can then verify that  $U(\mathcal{L}_{\{ u, v, w\}, \pi_{\pm}})\to BD(S_3)$ is surjective.  For example, one can access every element of $\C(S_3)\tens u$ via linear combinations of elements of the form $(r_{a}{}^{u}r_{b}{}^{c})r_{d}{}^{v}$ for $a, b, c, d \in \mathcal{C}$. \endproof
\end{example}

\subsection{Quotient (braided) Hopf algebras of $U(\CL)$ and $A(\mathcal{R})$}

 Finally, we would like to a construction for quotients of $U(\CL)$ and $A(\mathcal{R})$ which are (braided) Hopf algebras and through which the maps $U(\CL)\to BD(G)$ and $A(\mathcal{R})\to D^\vee(G)$ in Lemma~\ref{lem_map_from_enveloping_to_og} factor. These are not necessarily the largest such quotients (where one could use some form of braided or quantum determinant suitable to the $\CR$ matrix) on the other hand, our construction is quite explicit. The idea is, motivated by the case where the above maps are surjective, to `lift' the braided/usual antipodes on $BD(G)$ and $D^\vee(G)$ respectively. 
 
We focus on the case of a single block $\CL_{\CC,\pi}$, but the arguments also extend to the direct sum case. We first note that the braided antipode $\underline{S}$ of $BD(G)$ in Lemma~\ref{lem_braided_Hopf_DG} on the image of $\CL_{\CC,\pi}$ can be computed as
\begin{align*} \underline{S} r_{ai}{}^{bj} & =  \sum_{n \in C_{G}} \pi(n)^{j}{}_{i} \delta_{ab^{-1}q_{b}nq_{a}^{-1}} \otimes a = \sum_{m \in C_{G}} \pi(\zeta_{b}(ab^{-1})^{-1})^{j}{}_{k}\pi(m)^{k}{}_{i} \delta_{q_{ab^{-1}a^{-1} }mq_{a}^{-1}} \otimes a\\
& = \frac{1}{|C_{G}|}\pi(\zeta_{b}(ab^{-1})^{-1})^{j}{}_{k}\sum_{\pi' \in \textup{Irr}(C_{G})} \Big( \textup{dim}(V_{\pi'})\sum_{n\in C_{G}} \pi(n)^{k}{}_{i} \pi'(n)^{x'}{}_{y'} \Big) r_{ab^{-1}a^{-1}x'}{}^{a^{-1}y'}, \end{align*}
 where in the final equality we use the `grand orthonality relation' which gives that $\sum_{\pi' \in \textup{Irr}(C_{G})}\textup{dim}(V_{\pi'})\textup{Tr}_{\pi'} (p) = |C_{G}|\delta_{p, e}$ for all $p \in C_{G}$, and choosing $q_{a} = q_{a^{-1}}$ (which is possible by simply choosing $r^{-1}$ as the representative of $\mathcal{C}^{-1}$). If $\pi$ is real orthogonal, we have orthogonality relations
 \[ \sum_{n\in C_{G}} \pi(n)^{k}{}_{i} \pi'(n)^{x'}{}_{y'} = \sum_{n\in C_{G}} \pi(n)^{*}{}^{k}{}_{i} \pi'(n)^{x'}{}_{y'}  =\frac{|C_{G}|}{\textup{dim}(V_{\pi})} \delta_{\pi, \pi'} \delta_{k, x'}\delta_{i, y'}, \] 
and hence in this case $\underline{S}$ simplifies to a formula in the following proposition that sends elements corresponding to the pair $(\mathcal{C}, \pi)$ to elements corresponding to the pair $(\mathcal{C}^{-1}, \pi)$. This motivates our construction. 

\begin{figure}
\[\includegraphics[scale=1]{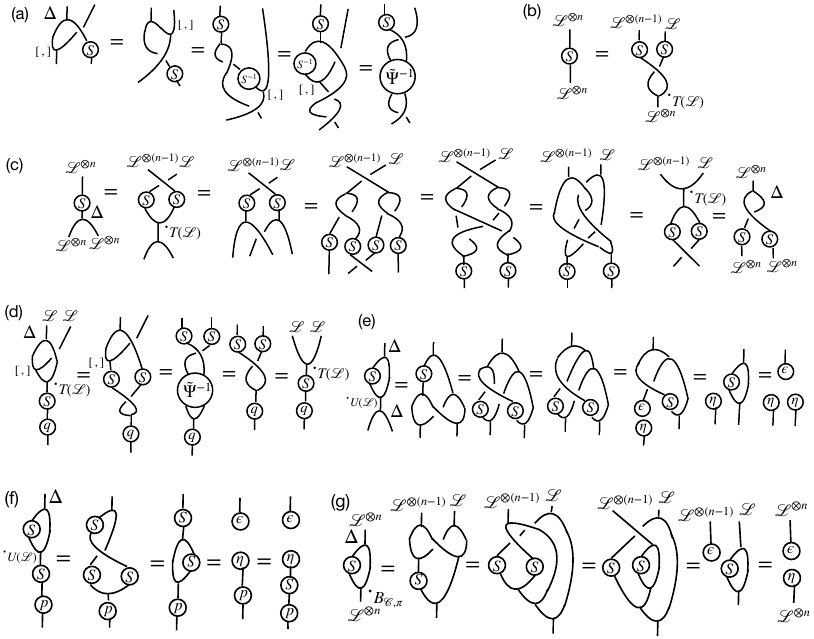}\]
 \vspace{-5mm}
\caption{\label{fig_UL_quotient} Proofs for Proposition~\ref{proposition_further_quotient} (a) Calculation in a braided Hopf algebra $B$ viewed as a braided-Lie algebra. (b) Inductive definition of $\underline{S}$ on the tensor algebra $T(\CL)$. (c) Proof of braided anticomultiplicativity of $\underline{S}$ on the tensor algebra. (d) Proof that $\underline S$ descends to a well-defined map on $U(\CL)$, where $q: T(\mathcal{L}_{\CC, \pi}) \twoheadrightarrow U(\mathcal{L}_{\CC, \pi})$ is the canonical projection.  (e) Proof that the coproduct descends to $B_{\CC, \pi}$.  (f) Proof that $\underline S$ descends to $B_{\CC,\pi}$, where $p: U(\mathcal{L}_{\CC, \pi}) \twoheadrightarrow B_{\CC, \pi}$ is the canonical projection.  (g) Proof by induction that $\underline{S}$ satisfies  the antipode axioms. Mirroring the proofs in (f),(g)  gives the other side.} \end{figure}

\begin{proposition} \label{proposition_further_quotient}
If $\pi$ is real orthogonal and $\mathcal{C}$ is stable under inversion there is a map $\underline{S}:U(\mathcal{L}_{\CC, \pi}) \to U(\mathcal{L}_{\CC, \pi})$ defined on generators
\[ \underline{S}E_{ai}{}^{bj} = \pi(\zeta_{b}(ab^{-1})^{-1})_{j}{}^{k} E_{ab^{-1}a^{-1}k}{}^{a^{-1}i}, \]
 and extended as a unital braided anti-algebra map. Then  
\[ B_{\CC,\pi} := U(\mathcal{L}_{\CC, \pi})/\langle \sum_{c \in \mathcal{C},k}E_{ai}{}^{ck} \underline{S}E_{ck}{}^{bj} - 1\delta_{a,b}\delta_{i,j}, \sum_{c \in \mathcal{C},k}(\underline{S}E_{ai}{}^{ck}) E_{ck}{}^{bj} -  1\delta_{a,b}\delta_{i,j}\ | \, a,b \in \mathcal{C} \rangle\] 
is a braided Hopf algebra in the category of $G$-crossed modules with invertible antipode, through which the map $U(\mathcal{L}_{\CC, \pi}) \to BD(G)$ factors. Moreover, $B_{\CC,\pi}$ can be obtained as the comodule transmutation of the coquasitriangular Hopf algebra 
\[ H_{\CC,\pi}:= A_{\CC,\pi}/\langle  \sum_{c \in \mathcal{C},k}t_{ai}{}^{ck} S t_{ck}{}^{bj} - 1\delta_{a,b}\delta_{i,j},  \sum_{c \in \mathcal{C},k}(St_{ai}{}^{ck})t_{ck}{}^{bj} - 1\delta_{a,b}\delta_{i,j} \, | \, a,b \in \mathcal{C} \rangle,  \] with $S t_{ai}{}^{bj}:= t_{b^{-1}j}{}^{a^{-1}i}$, and the map $A_{\CC, \pi} \to D^{\vee}(G)$ factors through $H_{\CC, \pi}$. 
\end{proposition}
\proof (i) We will use diagrammatic proofs read down the page, with braided coproduct $\Delta=\includegraphics{deltafrag.pdf}$ and  braided antipode $S$ and different `products' $\includegraphics{prodfrag.pdf}$ as marked. There is no ambiguity in the diagrams, but to avoid confusion we will as usual underline the braided antipodes the text.  

(a) Our  starting point is to define  $\underline S:\mathcal{L}_{\CC, \pi} \to \mathcal{L}_{\CC, \pi}$ as the restriction of the braided antipode of $BD(G)$, which applies under our assumptions and where we consider $\CL_{\CC,\pi}\subset BD(G)$ by Corollary~\ref{cor_form_of_inclusion_braided_lie}. From the stated formula for $\underline S$ on the basis of $\CL_{\CC,\pi}$,  it is clear that $\underline{S}^{-1} = \underline{S}$. Also, $\underline{S}$ being a $G$-crossed module map means that 
 \[ \textup{Ad}^{L}_{\delta_{g}h}(\underline{S}r_{ai}{}^{bj}) = \underline{S}\big(\textup{Ad}^{L}_{\delta_{g}h}(r_{ai}{}^{bj})\big),  \] 
hence there is analogous statement for $\underline{S}E_{ai}{}^{bj}$ and action of $D(G)$. Since the bracket of $\mathcal{L}_{\CC, \pi}$ was inherited from this adjoint action of $D(G)$ on itself, it follows that  
 \[ [E_{ai}{}^{bj}, \underline{S}E_{ck}{}^{dl}] = \underline{S}[E_{ai}{}^{bj}, E_{ck}{}^{dl}], \quad (\underline{S} \otimes \id) \circ \widetilde{\Psi}(E_{ai}{}^{bj} \otimes E_{ck}{}^{dl}) = \widetilde{\Psi}(E_{ai}{}^{bj} \otimes \underline{S}(E_{ck}{}^{dl})), \]
 with the second relation obtained by the definition of $\widetilde{\Psi}$ in Figure~\ref{figLie}. Figure~\ref{fig_UL_quotient}(a) deduces  another compatibility between $\widetilde{\Psi}$ and $\underline{S}$. We used $\widetilde{\Psi}^{-1}$ from Figure~\ref{figLie}(e)) restricted to $\CL_{\CC,\pi}$, the first equality is by (L2) of Figure~\ref{figLie}(a) and the second is the braided-anticommutativity of $\underline{S}$ as inherited from the braided antipode of $BD(G)$ (as this holds in any braided-Hopf algebra \cite{Ma}). 
 
Figure~\ref{fig_UL_quotient}(b) next extends $\underline{S}: \mathcal{L}_{\CC, \pi} \to \mathcal{L}_{\CC, \pi}$ as an braided anti-algebra isomorphism $T(\mathcal{L}_{\CC, \pi}) \to T(\mathcal{L}_{\CC, \pi})$ in the category, meaning the map is defined by $\underline{S}(1) = 1$ and on higher degree elements diagrammatically as shown. Figure~\ref{fig_UL_quotient}(c) then by induction that this extension of $\underline{S}$ remains anti-comultiplicative on higher degree elements. Note that  $\underline{S}^{-1}=\underline{S}$ degree 1 similarly extends to  $T(\mathcal{L}_{\CC, \pi})$ using braided anti-multiplicatively as in Figure~\ref{fig_UL_quotient}(b) but with reversed braid crossing. 

Figure~\ref{fig_UL_quotient}(d) then shows that $\underline{S}$ descends to a map on the quotient $U(\mathcal{L}_{\CC, \pi}) \to U(\mathcal{L}_{\CC, \pi})$, which we denote again by $\underline{S}$. Denoting the quotient map $q: T(\mathcal{L}_{\CC, \pi}) \twoheadrightarrow U(\mathcal{L}_{\CC, \pi})$, we show that $\textup{ker}(q) \subseteq \textup{ker}(q \circ \underline{S})$. The second equality uses the compatibilities between $\underline{S}$ and $\widetilde{\Psi}$, and the fourth uses the relation given by applying $\widetilde{\Psi}^{-1}$ to the defining relation \eqref{eq_defining_relation_UL} of $U(\mathcal{L}_{\CC, \pi})$. At this point we have a well-defined map $\underline{S}: U(\CL_{\CC,\pi})\to U(\CL_{\CC,\pi})$ and it inherits the braided antimultiplicativity and braided anti-comultiplicativity from the earlier parts. Moreover, it is again bijective, by analogous constructions for $\underline{S}^{-1}$.

We now show in Figure~\ref{fig_UL_quotient}(e) that quotienting by the two-sided ideal generated in $U(\mathcal{L}_{\CC, \pi})$ by the stated relations in the statement (where $1\delta_{a, b}\delta_{i, j} = \epsilon(E_{ai}{}^{bj})$) is indeed compatible with the coproduct of $U(\mathcal{L}_{\CC, \pi})$, and hence we obtain a bialgebra $B_{\CC, \pi}$. It then follows from Figure~\ref{fig_UL_quotient}(f) that $\underline{S}$ further descends to give a map $\underline{S}: B_{\CC, \pi} \to B_{\CC, \pi}$. This is again a bijection, since $\underline{S}^{-1}: U(\mathcal{L}_{\CC, \pi}) \to U(\mathcal{L}_{\CC, \pi}) $  descends by an analoguous proof. We denote by $p: U(\mathcal{L}_{\CC, \pi}) \twoheadrightarrow B_{\CC, \pi}$ the canonical quotient map and show that $\textup{ker}(p) \subseteq \textup{ker}(p\circ \underline{S})$. The final equality uses $\epsilon \circ \underline{S} = \epsilon$ which is trivial on degree 0, inherited on degree 1 and extends inductively to higher degrees. It just remains to show in  Figure~\ref{fig_UL_quotient}(g) that $\underline{S}$ on the quotient indeed satisfies the antipode axioms. Hence, $B_{\CC,\pi}$ is indeed a braided Hopf algebra in the category of $G$-crossed modules. It is clear that the map $U(\CL_{\CC,\pi}) \to BD(G)$ factors through $B_{\CC,\pi}$ since its additional relations already hold in $BD(G)$. 

(ii) Similarly, the map $S: A(\mathcal{R}) \to A(\mathcal{R})$ on degree 1 elements in the statement originates from \eqref{eq_relation_r_rho}, where 
\begin{align}S r_{ai}{}^{bj} & = \sum_{n \in C_{G}} \pi(n)^{j}{}_{i} \delta_{q_{b}nq_{a}^{-1}} \otimes a \nonumber\\
& = \sum_{\pi' \in \textup{Irr}(C_{G})}\Big(\frac{\textup{dim}(V_{\pi'})}{|C_{G}|} \sum_{n \in C_{G}} \pi(n)^{j}{}_{i}\pi'(n)^{x'}{}_{y'} \Big)r_{b^{-1}x'}{}^{a^{-1}y'} = r_{b^{-1}j}{}^{a^{-1}i}, \label{eq_antipode_on_generator}\end{align} 
again using the `grand orthogonality relation' for the second equality, and the assumption that $\pi$ is real for the final equality. Since this map $S$ is defined to commute with the inclusion $j: t_{ai}{}^{bj} \mapsto r_{ai}{}^{bj} \in D^{\vee}(G)$ and the antipode of the latter, we follow the same procedure and diagrams as above, but for ordinary Hopf algebras. This gives that $H_{\CC,\pi}$ is indeed a well-defined Hopf algebra through which $A_{\CC,\pi}\to D^\vee(G)$ factors for the same reasons as before. It inherits the coquasitriangular structure of $A(\mathcal{R})$ by elementary properties of the latter (and the fact that the map $A_{\CC, \pi} \to D^{\vee}(G)$ intertwined the coquasitriangular structures). It is also clear that the comodule transmutation $B(H_{\CC,\pi})$ is then $B_{\CC,\pi}$ since $\underline{S}$ was obtained as a transmutation of $S$, and so the relations give the same quotient object. \endproof 

We note that more generally for direct sum braided Lie algebras $\mathcal{L}$, we can obtain a quotient $B(\mathcal{L})$ in the same way as in the above if the braided antipode $\underline{S}$ of $BD(G)$ restricts to a map $\underline{S}: \mathcal{L} \to \mathcal{L}$, where $\mathcal{L}$ is viewed inside $BD(G)$ via $\rho^{\sharp}$ from Corollary~\ref{cor_form_of_inclusion_braided_lie}. Immediate from the definition of $B(\mathcal{L})$ is the following.

\begin{corollary}
Under the assumptions of Proposition~\ref{proposition_further_quotient}, $B(\mathcal{L})$ is universal amongst braided Hopf algebras $F$ (in the category of $G$-crossed modules) that factorise the map $U(\mathcal{L}) \to BD(G)$ in such a way that $S|_{\mathcal{L}} = \underline{S}$ with $S$ the antipode of $F$ and $\underline{S}: \mathcal{L} \to \mathcal{L}$.
\end{corollary}

\proof Denoting the braided bialgebra maps $U(\mathcal{L}) \xrightarrow{p} F \xrightarrow{q} BD(G)$, then $p \circ \cdot \circ (\underline{S} \otimes \id) \circ \Delta|_{\mathcal{L} \subset F} = \cdot \circ (S \otimes \id) \circ \Delta \circ p|_{\mathcal{L}} = \eta \circ \epsilon \circ p|_{\mathcal{L}} = p \circ \eta \circ \epsilon|_{\mathcal{L}}$ when viewed in $U(\mathcal{L})$. As such $\textup{ker}(U(\mathcal{L}) \twoheadrightarrow B(\mathcal{L})) \subseteq \textup{ker}(p)$. \endproof

Lastly, we note that we can extend Corollary~\ref{cor_sequence_braided_enveloping} by another splitting the lower row with 
\begin{equation*} 
          \begin{tikzcd}
                 U(\delta_{\CC, \{ e\}}\mathcal{L}_{\pi})  \arrow[hookrightarrow]{r} \arrow[twoheadrightarrow, swap]{d}& U(\mathcal{L}_{\CC, \pi} ) \arrow[twoheadrightarrow, swap]{r} \arrow[twoheadrightarrow, swap]{d}& U(\mathcal{L}_{\CC}) \arrow[twoheadrightarrow, swap]{d}\\[-10pt]
                  B(\delta_{\CC, \{ e\}}\mathcal{L}_{\pi}) \arrow[hookrightarrow]{r}  \arrow[d]& B_{\CC, \pi} \arrow[twoheadrightarrow, swap]{r} \arrow[d]& kG_{\CC} \arrow[d]\\[-10pt]
                  k(G) \arrow[hookrightarrow]{r}  & BD(G) \arrow[twoheadrightarrow, swap]{r} & kG
\end{tikzcd}
\end{equation*}
where $G_{\CC}$ and the right-hand column are as in  \cite{MaRie}, thereby generalising the point of view there. Indeed, one can check that $\textup{ker}(U(\mathcal{L}_{\CC, \pi}) \twoheadrightarrow B_{\CC, \pi}) \hookrightarrow \textup{ker}(U(\mathcal{L}_{\CC, \pi}) \twoheadrightarrow kG_{\CC}) $ by applying $i_{1}{}_{*}^{\sharp}$ to $E_{ai}{}^{ck} \underline{S}E_{ck}{}^{bj} - 1\delta_{a,b}\delta_{i,j}$, and similarly for the antipode on the other side.

\begin{example}\rm 
For $G=S_3$ and the 3 conjugacy classes, we first study $H_{\CC,\pi}$ as a quotient of $A_{\CC,\pi}$ in Example~\ref{example_rel_AR}.

(i) For $\mathcal{C} = \{ e\}$, we have that $S t_{i}{}^{j} = t_{j}{}^{i}$. For the one-dimensional representations we have $S t=t$ and hence  $H_{\CC,\pi} = \C[t]/(t^{2} - 1)$. For the two-dimensional representation (where we can indeed pick $\pi$ to be real orthogonal), the additional relations for  $H_{\CC,\pi}$ are then just $t_{i}{}^{1}t_{j}{}^{1} + t_{i}{}^{2}t_{j}{}^{2}  = \delta_{i, j}$ and similarly with upper and lower indices flipped. 

(ii) For $\mathcal{C} = \{ uv, vu\}$ with $\pi_j$, only $\pi_0$ is real as needed in our construction, but we can still consider the formula  \eqref{eq_antipode_on_generator}. Keeping track of which $\pi_j$ the generator is imaged in, 
\[S(r_{_{j}})_{a}{}^{b} = \frac{1}{3} \sum_{k = 0}^{2}\sum_{l = 0}^{2} q^{jl}q^{kl} (r_{_{k}})_{b^{-1}}{}^{a^{-1}} = \frac{1}{3} \sum_{k = 0}^{2}(1 + q^{j+k} + q^{-(j+k)})(r_{_{k}})_{b^{-1}}{}^{a^{-1}}. \] 
For $j=0$ we indeed recover $S(t_{a}{}^{b}) = t_{b^{-1}}{}^{a^{-1}}$ for $a,b\in \CC$ as labels. This gives  the additional relations 
\[ t_{a}{}^{1}t_{a^{-1}}{}^{2} +t_{a}{}^{2}t_{a^{-1}}{}^{1}  = 1, \quad t_{a}{}^{1}t_{a}{}^{2}+  t_{a}{}^{2}t_{a}{}^{1}= t_{1}{}^{a}t_{2}{}^{a}+  t_{2}{}^{a}t_{1}{}^{a}=0,\] 
for all $a\in \CC$, using the numbering in Example~\ref{example_rel_AR}(ii) where inversion in $\CC$ swaps $1,2$. These amount to the additional relations on $A_{\CC, \pi}$
\[ t_1{}^1 t_2{}^2+ t_1{}^2 t_2{}^1=1,\quad t_1{}^1 t_1{}^2=t_1{}^1 t_2{}^1= t_2{}^2 t_1{}^2=t_2{}^2 t_2{}^1=0, \]
which is a little unusual but compatible with the matrix coproduct (we already met the second set for $j=1,2$). Thus $H_{\CC,\pi_0}$  is the algebraic group of $\C^*\times \Z_2$ (with $\C^*$ diagonal matrices with product 1 and $\Z_2$ the transposition matrix). For the complex cases of $j = 1, 2$ (viewed mod 3),   $S (r_j)_{a}{}^{b} = (r_{-j})_{b^{-1}}{}^{a^{-1}}$ so that we would need to work with $A(\CR)$ for $\mathcal{L} = \mathcal{L}_{\mathcal{C}, \pi_{1}}\oplus \mathcal{L}_{\mathcal{C}, \pi_{2}}$ in order to be able to apply our construction. 

(iii) For $\mathcal{C} = \{ u,v, w\}$, we have $S(t_{a}{}^{b}) = t_{b}{}^{a}$ and the additional relations become 
\[t_{a}{}^{u}t_{b}{}^{u} +  t_{a}{}^{v}t_{b}{}^{v} + t_{a}{}^{w}t_{b}{}^{w} = \delta_{a, b}, \] and the same with upper and lower indices flipped. This amounts to the additional 5 relations in the case of $\pi_{+}$ (for simplicity)
\[ t_u^2+ (t_u{}^v)^2+(t_u{}^w)^2=t_v^2+ (t_v{}^u)^2+(t_v{}^w)^2=t_w^2+ (t_w{}^u)^2+(t_w{}^v)^2=1,\]
\[ t_ut_v{}^u+t_v t_w{}^v+t_w t_u{}^w= t_u t_w{}^u+ t_v t_u{}^v+t_w t_v{}^w=0, \]
and the same with upper and lower indices swapped, given  the existing commutation relations of $A_{\CC,\pi_+}$ from Example~\ref{example_rel_AR}(iii). Only two of the three new $(t_{a})^{2}$ relations are independent given the existing relations, so that $H_{\CC,\pi_+}$ is 24 dimensional in degree 2. By our construction, it is a coquasitriangular Hopf algebra. 
\end{example}

\begin{example}\rm 
For the braided case $B_{\CC,\pi}$ we just focus on the more interesting case (iii) with $\CC=\{u,v,w\}$ and, for simplicity, $\pi_+$. Then
\[ \underline{S} E_a{}^b=E_{\eps(a,b) }{}^a, \]
and one can check directly that this extends braided-antimultiplicatively to degree 2 respecting the relations. For example, applying $\underline{S}$ to the first of the off-diagonal relations in Example~\ref{example_UL}(iii),
\begin{align*}
\underline S(E_a{}^b E_a{}^{\eps(a,b)})=(\underline S E_{ab\la a}{}^{ab\la \eps(a,b)})\underline S E_a{}^b =E_{\eps(a,b)}{}^bE_{\eps(a,b)}{}^a,\\
\underline S(E_{\eps(a,b)}{}^a E_{\eps(a,b)}{}^b)=(\underline S E_{ab\la \eps(a,b)}{}^{ab\la b})\underline S E_{\eps(a,b)}{}^a=E_b{}^a E_b{}^{\eps(a,b)},
\end{align*}
using the $G$ crossed module braiding in expanding $\underline S$ on a product. We need these to be equal and indeed they are by using the same off-diagonal relation with $a,b$ swapped. Having a well-defined map, we then impose the antipode relations for $B_{\CC,\pi_+}$, 
\[ E_a{}^u E_{\eps(u, b)}{}^u+   E_a{}^v E_{\eps(v, b)}{}^v+ E_a{}^w E_{\eps(w, b)}{}^w=\delta_{a,b}=E_{\eps(a, u)}{}^{a}E_{u}{}^{b} +  E_{\eps(a, v)}{}^{a}E_{v}{}^{b} + E_{\eps(a, w)}{}^{a}E_{w}{}^{b}.\]
Given existing relations in $U(\CL_{\CC},\pi_+)$, the  first half of the antipode relations amount to the five relations
\[ e_u^2 + E_u{}^v E_w{}^v+ E_u{}^w E_v{}^w=e_v^2 +E_v{}^uE_w{}^u+E_v{}^wE_u{}^w=e_w^2+E_w{}^v E_u{}^v+E_w{}^u E_v{}^u=1,\]
\[ e_u E_w{}^u+e_vE_u{}^v+e_w E_v{}^w=e_u E_v{}^u+ e_v E_w{}^v+e_w E_u{}^w=0, \]
while the other half introduces two further relations (with a similar $e_w^2$ relation redundant)
\[ e_u^2+E_w{}^u E_v{}^u+E_v{}^u E_w{}^u=e_v^2+E_w{}^v E_u{}^v+ E_u{}^vE_w{}^v=1, \]
and two more relations which are the up-down swap of the stated diagonal $\times$ off-diagonal relations from the first half. This gives us $B_{\CC, \pi_+}$ as 24-dimensional in degree 2, as must be the case as the transmutation of $A_{\CC, \pi_{+}}$. Note that the surjection $B_{\CC, \pi_{+}} \twoheadrightarrow BD(G)$ cannot be an isomorphism since, for example, we know that $r_{u}{}^{u}r_{u}{}^{w} = 0$ in $BD(G)$  from Example~\ref{example_image_map}(iii) with no corresponding constraint on the $E_a{}^b$ side. So this is a genuine extension of $BD(G)$ in this example (in contract to the rack case \cite{MaRie}, where $G_{\CC} \to S_{3}$ is an isomorphism for this choice of $\CC$).
\end{example}

\subsection{The braided Killing form} Another general construction for braided-Lie algebras is the notion of a braided-Killing form $K:\CL\tens\CL\to \underline 1$  to the unit object of the braided category, in our case $k$. This is defined as a braided trace of $[\ ,[\ ,\ ]]$ as shown in Figure~\ref{figLie} part (c). This has already proven interesting for $k\CC\subseteq k G$ as a braided-Lie algebra with trivial braiding, see \cite{LMR}.

\begin{proposition}
The braided Killing forms of $\mathcal{L}_{BD(G)}$ and $\mathcal{L}_{\mathcal{C},\pi}$ are given by \[ K_{BD(G)}(\delta_{u}v,  \delta_{g}h) = \delta_{u, vgv^{-1}}\sum_{k , l \in Z(vu^{-1}(vhv^{-1}))} \delta_{u, |\delta_{k}l|},\]
\begin{align*} K_{\mathcal{C},\pi}(E_{ai}{}^{bj},  E_{ck}{}^{dl}) & =  \delta_{ab^{-1},dc^{-1}} \sum_{f, g \in \mathcal{C}}  \delta_{ab^{-1}, \{ gf^{-1}, d\}}\delta_{b^{-1}d^{-1}g, Z(f)}\delta_{dbf, Z(g)}  \\
& \hphantom{xxxxxx}\pi(\zeta_{c}(d^{-1}fg^{-1}d))^{l}{}_{k} \pi(\zeta_{a}(fg^{-1}))^{j}{}_{i} \textup{Tr}_{\pi}(\zeta_{b^{-1}d^{-1}fdb}(fg^{-1}db)^{-1})\textup{Tr}_{\pi}(\zeta_{g}(bf)), \end{align*} where $Z(h)\subseteq G$ is the centraliser of the element $h \in G$ and $a,b,c,d\in \CC$.
\end{proposition}
\proof We carry out the computation explicitly for $\mathcal{L}_{\mathcal{C},\pi}$: 
\begin{align*}
K_{\mathcal{C},\pi}&(E_{ai}{}^{bj} \otimes E_{ck}{}^{dl})  = \sum_{f, g \in \mathcal{C}}\textup{ev} \circ \Psi_{\mathcal{L}_{\mathcal{C},\pi}, \mathcal{L}_{\mathcal{C},\pi}^{\sharp}} \circ ([\_, \_] \otimes \id)( E_{ai}{}^{bj} \otimes [E_{ck}{}^{dl}, E_{fm}{}^{gn}] \otimes E_{gn}{}^{fm}\\
&= \sum_{f, g \in \mathcal{C}} \delta_{cd^{-1}fg^{-1}, d^{-1}fg^{-1}d}  \pi(\zeta_{c}(d^{-1}fg^{-1}d))^{l}{}_{k} \pi(\zeta_{f}(d^{-1}))^{x}{}_{m}\pi(\zeta_{g}(d^{-1})^{-1})^{n}{}_{y} \\
&\qquad\qquad \textup{ev} \circ \Psi_{\mathcal{L}_{\mathcal{C},\pi}, \mathcal{L}_{\mathcal{C},\pi}^{\sharp}} ( [E_{ai}{}^{bj}, E_{d^{-1}fd x}{}^{d^{-1}gdy}] \otimes E_{gn}{}^{fm})\\
&= \sum_{f, g \in \mathcal{C}} \delta_{cd^{-1}, \{d, gf^{-1}\}}  \pi(\zeta_{c}(d^{-1}fg^{-1}d))^{l}{}_{k} \pi(\zeta_{f}(d^{-1}))^{x}{}_{m}\pi(\zeta_{g}(d^{-1})^{-1})^{n}{}_{y}\\
&\qquad  \delta_{ab^{-1}d^{-1}fg^{-1}d, b^{-1}d^{-1}fg^{-1}db}\pi(\zeta_{a}(b^{-1}d^{-1}fg^{-1}db))^{j}{}_{i} \pi(\zeta_{d^{-1}fd}(b^{-1}))^{u}{}_{x}\pi(\zeta_{d^{-1}gd}(b^{-1})^{-1})^{y}{}_{v}\\
&\qquad \textup{ev} \circ \Psi(E_{b^{-1}d^{-1}fdbu}{}^{b^{-1}d^{-1}gdbv} \otimes E_{gn}{}^{fm})\\
&= \sum_{f, g \in \mathcal{C}} \delta_{cd^{-1}, \{d, gf^{-1}\}} \delta_{ab^{-1}, \{ b, d^{-1}gf^{-1}d\}} \pi(\zeta_{c}(d^{-1}fg^{-1}d))^{l}{}_{k} \pi(\zeta_{f}(d^{-1}))^{x}{}_{m}\pi(\zeta_{g}(d^{-1})^{-1})^{n}{}_{y} \\
&\qquad \pi(\zeta_{a}(b^{-1}d^{-1}fg^{-1}db))^{j}{}_{i} \pi(\zeta_{d^{-1}fd}(b^{-1}))^{u}{}_{x}\pi(\zeta_{d^{-1}gd}(b^{-1})^{-1})^{y}{}_{v}\\
&\qquad  \textup{ev} (b^{-1}d^{-1}fg^{-1}db \triangleright E_{gn}{}^{fm} \otimes E_{b^{-1}d^{-1}fdbu}{}^{b^{-1}d^{-1}gdbv}),
\end{align*}
where we recall that we had $\textup{End}(V_{\mathcal{C},\pi})^{\sharp} = \textup{End}(V_{(\mathcal{C}, \pi})$ as a $G$-crossed module with respect the linear evaluation extended to tensor products in the categorical way.  Therefore, $\textup{ev}(b^{-1}d^{-1}fg^{-1}db \triangleright E_{gn}{}^{fm} \otimes E_{b^{-1}d^{-1}fdbu}{}^{b^{-1}d^{-1}gdbv} ) $ equals
\begin{align*} & \pi(\zeta_{g}(b^{-1}d^{-1}fg^{-1}db))^{r}{}_{n} \pi(\zeta_{f}(b^{-1}d^{-1}fg^{-1}db)^{-1})^{m}{}_{s} \\
&\qquad\qquad \delta_{b^{-1}d^{-1}fg^{-1}dbfb^{-1}d^{-1}gf^{-1}db, b^{-1}d^{-1}fdb}\delta_{r}^{v} \delta_{b^{-1}d^{-1}fg^{-1}dbgb^{-1}d^{-1}gf^{-1}db, b^{-1}d^{-1}gdb} \delta^{s}_{u}\\
& =  \pi(\zeta_{g}(b^{-1}d^{-1}fg^{-1}db))^{r}{}_{n} \pi(\zeta_{f}(b^{-1}d^{-1}fg^{-1}db)^{-1})^{m}{}_{s} \delta_{g^{-1}dbfb^{-1}d^{-1}, fg^{-1}}\delta_{r}^{v} \delta_{g^{-1}dbgb^{-1}d^{-1}g, f^{-1}gf} \delta^{s}_{u}. \end{align*}
 Next, the first Kronecker $\delta$ above equivalently gives that $b^{-1}d^{-1}g \in Z(f)$, which substituting into the second delta gives that $dbf \in Z(g)$. Conversely, these two properties recover the two delta functions. These together also imply that $fg^{-1} \in Z(db)$. This simplifies the above computation to give
\begin{align*} 
& K_{\mathcal{C},\pi}(E_{ai}{}^{bj} \otimes E_{ci}{}^{dj})  \\
&= \sum_{f, g \in \mathcal{C}} \delta_{cd^{-1}, \{ d, gf^{-1}\}} \delta_{ab^{-1}, \{ gf^{-1}, d\}} \pi(\zeta_{c}(d^{-1}fg^{-1}d))^{l}{}_{k} \pi(\zeta_{a}(fg^{-1}))^{j}{}_{i} \delta_{b^{-1}d^{-1}gf, fb^{-1}d^{-1}g}\delta_{dbfg, gdbf}\\
& \qquad  \textup{Tr}_{\pi}(\zeta_{f}(d^{-1})\zeta_{f}(fg^{-1})^{-1}\zeta_{d^{-1}fd}(b^{-1}))\textup{Tr}_{\pi}(\zeta_{g}(d^{-1})^{-1}\zeta_{d^{-1}gd}(b^{-1})^{-1}\zeta_{g}(fg^{-1})),
\end{align*}
which, using \eqref{def_cocycle} and the Kronecker $\delta$s, gives the statement. One can check that the general formula \cite{Ma:lie} for the braided Killing form of a matrix braided Lie algebras  gives the same result using  the structure in Remark~\ref{rem_braided_matrix_lie}.  \endproof

For $\CC=\{e\}$,  the centralisers are $G$ and all the $\delta$-functions are automatic. Moreover $f=g=e$ so there is no sum and we have $\zeta_e(e)=e$ in the general expression for $K$. Hence we immediately obtain
\[ K(E_i{}^j\tens E_k{}^l)=\dim(V_\pi)^2\delta_{ij}\delta_{kl},\]
which is necessarily degenerate when $\pi$ is not 1-dimensional and otherwise has $K=1$. At the other extreme, for general $\CC$ and $\pi$ is 1-dimensional,  the general formula for the Killing form becomes
\begin{align*} K(e_a{}^b,e_c{}^d)& =  \delta_{ab^{-1},dc^{-1}} \sum_{f, g \in \mathcal{C}}  \delta_{ab^{-1}, \{ gf^{-1}, d\}}\delta_{b^{-1}d^{-1}g, Z(f)}\delta_{dbf, Z(g)}  \\
& \qquad \pi(\zeta_{c}(d^{-1}fg^{-1}d))  \pi(\zeta_{a}(fg^{-1})) \pi(\zeta_{b^{-1}d^{-1}fdb}(fg^{-1}db)^{-1}). \pi(\zeta_{g}(bf))\end{align*}

\begin{example}\rm For $G=S_3$, the case (i) where $\CC=\{e\}$ gives the trivial results already stated for the 3 different $\pi$. Two of these have $\pi$ 1-dimensional with $K=1$ and these can also be viewed as trivial examples of the second kind.

For case (ii) where $\CC=\{uv,vu\}$, we note that $\CC\subset\Z_3$ so all elements commute. Moreover, $Z(f)=Z(g)=\Z_3$ and $b^{-1}d^{-1}g, dbf \in \Z_3$ so there is no constraint on $f,g$ from these. Also the group commutator is trivial so the other $\delta$-function gives $\delta_{a,b}$. Using that $\zeta(a,b)=\zeta(b,a)$ in the present case and that the elements of $\CC$ and that $\zeta_a(\  )$ acts by group automorphisms in this case,  we then obtain
\[ K(e_a{}^b,e_c{}^d)=\delta_{a,b}\delta_{c,d}\sum_{f,g}{\pi(\zeta_f(g)^2)\over \pi(\zeta_d(g)\zeta_f(f))}=\delta_{a,b}\delta_{c,d}\sum_{f,g}\pi_j\left(r^{1+\delta_{g,f}+\delta_{g,d}}\right)=\delta_{a,b}\delta_{c,d} q^{2j}\]
in the $\pi_j$ representation. The Killing form here is highly degenerate. 

For case (iii) where $\CC=\{u,v,w\}$ and $\pi_\pm$, we note that the elements of $\CC$ are order 2,  $bdg,dbf\in \CC$ and $\CC\cap Z(f)=\{f\}$, $\CC\cap Z(g)=\{g\}$ so that the last two $\delta$-functions are the same and just set $g=dbf$. Then $\{gf^{-1},d\}=\{db,d\}=bdddbd=(bd)^2=db$ since $\CC^2=\{e,uv,vu\}$ has all elements cubing to $e$. Hence  the remaining $\delta$-function in the sum is $\delta_{a,d}$, independently of $f$. Making these substitutions,  we have
\[K(e_a{}^b,e_c{}^d) = \delta_{b,c}\delta_{a,d}  \pi(\zeta_{b}(d)\zeta_d(d)\zeta_b(b)) \sum_{f\in \mathcal{C}}  \pi(\zeta_{dbf}(bf)),\]
since 
\begin{align*} \zeta_{b}(dbd)) \zeta_{d}(bd)&=\zeta_b(bdb)\zeta_d(b)\zeta_d(d)=\zeta_b(bd)\zeta_b(b)\zeta_d(b)\zeta_d(d)=\zeta_{dbd}(b)\zeta_b(d)\zeta_b(b)\zeta_d(b)\zeta_d(d)\\
&=\zeta_{bdb}(b)\zeta_d(b) \zeta_b(d) \zeta_b(b) \zeta_d(d)=\zeta_{d}(b^2) \zeta_b(d) \zeta_b(b) \zeta_d(d)=\zeta_b(d) \zeta_b(b) \zeta_d(d), \end{align*}
using the cocycle properties and $dbd=bdb$ for all $b,d\in \CC$. In the $\pi_+$ case we just obtain $K_+(e_a{}^b,e_c{}^d) = 3\delta_{b,c}\delta_{a,d}$ while in the $\pi_-$ case we have
\[ K_-(e_a{}^b,e_c{}^d)=\delta_{b,c}\delta_{a,d}\begin{cases} \pi_{bd}(1-\pi_{bd}-\pi_{\eps(b,d)d}) &{\rm if\ }b\ne d\\
1 & {\rm if\ }b=d\ne v\\ 
-3 & {\rm if\ }b=d=v.\end{cases}\]
The only nonzero entries in the standard basis order are then
\[ K_-(e_d{}^b,e_b{}^d)=K_{bd}, \quad K_{\cdot \cdot}=\begin{pmatrix}1 & -3 & -1\\ -1 & -3 & -1\\ -1 & -3 & 1\end{pmatrix},\]
in contrast to $K_+$, where the corresponding $3\times 3$ matrix has all entries 3. In both cases, $K_\pm$ as a $9 \times 9 $ matrix is invertible, i.e. the braided-Killing form is non-degenerate.
 \end{example}

\section{Concluding remarks}\label{secrem}

We have studied the quantum differential geometry of $D^\vee(G)$ and its relationship to the construction of irreps of $D(G)$. Working first at the level of universal calculus or Hopf-Galois extensions in Section~\ref{secWigner}, we found that irreps of the form $\C \CC\tens V_\pi$ or more generally of certain subspaces of $\C \CC\tens W$ should be viewed as fields in $\C  G$ in the role of `spacetime' from the point of view of the Wigner construction for the Poincar\'e group. Here $\CC$ plays the role of a mass shell in Minkowski space and the construction initially gives the irrep as fields over this, namely the space of sections $E$ of a bundle. These are then converted  by Fourier transform/transfer map to a certain subspace of the sections $E'$ of a bundle over $\C  G$. In the Poincar\'e case the particular subspace is characterised by the Klein-Gordon equations, while in our case it ends up being the span of $\CC^{-1}$. This then left the question, which we addressed in Section~\ref{secPcalc}, of what should be the analogue of the Klein-Gordon equation? Such subspaces turned out to correspond to $D(G)$-covariant (or rather, $D^\vee(G)$-covariant in our comodule setting) maps $L:\C G \to \C G$ in $\C  G$ for a fixed bicovariant calculus there. The case that deserves to be called `2nd order' is characterised with respect to a $D(G)$-invariant bimodule inner product $(\ ,\ )$ on $\Omega^1$, justifying its role as `isometry quantum group'. We then asked if  there is a connection $\nabla$ on $\Omega^1$ that completes this quantum geometric picture and showed how this works in a nontrivial example with $G=S_3$. 

Strikingly, given $V_{\CC,\pi}$ obtained this way, there is also a calculus $\Omega_{\CC,\pi}$ on $D^\vee(G)$. Moreover it can be restricted to a $G$-module $\tilde \pi$ and used to define a calculus on $\C  G$. Thus, there is a kind of interplay between the calculus on $\C  G$ as `spacetime' and the quasiparticle or irrep which we are aiming to view as fields in this `noncommutative spacetime'. This seems to be a new phenomenon with no classical analogue, where classically the calculus on $\R^{1,3}$ is taken as the standard one prior to constructing irreps by the Wigner construction. This novel interplay was put into a coherent setting in Section~\ref{sec_duality} as a duality between conjugacy classes and irreps of $G$, where the dual model has $\C(G)$ as spacetime and has the roles of $\CC,\tilde\pi$  swapped. This is a novel phenomenon, building on \cite{MaTao}, which deserves to be studied further.  Similarly, in the last section, we looked, for completeness at the dual of the calculus $\Omega^1_{\CC,\pi}$ as braided-Lie algebras $\CL_{\CC,\pi}$ for $D(G)$. We also introduced associated coquasitriangular and braided matrix bialgebras $A_{\CC,\pi}$ and $U(\CL_{\CC,\pi})$ and, in nice cases, quotient Hopf algebra or braided-Hopf algebras, but without a fully general construction for these.

These are some directions for further work. In addition, many of the results generalise to the double $D(H)$ and $D^\vee(H)$ starting with a Hopf algebra $H$ in place of $\C (G)$. This and other, more algebraic, results will be a sequel \cite{MaMc} and may then relate back to the physics of 2+1 quantum gravity. Meanwhile, $D(G)$ has an immediate role in quantum computing, as in the Kitaev model\label{Kit,CowMa}, and its differential structures/braided Lie algebras as studied here should have a currently overlooked role in the interpretation of the model as some kind of (lattice) gauge theory, which currently is effectively limited to working with the universal calculus. This represents another direction for further work.

 \end{document}